%% file: thesis.tex
\documentclass[titlepage, 12pt, a4paper, twoside, openright]{report}
\usepackage{amssymb}
\usepackage{amsthm}
\usepackage{amsmath}
\usepackage{lscape}
\usepackage{amscd}
\usepackage{rotating}
\usepackage{epsfig}

\renewcommand{\abstractname}{Introduction}  
\renewcommand{\abstract}{\chapter*{\abstractname}
\addcontentsline{toc}{chapter}{\numberline{}\abstractname}}

\newcommand{\bib}{
\addtocounter{page}{1}
\addcontentsline{toc}{chapter}{\numberline{}\bibname}
\addtocounter{page}{-1}
\bibliographystyle{hplain}
\bibliography{gerbe}}

\def\bull{\vbox{\hrule\hbox{\vrule\kern3pt\vbox{\kern6pt}\kern3pt\vrule}\hrule}}

\def\dd{\text{dd}}

\def\L{\mathcal L}

\def\<#1,#2>{\langle #1,#2\rangle}
\def\Tr{\text{Tr}\,}
\def\dep(#1,#2){\text{det}_{#1}#2}
\def\norm(#1,#2){\parallel #1\parallel_{#2}}
\def\Z{{\mathbb Z}}

\def\R{{\mathbb R}}
\def\D{{\mathcal D}}
\def\C{{\mathbb C}}

\def\U{\mathcal U}

\def\hol{\text{hol}}
\def\Map{\text{Map}}
\def\Hom{\text{Hom}}
\def\Ext{\text{Ext}}
\def\text{\mathrm}
\def\isom{\cong}

\def\triv{\mathfrak{T}}
\def\V{\mathcal V}
\def\uz{\underline{\Z}}
\def\ur{\underline{\R}}
\def\uu{\underline{U(1)}}

\def\ua{\underline{\Omega}}
\def\P{\mathcal{P}}
\def\pa{\partial}
\def\A{\mathcal{A}}
\def\i{\iota}
\def\Ad{\text{Ad}}
\def\dg{dg g^{-1}}
\def\dd{\text{dd}}
\def\pfa{\text{Pfaff}}
\newcommand{\un}[1]{\underline{#1}}

\newtheorem{lemma}{Lemma}[chapter]
\newtheorem{proposition}{Proposition}[chapter]
\newtheorem{corollary}{Corollary}[chapter]

\theoremstyle{definition}
\newtheorem{definition}{Definition}[chapter]

\theoremstyle{remark}
\newtheorem{example}{Example}[chapter]

\begin{document}

\null
\pagestyle{empty}
\vspace{3 cm}
\begin{center}
    \Huge\bf Constructions with Bundle Gerbes
    \end{center}
    
\vspace{4 cm}
\begin{center}
  \bf \Large  Stuart Johnson
    \end{center}
    
\vspace{4 cm}

\begin{center}
\bf  Thesis submitted for the degree of \\
 Doctor of Philosophy \\
 in the School of Pure Mathematics \\
 University of Adelaide
    \end{center}
\vspace{1 cm} 
    \begin{center}
 \includegraphics[scale=0.3]{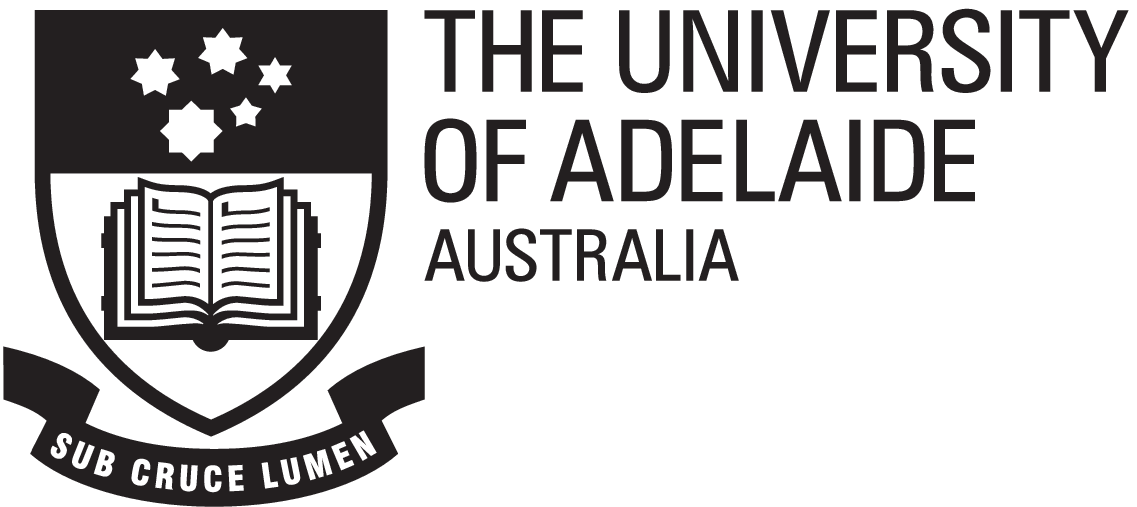}
 \end{center}

 \vspace{1 cm}
 
 \begin{center}
 \bf 13 December 2002
 \end{center}
 
\newpage
\null
\newpage
\pagestyle{plain}
\pagenumbering{roman}

\begin{center} 
   {\bf Abstract} 
\end{center} 

This thesis develops the theory of bundle gerbes and examines
a number of useful constructions in this theory. These allow
us to gain a greater insight into the structure of 
bundle gerbes and related objects. Furthermore they naturally
lead to some interesting applications in physics.

\newpage
\null
\newpage

\begin{center} 
   {\bf Statement of Originality} 
\end{center} 

This thesis contains no material which has been accepted for the 
award of any other degree or diploma at any other university 
or other tertiary institution and, to the best of my 
knowledge and belief, contains no material previously published 
or written by another person, except where due reference 
has been made in the text. 

I give consent to this copy of my thesis, when deposited in 
the University Library, being made available for loan and 
photocopying.  

\vspace{ 3 cm}
\begin{flushright}
    Stuart Johnson\\
    Adelaide, 13 December, 2002
    \end{flushright}
    
\newpage 
\null
\newpage

\begin{center} 
   {\bf Acknowledgement } 
\end{center} 

Firstly I would like to thank my supervisor Michael Murray. He 
has been extremely  helpful, insightful and encouraging at all times
and it
has been a great pleasure to work with him.
I would also like to thank Alan Carey for valuable 
assistance in the early stages of my work on this thesis 
and for his continuing support. 
Thanks are also due to Danny Stevenson for many interesting and helpful
conversations and to Ann Ross for always being of assistance 
with administrative matters.
I would also like to acknowledge the support of a University of Adelaide 
postgraduate scholarship. 

Finally I should acknowledge everyone, past and present, at Kathleen 
Lumley College as well as all of the other wonderful people that I 
have met in Adelaide and who have made my time here so enjoyable.

\newpage
\null

\newpage


\tableofcontents
\newpage

\pagenumbering{arabic}
\include{chapter1}

\include{chapter2}

\include{chapter3}

\include{chapter4}

\include{chapter5}

\include{chapter6}

\include{chapter7}

\include{chapter8}

\appendix

\bib

\end{document}

%% file: chapter1.tex
\chapter{Introduction}

Bundle gerbes were introduced by Murray \cite{mur} as geometric
realisations of classes in $H^3(M,\Z)$ on a manifold $M$. 
By geometric realisation we mean 
an equivalence class of geometric objects which is isomorphic to
$H^3(M,\Z)$. An example of this in lower degree is the relationship
between isomorphism classes of principal $U(1)$-bundles, or
equivalently line bundles, over $M$ and the Chern class in $H^2(M,\Z)$.
The interest in $H^3(M,\Z)$ was motivated by the appearance of
such integral cohomology classes in a number of situations 
including central extensions of structure 
groups of principal bundles and Wess-Zumino-Witten (WZW) theory. 
There were already a number of other realisations of $H^3(M;\Z)$,
of particular interest are the gerbes of Giraud as described in 
Brylinski's book \cite{bry}. These are, from a rather simplistic
view, defined as sheaves of groupoids. The idea of bundle gerbe theory
was to define realisations which did not involve sheaves. Instead
it is possible to build a realisation out of principal $U(1)$- bundles. 

Essentially
a bundle gerbe over a manifold $M$ consists of a submersion 
$Y \stackrel{\pi}{\rightarrow} M$ and a $U(1)$-bundle 
$P \rightarrow Y^{[2]}$ over the fibre product 
\[
Y^{[2]} = Y \times_\pi Y = \{ (y,y') \in Y^2 | \pi(y) = \pi(y') \}
\]
The fibres of $P$ are required to carry a certain associative 
product structure
which is called the bundle gerbe product. 

In \cite{mur} bundle gerbe connections and curving were defined. 
Together these form a higher analogue of connections on $U(1)$-bundles.
They also correspond to connective structures and curvings on
gerbes \cite{bry}. A bundle gerbe connection is a connection, $A$, 
on the bundle $P$ which is compatible with the bundle gerbe 
product. Denote the curvature of this connection by $F$. Let
$\pi_1$ and $\pi_2$ be the projections of each component
$Y^{[2]} \rightarrow Y$.
A curving is a 2-form, $\eta$, on $Y$ which  
satisfies $\delta(\eta) = F$ where $\delta(\eta) = \pi_2^* \eta - \pi_1^* \eta$.
A bundle gerbe with connection and curving defines a class in the 
Deligne cohomology group $H^3(M,\Z(3)_D)$ which may be thought of as the 
hypercohomology of a complex of sheaves, 
$H^2(M,\un{U(1)}_M \rightarrow
\Omega^1(M) \rightarrow \Omega^2(M))$.

The bundle gerbe construction has also been used to consider 
higher degree \v{C}ech and Deligne classes \cite{camuwa}. In particular
Stevenson has developed a theory of bundle 2-gerbes (\cite{ste},\cite{ste3})
which have an associated class in $H^4(M,\Z)$, and when given a higher 
analogue of connection and curving give rise to a class in
$H^4(M,\Z(4)_D)$.

There is a cup product in Deligne cohomology which was described by 
Esnault and Viehweg \cite{esvi} and which has been given a geometric
interpretation by Brylinski and McLaughlin (\cite{bry},\cite{brmc2}).

A generalisation of the concept of holonomy to bundle gerbes was first
considered in \cite{mur}. Just as the holonomy of a $U(1)$-bundle with 
connection associates an element of $U(1)$ to every loop in the base, the
holonomy of a bundle gerbe with connection and curving 
associates an element of $U(1)$ to every closed surface in the base. 
In one particular case bundle gerbe holonomy has been used to describe the 
WZW action 
\cite{camimu}. 

These considerations lead to the transgression formulae derived by
Gawedski \cite{gaw} for dealing with such actions in general settings.
These formulae generalise bundle holonomy and parallel transport to
higher degree Deligne classes. Completely general transgression formulae have
been given by Gomi and Terashima (\cite{gote2},\cite{gote}).

The relevance of \v{C}ech and Deligne cohomology classes to applications in
physics has been well established. For example Dijkgraaf and Witten \cite{diwi}
have used differential characters to find a general Chern-Simons Theory, 
Gomi \cite{gom} has considered the relation 
between gerbes and Chern-Simons theory and 
Freed and Witten \cite{frwi} have considered the role of Deligne classes
in anomaly cancellation in $D$-brane theory. As well as the WZW case we have
already alluded to there are a number of other examples discussed in 
\cite{camimu}.

The basic aim of this thesis is to provide further development of  
the theory of bundle gerbes. 
The goal has been to develop this theory in a way which
keeps in mind the need for a balance between an abstract approach which 
readily accommodates generalisation and an approach which more easily 
allows application of the theory and which could be of interest to a 
wider audience. For the first factor the most important feature is the
bundle gerbe hierarchy principle which is a guiding principle for
relating bundle gerbe type constructions corresponding to Deligne 
cohomology in various degrees. For the second factor we show how
various constructions may be described in geometric terms, often 
allowing manipulation of
diagrammatic representations of bundle gerbes to take the place of 
complicated calculations.

With these factors as a guide we have described a number
of constructions involving bundle gerbes. 
Some of these which have already been developed 
elsewhere
are given in different forms or with a different emphasis 
to demonstrate the hierarchy principle or to relate more easily to
our applications. We also describe some constructions which are new to
bundle gerbe theory. Finally we show how these constructions are 
useful in applications of bundle gerbe theory to physics. 

We begin with a review of the basic features of Deligne cohomology and
introduce the bundle gerbe family of geometric realisations via
bundle 0-gerbes as an alternative to $U(1)$-bundles. This would 
appear to be a complicated approach however it simplifies
the transition from bundles to bundle gerbes and allows us 
to develop some features of bundle gerbe theory in a setting which
is still relatively familiar. We then
define bundle gerbes and explain their role as  
representatives of degree 3 Deligne
cohomology. 

In Chapter 3 we consider some important examples of bundle gerbes. Tautological
bundle gerbes \cite{mur} are introduced by first defining a bundle 0-gerbe,
helping to gain a feel for the bundle gerbe hierarchy. Trivial bundle
gerbes are discussed in some detail since they play an important role in
many constructions. In particular we give a detailed account of
the distinction between trivial bundle gerbes which by definition
have trivial \v{C}ech class and $D$-trivial bundle gerbes which have trivial
Deligne class. We then consider torsion bundle gerbes which are defined
as bundle gerbes with a torsion \v{C}ech class. We describe bundle 
gerbe modules which were introduced in \cite{bcmms} and derive 
corresponding local data. We briefly describe the example of the lifting
bundle gerbe \cite{mur} and then describe bundle gerbes representing
cup products of Deligne classes. 
Each of these examples
is useful in subsequent constructions and applications.

In Chapter 4 the bundle gerbe hierarchy principle is 
introduced via comparison with
some other geometric realisations of Deligne cohomology. 
The correspondence between bundle gerbes and gerbes which has
previously been described in \cite{mur} and \cite{must} is put in the
context of the hierarchy. We then define bundle 2-gerbes 
following \cite{ste}, however we consider the product structures as members
of the hierarchy and define the Deligne class using 
the language of $D$-obstruction forms which we established in Chapter 3.
We prove the isomorphism between bundle 2-gerbes with connection
and curving and $H^4(M,\Z(4)_D)$. Next we go in the other direction and
define $\Z$-bundle gerbes which lie at the bottom of the hierarchy. These
would seem rather trivial however it is of interest to see how the
various aspects of the higher theory appear here, in particular the
$\Z$-bundle gerbe connection naturally motivates classifying theory, our 
next topic for consideration. This is a generalisation of classifying
theory for bundles. We present a number of results of Gajer \cite{Gaj}
relating to Deligne classes and of Murray and Stevenson relating 
specifically to bundle gerbes \cite{must}.
Chapter 4 concludes with a table which catalogues the various realisations
of Deligne cohomology which we have dealt with.

Chapter 5 begins with an account of the holonomy of $U(1)$-bundles
which differs somewhat from the standard treatment of the subject. 
The reason for this is that we need a theory of holonomy which relates
directly to the Deligne class rather than concepts such as horizontal
lifts which are not easily generalised to bundle gerbes and beyond.
We then define the holonomy of bundle gerbes with an explanation of 
how it relates to
the holonomy of $U(1)$-bundles and with details of how local
formulae are obtained. The concept of holonomy is extended to
bundle 2-gerbes and to general Deligne classes. 

In Chapter 6
we describe the extension of the notion of parallel transport from
$U(1)$-bundles to bundle gerbes, bundle 2-gerbes and general Deligne classes
and provide detailed derivations of local formulae. In particular we 
discuss how to obtain a $U(1)$-bundle on the loop space $LM$ from a
bundle gerbe on $M$ and bundle 2-gerbe generalisations of this 
construction.

The basic properties of bundle gerbe holonomy are described in
section 7.1. The
motivation for considering these particular 
properties comes from those of line bundles
obtained via transgression (\cite{bry},\cite{fre}). We consider 
the example of the tautological bundle gerbe which motivates holonomy 
reconstruction, that is, reconstructing a bundle gerbe with connection
and curving from its holonomy on closed surfaces. 
We show how the holonomy reconstruction for bundle gerbes relates 
to that of gerbes as described in \cite{mapi}. 
Transgression formulae are then used to give an alternative 
approach to holonomy reconstruction which allows us to
consider the case of bundle 2-gerbes.  We conclude the chapter with the 
gauge invariance properties of holonomy.

Finally in Chapter 8 we use bundle gerbe theory to examine applications in
physics.
Constructions in Wess-Zumino-Witten and Chern-Simons (CS) theories are
shown to follow naturally from various constructions in bundle gerbe theory.
In the WZW case we just interpret the standard results (see the Appendix A 
of \cite{fre}) in terms of bundle gerbes. In the case of Chern-Simons theory
we define a bundle 2-gerbe whose holonomy gives a general definition
of the Chern-Simons action. We can then interpret Chern-Simons lines,
gauge invariance and the relationship with the central extension of the
loop group in terms of bundle gerbe theory. All of these results are
quite straightforward from this point of view. 
Bundle gerbes also prove useful in studying anomaly cancellation in 
D-branes as described in \cite{cajomu}. Here we emphasise local aspects
which were not discussed in detail in that paper. We also add some comments on 
the potential for the application of bundle gerbes to the problem of 
anomalies involving $C$-fields and higher dimensional generalisations
of Chern-Simons theory. We conclude with comments on the relationship between 
bundle gerbes and the axiomatic approach to topological quantum field
theory.

It is necessary here for a brief comment on terminology. We shall
refer to line bundles and their associated principal bundles 
interchangeably. Also we shall usually work in the Hermitian setting, so
our bundles are $U(1)$-bundles. Since we deal almost exclusively with
these bundles we shall often simply refer to bundles, in situations
where we require a principal bundle with a more general structure group 
we refer to it as a $G$-bundle.

%% file: chapter2.tex
\chapter{Bundle Gerbes and Deligne Cohomology}

In this chapter we discuss a number of geometric realisations of
low degree Deligne cohomology, in particular bundle gerbes. 

\begin{section}{A Review of Sheaf Cohomology}
We begin with some background material regarding sheaf cohomology.
We mostly follow Brylinski's book \cite{bry}, though some material
is drawn from Bott and Tu \cite{botu}. 

We assume that the reader is familiar with the definitions of sheaves and
related concepts such as morphisms of sheaves. We provide the minimum
amount of detail necessary to define sheaf cohomology, for further details
see \cite{bry}.
Let $A$ be a sheaf of Abelian groups on a manifold $M$. 
Recall that this means that associated
with every open $U \subseteq M$ there is an Abelian group $A(U)$ which 
satisfies certain axioms with respect to restrictions. We always 
assume that manifolds are paracompact, that is, every open cover has
a locally finite subcover. 
We shall be interested in the following
examples:\\

\begin{tabular}{cl}
$\Z_M, \R_M, U(1)_M$: & sheaf of locally constant functions on a manifold $M$\\
$\un{\R}_M , \un{U(1)}_M$:  & sheaf of smooth $\R$ or $U(1)$-
valued functions
on a manifold $M$ \\
$\Omega^p_M$: & sheaf of real differential $p$-forms on $M$\\
&
\end{tabular}
\\
A {\it complex of sheaves} $K^\bullet$ is a sequence
\[
\cdots \stackrel{d^{n-1}}{\rightarrow} K^n \stackrel{d^n}{\rightarrow}
K^{n+1} \stackrel{d^{n+1}}{\rightarrow} \cdots
\]
where $n \in \Z$ and 
$d^n: K^n \rightarrow K^{n+1}$ are morphisms of sheaves of Abelian
groups satisfying $d^n \circ d^{n-1} = 0$. The map $d^n$ is called the
{\it differential} of the complex. We always assume that $K^p=0$ for
$p < 0$. A {\it morphism of complexes of sheaves} $\phi: K^\bullet
\rightarrow L^\bullet$ consists of a family of morphisms 
of sheaves $\phi^n : K^n \rightarrow L^n$ such that
$\phi^{n+1} \circ d^n_K = d^n_L \circ \phi^{n}$. Given two morphisms
$\phi$ and $\psi$ from $K^\bullet$ to $L^\bullet$ a {\it homotopy} H
from $\phi$ to $\psi$ consists of a series of morphisms
$H^n : K^n \rightarrow L^{n-1}$ such that $d^{n-}_L H^n + H^{n+1}
d^n_K = \phi^n - \psi^n$. A morphism of complexes of sheaves 
$\phi: K^\bullet \rightarrow L^\bullet$ is a {\it homotopy 
equivalence} if there exists a morphism $\psi : L^\bullet
\rightarrow K^\bullet$ and $\phi \psi$ and $\psi \phi$ are
both homotopic to the identity map.
A complex of sheaves is 
called {\it acyclic} if $Ker (d^n) = Im (d^{n-1})$ for all $n$. 
The {\it cohomology sheaves} $\un{H}^p(K^\bullet)$ are defined by the
presheaf $Ker(d^j) / Im(d^{j-1})$. For an acyclic complex all of the
cohomology sheaves are zero.
A sheaf
$I$ is called {\it injective} if given any morphism $f:A \rightarrow I$
and an injective morphism $i: A \rightarrow B$ there exists a morphism
$g:B \rightarrow I$ such that $g \circ i = f$. An {\it injective
resolution} of $A$ is a complex of injective sheaves $I^\bullet$
such that
$A \rightarrow I^\bullet$ is an acyclic complex of sheaves.
Injective
resolutions always exist and are unique up to homotopy equivalence.  
Let $\Gamma(M,.)$ be the functor of global sections which takes 
the sheaf $A$ to the Abelian group $\Gamma(M;A)$ defined by $A(M)$.
The {\it sheaf cohomology groups} $H^p(M,A)$ are defined as the 
$p$-th cohomology of the complex  
\[
\cdots \rightarrow \Gamma(M,I^j) \rightarrow \Gamma(M,I^{j+1}) \rightarrow
\cdots
\]
where $I^\bullet$ is an injective resolution of $A$.
Given a short exact sequence of sheaves
\[
0 \rightarrow A \rightarrow B \rightarrow C \rightarrow 0
\]
there is a long exact sequence in sheaf cohomology
\[
\cdots \rightarrow H^n(M,A) \rightarrow H^n(M,B) \rightarrow 
H^n(M,C) \rightarrow H^{n+1}(M,A) \rightarrow \cdots
\]
If
a morphism of complexes induces an isomorphism of
cohomology sheaves $\un{H}^n(K^\bullet) \isom
\un{H}^n(L^\bullet)$ then it is called a 
{\it quasi-isomorphism}. 

A useful example of a resolution is the {\it \v{C}ech 
resolution}. Let $\U$ be an open cover of $M$ and let
$U_{i_0,\ldots,i_p}$ denote an intersection 
$U_{\i_0} \cap \ldots \cap U_{i_p}$ of open sets in this cover.
The \v{C}ech resolution is a complex $\check{C}^\bullet(\U,A)$ which
is defined by
$\check{C}^p(\U,A) = \Pi_{i_0,\ldots,i_p} A(U_{i_0,\ldots,i_p})$.
and $\delta: \check{C}^p(\U,A) \rightarrow \check{C}^{p+1}(\U,A)$ is defined by
\[
\delta(\un{\alpha})_{i_0,\ldots,i_{p+1}} =
\sum_{j=0}^{p+1} (-1)^j (\alpha_{i_0,\dots,i_{j-1},i_{j+1},
\ldots,i_{p+1}})_{|U_{i_0,\ldots,i_{p+1}}}
\]
where we have introduced the notation $\un{\alpha} \in
\check{C}^p(\U,A)$ where the underline denotes a family $\alpha_{i_0,
\ldots,i_p} \in A(U_{i_0,\ldots,i_p})$. In general the resolution
and the resulting cohomology groups should depend on the choice of
open cover however
we shall be interested only in spaces which are 
manifolds and these always admit a good cover, that is, a cover in which all 
non-empty intersections are contractible. In this case the construction is
independent of the choice of cover and the \v{C}ech resolution computes
the sheaf cohomology. It is not an injective resolution, however there exists
a morphism with the \v{C}ech resolution which induces an isomorphism 
in cohomology. If $A$ is the sheaf $\Z$ then we recover the 
usual \v{C}ech cohomology groups $\check{H}^p(M,\Z)$, for example a class
$\check{H}^1(M,\Z)$ consists of a family of $\Z$ valued constants
$c_{ij}$ associated with double intersections $U_{ij}$ such that
$c_{jk} - c_{ik} + c_{ij} = 0$ on $U_{ijk}$ and under the equivalence
relation $c_{ij} \sim c_{ij} + b_j - b_i$ where $b_i$ and $b_j$ are from
a family of constants defined on single open sets.


A {\it double complex} of sheaves, $K^{\bullet\bullet}$, consists of 
sheaves $K^{p,q}$ and two differential maps $d: K^{p,q} \rightarrow
K^{p+1,q}$ and $\delta: K^{p,q} \rightarrow K^{p,q+1}$ such that
$dd=0$, $\delta \delta = 0$ and $d\delta = \delta d$. A double
complex may be represented diagrammatically as follows:
\[
\begin{CD}
& &  @AAdA @AAdA \\
@>\delta>> K^{p,q} @>\delta>> K^{p+1,q} @>\delta>> \\
& &  @AAdA @AAdA \\
@>\delta>> K^{p,q-1} @>\delta>> K^{p+1,q-1} @>\delta>> \\
& &  @AAdA @AAdA
\end{CD}
\]
Each row and column defines a complex of sheaves denoted $K^{\bullet,q}$
and $K^{p,\bullet}$ respectively. The {\it total complex} 
of a double complex $K^{\bullet\bullet}$ is an ordinary complex $K^\bullet$
which is defined by $K^n = \bigoplus_{p+q = n} K^{p,q}$ with differential
$D = \delta + (-1)^p d$.

Let $K^\bullet$ be a complex of sheaves with differential $d_K$. An
{\it injective resolution} of $K^\bullet$ is a double complex 
$I^{\bullet\bullet}$ with differentials $d$ and $\delta$ such that
for each $q$ the complex $I^{\bullet,q}$ with differential $\delta$ is an 
injective resolution of $K^q$, the complex $d(I^{\bullet,q-1})
\subseteq I^{\bullet,q}$ is an injective resolution of $d_K(K^{q-1})$,
the complex of sheaves $Ker(d) \subseteq I^{\bullet,q}$ is an
injective resolution of $Ker(d_K:K^{q-1} \rightarrow K^q)$ and
complex of cohomology sheaves of the rows, $\un{H}^{\bullet,q}$ is  
an injective resolution of $\un{H}^q(K^\bullet)$.

For all of our examples an injective resolution of a complex of sheaves 
exists and is unique up to homotopy. Given a complex $K^\bullet$ and 
an injective resolution $I^{\bullet\bullet}$ the {\it hypercohomology
group} $H^p(M,K^\bullet)$ is defined to be the $p$-th cohomology
of the double complex $\Gamma(M,I^{\bullet\bullet})$. Given 
a short exact sequence of complexes of sheaves there is 
a long exact sequence in hypercohomology. Quasi-isomorphisms $\phi:
K^\bullet \rightarrow L^\bullet$ induce 
isomorphisms in sheaf hypercohomology, $H^n(M,K^\bullet)
\isom H^n(M,L^\bullet)$. The \v{C}ech resolution 
may be extended to the case of a complex of sheaves by taking
the usual \v{C}ech resolution for each sheaf in the complex.

We shall describe a specific example of this.
Let $K^\bullet = \un{U(1)}_M \stackrel{d\log}{\rightarrow} \Omega^1_M$. 
The \v{C}ech
resolution looks like 
\[
\begin{CD}
 0 & & 0 & & 0 \\
@AAA @AAA @AAA \\
\Omega^1_M @>>> \check{C}^0(\U,\Omega^1) @>\delta>> \check{C}^1(\U,\Omega^1) 
@>\delta>> \\
@AAd\log A @AAd\log A @AAd\log A  \\
U(1)_M @>>> \check{C}^0(\U,\un{U(1)}) @>\delta>> \check{C}^1(\U,\un{U(1)})
@>\delta>>
\end{CD}
\]
A class in $H^0(M,\un{U(1)} \rightarrow \Omega^1)$ consists 
of  $\un{f} \in \check{C}^0(\U,\un{U(1)})$ such that
$f_\beta f^{-1}_\alpha = 1$ and $d\log f_\alpha = 0$. This is 
a locally constant $U(1)$-valued function on $M$.

A class in $H^1(M,\un{U(1)} \rightarrow \Omega^1)$ consists of a pair
$(\un{g},\un{A}) \in \check{C}^1(\U,\un{U(1)}) \oplus \check{C}^0(\U,\Omega^1)$
such that $g_{\beta \gamma} g^{-1}_{\alpha\gamma} g_{\alpha \beta} =1$ and
$d\log g_{\alpha \beta} = A_\beta - A_\alpha$ and is defined modulo
exact cocycles of the form $(h^{-1}_\alpha h_\beta, d\log h_\alpha)$ for
some $\un{h} \in \check{C}^0(\U,\un{U(1)})$. Classes of higher degree are 
defined in a similar way.

\end{section}
\begin{section}{$U(1)$-Functions}
\label{u1f}
We examine $U(1)$-valued functions as the starting point for our 
geometric objects corresponding to Deligne cohomology classes, focusing 
in particular on features which are of interest when we move on to 
geometric realisations of higher degree classes.

Let $M$ be a smooth manifold.
We shall consider smooth functions $f : M \rightarrow U(1)$. 
Such functions are elements of the sheaf cohomology group 
$H^0(M,\underline{U(1)}_M)$. 

Our interest in $U(1)$-valued functions is due to their role as representatives
of the smooth Deligne cohomology group $H^1(M,\Z(1)_D)$. 
\begin{definition} \label{delcoh} \cite{bry}
Let $\Z(p) = (2\pi \sqrt{-1})^p \cdot \Z$. Define a complex
of sheaves $\Z(p)_D$, for $p > 0$ by
\[
\Z(p)_M \stackrel{i}{\rightarrow} \ua^0_M 
\stackrel{d}{\rightarrow} \ua^1_M 
\stackrel{d}{\rightarrow} \cdots
\stackrel{d}{\rightarrow} \ua^{p-1}_M
\]
Where $\ua^k_M$ is the sheaf of real differential $k$-forms on $M$ and
the map $i$ is the inclusion. We define $\Z(0)_D$ to be $\Z_M$.
The Deligne cohomology groups of $M$ are defined as the hypercohomology
groups $H^q(M,\Z(p)_D)$. 
\end{definition}



Deligne classes are realised explicitly by first using a
quasi-isomorphism of sheaves \cite{bry}
\footnote{This quasi-isomorphism is derived from the 
exact sequence
\[
0 \rightarrow \uz_M \rightarrow \ur_M \rightarrow \uu_M \rightarrow 0
\] 
which may be replaced by 
\[
0 \rightarrow \uz_M \rightarrow \un{\C}_M \rightarrow 
\un{\C}^\times_M \rightarrow 0
\] 
giving an equivalent theory in terms of $\C^\times$ rather than 
$U(1)$.},
\begin{equation}\label{quasi}
\begin{array}{ccccccccccc}
0 & \rightarrow & \Z(p) &  
\rightarrow & \ur_M & \stackrel{d}{\rightarrow} &  
\ua^1_M & \rightarrow & \cdots & \rightarrow & \ua^{p-1} \\
& & & & \downarrow & & \downarrow & & & & \downarrow \\
& & & & \uu_M & \stackrel{d\log}{\rightarrow} &
\ua^1_M & \rightarrow & \cdots & \rightarrow & \ua^{p-1}
\end{array}
\end{equation}
which induces an isomorphism 
\begin{equation}\label{indisom}
H^q(M,\Z(p)_D) \isom H^{q-1}(M,
\uu_M  \stackrel{d\log}{\rightarrow} 
\ua^1_M  \rightarrow  \cdots  \rightarrow  \ua^{p-1})
\end{equation}
In general we shall denote the complex 
\[
\uu_M  \stackrel{d\log}{\rightarrow} 
\ua^1_M  \rightarrow  \cdots  \rightarrow  \ua^{p-1}
\]
by $\D^{p-1}$ so we can write \eqref{indisom} as
\[
H^q(M,\Z(p)_D) \isom H^{q-1}(M,\D^{p-1})
\]
We shall usually deal directly the groups $H^{q-1}(M,\D^{p-1})$ so we shall
also refer to these as Deligne cohomology. The original definition is still
required for certain purposes such as cup products.
To get concrete expressions for these sheaf
cohomology classes they are represented in terms of the \v{C}ech resolution
relative to a good open cover
as discussed in the example at the end of the previous section.

There are exact sequences \cite{Gaj}
\begin{equation} \label{gaj}
0 \rightarrow H^{p-1}(M,U(1))  \rightarrow 
H^p(M,\Z(p)_D)  \rightarrow \Omega^p_0(M) 
\rightarrow 0
\end{equation}
which for $p=1$ becomes
\[
0 \rightarrow H^0(M,U(1)) \stackrel{i}{\rightarrow} 
H^1(M,\Z(1)_D) \stackrel{d\log}{\rightarrow} \Omega^1_0(M) 
\rightarrow 0
\]
where $\Omega^p_0(M)$ denotes the set of closed 
$p$-forms on $M$ which have periods in $\Z(1)$ 
(that is, the integral
over a closed $p$-cycle is in $\Z(1)$)
. We shall refer to forms satisfying this integrality 
requirement as {\it $2\pi$-integral}.
Given $\{f_\alpha \} \in
H^0(M,\uu) \isom H^1(M,\Z(1)_D)$ we have $d\log f \in \Omega^1(M)_0$.
This is the globally defined since 
$d(\log_\alpha f - \log_\beta f) = 0$ and may be thought of as a 
lower dimensional version of the curvature of a connection.
If $\{ f_i \} \in \ker (d\log)$ then there are $U(1)$-valued
constants $c_i = f_i$ which are classes in $H^0(M,U(1))$.
\end{section}

\begin{section}{$U(1)$-Bundles}
We present basic material on the relationship between principal 
$U(1)$-bundles and Deligne cohomology to further develop the 
theory of geometric realisations of Deligne classes.

We follow the detailed treatment of the role of line bundles
as geometric realisations of degree 2 Deligne cohomology
in Brylinski's book \cite{bry}. Let $P$ denote a principal
$U(1)$ bundle over $M$. It is well known that the 
isomorphism classes of $U(1)$ bundles corresponds to the 
sheaf cohomology group $H^1(M,\uu_M)$. A representative,
$g_{\alpha \beta}$, of $H^1(M,\uu_M)$ corresponds to
the transition functions of the bundle. T
here is an isomorphism with \v{C}ech
cohomology  $H^1(M,\uu_M) \isom H^2(M,\Z)$. The 
image of $g_{\alpha \beta}$ under this isomorphism is
the Chern class, $n_{\alpha \beta \gamma} = 
-\log (g_{\beta \gamma}) + \log (g_{\alpha
\gamma}) - \log (g_{\alpha \beta})$. Note that
there is also an isomorphism with Deligne
cohomology $H^1(M,\uu_M) \isom 
H^2(M,\Z(1)_D)$, where the 
Deligne class corresponding to $g_{\alpha \beta}$ is
$(n_{\alpha \beta \gamma}, 
\log (g_{\alpha \beta}))$. 

It is also well known that isomorphism classes of
bundles with connection lie in the
hypercohomology group $H^1(M,\uu_M \rightarrow \ua^1_M) \equiv
H^1(M,\D^1)$. We use
the \v{C}ech resolution of the complex to produce explicit representatives
of these hypercohomology classes.
If $(g_{\alpha \beta},A_\alpha)$ is a class in 
$H^1(M,\D^1)$, then it represents a 
$U(1)$-bundle with transition functions $g_{\alpha \beta}$ and
local connection 1-forms $A_\alpha$.

The space of bundles with connection is related to the 
space of bundles via the exact sequence \cite{Gaj}
\begin{equation}\label{conseq1}
0 \rightarrow \Omega^1(M) / \Omega^1(M)_0 \rightarrow
H^2(M,\Z(2)_D) \rightarrow H^2(M,\Z) \rightarrow 0
\end{equation}
The quasi-isomorphism
of complexes of sheaves \eqref{quasi}
induces the usual isomorphism 
\[
H^1(M,\D^1)
\isom H^2(M,\Z(2)_D)
\]

Substituting $p=2$ into the exact sequence \eqref{gaj} gives the
exact sequence
\[
0 \rightarrow H^1(M,U(1)) \stackrel{i}{\rightarrow} 
H^1(M,\D^1) \stackrel{d}{\rightarrow} \Omega^2(M)_0 
\rightarrow 0
\]
where $d$ is the map which applies $d$ to the component
of $H^1(M,\D^1)$ with the highest $d$-degree. Geometrically,
$d$ maps a bundle with connection to its curvature 2-form.
This implies that $H^1(M,U(1))$ represents the set of flat
bundles on $M$. This can be seen explicitly in the following way
\cite{hit}, let $(g_{\alpha \beta}, A_\alpha)$ represent a flat 
bundle. Thus we have $d A_\alpha = 0$. Each element of a good cover is 
contractible so Poincare's Lemma applies and there exist $U(1)$-valued
functions 
$a_\alpha$ such that $d\log a_\alpha = A_\alpha$. Now we have
\begin{eqnarray*}
d\log g_{\alpha \beta} &=& A_\alpha - A_\beta \\
&=& d\log a_\alpha - d\log a_\beta \\
d\log (g_{\alpha \beta}\cdot a_\alpha^{-1}\cdot a_\beta) &=& 0 
\end{eqnarray*}
so we have constants $c_{\alpha \beta} = 
g_{\alpha \beta} \cdot a_\alpha^{-1} \cdot  a_\beta$ which
represent a cocycle in $H^1(M,U(1))$.
We shall refer to the cocycle $c_{\alpha \beta}$ as the flat holonomy of the
bundle represented by $(g_{\alpha \beta}, A_\alpha)$.  

The space of flat bundles may also be represented by the Deligne
cohomology groups $H^1(M,\D^p)$ for $p > 1$. To see why this is so,
consider what happens when the Deligne differential, $D$, is applied 
to a class $(g_{\alpha \beta},A_\alpha) \in H^1(M,\D^1)$.
\[
D(g_{\alpha \beta}, A_\alpha) = (\delta(g)_{\alpha \beta \gamma},
d\log(g_{\alpha \beta}) + \delta(A)_{\alpha \beta})
\]
This leads to the usual requirements for $(g_{\alpha \beta},A_\alpha)$
to represent a bundle with connection. If we truncate the Deligne 
complex at a higher value of $p$ then the third component
of $D(g_{\alpha \beta}, A_\alpha)$ will be $d A_\alpha$.
This means that a Deligne cycle will represent a flat bundle.
\end{section}

\begin{section}{Bundle 0-Gerbes}
These were introduced by Murray \cite{mur2}. The objects described here
should actually be called $U(1)$-bundle 0-gerbes, however since we only
use this type we omit the $U(1)$ prefix. Initially it may seem that 
bundle 0-gerbes are just a more complicated way of looking at line bundles,
certainly if one was interested only in line bundles then there would be little
point in studying them. Our motivation is that we are working towards 
bundle gerbes and bundle 2-gerbes. In this situation there are two 
advantages to considering bundle 0-gerbes. Several
properties of these higher objects also appear in the bundle 0-gerbe case
so it is useful to become familiar with them in a simpler setting. Secondly
we are interested in viewing all of these objects as part of a hierarchy and
it will become clear that the lower dimensional geometric realisation 
of Deligne cohomology should be a bundle 0-gerbe rather than a line bundle.
In this way bundle 0-gerbes will be useful in gaining an understanding of this 
hierarchy. 

\begin{definition} \label{b0g}
Let $Y \rightarrow M$ be a submersion. 
Let $Y^{[2]}$ denote the fibre product
\[
Y^{[2]} = Y \times_\pi Y = \{ (y,y') \in Y^2 | \pi(y) = \pi(y') \}
\]
and let 
\[
g : Y^{[2]} \rightarrow U(1)
\]
be a $U(1)$-function satisfying the cocycle identity
\[
g(y_1,y_2)g(y_2,y_3) = g(y_1,y_3).
\]
The triple $(g,Y,M)$ defines a {\it $U(1)$ bundle 0-gerbe}. 
\end{definition}
Note that the cocycle identity implies that $g(y,y)=1$ and
$g(y_1,y_2) = g^{-1}(y_2,y_1)$.

Recall that a submersion is an onto map with onto differential.
It admits local sections and all fibrations are submersions, however there
exist submersions which are not fibrations.

Bundle 0-gerbes may be represented diagrammatically in the following
way:
\[
\begin{array}{ccc}
& & U(1) \\
& \stackrel{g}{\nearrow} & \\
Y^{[2]} & \rightrightarrows & Y \\
& & \downarrow \\
& & M 
\end{array}
\]
A bundle 0-gerbe is called {\it trivial} if there exists a 
$U(1)$-function $h$ on $Y$ satisfying
\[
g(y_1,y_2) = {h(y_1)}^{-1}h(y_2).
\]
In this case we write $g = \delta(h)$.
The dual of a bundle 0-gerbe $(g,Y,M)$ is defined as $(g^{-1},Y,M)$.
Given bundle 0-gerbes $(g,Y,M)$ and $(g',Y',M)$ we can take
the product 
\[
(g,Y,M) \otimes (g',Y',M) = (g\cdot g',Y \times_\pi Y',M)
\]
which is easily verified to be a bundle 0-gerbe. 

A {\it bundle 0-gerbe morphism} is a smooth map $\phi: Y \rightarrow Y'$
such that $\pi' \circ \phi = \pi$ and $g =  g' \circ
\phi^{[2]}$ where $\phi^{[2]} : Y^{[2]} \rightarrow 
{Y'}^{[2]}$ is induced by the fibre product.

We say that two bundle 0-gerbes $(g,Y,M)$ and $(g',Y',M)$ are 
stably isomorphic if there exists a trivial bundle 0-gerbe $(\delta(h),X,M)$
and a bundle 0-gerbe morphism
\[
(g,Y,M) \isom (g',Y',M) \otimes (\delta(h),X,M).
\]

Since $(g,Y,M)\otimes (g^{-1},Y,M)$ is canonically trivial then this
condition is equivalent to requiring that $(g,Y,M) \otimes ({g'}^{-1},Y',M)$
is trivial.

An example of a stable isomorphism may be defined in the following
way.
Let $(g,Y,M)$ and $(g',Y',M)$ be two bundle 0-gerbes and
suppose there exists $\phi:Y' \rightarrow Y$ such that 
$\pi_{Y'} = \pi_Y \circ \phi$ and $g' = g \circ \phi^{[2]}$.
Then $(g,Y,M)$ and $(g',Y',M)$ are stably isomorphic.

To see this
consider the product bundle 0-gerbe
\[
\begin{array}{ccc}
&  & U(1) \\
& \stackrel{{g'}^{-1}g}{\nearrow}  & \\
(Y' \times Y)^{[2]} & \rightrightarrows & Y' \times Y \\
& & \downarrow \\
& & M
\end{array}
\]
The function ${g'}^{-1}g : (Y' \times Y)^{[2]}
\rightarrow U(1)$ is defined by
\begin{eqnarray*}
{g'}^{-1}g (y'_1,y'_2,y_1,y_2) &=&
{g'}^{-1}(y'_1,y'_2) g (y_1,y_2) \\
&=& g^{-1}(\phi(y'_1),\phi(y'_2))g(y_1,y_2) \\
&=& g^{-1}(\phi(y'_1),y)g^{-1}(y,\phi(y'_2))
g(y_1,y)g(y,y_2) \\
&=& g^{-1}(\phi(y'_1),y_1)g(\phi(y'_2),y_2) \\
&=& \delta (g(\phi(y'),y))
\end{eqnarray*}
and hence the two bundle 0-gerbes are stably isomorphic.

\begin{lemma}
The set of stable isomorphism classes of bundle 0-gerbes 
forms a group.
\end{lemma} 
The associativity of the product is clear. The identity
element is the equivalence class of trivial bundle 0-gerbes and 
the inverse of $(g,Y,M)$ is $(g^{-1},Y,M)$.
\begin{proposition} \label{p1}
The group of stable isomorphism classes of bundle 0-gerbes over 
$M$ is isomorphic to $H^1(M,\uu)$.
\end{proposition}
\begin{proof}
Let $(g,Y,M)$ be a bundle 0-gerbe. Define a \v{C}ech cycle on $M$ by
\[
g_{\alpha \beta}(x) = g(s_\alpha,s_\beta)
\]
where $s_\alpha$ and $s_\beta$ are local sections on $M$. Independence
of the choice of cover follows from the standard argument in the
case of the Chern class of a bundle, for example see Theorem 2.1.3 of
Brylinski \cite{bry}. 

Suppose we choose different sections and define
\[
g'_{\alpha \beta}(m) = g(s'_\alpha(m),s'_\beta(m)).
\]
Using the cocycle identity for $g$,
\begin{eqnarray*}
g_{\alpha \beta}(m) &=& g(s_\alpha(m),s_\beta(m)) \\
&=& g(s_\alpha(m),s'_\alpha(m))g(s'_\alpha(m),s'_\beta(m))
g(s'_\beta(m),s_\beta(m)) \\
&=&  g'_{\alpha \beta}(m) \delta(h(m))_{\alpha \beta}
\end{eqnarray*}
where $h_\alpha(m) = g(s'_\alpha(m),s_\alpha(m))$. Therefore
$g_{\alpha \beta}$ is a well defined \v{C}ech cochain on $M$.
Furthermore the cocycle identity on $g$ makes $g_{\alpha \beta}$
a \v{C}ech cocycle. 

Consider the cocycle corresponding to the product bundle 0-gerbe
$g \otimes g'$. Clearly
\[
(g \otimes g')_{\alpha \beta} (m) =
g_{\alpha \beta} (m) g'_{\alpha \beta}(m)
\] 
Thus we have a homomorphism from the set of bundle 0-gerbes to 
$Z^1(M,\uu)$.
Suppose $(g,Y,M)$ is a trivial bundle 0-gerbe with trivialisation $h :
Y \rightarrow U(1)$. In this case
\begin{eqnarray*}
g_{\alpha \beta}(m) &=& h^{-1}(s_\alpha (m))h(s_\beta (m)) \\
&=& \delta (h(m))_{\alpha \beta}
\end{eqnarray*}
where $h_\alpha(m) = h(s_\alpha (m))$. 
This ensures that stable equivalence classes map to \v{C}ech cohomology 
classes and hence we have a homomorphism from the 
set of stable equivalence classes of bundle 0-gerbes to 
$H^1(M,\uu)$. 

Suppose that for a bundle 0-gerbe $(g,Y,M)$, $g_{\alpha \beta}$ is
trivial. Then there is a bundle 0-gerbe trivialisation given by
\[
h(y) = g(s_\alpha(\pi(y)),y)h_\alpha (\pi(y)).
\]
It is easily verified that this is independent of $\alpha$. 
This proves injectivity of the homomorphism. To prove 
surjectivity we construct a bundle 0-gerbe corresponding
to a class in $H^2(M,\Z)$ by following the method of
theorem 2.1.3 of \cite{bry}. 

Let $g_{\alpha \beta} \in H^1(M,\uu)$. We define a bundle 0-gerbe
$(g,Y,M)$ where $Y = \coprod_{\alpha \in A} U_\alpha$
and $g(y_1,y_2) = g_{\alpha \beta}(\pi(y_1))$ where
$y_1 \in U_\alpha \subset Y$ and $y_2 \in U_\beta \subset Y$.
Since $g(s_\alpha(m),s_\beta(m)) = g_{\alpha \beta}(m)$ this 
construction proves surjectivity.
\end{proof}
\begin{corollary}
The group of stable isomorphism classes of bundle 0-gerbes on $M$ is isomorphic
to $H^2(M,\Z)$.
\end{corollary}
\begin{corollary}
The group of stable isomorphism classes of bundle 0-gerbes on $M$ is 
isomorphic to the group of isomorphism classes of bundles on $M$.
\end{corollary}
The $U(1)$ bundle corresponding to a bundle 0-gerbe is defined by
letting the total space be $Y \times S^1$ with the equivalence relation
\[
(y_1,g(y_1,y_2)) \sim (y_2,1).
\]
Conversely, given a bundle $(P,M)$ the corresponding bundle 0-gerbe is
$(g,P,M)$ where $g$ is defined by 
$p_1 g(p_1,p_2) = p_2$. 

We can describe an explicit correspondence between bundle isomorphisms
and bundle 0-gerbe stable isomorphisms.
Let $(g,Y,M)$ and $(h,X,M)$ be two bundle 0-gerbes and suppose
$\phi: Y \rightarrow X$ is a stable isomorphism. We claim that
$\tilde{\phi}: Y \times S^1 / \sim \rightarrow X \times S^1 / \sim$ defined
by $\tilde{\phi}([y,\theta]) = [\phi(y),\theta]$ is an isomorphism
of the corresponding bundles. First we check that this map is well defined on 
equivalence classes. Consider 
\begin{equation*}
\begin{split}
\tilde{\phi}([y_1,g(y_1,y_2)]) &= [\phi(y_1),g(y_1,y_2)] \\
&= [\phi(y_1),h(\phi(y_1),\phi(y_2))] \\
&= [\phi(y_2),1] \\
&= \tilde{\phi}([y_2,1])
\end{split}
\end{equation*}
Clearly the $S^1$ action is preserved by this map. Now suppose that
$\tilde{\phi}([y_1,\theta_1]) = \tilde{\phi}([y_2,\theta_2])$. Then
\begin{equation*}
\begin{split}
[\phi(y_1),\theta_1] &= [\phi(y_2),\theta_2] \\
&= [\phi(y_1),h(\phi(y_1),\phi(y_2))\theta_2] \\
&= [\phi(y_1),g(y_1,y_2)\theta_2]
\end{split}
\end{equation*}
so $\theta_1 = g(y_1,y_2)\theta_2$ and
thus $[y_1,\theta_1] = [y_1, g(y_1,y_2)\theta_2] = [y_2,\theta_2]$.
We define an inverse of $\phi$ by
\begin{equation}
\tilde{\phi}^{-1}([x,\theta]) = [y,h(\phi(y)x)\theta]
\end{equation}
where $y$ is any element of $Y$. It is independent of this 
choice since given $y_1,y_2 \in Y$ we have
\begin{equation}
\begin{split}
[y_2,h(\phi(y_2),x)\theta] &= [y_1,g(y_1,y_2)h(\phi(y_2),x)\theta] \\
&= [y_1,h(\phi(y_1),\phi(y_2))h(\phi(y_2),x)\theta] \\
&= [y_1,h(\phi(y_1),x)\theta]
\end{split}
\end{equation}
Thus the map $\tilde{\phi}$ is a bundle isomorphism. It is easy to check that
a bundle isomorphism defines a stable isomorphism on the associated bundle
0-gerbe.

Let $\pi_1$ and $\pi_2$ be the maps from $Y^{[2]}$ to $Y$ defined by
\begin{eqnarray*}
\pi_1(y_1,y_2) &=& y_2 \\
\pi_2(y_1,y_2) &=& y_1 
\end{eqnarray*}
This notation may appear counter intuitive. The idea is that the subscript on 
$\pi$ indicates which element will be omitted. This allows the maps
$\pi_i$ to be generalised to $p_i : Y^{[p]} \rightarrow Y^{[p-1]}$ for
$i = 1 \ldots p$.
 
Let $\delta$ be the pull back $\pi_1^* - \pi_2^*$.
\begin{definition}
A {\it bundle 0-gerbe connection}, $A$, is a 1-form on
$Y$ satisfying 
\[
\delta(A) = g^{-1}dg.
\]
\end{definition}
A bundle 0-gerbe $(g,Y,M)$ with connection $A$ may be 
written as $(g,Y,M;A)$ or $(g;A)$ where there is no 
ambiguity.

The existence of bundle 0-gerbe connections is established by
considering the following complex which has no cohomology
\cite{mur}
\begin{equation} \label{nocoh}
\Omega^q(M) \stackrel{\pi^*}{\rightarrow} \Omega^q(Y) \stackrel{\delta}{\rightarrow}
\cdots \stackrel{\delta}{\rightarrow} \Omega^q(Y^{[p-1]})
\stackrel{\delta}{\rightarrow} \Omega^q(Y^{[p]})
\stackrel{\delta}{\rightarrow} \Omega^q(Y^{[p+1]})
\stackrel{\delta}{\rightarrow} \cdots
\end{equation}
The general $\delta$ map from $\Omega^q(Y^{[p]})$ to $\Omega^q(Y^{[p+1]})$
is defined by $\delta(f) = \sum_{i=1}^{p} \pi_i^* f$.
The cocycle condition on $g$ implies that
$\delta(g) = 0 \in \Omega^0(Y^{[3]}) $. Since $d$ and $\delta$ commute
this means that $\delta(d\log(g)) = 0 \in \Omega^1(Y^{[3]})$ and so the
exactness of the complex \ref{nocoh} implies the existence
of $A \in \Omega^1(Y)$ such that $\delta(A) = d\log(g)$. 

It was established in $\cite{mur}$ that if we have
$\delta(dA) = 0$ for $A \in \Omega^q(Y)$ then there exists a unique 
$F \in \Omega^{q+1}(M)$ satisfying $dA = \pi^*(F)$.
In this case we have a two-form $F$ which we call the 
{\it bundle 0-gerbe curvature}. It is easily shown that changing the choice
of connection does not change the de Rham class of the curvature.

All of the operations which we have described on bundle 0-gerbes are possible
for bundle 0-gerbes with connection as well, 
\begin{eqnarray*}
(g,A)^* &=& (g^{-1},-A) \\
(g_1,A_1) \otimes (g_2,A_2) &=& (g_1g_2,A_1 + A_2)
\end{eqnarray*}
Now we show that corresponding to the bundle 0-gerbe connection is a 
connection for the corresponding bundle. Define 
a one-form on $Y \times S^1$ by 
\[
\tilde{A} = A + \theta^{-1}d\theta.
\]
Consider this form at two equivalent points
$(y_1,g(y_1,y_2))$ and $(y_2,1)$. The difference is given by
\[
A_{y_1} + d\log (g(y_1,y_2)) - A_{y_2}
\]
which is equal to $d\log(g) - \delta(A)$ on $Y^{[2]}$. Since this
is zero by the definition of the bundle 0-gerbe connection 
this 1-form is well defined on equivalence classes. Furthermore it can 
be easily shown that it satisfies the conditions for a connection 
1-form. 

Suppose $A$ is a connection 1-form on a bundle $L \rightarrow M$. We claim
that $A$ is also a connection on the corresponding bundle 0-gerbe.
This is true because
\begin{eqnarray*} 
\delta(A)_{(y_1,y_2)} &=& A_{y_2} - A_{y_1} \\
&=& A_{(y_1 g(y_1,y_2))} - A_{y_1} \\
&=& A_{y_1} + d\log(g)_{(y_1,y_2)} - A_{y_1} \\
&=& d\log(g)_{(y_1,y_2)}
\end{eqnarray*}
Thus we have a correspondence between bundle 0-gerbes with connection
and bundles with connection, however it is not yet clear whether this 
carries over to an isomorphism with Deligne cohomology. Recall that 
in the case without connections the role of isomorphism classes for
bundles was taken by stable isomorphism classes for bundle 0-gerbes.
We must define a slightly different notion of stable isomorphism
for bundle 0-gerbes with connection. This is because a trivialisation
of a bundle 0-gerbe corresponds to a trivialisation in \v{C}ech cohomology 
of the Chern class. Explicitly this is a cochain $\underline{h}$
satisfying $\delta (\underline{h}) = \underline{g}$.
When we have a choice of connection there is a further requirement on
$\underline{h}$ since we consider it as a Deligne cochain and require 
that $D(\underline{h}) = (\underline{g},\underline{A})$. This means 
that $h$ must satisfy
\begin{eqnarray}
\delta(\underline{h}) &=& \underline{g} \qquad \text{and} 
\label{0triv1}\\
d\log (\underline{h}) &=& \underline{A}. \label{0triv2}
\end{eqnarray}
The geometric realisation of the cochain $\underline{h}$ is
a function $h : Y \rightarrow S^1$ such that $\delta(h)
= g$ and $d\log h = A$. 
We shall refer to a cochain satisfying \eqref{0triv1} as a 
trivialisation and to one satisfying both \eqref{0triv1} and
\eqref{0triv2} as a {\it $D$-trivialisation} and a bundle 0-gerbe
with connection which has a $D$-trivialisation is called {\it $D$-trivial}.
\begin{definition}
Let $(g_1;A_1)$ and $(g_2;A_2)$ be bundle 0-gerbes with 
connection. We say that they are {\it $D$-stably isomorphic}
if there exists a $D$-trivial bundle gerbe with connection $(\tau;C)$
and an isomorphism
\[
(g_1;A_1) \isom (g_2;A_2) \otimes (\tau;C).
\] 
\end{definition}
It is easy to verify that the set of $D$-stable isomorphism classes of
bundle 0-gerbes with connection forms a group. 
\begin{proposition}\label{Dclass1}
The group of $D$-stable isomorphism classes of bundle 0-gerbes with
connection is isomorphic to $H^1(M,\D^1)$.
\end{proposition}
\begin{proof}
The proof is a simple extension of that for Proposition \ref{p1}.
\end{proof}
\end{section}

\begin{section}{Bundle Gerbes}\label{bndlgbs}
Bundle gerbes were introduced in $\cite{mur}$ as a geometric
realisation of classes in $H^3(M,\Z)$. They are the key object of 
interest in this thesis, here we present the basic theory.

\begin{definition}
Let $Y \stackrel{\pi}{\rightarrow} M$ be a submersion and
let $P \stackrel{\pi_P}{\rightarrow} Y^{[2]}$ be a $U(1)$ bundle. 
A {\it $U(1)$-bundle gerbe } is a triple $(P,Y,M)$ together
with a $U(1)$-bundle isomorphism $P_{(y_1,y_2)} \otimes P_{(y_2,y_3)}
\rightarrow P_{(y_1,y_3)}$ which is called the bundle gerbe product.
Associativity is required whenever triple products
are defined.
\end{definition}
The bundle gerbe $(P,Y,M)$ is represented diagrammatically by
\[
\begin{array}{ccc}
P & & \\
\downarrow & & \\
Y^{[2]} & \rightrightarrows & Y \\
& & \downarrow \\
& & M 
\end{array}
\]
Since we only deal with $U(1)$-bundle gerbes we shall refer to them
simply as bundle gerbes. 
Often we will say that $(P,Y)$ or $P$ is a bundle gerbe over $M$ when
there is no ambiguity. Given a map $\phi : N \rightarrow M$ we 
may define the pullback $\phi^{-1} P$ which is a bundle gerbe
on $M$. Given two bundle gerbes $(P,Y,M)$ and $(Q,X,M)$ there
is a product bundle gerbe $(P \otimes Q, Y \times_M X,M)$.
For any bundle gerbe $P$ there exists a dual bundle gerbe $P^*$.
For details of these constructions see \cite{mur}.\\

In the definition of a bundle gerbe the bundle over
$Y^{[2]}$ may be replaced with a bundle 0-gerbe $(\rho,X,Y^{[2]})$.
In this case the product is no longer a morphism since $X$ is not
acted on by $S^1$. Since we are dealing with bundle 0-gerbes rather
than bundles it is not surprising that morphisms should be 
replaced by stable morphisms. A choice of stable morphism
$\rho_{(y_1,y_2)} \otimes \rho_{(y_2,y_3)} \rightarrow \rho_{(y_1,y_3)}$
is equivalent to a choice of trivialisation 
\[
\rho_{(y_1,y_2)} \otimes \rho_{(y_2,y_3)} \otimes \rho_{(y_1,y_3)}^* 
\isom \delta(m_{123})
\]
This trivialisation represents a bundle gerbe product if it 
satisfies the associativity condition
\[
m_{123} \cdot m_{134} = m_{124} \cdot m_{234}.
\] 
\begin{definition}
A bundle gerbe $(P,Y,M)$ is called ${\it trivial}$ if there exists
a bundle $J \rightarrow Y$ such that there is a bundle isomorphism
\[
P \isom \pi_1^{-1} (J) \otimes \pi_2^{-1} (J)^*
\]
where $\pi_1$ and $\pi_2$ are the projections of each 
component of $Y^{[2]}$ onto $Y$. The product 
$\pi_1^{-1} (J) \otimes \pi_2^{-1} (J)^*$ is also denoted by $\delta(J)$.
\end{definition}
\begin{definition}
A {\it bundle gerbe morphism} between $(P,Y,M)$ and $(Q,X,N)$ is a triple of 
maps $(\alpha,\beta,\gamma)$
where $\beta : Y \rightarrow X$ is a fibre preserving map covering 
$\gamma: M \rightarrow N$ and $\alpha :P \rightarrow Q$ is a bundle
morphism covering the induced map $\beta^{[2]}: Y^{[2]} \rightarrow 
X^{[2]}$. Furthermore $\alpha$ must commute with the bundle gerbe
product. An {\it isomorphism of bundle gerbes} is a bundle
gerbe morphism
with $M=N$ and where $\gamma$ is the identity map. Two bundle 
gerbes $P$ and $Q$ are {\it stably isomorphic} if $P \isom Q \otimes 
\delta(J)$.
\end{definition}
\begin{proposition}\cite{mur}
The set of stable isomorphism classes of bundle gerbes over $M$ is
isomorphic to $H^3(M,\Z)$.
\end{proposition}
We construct a class $g_{\alpha \beta \gamma} \in H^2(M,\un{U(1)})$ 
corresponding to a bundle gerbe $(P,Y,M)$. The standard isomorphism gives
a corresponding class in $H^3(M,\Z)$
which is known as the Dixmier-Douady class. We shall also refer to 
$ g_{\alpha \beta \gamma}$ as the Dixmier-Douady class, or by analogy with
the local data associated with a bundle we shall also call these 
transition functions. 
Let $s_\alpha$ and $s_\beta$ be two local sections of $Y \rightarrow M$
defined on $U_\alpha \subset M$ and $U_\beta \subset M$ 
respectively. These define a section $(s_\alpha,s_\beta): U_{\alpha
\beta} \rightarrow Y^{[2]}$. Use this section to form the pull-back
bundle $P_{\alpha \beta} = (s_\alpha , s_\beta)^*P$ over 
$U_{\alpha \beta}$. Since $U_{\alpha \beta}$ is contractible 
$P_{\alpha \beta}$ is trivial and so admits a global section which we 
shall denote by $\sigma_{\alpha \beta} : U_{\alpha \beta}
\rightarrow P_{\alpha \beta}$. Over the triple intersection
$U_{\alpha \beta \gamma}$ the bundle gerbe product gives a bundle
isomorphism $P_{\alpha \beta} \otimes P_{\beta \gamma} \isom
P_{\alpha \gamma}$. Thus we can define $g_{\alpha \beta \gamma} :
U_{\alpha \beta \gamma} \rightarrow U(1)$ by 
\[
\sigma_{\alpha \beta} \otimes \sigma_{\beta \gamma}
= \sigma_{\alpha \gamma} g_{\alpha \beta \gamma}.
\]
To get a class in Deligne cohomology we will also need to 
define connections and curvings on bundle gerbes. 
\begin{definition}
Let $(P,Y,M)$ be a bundle gerbe. A {\it bundle gerbe connection}, $A$,
is a connection on the bundle $P \rightarrow Y^{[2]}$ which 
commutes with the bundle gerbe product.
\end{definition}
\begin{definition}
Let $(P,Y,M)$ be a bundle gerbe with connection $A$. Let 
$F_A \in \Omega^2(Y^{[2]})$ be
the curvature of $A$ considered as a bundle connection on 
$P \rightarrow Y^{[2]}$. A {\it curving} is a 2-form $\eta$ on $Y$ 
satisfying $\delta(\eta) = F_A$.
\end{definition} 
A bundle gerbe $(P,Y,M)$ with connection $A$ and curving $\eta$ may
also be referred to as $(P,Y,M;A,\eta)$ or $(P;A,\eta)$.
We may now define the Deligne class associated to a 
bundle gerbe $(P,Y,M)$ with connection, $A$, and curving,
$\eta$. Given local sections $s_\alpha : U_\alpha \rightarrow
Y$ we may define the {\it local curvings}
\[
\eta_\alpha = s_\alpha^* \eta.
\]
We have already defined the bundles $P_{\alpha \beta}$. 
The pull back by $(s_\alpha, s_\beta)$ induces connections on
each of these bundles which may be pulled back to $U_{\alpha \beta}$
using the sections $\sigma_{\alpha \beta}$ to give a 
collection of 1-forms $A_{\alpha \beta}$ on double intersections of
open sets on $M$. We call these {\it local connections}. 
\begin{proposition} \cite{mur}
Let $(P,Y,M)$ be a bundle gerbe with connection and curving. Let 
$g_{\alpha \beta \gamma}$ be the Dixmier-Douady class, $A_{\alpha
\beta}$ be the local connections and $\eta_\alpha$ be the local curvings.
Then $(g_{\alpha \beta \gamma},A_{\alpha \beta},\eta_\alpha)$ defines
a class in the sheaf cohomology group $
H^2(M,\D^2)$.
\end{proposition}
For $(g_{\alpha \beta \gamma}, A_{\alpha \beta}, \eta_\alpha)$
to be a class in  $H^2(M,\D^2)$ it must satisfy
\begin{eqnarray}
g_{\beta \gamma \delta} - g_{\alpha \gamma \delta}
+ g_{\alpha \beta \delta} - g_{\alpha \beta \gamma} &=& 0 
\label{bgdc1} \\
A_{\beta \gamma} - A_{\alpha \gamma} + A_{\alpha \beta}
&=& d\log(g_{\alpha \beta \gamma}) \label{bgdc2} \\
\eta_\alpha - \eta_\beta &=& dA_{\alpha \beta} \label{bgdc3}
\end{eqnarray}
As in the previous cases there is an isomorphism 
\[
H^3(M,\Z(3)_D) \isom H^2(M,\D^2)
\]
so each bundle gerbe with connection and curving gives rise to an
element of $H^3(M,\Z(3)_D)$.
Explicitly the Deligne class is given by 
\[
(n_{\alpha \beta \gamma \delta},\log(g_{\alpha \beta \gamma}),
A_{\alpha \beta}, \eta_\alpha)
\]
where $n_{\alpha \beta \gamma \delta} = \delta (\log(g))_{\alpha
\beta \gamma \delta}$ though we will often refer to the 
class $(g_{\alpha \beta \gamma}, A_{\alpha \beta}, \eta_\alpha)$
as the Deligne class. 

As with the case of bundle 0-gerbes it is necessary to introduce
$D$-trivialisations for bundle gerbes with connection and curving.
A $D$-trivialisation of a Deligne class $(\underline{g},\underline{A},
\underline{\eta})$ is a cochain $(\underline{h},\underline{B})$ which
satisfies 
\begin{eqnarray}
\delta(\un{h}) &=& \un{g} \label{1triv1}\\
d\log(\un{h}) - \delta(\un{B}) &=& \un{A} \label{1triv2}\\
d\un{B} &=& \un{\eta} \label{1triv3}
\end{eqnarray}
Geometrically a $D$-trivialisation of $(P;A,\eta)$ is a bundle $J$ with 
connection $B$ such that 
$\delta(J;B) \isom (P;A)$ as bundle gerbes with connection, where $\delta(J;B)$
is the bundle $\delta(J)$ with connection induced from $B$ by $\delta$.  
Furthermore, in order to satisfy \eqref{1triv3} the curvature of $(J;B)$ 
must be equal to the curving $\eta$. We may define $D$-stable isomorphisms
in the obvious way and state a bundle 
gerbe version of Proposition \ref{Dclass1}:
\begin{proposition}\label{Dclass2}
The group of $D$-stable isomorphism classes of bundle gerbes with 
connection and
curving are isomorphic to $H^2(M,\D^2)$.
\end{proposition}

\begin{proof}
First we show independence of the choice of sections. There are two
different types of section involved in the construction of the 
Deligne class. Suppose the sections $\sigma_{\alpha\beta}$ are replaced
by $\tilde{\sigma}_{\alpha\beta}$. We have two choices of section of 
a principal bundle so they differ by functions $f_{\alpha\beta}$ and
the corresponding change in transition functions is 
\[
\tilde{g}_{\alpha\beta\gamma} = g_{\alpha\beta\gamma}f_{\alpha\beta}
f_{\beta\gamma}f^{-1}_{\alpha\gamma}
\]
The local connections are related by the usual change of connection 
formula
\[
\tilde{A}_{\alpha \beta} = A_{\alpha\beta} + d\log f_{\alpha\beta}
\]
and the local curvings are unaffected so the overall contribution is
the trivial cocycle $D(f_{\alpha \beta},0)$.

Now suppose we change the sections $s_\alpha$ to $s'_\alpha$. In general 
these are not sections of a principal bundle so they do not differ by a 
function. Using the bundle gerbe product we have an isomorphism
\begin{equation}\label{chsec}
P_{(s'_\alpha,s'_\beta)} = P_{(s'_\alpha,s_\alpha)} \otimes 
P_{(s_\alpha,s_\beta)} \otimes P_{(s_\beta,s'_\beta)}
\end{equation}
Let $\sigma_{\alpha \beta}$, $\sigma'_{\alpha\beta}$, $\delta_\alpha$
and $\delta_\beta$ be sections of the trivial bundles 
$P_{(s_\alpha,s_\beta)}$, $P_{(s'_\alpha,s'_\beta)}$, 
$P_{(s'_\alpha,s_\alpha)}$ and $P_{(s'_\beta,s_\beta)}$ respectively.
We have two sections, $\sigma'_{\alpha \beta}$ and 
$\delta_\alpha \sigma_{\alpha\beta} \delta^{-1}_{\beta}$ of isomorphic 
bundles so they differ by functions $h_{\alpha\beta}$. When comparing
the transition functions defined using $\sigma_{\alpha \beta}$ or 
$\sigma'_{\alpha \beta}$ the $\delta$ sections all cancel out 
and we have essentially the previous case. 
Equation \eqref{chsec} also leads to an equation involving local
connections,
\begin{equation}
A'_{\alpha \beta} = k_\alpha + A_{\alpha \beta} - k_\beta
\end{equation}
where $k_\alpha$ is defined by pulling back the bundle gerbe 
connection to $P_{(s'_\alpha,s_\alpha)}$ and then 
pulling this connection back to $U_\alpha$ using the section $\delta_\alpha$.
Consider what happens to the local curvings. Since $\eta$ satisfies
$\delta(\eta) = F$ then ${s'_\alpha}^* \eta - {s_\alpha}^* \eta$ is equal to
the curvature of $P_{(s'_\alpha,s_\alpha)}$ which is $dk_\alpha$, so
\begin{equation}
\eta'_\alpha = \eta_\alpha + dk_\alpha
\end{equation}
so we have added a trivial cocycle $D(1,k_\alpha)$. Hence the Deligne 
class is independent of all choices of sections.

The homomorphism property is a straightforward consequence of 
the definition of the tensor product of bundle gerbes and the 
Deligne class so we omit details.

The result that a bundle gerbe is trivial if and only if
it has a trivial \v{C}ech class has been discussed in detail 
elsewhere (\cite{mur},\cite{ste}). Essentially it comes down to the 
fact that for a trivial bundle gerbe the sections $\sigma_{\alpha\beta}$
are of the form $\delta_\alpha^* \otimes \delta_\beta$. The inclusion
of connections and curvings does not add any significant complications.

Finally we need to describe a bundle gerbe which is classified by a 
particular Deligne class $(g_{\alpha \beta\gamma},A_{\alpha \beta},\eta_\alpha)$.
Let $Y = \amalg_\alpha U_\alpha$, the disjoint product of 
all of the elements of the open cover of $M$. Let $P \rightarrow Y^{[2]}$
be the trivial bundle. An element of $Y^{[2]}$ is of the form
$(m_\alpha,m_\beta)$ where $m \in U_{\alpha\beta}$ and $m_{\alpha}$ is
$m$ considered as an element of $U_\alpha \in Y$. The define the 
bundle gerbe product by
\begin{equation}
(m_\alpha,m_\beta,z_1)\cdot (m_\beta,m_\gamma,z_2) = (m_\alpha,m_\gamma,
z_1 z_2 g_{\alpha\beta\gamma})
\end{equation}
where $z_1, z_2 \in U(1)$.
Since $P$ is trivial then we can define the connection as a 1-form
on $Y^{[2]}$. At $(m_\alpha,m_\beta) \in Y^{[2]}$ the connection 
1-form is given by $A_{\alpha \beta}$ at $m$. Define the curving on
$U_\alpha \in Y$ by $\eta_\alpha$. 

\end{proof}

We have only considered bundle gerbes with connection and curving.
It is easily seen that bundle gerbes with a choice of connection
but no choice of curving are classified by $H^3(M,\Z(2)_D) \equiv
H^2(M,\D^1)$.
It is a standard result (see \cite{bry}) that 
$H^p(M,\Z(q)_D) \isom H^p(M,\Z(1)_D)$ whenever
$p > q$, thus the stable isomorphism class of bundle gerbe with
connection is invariant under a change of connection.

As in the previous cases the exact sequence \eqref{gaj} 
\[
0 \rightarrow H^2(M,U(1)) \stackrel{i}{\rightarrow} 
H^2(M,\D^2) \stackrel{d}{\rightarrow} \Omega^3(M)_0 
\rightarrow 0
\]
gives the curvature and flat holonomy. Explicitly the 
curvature is $\omega \in \Omega^3(M)$ satisfying 
\[
\pi^* \omega = d\eta
\]
and is guaranteed to exist since $\delta(d\eta)=0$. In terms
of local curvings the 3-curvature is defined in terms of
a collection of local 3-forms $\omega_\alpha = d\eta_\alpha$ 
which agree on overlaps since $\delta(d\eta_\alpha) = 0$. \\
The flat holonomy is calculated in the following way \cite{hit}.
Suppose $\omega = 0$. Then $d\eta_\alpha = 0$ so there
exist local 1-forms  $B_\alpha$ satisfying $dB_\alpha = \eta_\alpha$.
Furthermore 
\[
\eta_\beta - \eta_\alpha = dA_{\alpha \beta} = d(B_\beta - B_\alpha)
\]
so there exists functions $a_{\alpha \beta}$ which are defined on
double intersections and satisfy
\[
A_{\alpha \beta} - B_\beta + B_\alpha = d\log(a_{\alpha \beta}).
\]
It follows that
\[
d\log(a_{\alpha \beta} \cdot a_{\beta \gamma} \cdot a_{\alpha \gamma}^{-1}
\cdot  g_{\alpha \beta \gamma}^{-1}) = 0
\]
and the flat holonomy is 
\[
c_{\alpha \beta \gamma}  = a_{\alpha \beta} \cdot a_{\beta \gamma} 
\cdot a_{\alpha \gamma}^{-1} \cdot 
 (g_{\alpha \beta \gamma}^{-1}
)
\]
We conclude our discussion of flat bundle gerbes with the observation
that the Deligne cohomology groups $H^2(M,\Z(p)_D)$ represent flat 
bundle gerbes for any $p > 2$. 
A class in this cohomology group is the same as a class
\[
(g_{\alpha \beta \gamma},A_{\alpha \beta},\eta_\alpha) \in 
H^2(M,\D^2)
\]
with the additional condition that $d\eta_\alpha = 0$, therefore the class 
represents a flat bundle gerbe. 
\end{section}

%% file: chapter3.tex
\chapter{Examples of Bundle Gerbes}

We define and present the basic properties of a number of examples
of bundle gerbes which shall be of use to us. 

\begin{section}{Tautological Bundle Gerbes}
The tautological bundle gerbe was introduced in \cite{mur} as a way
to construct a bundle gerbe with any given closed, $2\pi$-integral 3-form 
as its 3-curvature. Our approach will be similar to that in \cite{ste}
however we use bundle 0-gerbes rather than bundles.

Let $M$ be a 1-connected manifold with distinguished base point $m_0$. 
Denote by 
$\P_0 M$ the space of paths in $M$ which are based at $m_0$. 
An element of $\P_0 M$ is a map $\mu: [0,1] \rightarrow M$ such that
$\mu(0) = m_0$. 
There is a fibration $\P_0 M \rightarrow
M$ defined by the projection $\pi: \mu \mapsto \mu(1)$.
The fibre product ${\P_0 M}^{[2]}$ over $m \in M$ consists of pairs 
of paths between $m_0$ and $m$. By reversing the orientation of one of the
paths
this pair may be identified with a loop based at $m_0$. Thus
we can identify ${P_0 M}^{[2]}$ with $\L_0 M$, the space of
smooth loops in $M$ which are based at $m_0$. There is a technical
point that needs to be dealt with here. When two paths are joined together
the resultant loop may not be smooth at the two points where the 
paths are joined. To overcome this problem we follow 
Caetano and Picken \cite{capi} and re-parametrise the
paths around these points such that there is a sitting instant at each of 
these points, that is,
an interval of length $\epsilon$ around a point $t_0 \in [0,1]$ 
such that the loop is constant in the interval $(t_0-\epsilon,t_0+
\epsilon)$. The obvious adjustment is made when $t_0 = 0$ (or equivalently
$t_0=1$). The structure above $\L_0 M$ is defined in terms of an integral
which is invariant under such reparametrisations.

Let $F$ be
a closed, $2\pi$-integral 2-form.
Let $\rho : \L_0 M \rightarrow U(1)$ be defined by
\begin{equation} \label{rsig}
\rho (\gamma) = \exp (\int_\Sigma F)
\end{equation}
where $\Sigma$ is any surface such that $\partial \Sigma = 
\gamma$. Equivalently we may write
\begin{equation} \label{rhom}
\rho (\mu_1,\mu_2) = \exp (\int_{I^2} H^*F)
\end{equation}
where $\mu_1, \mu_2 \in {\P M}^{[2]}$ and $H$ is a homotopy
between $\mu_1$ and $\mu_2$. To see that $\rho$ is independent
of the choice of $\Sigma$ note that if we choose a different
surface, $\bar{\Sigma}'$, where the bar indicates that the opposite 
orientation is induced on the boundary, and let 
\[
\rho' (\gamma) = \exp (\int_{\Sigma'} F)
\]
then we have
\begin{eqnarray*}
\rho (\gamma) / \rho' (\gamma) &=&
\exp (\int_{\Sigma \cup \bar{\Sigma}'} F) \\
&=& 1
\end{eqnarray*}

Let $\mu_1,\mu_2,\mu_3 \in \pi^{-1}(m) \subset 
\P_0 M$ and let $\gamma_{ij}$ denote
the loop identified with $(\mu_i,\mu_j) \in {\P_0 M}^{[2]}$. Then
\begin{eqnarray*}
\rho(\gamma_{12})\rho(\gamma_{23}) &=& 
\exp(\int_{\Sigma_{12}} F + \int_{\Sigma_{23}}F) \\
&=& \exp(\int_{\Sigma_{12}\cup\Sigma_{23}} F). 
\end{eqnarray*}
Note that $\gamma_{12}$ and $\gamma_{23}$ are connected along 
$\mu_2$ and hence the surface $\Sigma_{12} \cup \Sigma_{23}$ has
boundary $\gamma_{13}$ and the cocycle condition
\[
\rho(\gamma_{12})\rho(\gamma_{23}) = \rho(\gamma_{13}) 
\]
is satisfied. Therefore $(\rho,\P_0M,M)$ defines a
bundle 0-gerbe. Furthermore if we let the connection 
form on $\P_0M$ be given by
\[
A = \int_I ev^* F
\]
where $ev$ is the evaluation map $ev: \P_0M \times I \rightarrow M$
then it may be shown that $A$ satisfies
\begin{eqnarray*}
\delta(A) &=& d\log \rho \\
dA &=& \pi^* F
\end{eqnarray*}
If $M$ is not connected then we may carry out this construction on each
connected component.

\begin{lemma}\label{indbp}
The tautological bundle 0-gerbe is independent (up to stable isomorphism) 
of the choice of base point in $M$.
\end{lemma}
\begin{proof}
Suppose we have a curvature form $F$ and two choices of base point,
$m_0$ and $m_1$. We shall show that the resulting tautological
bundle 0-gerbes are stably isomorphic. Over $M$ we can form two 
different path fibrations, $\P_0M$ and $\P_1M$ using the two choices
of base point. Form the corresponding tautological bundle 0-gerbes
and take the fibre product,
\begin{equation}
\begin{array}{ccc}
& & S^1 \\
& \stackrel{\rho_0^{-1}\rho_1}{\nearrow} & \\
\L_0M \times_\pi \L_1M & \rightrightarrows & \P_0 \times_\pi
\P_1 \\
& & \downarrow \\
& & M 
\end{array}
\end{equation}
An element of $\pi^{-1} (m) \subset \P_0 \times_\pi \P_1$ is a path
from $m_0$ to $m_1$ passing through $m$. An element of 
$\L_0M \times_\pi \L_1M$ is a figure eight with $m$ at the centre 
with each loop passing through either $m_0$ or $m_1$. To define a 
trivialisation of this bundle 0-gerbe we need to choose a path $q$ from
$m_0$ to $m_1$. The trivialisation is then given by the function
$h(\mu,\eta) = \exp \int_{\Sigma} F$ where $\Sigma$ is a surface bounded by 
$q^{-1} \star \mu^{-1} \star \eta $. Taking $\delta$ of $h$ gives the integral
of $F$ over a surface with boundary $\eta_1^{-1} \star 
\mu_1 \star q \star q^{-1} \star \mu_2^{-1} \star \eta_2$. After eliminating
the $q \star q^{-1}$ component this
is equal to $\rho_0^{-1}\rho_1$, therefore the two bundle 0-gerbes are stably
isomorphic. Calculation of $d\log h$ at $(X_0,X_1) \in T(\P_0 M \times_\pi
\P_1 M)$ gives 
$\int_{\mu_1} F(\mu_1',X_1) - \int_{\mu_0} F(\mu_0',X_0)$ which is 
equal to $A_1 - A_0$, the difference of the connections corresponding 
to each choice of base point, so $h$ defines a $D$-stable morphism.
Since the construction depends on the choice of $q$ 
this is not a canonical stable 
isomorphism.
\end{proof}

\begin{example}\label{egLG}
Let $G$ be a 
compact simply connected semisimple Lie group.
Let $M$ be the loop group 
$\L G$ and let the curvature 2-form be $\int_{S^1} ev^* 
<g^{-1}dg \wedge [g^{-1}dg \wedge g^{-1}dg]>$ where $ev$ is the
evaluation map $ev: \L G \times S^1 \rightarrow G$, $<,>$ is the 
Killing form and $[,]$ is the Lie bracket. We may then 
construct a tautological bundle 0-gerbe. Since $G$ is simply
connected, discs in $G$ may
be thought of as paths of loops based at a constant loop and 
may be recentred as in lemma \ref{indbp}. This means that we 
may consider the fibre over $\gamma$ to consist of discs bounded 
by $\gamma$.
The bundle obtained by the standard construction from this 
tautological bundle 0-gerbe is 
the central extension of the loop group $\widetilde{\L G} \rightarrow \L G$
as described 
by Mickelsson \cite{mic}.
\end{example}

Now suppose that we have a closed, $2\pi$-integral 3-form, $\omega$ on a 
2-connected manifold $M$.
Let $Q[F] \rightarrow {(\P_0M)}^{[2]}$ be the tautological bundle over 
${(\P_0M)}^{[2]}$ with curvature $F = \int_{S^1} ev^* \omega$. Here we
have identified $(\P_0M)^{[2]}$ with $\L_0 M$ and used the 
evaluation map $ev : \L_0 M \times S^1 \rightarrow M$. The tautological
bundle on $\L_0 M$ may be defined since $M$ is 2-connected. 
To avoid the need for a base point in $\L_0 M$  we shall use a slightly
different definition of tautological bundle. In fact the tautological 
construction is more natural over a fibre product space, the introduction
of a base point when the base is not a fibre product compensates for this.
This construction of the tautological bundle follows the approach of 
\cite{camuwa}.
Over $\L_0M$ we have the space $\Sigma^\partial M$ of 2-surfaces, 
such that the fibre 
over $\gamma$ is a surface with $\gamma$ as its boundary. 
Elements of the fibre product may be considered as elements of $\Sigma M$, 
the space of smooth maps of closed 
2-surfaces into M (with possible reparametrisation to deal with
any problems with smoothness along $\gamma$) and 
we may define the tautological
function in the usual way to give a tautological bundle $Q[F]\rightarrow \L_0
 M$.
We 
can now construct a bundle gerbe over $M$,
\[
\begin{array}{ccc}
Q[F] & & \\
\downarrow & & \\
\L_0M & \rightrightarrows & \P_0M \\
& & \downarrow \\
& & M \\
\end{array}
\]
To define the product we observe that for any $Y$ the fibration
$\P_0(Y^{[2]}) \rightarrow Y^{[2]}$, where the base
point lies in the diagonal subset of $Y^{[2]}$,  admits a product 
covering $(y_1,y_2) \times (y_2,y_3)
\rightarrow (y_1,y_3)$ which is defined by 
composition of paths. Strictly speaking the composition of paths is 
not associative, however we do have associativity up to 
reparametrisations which do not affect the overall structure of the
bundle gerbe. In general this will be the bundle gerbe product for 
any bundle gerbes which we define in terms of a bundle on the loop
space. 
The connection is given by the connection on the tautological bundle,
which in this case may be written as
\[
A = \int_{\Sigma^\partial} ev^* \omega
\]

If the curving is defined by
\[
\eta = \int_I ev^* \omega
\]
then the curvature is $\omega$. 

As with the tautological bundle 0-gerbe, the two tautological
bundle gerbes obtained by a change of base point are stably isomorphic.
The trivialisation over $\P_0 \times_\pi \P_1$ is defined by
$J_{(\mu,\eta)} = Q[F]_{(\mu \star q^{-1},\eta)}$ where $\mu \in \P_0M$, 
$\eta \in \P_1M$, $q$ is a path from $m_0$ to $m_1$ and $Q[F]$ is the 
tautological bundle over $\L_1M$. Using the product on $Q[F]$ it can be
shown that
\begin{equation}
\delta(J)_{(\mu_1,\eta_1,\mu_2,\eta_2)} = 
Q[F]^*_{(\mu_1 \star q^{-1},\mu_2 \star q^{-1})} \otimes 
Q[F]_{(\eta_1,\eta_2)}
\end{equation}
Consider the fibre in $Q[F]$ over $(\mu_1 \star q^{-1},\mu_2 \star q^{-1})$.
This consists of surfaces bounded by $\mu_2 \star q^{-1} \star q \star \mu_1$
and may be identified with $Q[F]_{(\mu_1,\mu_2)}$ over $\L_0M$. Thus we see
that $J$ is a trivialisation.

We sometimes abbreviate the tautological bundle as $Q[F] \rightarrow M$ 
and the tautological bundle gerbe as $Q[\omega] \Rightarrow M$.
\end{section}

\begin{section}{Trivial Bundle Gerbes}\label{trivial}
In the previous chapter we defined what it means for a bundle gerbe to be
trivial or $D$-trivial. In this section we examine the properties of 
these classes of bundle
gerbe.
\begin{lemma} {\bf (\cite{mur})}\label{striv}
Let $(P,Y,M)$ be a bundle gerbe. Suppose the projection
$Y \stackrel{\pi}{\rightarrow} M$ admits a global section, $s$. 
Then $(P,Y,M)$ is a trivial bundle gerbe.
\end{lemma}
The trivialisation is $(s\circ \pi , 1)^{-1} P$. The converse of this 
proposition is not true. To see this, consider the following counterexample.
Let $Y \rightarrow M$ be a projection which admits local sections, but has
no global section. Let $Y^{[2]} \times S^1 \rightarrow Y^{[2]}$ be the 
trivial bundle. We make $(Y^{[2]} \times S^1,Y,M)$ into a bundle gerbe 
with the product $(y_1,y_2,\theta) \times (y_2,y_3,\phi) =
(y_1,y_3,\theta \phi)$. There are sections 
\[
\sigma_{\alpha \beta} :
U_{\alpha \beta} \rightarrow (s_\alpha,s_\beta)^{-1}(Y^{[2]} \times S^1)
\]
which are given by 
\[
\sigma_{\alpha \beta} (m) = (s_\alpha(m),s_\beta(m),1)
\]
and which clearly satisfy the cocycle identity
\[
\sigma_{\alpha \beta} \sigma_{\beta \gamma} = \sigma_{\alpha \gamma}.
\]
Thus the Dixmier-Douady class is 1 and the bundle gerbe is trivial.

As a special case of lemma \ref{striv} we may consider restricting a 
bundle gerbe over $M$ to an open set $U_\alpha$. This admits a 
section $s_\alpha$ and so we may construct a trivialisation 
as described above.

We now review the geometric realisation of a trivial Dixmier-Douady
class which was originally given in \cite{mur}, and described in greater
detail in \cite{ste}.
Let $(P,Y,M)$ be a bundle gerbe with Dixmier-Douady class $\underline{g}$
and let $\un{h}$ be a trivialisation. Let $J_\alpha$ be the bundle
on $s_\alpha (U_\alpha) \subset Y$ defined by
\[
J_{\alpha} = (1,s_\alpha \circ \pi)^{-1} P.
\]
Define isomorphisms $\phi_{\alpha \beta} : J_\alpha \rightarrow
J_{\beta}$ by
\[
\phi_{\alpha \beta} (u) = m(\sigma_{\alpha \beta} h_{\alpha \beta}^{-1}
\otimes u)
\]
where $u \in J_\alpha$. The bundle $J$ obtained by gluing together 
the $J_\alpha$ using the standard clutching construction 
with the isomorphisms $\phi_{\alpha \beta}$ is a
trivialisation of $P$. 

Conversely, if we are given a trivialisation $J$ then we can recover
the trivialisation of the Dixmier-Douady class in the following way.
Let $J_\alpha$ be defined by $s_{\alpha}^{-1} J$. Since this is a bundle
over $U_\alpha$ it must be trivial and admits a global 
section $\delta_\alpha$. Since $\delta(J) \isom P$ then there exist
functions $h_{\alpha \beta}: U_{\alpha \beta} \rightarrow S^1$ such that
\begin{equation}\label{hdef}
\sigma_{\alpha \beta}(m) = (\delta_\alpha^{-1}(m) \otimes 
\delta_\beta(m))h_{\alpha \beta}(m).
\end{equation}
It may be shown that the $h_{\alpha \beta}$ trivialise the Dixmier-Douady
class. 

As an example we may calculate the local data for the canonical 
trivialisation over an open set $U_0$. The trivialisation is defined 
by $J^0_y = P_{(s_0(\pi(y)),y)}$ where $s_0 : U_0 \rightarrow Y$ is a section.
Over any $U_\alpha$ restricted to $U_0$ we can pull back $J^0$ by a section
$s_\alpha$ to get $s_\alpha^{-1} J^0_m = P_{(s_0(m),s_\alpha(m))}$.
These have sections $\delta_\alpha = \sigma_{0\alpha}$. The local
data for the trivialisation, $h_{\alpha\beta}$ is then defined by
\begin{equation}
\sigma_{\alpha\beta} = \sigma^{-1}_{0\alpha} \otimes \sigma_{
0\beta} h_{\alpha\beta}
\end{equation}
so $h_{\alpha\beta} = g_{0\alpha\beta}$. 

Next we consider the relationship between trivial bundle gerbes
and $D$-trivial bundle gerbes. To do this we first consider 
the relationship between bundle 0-gerbes and bundle 0-gerbes with
connection which
is given by the exact sequence \ref{conseq1},
\[
0 \rightarrow \Omega^1(M) / \Omega^1_0(M) \rightarrow 
H^1(M,\D^1) \rightarrow \check{H}^2(M,\Z) \rightarrow 0.
\]
The space $\Omega^1(M) / \Omega^1_0(M)$ may be interpreted as equivalence
classes of 
connections on the trivial bundle 0-gerbe.
This implies that a trivial bundle gerbe with connection
$A \in \Omega^1_0(M)$ is $D$-trivial. We can define the 
$D$-trivialisation in the following way. Let $h = \pi^* \rho$,
where
\[
\rho(m) = \exp (\int_\mu A)
\]
where $\mu \in \P_0 M$ for some base point $m_0$ and
$\mu(1) = m$. The $2\pi$-integrality of $A$ ensures that $\rho$
is independent of the choice of path and applying the
Deligne differential to the class representing $\rho$ gives
the bundle 0-gerbe $(1,A)$. This construction is essentially
the same as that used for the tautological bundle 0-gerbe and 
bundle gerbe, so we may refer to $\rho$ as the tautological
function. In fact the function defining the tautological
bundle 0-gerbe is the tautological function. Note that the
construction relies on the assumption that $M$ is connected.
If $M$ is not connected then the construction may be repeated
for each connected component.

Suppose the bundle 0-gerbe with connection $(g;A)$ is 
$D$-trivial. Furthermore suppose we have a particular choice
of trivialisation, $\un{h}$, which is not necessarily a $D$-trivialisation.
We would like to see how this trivialisation differs from a 
$D$-trivialisation. Using the fact that $\delta(\un{h}) = \un{g}$ 
and applying $d\log$ gives
\begin{eqnarray*}
d\log(\delta (\un{h})) &=& d\log(\un{g}) \\
\delta(d\log(\un{h})) &=& \delta (A) 
\end{eqnarray*}
so there exists a 1-form $\chi$ such that
\[
d\log(h_\alpha) = A_\alpha - \chi
\]
If we had not assumed that $(g;A)$ is $D$-trivial then $\chi$ would 
represent the obstruction in $\Omega^1(M)/\Omega^1_0(M)$ to a trivial 
bundle 0-gerbe being $D$-trivial as well. This is true since
a change in choice of trivialisation changes
$\chi$ by an element of $\Omega^1_0(M)$ which is the 1-curvature of 
the function defined by the difference between two trivialisations.
Furthermore
\begin{eqnarray*}
d\chi &=& dA_\alpha \\
&=& F
\end{eqnarray*}
where $F$ is the bundle 0-gerbe curvature. We shall refer to
$\chi \in \Omega^1(M)/\Omega^1_0(M)$ as the {\it $D$-obstruction form}.

Now we return to the case where $(g;A)$ is $D$-trivial, hence 
it is flat and $A$ is locally exact. 
Thus $\chi$ is closed. If it is also $2\pi$-integral then we may construct the
tautological function $\rho$ on $M$ with curvature $\chi$. 
Finally we define a $D$-trivialisation by the product 
$h \cdot \pi^*\rho$. To check that it is indeed a $D$-trivialisation
observe that 
\begin{eqnarray*}
D(h_\alpha \cdot \rho) &=& (\delta(h)_{\alpha \beta},d\log(h_\alpha)) 
+ (1,\chi) \\
&=& (\delta(h)_{\alpha \beta},d\log(h_\alpha) + \chi) \\
&=& (g_{\alpha \beta},A_\alpha)
\end{eqnarray*}
The bundle gerbe case is very similar to that for bundle 0-gerbes and was 
described in \cite{must}. Let $(\un{g},\un{A},\un{\eta})$ be the Deligne
class of a bundle gerbe with connection and curving,
$(P;A,\eta)$.
Suppose we have a trivialisation $J$ which is represented
by a Deligne cochain $\un{h}$. Since there is an isomorphism between
bundle gerbes and bundle gerbes with connection we may choose
a connection on $B$ such that $(J;B)$ trivialises $(P;A)$, however
we may not assume that it trivialises $(P;A,\eta)$. In terms 
of cochains this means that we have
\begin{eqnarray}
\un{g} &=& \delta(\un{h}) \label{Dtriv1}\\
\un{A} &=& d\log(\un{h}) - \delta(\un{B}) \label{Dtriv2}
\end{eqnarray}
but it is not true that $\un{\eta} = d\un{B}$. We can, however, deduce from
\eqref{Dtriv2} that 
\[
\eta_\alpha - dB_\alpha = \chi
\]
If this $D$-obstruction form is closed and $2\pi$-integral then we 
can construct the tautological bundle $(Q[\chi];\int_I ev^*\chi)$ with 
curvature $\chi$ and
define a $D$-trivialisation by 
\[
(P;A,\eta) = D((J;B)\otimes \pi^{-1}(Q[\chi];\int_I ev^* \chi)).
\]

An extension of $D$-obstruction theory which is useful is 
the situation where we have two trivialisations $\delta(L)$ and
$\delta(J)$ of the same bundle gerbe. If they have $D$-obstruction
forms $\chi_L$ and $\chi_J$ respectively which satisfy
$d\chi_L = d\chi_J$ then
$D(L) = D(J \otimes \pi^{-1}K)$ where $K$ is the tautological bundle 
with curvature $\chi_L - \chi_J = dB_J - dB_L$ which is closed and
$2\pi$-integral since it is the difference of two curvatures.
We now consider the situation of a bundle gerbe with two different
trivialisations. 
\begin{proposition}\cite{must}
Let $(P,Y,M)$ be a bundle gerbe and let $L$ and 
$J$ be two trivialisations. Then there exists a bundle $(K,M)$
such that $L = J \otimes \pi^{-1}K$ as bundles over $Y$.
\end{proposition}
\begin{proposition}
Let $(P,Y,M;A,\eta)$ be a bundle gerbe and let $L$ and $J$ be
two $D$-trivialisations. Then there exists a flat bundle $(K,M)$
such that $L = J \otimes \pi^{-1}K$ as bundles with connection 
over $Y$.
\end{proposition}
This follows from the previous proposition together with the 
observation that the curvatures of $L$ and $J$ must both be 
equal to $\eta$.
\begin{proposition}\label{dobs}
Suppose $\delta(L) = \delta(K)$ and $F_L = F_K$. Furthermore 
suppose that $\chi_K - \chi_L$ is closed and $2\pi$-integral. Then there exists 
a bundle $J$ with curvature $\chi_K - \chi_L$ such that
$L = K \otimes \pi^{-1} J$. 
\end{proposition}
This result shall be useful for studying bundle 2-gerbes in the next 
chapter.
\begin{proof}
Suppose $\delta(L) = \delta(K) = P$. Then if $T$ is the canonically
trivial bundle $P^* \otimes P$ then $T = \delta(L) \otimes \delta(K)$.
Since the $D$-obstruction of the left hand side is trivial we have
$T = D(L^* \otimes K \otimes \pi^{-1}E)$ where $E$ has curvature 
$\chi_K - \chi_L$. There exists a flat bundle, $R$ such that 
$\pi^{-1} R = L^* \otimes K \otimes \pi^{-1}E$, so 
$L = K \otimes \pi^{-1}J$ where $\pi^{-1} J$ is the bundle $\pi^{-1}E \otimes 
\pi^{-1}R$ which has curvature $\chi_K - \chi_L$.
\end{proof}

There is also a local theory of trivialisations with 
connection. First consider a bundle gerbe that is $\delta$-trivial
but not necessarily $D$-trivial. Given a trivialisation 
with connection $(J;B)$ define 1-forms $k_\alpha$ by 
$\delta_\alpha^* B_\alpha$ where $\delta_\alpha$ is a section of
$s_\alpha^{-1} J_\alpha$ and $B_\alpha$ is the pulled back connection
on $J_\alpha$. Using the definition of $h_{\alpha\beta}$ \eqref{hdef}
it immediately follows that
\begin{equation}
d\log h_{\alpha\beta} - k_\beta + k_\alpha = A_{\alpha\beta}
\end{equation}
In the example of the bundle gerbe over $U_0$ we have $k_\alpha =
A_{0\alpha}$. To be $D$-trivial there is the additional requirement
that $dk_\alpha = \eta_\alpha$ which is satisfied if $F_{J^0} = \eta$. Note that
in the $U_0$ example we have $dk_\alpha = dA_{0\alpha} = \eta_\alpha - \eta_0$
so the $D$-obstruction is $\eta_0$.

\end{section}

\begin{section}{Lifting Bundle Gerbes}\label{lifting}
The lifting bundle gerbe was introduced in \cite{mur} as one of 
the motivating examples of the theory of bundle gerbes.
Let 
\[
0 \rightarrow U(1) \rightarrow \tilde{G} \stackrel{p}{\rightarrow} G
\rightarrow 0
\]
be a central extension of groups and let $P_G \rightarrow M$ be a principal
$G$ bundle. The lifting bundle gerbe is defined by the following diagram:
\[
\begin{array}{ccc}
& & \tilde{G} \\
& & \downarrow \\
& & G \\
& \stackrel{g}{\nearrow} & \\
P_G^{[2]} & \rightrightarrows & P_G \\
& & \downarrow \\
& & M 
\end{array}
\] 
The map $g:P_G^{[2]} \rightarrow G$ is defined such that $g(p_1,p_2)$
is the element of $G$ which satisfies $p_2 = p_1 g(p_1,p_2)$.
Alternatively $P_G^{[2]}$ may be identified with $P \times G$ via
$(p_1,p_2) \mapsto (p_1,g(p_1,p_2))$ and $(p,g) \mapsto (p,pg)$ in 
which case let $g: P_G \times G \rightarrow G$ be the projection
of the second factor.
It is to be understood that $\tilde{G} \rightarrow G$ is pulled back to
a bundle over $P_G^{[2]}$ by $g$ and the bundle gerbe product is
induced from the group product on $\tilde{G}$.
\begin{proposition}\cite{mur}
The lifting bundle gerbe associated with a $G$ bundle $P_G \rightarrow M$
and a central extension $\tilde{G}$ is trivial if and only if 
$P_G$ lifts to a $\tilde{G}$ bundle.
\end{proposition}
In general a connection $A$ on $\tilde{G} \rightarrow G$ does not
define a bundle gerbe connection. This is because the corresponding curvature
form, $g^* F_A$, for the bundle $g^{-1} \tilde{G}$ may not satisfy 
the condition $\delta( g^* F_A)$ on $P^{[3]}$. It is shown 
in \cite{must2} that in general there exists a 1-form $\epsilon$ on
$P_G^{[2]}$ such that $g^*(A) - \epsilon$ is a bundle gerbe connection.
\end{section}

\begin{section}{Torsion Bundle Gerbes}
If the Dixmier-Douady class of a bundle gerbe is torsion then we 
refer to it as a {\it torsion bundle gerbe}. These bundle gerbes 
naturally arise in applications to physics and there are two
particular aspects which are of interest: the canonical bundle gerbe
of a class in $H^2(M,\Z_p)$ and bundle gerbe modules.

Associated to the short exact sequence
\[
\Z \stackrel{p \times}{\rightarrow} \Z \rightarrow \Z_p 
\]
is a Bockstein operator $\beta: H^2(M,\Z_p) \rightarrow H^3(M,Z)$.
This indicates that given a class $w \in H^2(M,\Z_p)$ we 
may define a bundle gerbe with Dixmier-Douady class $\beta(w)$. 
We would like to demonstrate that there is a canonical choice of
Deligne class arising from $w$. 

The class $\beta(w)$ must satisfy $p\cdot\beta(w) = 0$ so we
have a torsion bundle gerbe. Consider the exact sequence
\[
0 \rightarrow H^2(M,U(1)) \rightarrow H^2(M,\D^2)
\rightarrow \Omega_0^3(M) \rightarrow 0
\]
Recall that this may be interpreted as 

\begin{center}
{\it flat holonomy class $\rightarrow$ bundle gerbe with connection and curving
$\rightarrow$ curvature }\\
\end{center}

\noindent 
If the bundle gerbe is torsion then since the curvature is the image of
the Dixmier-Douady class in de Rham cohomology it must be an exact form. 
The possible $D$-stable isomorphism class of bundle gerbes with a particular
choice of curvature are given by flat holonomy classes. The map
$H^2(M,\Z_p) \rightarrow H^2(M,U(1))$ allows us to consider $w$ as
a flat holonomy class. This in turn defines a bundle gerbe with
Deligne class $(\un{w},0,\un{dB})$. The fact that
the transition functions are constant and the curvature is exact mean
that this is a Deligne cocycle. The Dixmier-Douady class of this 
bundle gerbe is $w_{\alpha\beta\gamma} =
-\log w_{\beta \gamma} + \log w_{\alpha \gamma} 
- \log w_{\alpha \beta}$. The class $pw_{\alpha \beta \gamma}$
is trivial in $H^3(M,\Z)$, with trivialisation 
$p\log w_{\alpha \beta}$. The mod $p$ reduction then gives 
$w_{\alpha \beta} \in H^2(M,\Z)$, so the transition functions are
given by $\beta(w)$ as desired.
Trivially the flat holonomy class of
this bundle gerbe is $\un{w}$. A canonical choice of such a bundle
gerbe is given by setting $\un{dB} = 0$. Thus associated with 
a torsion class $w \in H^2(M,\Z_p)$ we have a canonical
torsion Deligne class $(\un{w},0,0)$. This construction 
extends to Deligne classes of arbitrary degree.

An interesting class of examples of torsion bundle gerbe
is given by the lifting bundle gerbes associated with a central extension
\[
\Z_p \rightarrow \tilde{G} \rightarrow G 
\]
We may use the map
$H^2(M,\Z_p) \rightarrow H^2(M,U(1))$ and the discussion above
to see that this is a torsion bundle gerbe. A particular example
is given by the obstruction to lifting a projective unitary bundle to a 
unitary bundle \cite{cajomu}, which is defined in terms 
of the central extension
\[
\Z_n \rightarrow U(n) \rightarrow PU(n) 
\]

Next we define bundle gerbe modules.
\begin{definition}\cite{bcmms}
Let $(P,Y,M)$ be a bundle gerbe. Let $E \rightarrow Y$ be a finite
rank hermitian vector bundle such that there exists a hermitian
bundle isomorphism 
\[
\phi: P \otimes \pi_1^{-1} E \isom \pi_2^{-1} E
\]
We require that this isomorphism is compatible with the bundle 
gerbe product in the sense that the maps
\[
P_{(y_1,y_2)} \otimes (P_{(y_2,y_3)} \otimes E_{y_3}) \rightarrow
P_{(y_1,y_2)} \otimes E_{y_2} \rightarrow E_{y_1}
\]
and
\[
(P_{(y_1,y_2)}\otimes P_{(y_2,y_3)}) \otimes E_{y_3} \rightarrow
P_{(y_1,y_3)} \otimes E_{y_3} \rightarrow E_{y_1}
\]
are the same. Call $E$ a {\it bundle gerbe module} and 
say that the bundle gerbe $P$ acts on $E$.
\end{definition}
A rank one bundle gerbe module is a trivialisation. Given
a rank $r$ bundle gerbe module the product $P^r$ acts on the
rank one bundle $\Lambda^r(E)$ and hence we have
\begin{proposition}\cite{bcmms}
Let $(P,Y,M)$ be a bundle gerbe with Dixmier-Douady class $\dd(P)$. 
Suppose $(P,Y,M)$ has a bundle gerbe module $E$ of rank $r$. Then
$P$ is a torsion bundle gerbe with $r\dd(P) = 0$.
\end{proposition} 
A connection on a bundle gerbe module is called a 
{\it bundle gerbe module connection} if the isomorphism
$\phi$ is an isomorphism of bundles with connection. 

A bundle gerbe module with connection may also be defined in terms of
local data. Suppose we have a bundle gerbe represented locally in
Deligne cohomology by $(\un{g},\un{A},\un{f})$. Let $E$ be a bundle
gerbe module, and define a set of local bundles on $M$ by
$E_\alpha = s_\alpha^{-1} E$. These bundles are trivial with
sections $\delta_\alpha$. We consider the isomorphism
$P \otimes \pi^{-1}_1 E \isom \pi^{-1}_2 E$ at the local
level. In terms of sections we may define local matrix valued 
functions $h_{\alpha \beta}$ such that
\begin{equation}
\sigma_{\alpha \beta} \otimes \delta_\beta = \delta_\alpha h_{\alpha \beta}
\end{equation}
Consider the section $\sigma_{\alpha \beta} \otimes \sigma_{\beta
\gamma} \otimes \delta_\gamma$  associated with $P_{\alpha \beta} \otimes 
P_{\beta \gamma} \otimes E_\gamma$. This can be simplified in two different
ways,
\begin{equation}
\begin{split} 
\sigma_{\alpha \beta} \otimes \sigma_{\beta \gamma} \otimes 
\delta_\gamma &= \sigma_{\alpha \gamma} g_{\alpha \beta \gamma} \otimes
\delta_\gamma \\
&= \delta_\alpha h_{\alpha \gamma} g_{\alpha \beta \gamma}1
\end{split}
\end{equation}
where $1$ is the identity matrix of the same rank as $E$, or
\begin{equation}
\begin{split}
\sigma_{\alpha \beta} \otimes \sigma_{\beta \gamma} \otimes 
\delta_\gamma &= \sigma_{\alpha\beta} \otimes \delta_{\beta} h_{\beta \gamma}
\\
&= \delta_\alpha h_{\alpha \beta} h_{\beta\gamma} 
\end{split}
\end{equation}
and hence we have 
\begin{equation}
h_{\alpha \beta} h_{\beta \gamma} = h_{\alpha \gamma} g_{\alpha \beta \gamma}1
\end{equation} 
\end{section}
To get a local expression for the connection let $\nabla_P$ and 
$\nabla_E$ be the connections on $P$ and $E$ respectively. Using
the sections $\sigma_{\alpha\beta}$, $\delta_\alpha$ and $\delta_\beta$
we have two choices of connection on bundles $P_{\alpha \beta} 
\otimes E_{\alpha}$ and $E_\beta$ over $U_\alpha$. These are
$\sigma_{\alpha\beta}^{-1}\nabla_P + \delta_\alpha^{-1} \nabla_E$ 
and $\delta_\beta^{-1} \nabla_E$. The isomorphism of bundles is 
given by the local functions $h_{\alpha \beta}$ and so, 
letting $\delta_\alpha^{-1} \nabla_E = k_\alpha$, the 
two choices of connection are related by the usual formula
for connections under a change of section,
\begin{equation}
A_{\alpha \beta}1 + k_\alpha = h_{\alpha \beta}^{-1}
k_\beta h_{\alpha \beta} + h^{-1}_{\alpha \beta} 
dh_{\alpha\beta}
\end{equation}

\begin{section}{Cup Product Bundle Gerbes}\label{cupsect}
There is a cup product in Deligne cohomology (see \cite{esvi} or
\cite{bry}). The correspondence between Deligne cohomology and geometric
objects implies that the cup product may be used to construct new 
examples of
geometric objects which realise Deligne classes \cite{brmc2}. We shall
demonstrate how to construct bundle gerbes corresponding to 
various cup products. First we consider the bundle obtained by
taking the cup product of two functions (\cite{bry},\cite{esvi}) as 
a bundle 0-gerbe. This helps us to find geometric realisations of 
various bundle gerbes which may be obtained by taking cup products. Let 
$f$,$g$ and $h$ be $U(1)$-valued functions on $M$ and let $L \rightarrow M$ be
a $U(1)$-bundle. We consider the cup products $f\cup L$, $L\cup f$ and
$f \cup g\cup h$.

The cup product is induced by a product in the Deligne complex
\[
\cup : \Z(p)_D \otimes \Z(q)_D \rightarrow \Z(p+q)_D
\] 
which is defined by
\begin{equation}\label{defcp}
x \cup y = \begin{cases}
x\cdot y & \mbox{if deg $x =0$},\\
x\wedge dy & \mbox{if deg $x > 0$ and deg $y = q$},\\
0 & \mbox{otherwise}.
\end{cases}
\end{equation}
It is a standard result that this product is associative and 
induces a product of Deligne cohomology groups. Furthermore it is 
anticommutative, that is, for $\alpha \in H^q(M,\Z(p)_D)$ and
$\beta \in H^{q'}(M,\Z(p')_D)$ the cup product satisfies
$\alpha \cup \beta = (-1)^{qq'} \beta \cup \alpha$. We shall
calculate some specific examples and construct corresponding geometric
objects. 

The cup product of two functions was described in \cite{esvi} and
\cite{bry}. We review it in detail as it provides the basis
for all of our subsequent examples.
Let $f$ and $g$ be $U(1)$-valued functions on $M$.
As we have seen $f$ and $g$ may be represented by the 
Deligne classes $(n_{\alpha \beta},\log_\alpha (f)) \in H^1(M,\Z(1)_D)$ and
$(m_{\alpha \beta},\log_\alpha (g)) \in H^1(M,\Z(1)_D)$ respectively, where 
$\log_\alpha$ is a branch of the logarithm function which
is defined on $U_\alpha \subset M$, and the integers
$n_{\alpha \beta}$ and $m_{\alpha \beta}$ are defined
by the differences $\log_\beta (f) - \log_\alpha(f)$ and
$\log_\beta(g) - \log_\alpha (g)$ respectively. The cup
product $f \cup g$ is the Deligne class $(n_{\alpha \beta}
m_{\beta \gamma},n_{\alpha \beta}\log_\beta(g),\log_\alpha(f) d\log(g))
\in H^2(M,\Z(2)_D)$.
Under the isomorphism $H^2(M,\Z(2)_D) \isom H^1(M,\D^1)$
this becomes $(g^{n_{\alpha \beta}},\log_{\alpha}(f) d\log(g))$ which
represents a bundle 0-gerbe with connection which may be described
explicitly. 

We construct a bundle 0-gerbe over $S^1 \times S^1$ and pull it back to
$M$ via the map $(f,g): M \rightarrow S^1 \times S^1$. Let the bundle 0-gerbe
over $S^1 \times S^1$ be defined by the following diagram:
\[
\begin{array}{ccc}
& & S^1 \\
& \stackrel{\rho_\cup}{\nearrow} & \\
(\R \times \Z) \times (\R \times \Z) 
& \rightrightarrows & \R \times \R \\
& & \downarrow \\
& & S^1 \times S^1 
\end{array}
\]
where the projection to the base is given by two copies of the
exponential and the map $\rho_\cup$ is defined by
$\rho_\cup(r,n,s,m) = \exp (sn)$. We have used the identification 
$\R^{[2]} = \R \times \Z$ for the $\Z$-bundle $\R \rightarrow S^1$ as
discussed for a general $G$-bundle in \S \ref{lifting}.
It is 
easily shown that the 1-form
$rds$ is a connection
for this bundle 0-gerbe.
\begin{proposition}
The bundle 0-gerbe $(f,g)^{-1}(\rho_\cup,\R \times \R,S^1 \times S^1)$ with 
connection $rds$ has Deligne class 
$(g^{n_{\alpha \beta}},\log_{\alpha}(f) d\log(g))$.
\end{proposition}
\begin{proof}
Define local sections $s_{\alpha}$ of $ \R \times \R \rightarrow
M$
by $s_\alpha (\theta,\phi) = (\log_\alpha (\theta),\log_\alpha(\phi))$. Then 
$(s_\alpha,s_\beta) = (\log_\alpha (\theta), \log_\beta(\theta) - \log_\alpha
(\theta),
\log_\beta(\phi),
\log_\beta (\phi) - \log_\alpha(\phi))$ and
$\rho(s_\alpha,s_\beta) = \exp(\log_\beta(\phi)(\log_\beta(\theta) - \log
_\alpha(\theta)))$. To get the pull back to $M$ we simply replace $\theta$ 
and $\phi$ with $f(m)$ and $g(m)$ respectively to get the required 
transition functions,
$g^{n_\alpha \beta}$. To complete
the proof observe that $s_\alpha^*(rds) = \log_\alpha(\theta)
d\log(\phi)$, so under the pull back we get $\log_\alpha(f)d\log(g)$.
\end{proof}
Clearly the bundle 0-gerbe $g \cup f$ is obtained by replacing $\rho_\cup$ with
$\rho_\cup^*: (r,n,s,m) \mapsto \exp(rm)$. By anticommutativaty the
product bundle 0-gerbe $\rho_\cup \otimes \rho_\cup^*$ should be trivial.
We may demonstrate this directly by defining
$r_1,r_2,s_1$ and $s_2$ such that
$(r,n,s,m)=(r_1,s_1,r_2-r_1,s_1,s_2-s_1)$ and considering
\begin{eqnarray*}
(\rho_\cup\rho_\cup^*)(r,n,s,m) &=& \exp(sn+rm) \\
&=& \exp(s_1(r_2-r_1) + (r_2 - n)m) \\
&=& \exp(s_1r_2 -s_1r_1 + r_2(s_2 - s_1) - nm) \\
&=& \exp(s_1r_2 -s_1r_1 + r_2s_2 - r_2s_1) \\
&=& \exp(r_2s_2 - r_1s_1) \\
&=& \delta(\exp(rs))
\end{eqnarray*}

There are three ways of obtaining a bundle gerbe via cup products. We shall
calculate the Deligne class and provide a geometric construction for each one.

\subsubsection{The Bundle Gerbe $ f \cup L$}
Let $f$ be a $U(1)$-valued function on $M$ and let $L$ be a
bundle 0-gerbe over $M$. 
Let $f$ and $L$ have Deligne class 
$(n_{\alpha \beta},\log_\alpha (f)) \in H^1(M,\Z(1)_D)$ and 
$(m_{\alpha \beta \gamma},\log (g_{\alpha \beta})) \in H^2(M,\Z(1)_D)$ 
respectively. Then $m_{\alpha\beta\gamma} 
= -\log (g_{\beta \gamma}) + \log  
(g_{\alpha \gamma}) - \log (g_{\alpha \beta})$.
The product $f \cup L$ is 
\[
(n_{\alpha \beta}m_{\beta \gamma \delta}, n_{\alpha \beta}
\log (g_{\beta \gamma}), \log_\alpha (f) d\log (g_{\alpha \beta})) \in 
H^3(M,\Z(2)_D)
\]
which, under the usual isomorphism, becomes $(g_{\beta \gamma}^{n_{\alpha 
\beta}}, \log_\alpha (f) d\log (g_{\alpha \beta}))$. We define the 
corresponding bundle gerbe with the following diagram:
\[
\begin{array}{ccc}
& & K_\cup \\
& & \downarrow \\
& & S^1 \times S^1 \\
& \stackrel{\rho}{\nearrow} & \\
\R \times \Z \times L \times S^1 & \rightrightarrows 
& \R \times L \\
& & \downarrow \\
\mspace{90mu} M & \stackrel{(f,m)}{\longrightarrow} & S^1 \times M 
\end{array}
\]
where $K_\cup$ is the cup product bundle described in the previous section,
$m:M\rightarrow M$ is the identity map
and the pullback to $M$ by $(f,m)$ is implied.
Local sections are defined by $(\log_\alpha f , s_\alpha)$ where 
$s_\alpha$ is a local section of $L$. The sections $\sigma_{\alpha
\beta}$ are given by sections of $\rho^{-1}K_\cup$ over 
$(\log_\alpha f,n_{\alpha\beta},s_\alpha, g_{\alpha\beta})$. 
Essentially the fibres of
$\rho^{-1}K_\cup$ look like $\R \times \R \times S^1$ with an equivalence 
relation $(a+n,b+m,z) \sim (a,b,z e^{nb})$ where $a, b \in \R$, 
$n,m \in \Z$ and $z \in S^1$. In the fibre there are also copies
of $\Z$ and $L$ but we omit these to simplify the expressions.
The projection takes $(a,b,z)$ to
$(a,n,l,e^b) \in \R \times \Z \times L \times S^1$. We need an
expression for the bundle gerbe product. It must satisfy two conditions:
it must cover a particular product on the base and it must respect the
equivalence relation in the definition of $\rho^{-1} K_\cup$.
To find the map on the base which must be covered by the bundle gerbe
product consider it in the form $(\R \times L)^
{[2]}$ in which case the product is
\begin{equation}
(r_1,r_2,l_1,l_2) \times (r_2,r_3,l_2,l_3) \rightarrow (r_1,r_3,l_1,l_3)
\end{equation}
Under the standard identification with $\R \times \Z \times 
L \times S^1$ given by $(r,n,l,\theta) = (r,r+n,l,l\theta)$ this becomes
\begin{equation}
(r_1,n_1,l_1,\theta_1) \times (r_2,n_2,l_2,\theta_2)
\rightarrow (r_1,n_1 + n_2,l_1,\theta_1\theta_2)
\end{equation}
so the bundle gerbe product must be of the form
\begin{equation}
(a_1,b_1,z_1) \times (a_2,b_2,z_2) = (a_1,b_1+b_2,z_1z_2 \Pi)
\end{equation}
where $\Pi$ is some function $\Pi (a_1,b_1,z_1,a_2,b_2,z_2)$.

To determine an expression for $\Pi$ we consider what happens to this 
product under the equivalence relation. First we replace $(a_1,b_1,z_1)$
by $(a_1+n,b_1+m,z_1 e^{-nb})$.
\begin{equation}
\begin{split}
(a_1+n,b_1+m,z_1 e^{-nb_1}) \times (a_2,b_2,z_2)
&= (a_1+n, b_1 + b_2 + m, z_1z_2e^{-nb_1} \Pi) \\
&= (a_1, b_1 + b_2,z_1z_2e^{-nb_1}e^{n(b_1+b_2)}\Pi) \\
&= (a_1, b_1 + b_2,z_1z_2e^{nb_2}\Pi)
\end{split}
\end{equation}
so $\Pi(a_1+n,b_1+m,z_1 e^{-nb_1},a_2,b_2,z_2) =
e^{-nb_2}\Pi(a_1,b_1,z_1,a_2,b_2,z_2)$. Now we consider the 
second factor,
\begin{equation}
\begin{split}
(a_1,b_1,z_1) \times (a_2+n,b_2+m,z_2 e^{-nb_2})
&= (a_1,b_1+b_2+m,z_1z_2 e^{-nb_2}\Pi) \\
&= (a_1,b_1+b_2,z_1z_2 e^{-nb_2}\Pi) 
\end{split}
\end{equation}
so $\Pi(a_1,b_1,z_1,a_2+n,b_2+m,z_2 e^{-nb_2}) = e^{nb_2}\Pi
(a_1,b_1,z_1,a_2,b_2,z_2)$.

Let $\Pi(a_1,b_1,z_1,a_2,b_2,z_2) = e^{b_2(a_2-a_1)}$. 
Under the transformation $a_1 \rightarrow a_1+n$ we have
$\Pi \rightarrow e^{b_2(a_2-a_1-n)} = \Pi e^{-nb_2}$. Under
the transformation $(a_2,b_2) \rightarrow (a_2+n,b_2+m)$ we
have $\Pi \rightarrow e^{b_2(a_2-a_1)}e^{nb_2}e^{m(a_2-a_1)}
= \Pi e^{nb_2}$ since $m(a_2-a_1) \in \Z$. 

Define sections $\sigma_{\alpha\beta} = (\log_\alpha f,\log g_{\alpha
\beta},1)$.
We may now calculate the product $\sigma_{\alpha\beta}\sigma_{\beta\gamma}$,
\begin{equation}
\begin{split}
\sigma_{\alpha\beta}\sigma_{\beta\gamma} &=
(\log_\alpha f , \log g_{\alpha\beta}, 1) \times
(\log_\beta f , \log g_{\beta\gamma}, 1) \\
&= (\log_\alpha f , \log g_{\alpha \beta} + \log g_{\beta\gamma}, 
e^{n_{\alpha\beta}\log g_{\beta\gamma}} ) \\
&= (\log_\alpha f, \log g_{\alpha \gamma} + m_{\alpha\beta\gamma},
g_{\beta\gamma}^{n_{\alpha\beta}})\\
&= (\log_\alpha f, \log g_{\alpha \gamma},g_{\beta\gamma}^{n_{\alpha\beta}})\\
&= \sigma_{\alpha\gamma} g_{\beta\gamma}^{n_{\alpha\beta}}
\end{split}
\end{equation}
This gives the required transition functions.


The local
connections $\log_\alpha (f) d\log (g_{\alpha \beta})$
are induced by a bundle gerbe connection $adb$ at $(a,b,z) \in \R \times \R
\times S^1 / \sim $, the total space of $K_\cup$. 

We may also consider the cup product bundle gerbe $f \cup L$ where
$L$ is a bundle with connection $A$. The Deligne
class of the product is $(g_{\beta \gamma}^{n_{\alpha \beta}},
n_{\alpha \beta}A_\beta,\log_\alpha (f) dA)$. The 
Dixmier-Douady class is the same as in the previous case however
the connection is different and we also have a choice of curving.
The connection form is $n(A_l + 
db)$ at a point $(r,n,l,\theta,z)$ in the total space of $\rho^{-1}K$, 
$\R \times \Z \times L \times \R \times S^1$.
The curving is $(r dA)$ on $\R \times L$. It is not difficult to see that
these give the appropriate local expressions. \\
 
There is an alternate geometric representation of the bundle gerbe 
$f \cup L$ which is of interest. Define it with the following 
diagram:
\[
\begin{array}{ccc}
L^{-n} & & \\
\downarrow & & \\
\R \times \Z \times M & \rightrightarrows & \R \times M \\
& & \downarrow \\
M & \stackrel{(f,m)}\rightarrow & S^1 \times M  
\end{array}
\]
where $L^{-n}$ is the bundle with fibre at $(r,n,m)$ equal
to the fibre of the $n$-th tensor product bundle of $L^*$ at $m$.
Sections of $Y^{[2]}$ over $U_{\alpha\beta}$ are given by 
$(\log_\alpha f , n_{\alpha \beta}, m)$, so we define the 
sections $\sigma_{\alpha\beta}$ by
\begin{equation}
\sigma_{\alpha\beta} = s_{\beta}^{-n_{\alpha\beta}}
\end{equation}
where $s_\alpha$ is a local section of $L$. Using these to calculate
the transition function we get
\begin{equation}
\begin{split}
\sigma_{\alpha\beta} \sigma_{\beta\gamma} &=
s_{\beta}^{-n_{\alpha\beta}}s_{\gamma}^{-n_{\beta\gamma}} \\
&= s_\gamma^{-n_{\alpha\beta}} g_{\beta\gamma}^{n_{\alpha\beta}}
s_{\gamma}^{-n_{\beta\gamma}} \\
&= s_{\gamma}^{-n_{\alpha\gamma}}g_{\beta\gamma}^{n_{\alpha\beta}}\\
&= \sigma_{\alpha\gamma} g_{\beta\gamma}^{n_{\alpha\beta}}
\end{split}
\end{equation}
The connection at $(r,n,m)$ is given by $-nA$ at $m$, 
the connection on $L^{-n}$ induced in the natural way from that of $L$, 
and the curving at $(r,m)$ is $rF$ where $F$ is the curvature of $L$.

This representation allows for a much simpler calculation of the Deligne
class, 
however from the
point of view of bundle gerbe theory it is 
not as general as the previous case since it depends on $L$ as a bundle
rather than a bundle 0-gerbe. The significance of the first method 
is that the cup product bundle appears in the definition of the
cup product bundle gerbe. This is related to the bundle gerbe
hierarchy which we shall consider in the next chapter. We shall find
both methods useful in considering cup products and bundle 2-gerbes in
\S \ref{bun2ger}. Also
it is not obvious how to approach the bundle gerbe $L \cup f$ using
the second method.

\subsubsection{The Bundle Gerbe $L \cup f$}
The Deligne class for this bundle gerbe with connection is
$(f^{m_{\alpha\beta\gamma}},\log (g_{\alpha\beta}) d\log f)$. By the 
commutativity of the cup product this bundle gerbe should be stably
isomorphic to $f \cup L$. We define it by replacing $K$ with the
dual bundle obtained by swapping the two functions in the 
cup product.
The result is that the equivalence
relation on $\R \times \R \times S^1$ becomes $(a,b,z) \sim
(a+n,b+m,ze^{ma})$. The product is still of the form
\begin{equation*} 
(a_1,b_1,z_1) \times (a_2,b_2,z_2) = (a_1,b_1+b_2,z_1z_2\Pi)
\end{equation*}
Changing representatives of the equivalence class gives
\begin{equation}
\begin{split}
(a_1+n,b_1+m,z_1 e^{-ma_1}) \times (a_2,b_2,z_2)
&= (a_1 +n ,b_1+b_2+m,z_1z_2 e^{-ma_1} \Pi) \\
&= (a_1, b_1+b_2, z_1z_2 \Pi) \\
(a_1,b_1,z_1) \times (a_2+n,b_2+m,z_2 e^{-ma_2})
&= (a_1,b_1 + b_2, z_1z_2 e^{m(a_1-a_2)}\Pi) \\
&= (a_1,b_1 + b_2, z_1z_2 \Pi)
\end{split}
\end{equation}
since $m(a_1 - a_2) \in \Z$. This means we may define a 
bundle gerbe product by $\Pi = 1$. The transition 
functions may now be calculated,
\begin{equation}
\begin{split} 
\sigma_{\alpha \beta} \sigma_{\beta\gamma} &= 
(\log_\alpha f, \log g_{\alpha \beta} , 1) 
\times (\log_\beta f , \log g_{\beta \gamma} , 1) \\
&= (\log_\alpha f, \log g_{\alpha \beta} + \log g_{\beta \gamma},
1) \\
&= (\log_\alpha f, \log g_{\alpha \gamma} - m_{\alpha\beta\gamma},
1) \\
&= (\log_\alpha f, \log g_{\alpha \gamma}, e^{m_{\alpha\beta\gamma}
\log_\alpha f}) \\
&= \sigma_{\alpha\gamma} f^{m_{\alpha \beta\gamma}}
\end{split}
\end{equation}
The connection is given by the 1-form $bda$ at 
$(a,b,z) \in \R \times \R \times S^1 /\sim$.

When the line bundle $L$ has connection $A$ the Deligne class of the cup
product is $(f^{m_{\alpha\beta\gamma}},\log (g_{\alpha\beta}) d\log f,
A_\alpha \wedge d\log(f))$. The connection is still the same and the 
curving is defined by $A \wedge dr$.

\subsubsection{The Triple Cup Product Bundle Gerbe}
The triple cup product bundle gerbe is defined by $f\cup g\cup h$ where
$f$, $g$ and $h$ are all $U(1)$ valued functions on $M$. 
The Deligne class of this cup product is 
\[
(h^{n_{\alpha\beta}m_{\beta\gamma}},n_{\alpha\beta}\log_\beta(g)
d\log h, \log_\alpha (f) d\log g \wedge d\log h) \in H^3(M,\Z(3)_D)
\]
where
$n_{\alpha\beta} = \log_\beta(f) - \log_\alpha(f)$ and
$m_{\beta\gamma} = \log_\gamma(g) - \log_\beta(g)$.
This bundle gerbe
could be represented geometrically by a combination of the cup 
product bundle and either of the cup product bundle
gerbes already discussed. There also a simpler representation which
we discuss here. Define a bundle gerbe by the following 
diagram:
\[
\begin{array}{ccc}
T & & \\
\downarrow & & \\
(\R\times\Z)^3 & \rightrightarrows & \R^3 \\
& & \downarrow \\
\mspace{40mu} M & \stackrel{(f,g,h)}{\rightarrow} & S^1 \times S^1 \times S^1
\end{array}
\]
where $T$ is the trivial bundle. We define a bundle 
gerbe product on $T$ by
\begin{equation}
\begin{split}
(r_1,n_1,s_1,m_1,t_1,&k_1,z_1) \times
(r_2,n_2,s_2,m_2,t_2,k_2,z_2) \\ &=
(r_1,n_1 + n_2, s_1, m_1+m_2, t_1, k_1 + k_2,z_1z_2e^{t_1 n_1 m_2})
\end{split}
\end{equation}
where for $i = 1,2$, $r_i,s_i,t_i \in \R$, $n_i,m_i,k_i \in \Z$ and
$z_i \in S^1$ for $i = 1,2$. The sections $\sigma_{\alpha\beta}$ may
be defined by
\begin{equation}
\sigma_{\alpha\beta} = (\log_\alpha (f), n_{\alpha\beta},
\log_\alpha(g), m_{\alpha\beta}, \log_\alpha (h), k_{\alpha\beta},1)
\end{equation}
Using $n_{\alpha\beta} + n_{\beta\gamma} = n_{\alpha\gamma}$ and
similar results for $m$ and $k$ we calculate
the product $\sigma_{\alpha\beta}\sigma_{\beta\gamma}$,
\begin{equation}
(\log_\alpha (f), n_{\alpha\gamma},\log_\alpha (g), m_{\alpha\gamma},
\log_\alpha(h), k_{\alpha\gamma}, e^{\log_\alpha (h) n_{\alpha\beta}
m_{\alpha\beta}}) = \sigma_{\alpha\gamma} h^{n_{\alpha\beta}
m_{\alpha\beta}}
\end{equation}
giving the required transition functions. Observe that this bundle 
has been constructed using a similar method to the canonical bundle 
associated with a Deligne class which was described in the 
proof of proposition \ref{Dclass2}. 
The connection at the point $(r,n,s,m,t,k,z) \in T$ is $n s dt$ and the 
curving at $(r,s,t) \in \R^3$ is $rds\wedge dt$.
\end{section}

%% file: chapter4.tex
\chapter{Other Geometric Realisations of Deligne Cohomology}

In this chapter we consider the bundle gerbe hierarchy of 
geometric realisations of Deligne cohomology. First we review what we
have considered so far, then we consider the relationship of bundle gerbes
with gerbes to clarify this picture of the hierarchy. We then extend the
hierarchy by considering bundle 2-gerbes. Following this we use 
$\Z$-bundle 0-gerbes to complete our catalogue of realisations and to
motivate a comparison with the theory of $B^pS^1$-bundles.

\begin{section}{The Bundle Gerbe Hierarchy}
We summarise the results on geometric realisations
of Deligne cohomology with a table:
\newline
\begin{table}[h] \label{t1}
\begin{center}
\caption{Low Dimensional Realisations of Deligne Cohomology}
\begin{tabular}{ll}
\\
\underline{Deligne Cohomology Group} & 
\underline{Geometric Realisation} \\
$H^0(M,\un{U(1)})$ & $U(1)$-functions \\
$H^0(M,\D^p), p > 0$ & constant $U(1)$-functions \\
$H^1(M,\un{U(1)})$ & $U(1)$-bundles \\
& $U(1)$-bundle 0-gerbes\\
$H^1(M,\D^1)$ & $U(1)$-bundles with connection \\
& $U(1)$-bundle 0-gerbes with connection \\
$H^1(M,\D^p), p > 1$ & flat $U(1)$-bundles\\
& flat $U(1)$-bundle 0-gerbes \\
$H^2(M,\un{U(1)})$ & $U(1)$-bundle gerbes \\
$H^2(M,\D^1)$ & $U(1)$-bundle gerbes with connection \\
$H^2(M,D^2)$ & $U(1)$-bundle gerbes with connection and curving \\
$H^3(M,\D^p), p>2$ & flat $U(1)$-bundle gerbes \\
\end{tabular}
\end{center}
\end{table}
\newline
It is to be understood that the right hand column of the above
table refers to equivalence classes of geometric object, that is,
isomorphism classes in the case of bundles and stable isomorphism classes 
in the case of bundle 0-gerbes and bundle gerbes.
\newpage
Recall that bundle 0-gerbes are defined as functions on $Y^{[2]}$ and
bundle gerbes are defined as bundles over $Y^{[2]}$. Since there
is an equivalence between bundles and bundle 0-gerbes we could
also think of bundle gerbes as bundle 0-gerbes over $Y^{[2]}$.
Thus there is a hierarchy
\newline 
\begin{center}
\begin{tabular}{c}
functions \\
$\downarrow$ \\
bundle 0-gerbes \\
$\downarrow$ \\
bundle gerbes 
\end{tabular} 
\end{center} 
where each object is built out of the preceding object and a 
submersion $Y \rightarrow M$.

Note also the similarity in the definitions of bundle 0-gerbe 
connections and bundle gerbe curvings. In both cases we have
a differential form on $Y$ satisfying the condition that applying 
$\delta$ gives the curvature of the object over $Y^{[2]}$. 
\end{section}

\begin{section}{Gerbes}
Gerbes are the most well known geometric realisation of $H^3(M,\Z(p)_D)$.
We shall review some relevant results about gerbes, 
for a detailed account
see \cite{bry}. It will suffice for us to think of a gerbe as
a sheaf of groupoids. Isomorphism classes of gerbes are represented by
classes in $H^3(M,\Z(1)_D)$. To each bundle gerbe there is an 
associated gerbe and equivalence classes of gerbes are in bijective
correspondence with stable isomorphism classes of bundle gerbes
\cite{must}. It is possible to define certain differential geometric 
structures on gerbes which are called a connective structure and
a choice of curving. Under the bijective correspondence these are 
equivalent to a connection 
and curving on a bundle gerbe. Gerbes with connective structure are
classified by $H^3(M,\Z(2)_D)$ and gerbes with connective
structure and a choice of curving are classified by $H^3(M,\Z(3)_D)$.\\

We wish to demonstrate that gerbes are to bundle gerbes what
bundles are to bundle 0-gerbes and thus remove the 
ambiguity in our bundle gerbe hierarchy. We briefly review the construction 
of a gerbe from a bundle gerbe in \cite{must}. Since a  
gerbe is a sheaf of groupoids we need to define a category over
each open set. The objects corresponding to $U \subset M$ are
bundle gerbe trivialisations over $U$. The morphisms
are morphisms of bundle gerbe trivialisations 
over $Y_U = \pi^{-1}(U)$. A gerbe is 
a bundle of groupoids in the sense that over each $m$ we have
the fibre of a bundle gerbe which is a groupoid with objects 
defined by trivialisations.

Consider a bundle 0-gerbe $(g,Y,M)$. Over $U \subset M$ there 
exists a trivialisation $h : \pi^{-1} (U) \rightarrow S^1$.
Given $h'$, a second trivialisation over $U$, there exists
a function $q_U: U \rightarrow S^1$ such that $h'=h\cdot
\pi^*q_U$.
We may define a bundle with trivialisations over $U$ given 
by $\sigma_U(x) = 
(x,h(s_u(x))$. The transition functions will be identical
to those of the original bundle 0-gerbes. If we replace
$h$ with $h'$ it is clearly seen that on overlaps
$q_U^{-1}q_V = 1$ and so we have a function
$q : M \rightarrow S^1$ which defines an automorphism of bundles. 
The fibre of the bundle may be considered as made up of bundle
0-gerbe trivialisations with any two differing by $q(m) \in S^1$.

This analysis suggests a refinement of the bundle gerbe hierarchy
described above.
\newline
\begin{center}
\begin{tabular}{ccccc}
& & functions & & \\
& $\swarrow$ & & $\searrow$ & \\
bundle 0-gerbes & & & & bundles \\
$\downarrow$ & & & & $\downarrow$ \\
bundle gerbes & & & & gerbes   
\end{tabular}
\end{center}
There is no simple diagrammatic representation of gerbes such as
we have for bundle gerbes.
\end{section}

\begin{section}{Bundle 2-Gerbes}\label{bun2ger}

We would like to construct a geometric realisation of $H^3(M,\D^p)$. 
This leads to the notion of a bundle 2-gerbe which has been developed
in \cite{ste}. We use a slightly different approach which is 
more suited to the bundle gerbe hierarchy. We shall require bundle 
2-gerbes for some applications in Chapter 8.

First we must deal with a matter of notation. Consider
the fibre product spaces $X^{[2]}$ and $X^{[3]}$ associated with
a submersion $X \rightarrow M$. We may define three projection maps 
$\pi_i : X^{[3]} \rightarrow X^{[2]}, i \in \{1,2,3\}$ by omission of the 
$i$th component of $X^{[3]}$. For example $\pi_1(x_1,x_2,x_3) = (x_2,x_3)$.
If there is a bundle 0-gerbe or bundle gerbe $P$ over $X^{[2]}$ then this 
may be pulled back by each of these projections. We will use the 
notation 
\begin{eqnarray*}
P_{12} &=& \pi_3^{-1}P \\
P_{13} &=& \pi_2^{-1}P \\
P_{23} &=& \pi_1^{-1}P \\
\end{eqnarray*}
Using this notation the bundle gerbe product can be written as 
a stable morphism of bundle 0-gerbes 
\begin{equation}\label{bgpr}
P_{12} \otimes P_{23} \isom P_{13}.
\end{equation}
This notation extends to projections
$X^{[p+1]} \rightarrow X^{[p]}$ for any positive integer $p$.

We now 
examine what happens if we replace the bundle 0-gerbes in
\eqref{bgpr} by bundle gerbes.
A choice of such a stable isomorphism is equivalent to a choice of
bundle gerbe trivialisation such that there is a bundle gerbe
isomorphism 

\begin{equation}\label{J123}
P_{12} \otimes P_{23} = P_{13} \otimes \delta(J_{123})
\end{equation}
The collection of trivialisations $J_{123}$ will be referred to as the 
{\it bundle 2-gerbe product}, $J$.
Observe that over $X^{[4]}$
\begin{equation}\label{PJJ}
P_{12} \otimes P_{23} \otimes P_{34} =
P_{14} \otimes \delta(J_{123}) \otimes \delta(J_{134}) 
= P_{14} \otimes \delta (J_{124}) \otimes \delta(J_{234}).
\end{equation}
We would like to write $\delta(J_{123}) \otimes \delta(J_{134}) 
= \delta( J_{123} \otimes J_{134})$ but this is not true
since the symbol $\otimes$ represents the bundle 0-gerbe
product which is a contracted tensor product of bundle 0-gerbes which have 
the same base space. Instead, given bundle 0-gerbes $(L,X)$ and
$(J,Y)$, with projections $X \rightarrow M$ and 
$Y \rightarrow M$  we define a product $(L \otimes_\delta J, 
X \times_M Y)$ as the bundle 0-gerbe with fibre over $(x,y)$ given
by $L_{x} \otimes J_{y}$. It is easy to show that given trivial
bundle gerbes $(\delta(L),X,M)$ and $(\delta(J),Y,M)$ there 
is an isomorphism $\delta(L) \otimes \delta(J) = \delta(
L \otimes_\delta J)$. We shall refer to this product as 
the {\it trivialisation product}.

We can now express \eqref{PJJ} in terms of trivialisation products as
\begin{equation}
P_{12} \otimes P_{23} \otimes P_{34} =
P_{14} \otimes \delta(J_{123} \otimes_\delta J_{134}) 
= P_{14} \otimes \delta (J_{124} \otimes_\delta J_{234}).
\end{equation}
Since we have two trivialisations of the same bundle gerbe there exists
a bundle 0-gerbe $A_{1234}$ on 
$X^{[4]}$,
called the {\it associator bundle 0-gerbe},
satisfying
\[
\pi^{-1}A_{1234} \otimes  (J_{123} \otimes_\delta J_{134}) =
(J_{124} \otimes_\delta J_{234})
\]
There is a technical point to be dealt with here. Up to this point we have
not needed to know anything about the base spaces for the trivialisations 
when considered as bundle 0-gerbes. For the formula above to make sense we 
need the bundle 0-gerbes on both sides to have the same base. There is no
reason for this to be true in general, however since we are really only
interested in $A_{1234}$ on $X^{[4]}$ we can get around this problem easily.
We just take the fibre product over $X^{[4]}$ of the base spaces from each
side and assume that we are actually dealing with the pullbacks to this 
product by the appropriate projection maps. The resultant bundle 0-gerbes
still define trivialisations and $A_{1234}$ is well defined. Throughout the
rest of our definition of bundle 2-gerbes we shall assume this construction
is used and will not specify base spaces for trivialisations.
  
Now suppose  
that there is a trivialisation $a_{1234}$ of $A_{1234}$ which
we call the {\it associator function}. Recall that if $A_{1234}$ is 
a bundle 0-gerbe over $X^{[4]}$,
\[
\begin{array}{ccc}
& & S^1 \\
& \stackrel{g}{\nearrow} & \\
A_{1234}^{[2]} & \rightrightarrows & A_{1234} \\
& & \downarrow \\
& & X^{[4]} 
\end{array}
\]
then $a_{1234}$ is a function $A_{1234} \rightarrow S^1$ satisfying 
$\delta(a_{1234}) = g$.
Furthermore we require that $a_{1234}$ satisfies a coherency condition
over $X^{[5]}$. Consider the series of bundle 0-gerbe isomorphisms given 
by each of the embeddings $X^{[4]} \rightarrow X^{[5]}$
\begin{eqnarray}
\pi^{-1}A_{1234} \otimes (J_{123} \otimes_\delta J_{134})
&=& (J_{124} \otimes_\delta J_{234})  \label{AJ1}\\
\pi^{-1}A_{1235} \otimes (J_{123} \otimes_\delta J_{135})
&=& (J_{125} \otimes_\delta J_{235})  \label{AJ2}\\
\pi^{-1}A_{1245} \otimes (J_{124} \otimes_\delta J_{145})
&=& (J_{125} \otimes_\delta J_{245})  \label{AJ3}\\
\pi^{-1}A_{1345} \otimes (J_{134} \otimes_\delta J_{145})
&=& (J_{135} \otimes_\delta J_{345})  \label{AJ4}\\
\pi^{-1}A_{2345} \otimes (J_{234} \otimes_\delta J_{245})
&=& (J_{235} \otimes_\delta J_{345})  \label{AJ5}
\end{eqnarray}
Consider the trivialisation product of
\eqref{AJ2} and \eqref{AJ4}
\[
(\pi^{-1}A_{1235} \otimes (J_{123} \otimes_\delta J_{135})) \otimes_{\delta}
(\pi^{-1}A_{1345} \otimes (J_{124} \otimes_\delta J_{145})).
\] 
This is isomorphic to $(J_{125} \otimes_\delta J_{235}) \otimes_{\delta}
(J_{135} \otimes_\delta J_{345})$ 
\[
J_{125} \otimes_\delta J_{235} \otimes_\delta J_{135} 
\otimes_{\delta} J_{345}
\]
and using \eqref{AJ5} this is isomorphic
to
\[
J_{125} \otimes_{\delta} J_{135} \otimes_{\delta}
(\pi^{-1}A_{2345} \otimes (J_{234} \otimes_\delta J_{245}))
\]
If we continue this process using the remaining isomorphisms
\eqref{AJ1} and \eqref{AJ3} we end up with 
\[
\pi^{-1}A_{2345} \otimes_{\delta} \pi^{-1}A_{1245} \otimes_{\delta} 
\pi^{-1}A_{1234} \otimes_{\delta} (J_{123} \otimes_\delta J_{135})
\otimes_{\delta} (J_{134} \otimes_\delta J_{145}).
\]
Using the trivialisations of the associator bundle 0-gerbes
this implies that there exists a function, $f_{12345}$, on $X^{[5]}$
such that
\[
a_{1234} \otimes a_{1245} 
\otimes a_{2345} = a_{1235} \otimes a_{1345} \otimes \pi^{-1}f_{12345}
\]
We call $f_{12345}$ the {\it coherency function}.
\\
We now return to our definition of a higher bundle gerbe.
Consideration of the bundle gerbe hierarchy leads to the following
\begin{definition}\cite{ste}
Let $X \rightarrow M$ be a submersion. Let
$(P,Y,X^{[2]})$ be a bundle gerbe. Then the quadruple $(P,Y,X,M)$ is
a {\it bundle 2-gerbe} if there is a bundle gerbe stable isomorphism
\[
P_{12} \otimes P_{23} \isom P_{13}.
\]
such that the corresponding associator bundle 0-gerbe is trivial, and 
the coherency function is identically 1. These last two conditions are
called the associator trivialisation and the coherency condition respectively.
The stable isomorphism together with the associator trivialisation and
the coherency condition is called the bundle 2-gerbe product.
\end{definition}
This definition corresponds to Stevenson's definition of a 
stable bundle 2-gerbe \cite{ste}.
The bundle 2-gerbe $(P,Y,X,M)$ may be represented diagrammatically
in the following way:
\[
\begin{array}{ccccc}
P & & & & \\
\downarrow & & & & \\
Y^{[2]} & \rightrightarrows & Y & & \\
& & \downarrow & & \\
& & X^{[2]} & \rightrightarrows & X \\
& & & & \downarrow \\
& & & & M 
\end{array}
\]
We may define pullbacks, products, duals, morphisms and trivial
bundle 2-gerbes by analogy with the definitions for 
bundle 0-gerbes and bundle gerbes.

By following the lower dimensional cases we may define 
connections and curvings by choosing a bundle gerbe 
connection and curving on $(P,Y,X^{[2]})$ and a 3-form
$\nu \in \Omega^3(X)$ such that $\delta(\nu) = \omega$ where
$\omega$ is the 3-curvature of $(P,Y,X^{[2]})$. 
The 3-form is called the 3-curving, the curving on $(P,Y,X^{[2]})$ is 
called the 2-curving and the connection on $(P,Y,X^{[2]})$ is also
referred to as the connection on the bundle 2-gerbe. We shall sometimes
refer to a bundle 2-gerbe with connection, 2-curving and 3-curving simply
as a bundle 2-gerbe with curvings. There is essentially no difference 
between the connection and the curvings, a connection could be referred to
as a 1-curving however we shall continue to use the familiar terminology.
The 
{\it 4-curvature} of a bundle 2-gerbe is a 4-form 
$\Theta \in \Omega^4(M)$ satisfying
$\pi^* \Theta = d\nu$. A bundle 2-gerbe is {\it flat} if the
curvature is zero. For bundle 2-gerbes with curvings we also
require that the bundle 2-gerbe product $J$ and the associator trivialisation
 $a$ 
are $D$-trivialisations.

\begin{proposition}
Associated to every bundle 2-gerbe with curvings is a 
class in \\$H^3(M,\D^3)$.
\end{proposition}
\begin{proof}
See \cite{ste}.
\end{proof}

All constructions and operations involving bundle gerbes can be 
carried out for bundle 2-gerbes. We describe some which are
relevant here (see \cite{ste} for more detail).

If there exists a bundle gerbe $R \Rightarrow X$ such that
      there is a bundle gerbe morphism $\delta(R) \isom P$ over
      $X^{[2]}$ which is compatible with the bundle 2-gerbe 
      product and associator function 
       then the bundle 2-gerbe is called {\it trivial}.

The set of $D$-stable isomorphism classes of  
      bundle 2-gerbes with 2-curving
      form a group under the tensor product, which is defined by analogy
      with the bundle gerbe case.

A flat bundle 2-gerbe has a flat holonomy which is 
a class in $H^3(M,U(1))$.

A $D$-trivial bundle 2-gerbe with curvings has a trivialisation 
with connection
and curving such that the curvature 3-form of the trivialisation
is equal to the 3-curving of the bundle 2-gerbe.

A trivial bundle 2-gerbe with curvings, $(P,Y,X,M;A,\eta,\nu)$ has a 
$D$-obstruction 3-form $\chi$. If $\chi \in \Omega^3_0(M)$ then for any
bundle gerbe $(Q,X;B,\zeta)$ which trivialises $(P;A,\eta)$ there exists
a bundle gerbe $(R,M;C,\mu)$ such that $Q \otimes \pi^{-1} R$ is 
a $D$-trivialisation of $(P;A,\eta,\nu)$.

\begin{example}\cite{ste}
There exists a tautological bundle 2-gerbe associated with any 
closed, $2\pi$-integral 4-form $\Theta$ on a 3-connected base $M$. This is defined
by the following diagram
\[
\begin{array}{ccc}
Q[\omega] & & \\
\Downarrow & & \\
\L_0M & \rightrightarrows & \P_0M \\
& & \downarrow \\
& & M
\end{array}
\]
where $Q[\omega] \Rightarrow \L_0$ is the tautological 
bundle gerbe with curvature 3-form $\omega = \int_{S^1} ev^* \Theta$. The
3-curving is $\int_I ev^* \Theta$ and the product is defined by
composition of paths in $\Omega_0M$ as in the bundle gerbe case.
We need $M$ to be 3-connected so that $\Omega_0(M)$ is 2-connected and
hence the tautological bundle gerbe is well defined. 
\end{example}

\begin{example}\cite{ste}
The bundle 2-gerbe associated with a principal $G$-bundle
$P_G \rightarrow M$, where $G$ is a compact, simply connected, simple
Lie group, is defined by the following diagram:
\[
\begin{array}{ccc}
& & Q[\omega] \\
& & \Downarrow \\
& & G \\
& \stackrel{g}{\nearrow} & \\
P_G^{[2]} & \rightrightarrows & P_G \\
& & \downarrow \\
& & M
\end{array}
\]
where $g$ is defined as in the bundle gerbe case (see \S\ref{lifting}) and
$Q \Rightarrow G$ is the tautological bundle gerbe with Dixmier-Douady class
given by the
canonical generator of $H^3(G,\Z) = \Z$. The main 
result regarding such bundle 2-gerbes is that the \v{C}ech 4-class
is equal to the first Pontryagin class of the bundle $P$.
\end{example}

Classes in $H^3(M,\D^3)$ may also be represented by
      { \it 2-gerbes }. We shall not give a proper definition
      of these here since it is quite complicated and is not
      of direct relevance. Full definitions may be found in \cite{brmc1}
      or \cite{ste}. Essentially if we think of a gerbe as a sheaf of 
groupoids then a 2-gerbe is a sheaf of 2-groupoids. These are defined 
in terms of higher categories. A 2-category consists of objects, 1-arrows
(morphisms) and 2-arrows (transformations between morphisms) with a
number of axioms relating to composition, associativity and identity.
A 2-groupoid is a 2-category with invertible 2-arrows and 1-arrows
which are invertible up to 2-arrows. A 2-gerbe is a sheaf of 
2-groupoids with a number of gluing and descent axioms. 
Given a bundle 2-gerbe the objects of the 2-groupoid 
associated with an open set are defined by the trivialisations of
the bundle 2-gerbe over the set. The 1-arrows are morphisms between the
trivialisations. Since trivialisations of bundle 2-gerbes are bundle 
gerbes over some space their morphisms may be thought of as bundle
gerbe trivialisations. The 2-arrows are morphisms of these. It is 
a result of Stevenson \cite{ste} that this construction 
gives a 2-gerbe with the same class in $H^4(M,\Z)$ as the original
bundle 2-gerbe.
There are differential geometric structures on 2-gerbes which
may be used to obtain a full Deligne class in $H^3(M,\D^3)$
\cite{bry2}. We do not have a direct relationship between these and 
bundle 2-gerbes with connection and curvings however we shall
see that 
such a relationship may
be established indirectly via the Deligne class.

\begin{subsubsection}{Bundle 2-gerbes and Deligne Cohomology}

We prove here the main result on bundle 2-gerbes 
which places them in the bundle gerbe hierarchy.

\begin{proposition}\label{delb2g}
The group of $D$-stable isomorphism classes of bundle 
2-gerbes with 2-curving is isomorphic to $H^3(M,\D^3)$.
\end{proposition}
This extends the results of Stevenson \cite{ste} which state that
a bundle 2-gerbe with connection and curving defines a Deligne 
class and that a trivial bundle 2-gerbe has a trivial \v{C}ech class.
\begin{proof}
We shall first describe the element of $H^3(M,\D^3)$ representing a 
bundle 2-gerbe. Suppose we have a bundle 2-gerbe $(P,Y,X,M)$,
with connection, $A$, 2-curving $\eta$ and 3-curving $\nu$.
On $U_\alpha \subset M$ define
\[
\nu_\alpha = s_\alpha^* \nu.
\]
Now consider the family of bundle gerbes obtained by pulling back
the bundle gerbe $(P,Y,X^{[2]})$ with 
\[
(s_\alpha,s_\beta) : U_{\alpha\beta} \rightarrow X^{[2]}.
\] 
Denote these pullback bundle gerbes by $(P_{\alpha\beta},
Y_{\alpha \beta},U_{\alpha \beta})$. They have induced connection and
curving $A_{\alpha \beta}$ and $\eta_{\alpha \beta}$ respectively. 
Since each base space
$U_{\alpha \beta}$ is contractible, each bundle gerbe $P_{\alpha
\beta}$ is trivial. Thus associated with each $P_{\alpha \beta}$ is
a $D$-obstruction 2-form $\chi_{\alpha \beta}$. 
Recall that if we 
choose trivialisations with connections, $L_{\alpha \beta} 
\rightarrow Y_{\alpha \beta}$ and
let $F_{L_{\alpha \beta}}$ denote
the curvature of the connection on $L_{\alpha \beta}$, then 
$\chi_{\alpha \beta}$ is defined by
$F_{L_{\alpha \beta}} = \eta_{\alpha \beta} - \pi^* \chi_{\alpha \beta}$.
Also
recall from \S \ref{trivial} that the
$D$-obstruction form satisfies 
$d\chi_{\alpha \beta} = \omega_{\alpha \beta}$ where  
$\omega_{\alpha \beta}$ is the
3-curvature of $P_{\alpha \beta}$. The 3-curving $\nu$ is defined
such that $\delta(\nu) = \omega$ and it follows that
\[
d\chi_{\alpha \beta} = \omega_{\alpha \beta} = \nu_\beta - \nu_\alpha
\]
Consider the isomorphism over $U_{\alpha \beta \gamma}$,
\[
P_{\alpha \beta} \otimes P_{\beta \gamma} =
P_{\alpha \gamma} \otimes D(J_{\alpha \beta \gamma})
\]
Using the trivialisations $L_{ij}$ this becomes
\[
\delta(L_{\alpha \beta}) \otimes \delta(L_{\beta \gamma}) =
\delta(L_{\alpha \gamma}) \otimes D(J_{\alpha \beta \gamma}).
\]
The $D$ obstruction form for the left hand side is 
$\chi_{\alpha \beta} + \chi_{\beta \gamma}$, and for the 
right hand side is $\chi_{\alpha \gamma}$. We have assumed without
loss of generality that the $D$-obstruction form of $D(J_{\alpha\beta\gamma})$
is zero rather than a general closed, $2\pi$-integral form. If this were not the
case then we could redefine it as described in \S \ref{trivial}. 
Comparing the curvatures of both sides we have 
\begin{equation}
d\chi_{\alpha \beta} + d\chi_{\beta \gamma} = 
d\chi_{\alpha \gamma} = \nu _\gamma - \nu_\alpha
\end{equation}
Also in terms of the definition of $\chi$ we have 
\begin{equation}
\begin{split}
\chi_{\alpha \beta} + \chi_{\beta \gamma} - \chi_{\alpha \gamma}
&= \eta_{\alpha\beta} + \eta_{\beta \gamma} - \eta_{\alpha \gamma} 
 - F_{L_{\alpha\beta}}- F_{L_{\beta \gamma}} 
+ F_{L_{\alpha\gamma}} \\
&= F_{J_{\alpha\beta\gamma}} - F_{L_{\alpha\beta}}- F_{L_{\beta \gamma}} 
+ F_{L_{\alpha\gamma}}
\end{split}
\end{equation}
so this difference is closed and $2\pi$-integral. Hence we may apply 
proposition \ref{dobs} and there
exists a bundle 0-gerbe $K_{\alpha \beta \gamma}$ with curvature
$-\chi_{\alpha \beta} - \chi_{\beta \gamma}
+ \chi_{\alpha \gamma}$ such that
\begin{equation}\label{defK}
L_{\alpha \beta} \otimes_\delta L_{\beta \gamma} = L_{\alpha \gamma}
\otimes_\delta J_{\alpha \beta \gamma} \otimes_\delta 
\pi^{-1}K_{\alpha \beta
\gamma}
\end{equation}
Since $K_{\alpha \beta \gamma}$ is 
a bundle 0-gerbe on $U_{\alpha \beta \gamma}$ it is trivial and has a 
$D$-obstruction 1-form $\kappa_{\alpha \beta \gamma}$ which satisfies
\[
d \kappa_{\alpha \beta \gamma} = -\chi_{\alpha \beta} - 
\chi_{\beta \gamma} + \chi_{\alpha \gamma}.
\]
Using
\eqref{defK} we get
\begin{eqnarray*}
L_{\alpha \beta} \otimes_\delta L_{\beta \gamma} \otimes_\delta 
L_{\gamma \delta} &=& L_{\alpha \delta} \otimes_\delta J_{\beta \gamma \delta}
\otimes_\delta \pi^{-1}K_{\beta \gamma \delta} \otimes_\delta 
J_{\alpha \beta \delta}
\otimes_\delta \pi^{-1}K_{\alpha \beta \delta} \\
&=& L_{\alpha \delta} \otimes_\delta J_{\alpha \gamma \delta} 
\otimes_\delta \pi^{-1}
K_{\alpha \gamma \delta} \otimes_\delta J_{\alpha \beta \gamma} 
\otimes_\delta
\pi^{-1} K_{\alpha \beta \gamma}.
\end{eqnarray*}
Furthermore, using the definition of the associator bundle we have
\begin{multline}
L_{\alpha \delta} \otimes_\delta J_{\alpha \gamma \delta} 
\otimes_\delta 
J_{\alpha \beta \gamma} \otimes_\delta \pi^{-1} A_{\alpha \beta \gamma \delta}
\otimes_\delta \pi^{-1}K_{\beta \gamma \delta} \otimes_\delta 
\pi^{-1}K_{\alpha \beta \delta}
= \\
L_{\alpha \delta} \otimes_\delta J_{\alpha \gamma \delta} 
\otimes_\delta 
J_{\alpha \beta \gamma} \otimes_\delta \pi^{-1}K_{\alpha \gamma \delta} 
\otimes_\delta 
\pi^{-1} K_{\alpha \beta \gamma}.
\end{multline}
Thus over $U_{\alpha \beta \gamma \delta}$
\begin{equation}\label{AK}
A_{\alpha \beta \gamma \delta} \otimes K_{\beta \gamma \delta} 
\otimes K_{\alpha \beta \delta}
=  K_{\alpha \gamma \delta} \otimes 
K_{\alpha \beta \gamma}
\end{equation}
Let $h_{\alpha\beta\gamma}$ be the trivialisation of 
$K_{\alpha\beta\gamma}$. Using these trivialisations
together with $A_{\alpha \beta \gamma \delta} = D (a_{\alpha
\beta \gamma \delta})$, equation \eqref{AK} becomes
\begin{equation}
D(a_{\alpha \beta \gamma \delta}) \otimes \delta(h_{\beta \gamma \delta})
\otimes \delta(h_{\alpha \beta \delta})
= \delta(h_{\alpha \gamma \delta}) \otimes \delta(h_{\alpha \beta \gamma}) 
\end{equation}
The $D$-obstruction of the right hand side is $\kappa_{\alpha \gamma \delta} 
+ \kappa_{\alpha \beta \gamma}$ and for the left hand side is
$\kappa_{\beta \gamma \delta} + \kappa_{\alpha \beta \delta}$.
The curvature of each of these is $-\chi_{\alpha \beta} - \chi_{\beta
\gamma} - \chi_{\gamma \delta} + \chi_{\alpha \delta}$, and by 
a version of proposition \ref{dobs} for bundles there
is a function $g_{\alpha \beta \gamma \delta}$ on $U_{\alpha \beta \gamma
\delta}$ which satisfies 
\begin{equation}
d\log(g_{\alpha \beta \gamma \delta}) = -\kappa_{\alpha \beta \gamma} 
+ \kappa_{\alpha \beta \delta} - \kappa_{\alpha \gamma \delta}
+ \kappa_{\beta \gamma \delta}
\end{equation}
and
\begin{equation}\label{def4g}
a_{\alpha \beta \gamma \delta}
\otimes_\delta h_{\beta \gamma \delta}
\otimes_\delta h_{\alpha \beta \delta} = \pi^*g_{\alpha \beta \gamma \delta}
 \otimes_\delta 
h_{\alpha \gamma \delta} \otimes_\delta h_{\alpha \beta \gamma}
\end{equation}
on $A_{\alpha \beta \gamma \delta} \times_M K_{\beta \gamma \delta}
\times_M K_{\alpha \beta \delta} \times_M K_{\alpha \gamma \delta}
\times_M K_{\alpha \beta \gamma}$ as a fibre product of total spaces
of the respective bundle 0-gerbes. Furthermore it may be shown that the
coherency condition implies that
\begin{equation}
g_{\beta \gamma \delta \epsilon} \cdot g_{\alpha \beta \delta \epsilon}
\cdot g_{\alpha \beta \gamma \delta} = g_{\alpha \gamma \delta \epsilon}
\cdot g_{\alpha \beta \gamma \epsilon}.
\end{equation}
on $U_{\alpha \beta \gamma \delta \epsilon}$.
\\
The Deligne class is given by 
\[
(g_{\alpha \beta \gamma \delta},\kappa_{\alpha \beta \gamma},
\chi_{\alpha \beta}, \nu_\alpha).
\]
In trivialisations $h$ could be replaced by sections since these 
notions are equivalent for bundles. Similarly the $D$-obstruction form
$\kappa_{\alpha\beta\gamma}$ could be replaced by the pull back of the 
connection on $K_{\alpha\beta\gamma}$ by this section. We have used
the more general terms above to highlight the role of the hierarchy, 
and because we believe that the language of trivialisations and 
$D$-obstructions has potential use in dealing with higher objects.

It appears that the connection form $A$ was not used in this derivation,
while it does not appear explicitly it is involved. When we trivialise
$P_{\alpha\beta}$ information about $A$ is carried in the connections on
the trivialisations $L_{\alpha\beta}$. The local one forms $\kappa$ depend
on both the connection of the bundle gerbe and the connections on the
bundles $J_{123}$.

Suppose we have a bundle 2-gerbe $(P,Y,X,M;A,\eta,\nu)$ represented
by $(g_{\alpha \beta \gamma \delta},\kappa_{\alpha \beta \gamma},
\chi_{\alpha \beta},\nu_\alpha)$.
To prove that this gives an isomorphism we first show that 
the Deligne class is independent of all choices in the construction.
Suppose we were to replace $h_{\alpha \beta \gamma}$ by $\tilde{h}_
{\alpha\beta\gamma}$. These are both trivialisations of the bundle 
0-gerbe  
$K_{\alpha\beta\gamma} \rightarrow U_{\alpha\beta\gamma}$ 
so they differ by a function $p_{\alpha\beta\gamma}$. Comparing
the two versions of equation
\eqref{def4g} obtained from the two choices we have 
\begin{equation}
 a_{\alpha\beta\gamma\delta}
\otimes_\delta a^{-1}_{\alpha\beta\gamma\delta}
\otimes_\delta \pi^* p_{\beta \gamma \delta}
\otimes_\delta \pi^* p_{\alpha \beta \delta} = 
\pi^*\tilde{g}_{\alpha\beta\gamma\delta}  \otimes_\delta 
\pi^*g^{-1}_{\alpha\beta\gamma\delta}  \otimes_\delta
\pi^* p_{\alpha \gamma \delta} \otimes_\delta \pi^* p_{\alpha \beta \gamma}
\end{equation}  
and so we have
\begin{equation}\label{zz1}
\tilde{g}_{\alpha\beta\gamma\delta} = g_{\alpha\beta\gamma\delta}
\delta(p)_{\alpha\beta\gamma\delta}
\end{equation}
Recall that the local connections $\kappa_{\alpha\beta\gamma}$ were
defined as $D$-obstruction forms for the bundle 0-gerbes $K_{\alpha
\beta \gamma}$. Changing the choice of trivialisation changes this 
$D$-obstruction form by the 1-curvature of the function defined
by the difference of the two trivialisations so we have
\begin{equation}\label{zz2}
\tilde{\kappa}_{\alpha \beta \gamma} = \kappa_{\alpha \beta \gamma}
+ dp_{\alpha\beta\gamma}
\end{equation}
Together equations \eqref{zz1} and \eqref{zz2} change the Deligne 
class by a trivial cocycle.

Now suppose that we change $L_{\alpha \beta}$ to $\tilde{L}_{\alpha\beta}$.
These differ by a bundle with connection $T_{\alpha \beta} \rightarrow
U_{\alpha \beta}$. Comparing the two versions of equation \eqref{defK} 
we have  
\begin{equation}
\pi^{-1}T_{\alpha \beta} \otimes_\delta \pi^{-1}T_{\beta \gamma} = 
\pi^{-1}T_{\alpha \gamma} \otimes_\delta J_{\alpha \beta \gamma} 
\otimes_\delta \pi^{-1}\tilde{K}_{\alpha \beta
\gamma} \otimes_\delta J^*_{\alpha \beta \gamma} 
\otimes_\delta \pi^{-1}K^*_{\alpha \beta
\gamma}
\end{equation}
so the connection on $K_{\alpha\beta\gamma}$ changes by 
$B_{\alpha \beta} + B_{\beta \gamma} - B_{\alpha\gamma}$, where
$B_{\alpha \beta}$ is the connection on $T_{\alpha\beta}$ and
hence the $D$-obstruction form $\kappa_{\alpha\beta\gamma}$ is 
changed in the same way,
\begin{equation}\label{zz3}
\tilde{\kappa}_{\alpha\beta\gamma} = \kappa_{\alpha\beta\gamma}
+ B_{\alpha\beta} + B_{\beta\gamma} - B_{\alpha \gamma}
\end{equation}
Under the change from $L_{\alpha \beta}$ to
$\tilde{L}_{\alpha\beta}$ the $D$-obstruction forms $\chi_{\alpha\beta}$
will change by the curvature of $T_{\alpha\beta}$,
\begin{equation}\label{zz4}
\tilde{\chi}_{\alpha \beta} = \chi_{\alpha \beta} + dB_{\alpha\beta}
\end{equation}
Together equations \eqref{zz3} and \eqref{zz4} change the Deligne 
class by a trivial cocycle.

The final choice that we have made is of the sections $s_\alpha$. 
For each choice of section there is a trivialisation $R_\alpha$, 
a bundle gerbe 
over $\pi^{-1} (U_\alpha)$ which is defined by $R_{\alpha} = (1,
s_{\alpha})^{-1} P$, where $P$ is considered as a bundle gerbe on 
$Y^{[2]}$. Then $\pi^{-1} P_{\alpha \beta} = R_{\alpha}^* \otimes
R_\beta$ and a different choice of section, $\tilde{s}_\alpha$
defines a bundle gerbe $\xi_{\alpha}$ on $U_\alpha$ satisfying
$\tilde{R}_\alpha = R_{\alpha} \otimes \pi^* \xi_\alpha$. 
Thus a change of section changes $P_{\alpha\beta}$ by
\begin{equation}
\tilde{P}_{\alpha\beta} = P_{\alpha\beta} \otimes \xi^*_\alpha
\otimes \xi_\beta
\end{equation}
If $\xi_\alpha$ has curving $\mu_\alpha$ then the 2-forms $f_{\alpha
\beta}$ and hence the $D$-obstruction forms $\chi_{\alpha\beta}$ 
change by $\mu_\beta - \mu_\alpha$. The local 3-curvings may be thought
of as a $D$-obstruction form for the trivialisation $R_\alpha$ so they
change by $d\mu_\alpha$ and once again we have
a trivial contribution to the Deligne class.


Next we claim that this assignment of a Deligne class to
a bundle 2-gerbe is a homomorphism. Since the 3-curving 
of the tensor product of two bundle 2-gerbes is the sum of their 
respective 3-curvings, and the local 3-curvings are defined by pullback,
then it is clear that the local two curvings will be 
additive under tensor products. The local 2-curvings and connections are 
both defined as $D$-obstruction forms for a bundle gerbe and bundle 0-gerbe
respectively. Since $D$-obstructions are additive under tensor products
in both of these cases then so will the local 2-curvings and connections.
The transition functions may be thought of in similar terms as a 
$D$-obstruction defined in terms of two trivialisations of a bundle 0-gerbe
and hence the assignment of a Deligne class preserves the tensor 
product of bundle 2-gerbes with curvings. To see that this gives a 
homomorphism between equivalence classes we 
show that a D-trivial bundle 2-gerbe has a trivial 
Deligne class. Let $Q \rightarrow Y$ be a $D$-trivialisation
of the bundle 2-gerbe $(Q,X,Y,M)$. 
The 3-curvature 
$\omega_Q$ of $Q$ satisfies $\omega_Q = \nu$, where $\nu$ is the
3-curving of $P$. Using a section 
$s_\alpha : U_\alpha \rightarrow Y$ we pull back
$Q$ to a bundle gerbe $Q_\alpha$. This bundle gerbe must be 
trivial so let $q_\alpha$ be the $D$-obstruction form. This satisfies
\begin{equation}
dq_\alpha = s_\alpha^* \omega_Q = s_\alpha^* \nu = \nu_\alpha
\end{equation}
If $R_\alpha$ is a trivialisation of $Q_\alpha$ then
the isomorphism 
\begin{equation}\label{PQQ}
P_{\alpha\beta} = Q_\alpha^* \otimes Q_\beta
\end{equation}
induces 
an isomorphism $\delta(L_{\alpha\beta}) = \delta(R^*_\alpha) \otimes
\delta(R_\beta)$. This means we may define bundles $N_{\alpha \beta}$
which satisfy 
\begin{equation}\label{LRRN}
L_{\alpha\beta} = R^*_\alpha \otimes R_\beta \otimes 
\pi^{-1} N_{\alpha\beta}
\end{equation} 
We may now find an expression for the
$D$-obstruction form for $P_{\alpha\beta}$ which is defined by
$\pi^* \chi_{\alpha\beta} = f_{\alpha\beta} - F_{L_{\alpha\beta}}$. Since
$Q$ is a $D$-trivialisation then from equation \eqref{PQQ} we 
have $f_{\alpha\beta} = f_{Q_\beta} - f_{Q_\beta}$ where the 
terms on the right are the curvings induced on $Q_\alpha$ and
$Q_\beta$ from that of $Q$. Equation \eqref{LRRN} gives an
equation for the curvature of $L_{\alpha\beta}$, 
$F_{L_{\alpha\beta}} = F_{Q_\beta} - F_{Q_\alpha} + \pi^{-1} 
F_{N_{\alpha\beta}}$. Since $N_{\alpha\beta}$ is trivial it has
a $D$-obstruction form $n_{\alpha\beta}$ which satisfies 
$dn_{\alpha\beta} = F_{N_{\alpha\beta}}$. Putting these together 
we get 
\begin{equation}
\begin{split}
\pi^* \chi_{\alpha\beta} &= f_{Q_\beta} - f_{Q_\beta} - F_{Q_\beta} 
+ F_{Q_\alpha} - \pi^* dn_{\alpha\beta} \\
&= \pi^* q_\beta - \pi^* q_\alpha - \pi^* dn_{\alpha\beta}
\end{split} 
\end{equation}
and so 
\begin{equation}
\chi_{\alpha\beta} = q_\beta - q_\alpha - dn_{\alpha\beta}
\end{equation}

Next we need a $D$-obstruction form for $K_{\alpha\beta\gamma}$. First
observe that
we may express the bundle 
2-gerbe product in terms of $Q$, using 
$P_{\alpha\beta} \otimes P_{\beta\gamma} = 
P_{\alpha\gamma} \otimes D(J_{\alpha\beta\gamma})$ to get
$Q^*_\alpha \otimes Q_\beta \otimes Q^*_\beta \otimes Q_\alpha = Q^*_\alpha 
\otimes
Q_\gamma \otimes D(J_{\alpha\beta\gamma})$. Using the trivialisations
$R_\alpha$ we may express this as a difference of two trivialisation and
define a bundle $M_{\alpha\beta\gamma}$ such that
\begin{equation}
R^*_\alpha \otimes R_\beta \otimes R^*_\beta \otimes R_\gamma = R^*_\alpha 
\otimes
R_\gamma \otimes J_{\alpha\beta\gamma} \otimes \pi^{-1} M_{\alpha\beta\gamma}
\end{equation}
The combination of $R$ terms is actually a $D$-trivialisation since
the sum of the $D$-obstructions cancels so we may assume that 
$M_{\alpha\beta\gamma}$ is flat.

Substituting into equation \eqref{defK} gives
\begin{equation}\label{KNM}
K_{\alpha \beta\gamma} = 
N_{\alpha\beta} \otimes N_{\beta\gamma} \otimes N^*_{\alpha\beta} \otimes
M_{\alpha\beta\gamma}
\end{equation}
The $D$-obstruction form is defined by $\pi^* \kappa_{\alpha\beta\gamma}
= A_{K_{\alpha\beta\gamma}} - d\log h_{\alpha\beta\gamma}$.
A formula for the connection may be obtained from equation 
\eqref{KNM},
\begin{equation}
A_{K_{\alpha\beta\gamma}} = A_{N_{\alpha\beta}} + A_{N_{\beta\gamma}} 
- A_{N_{\alpha\beta}} +
A_{M_{\alpha\beta\gamma}}
\end{equation}
where the terms on the right are the connections on the respective bundles.
Let $\zeta_{\alpha\beta}$ be a trivialisation of $N_{\alpha\beta}$ and
let $\epsilon_{\alpha\beta\gamma}$ be a trivialisation of 
$M_{\alpha\beta\gamma}$. Then using \eqref{KNM} we may define 
two trivialisations and hence define functions $\rho_{\alpha\beta\gamma}$ 
on $U_{\alpha\beta\gamma}$ which satisfy
\begin{equation}\label{hzer}
h_{\alpha\beta\gamma} = \zeta_{\alpha\beta} \zeta_{\beta\gamma}
\zeta^{-1}_{\alpha\gamma} \epsilon_{\alpha\beta\gamma} \pi^* 
\rho_{\alpha\beta\gamma}
\end{equation}
Thus we have
\begin{equation}
d\log h_{\alpha\beta\gamma} = d\log \zeta_{\alpha\beta} 
+ d\log \zeta_{\beta\gamma} - d\log \zeta_{\alpha\gamma} 
+ d\log \epsilon_{\alpha\beta\gamma} - \pi^* 
d\log \rho_{\alpha\beta\gamma}
\end{equation}
Observe that $A_{N_{\alpha\beta}} - d\log \zeta_{\alpha\beta} = \pi^* 
n_{\alpha
\beta}$ by definition. Also, since $M_{\alpha\beta\gamma}$ is flat and
is defined on a contractible set then we may assume that we have a 
flat trivialisation so $A_{M_{\alpha\beta\gamma}} = d\log \epsilon_
{\alpha\beta\gamma}$. Combining all of this we have
\begin{equation}
\kappa_{\alpha\beta\gamma} = n_{\alpha\beta} + n_{\beta\gamma}
- n_{\alpha\gamma} + d\log \rho_{\alpha\beta\gamma}
\end{equation}
The final step follows the same pattern however the arguments 
will be simpler as we are dealing with functions. First use 
\eqref{AK} to compare $a$ with $\zeta$ and $\epsilon$. Substituting
this expression for $a$ and the expressions \eqref{hzer} for the $h$
terms into \eqref{def4g} will lead to the cancellation of all
$\zeta$ and $\epsilon$ terms leaving
\begin{equation}
g_{\alpha\beta\gamma\delta} = \delta(\rho)_{\alpha\beta\gamma\delta}
\end{equation}
thus we have a trivial Deligne class
\begin{equation}
(g_{\alpha\beta\gamma\delta},\kappa_{\alpha\beta\gamma},
\chi_{\alpha\beta},\nu_\alpha) =
(\delta(\rho)_{\alpha\beta\gamma\delta},\delta(n)_{\alpha\beta\gamma}
+ d\log \rho_{\alpha\beta\gamma}, \delta(q)_{\alpha\beta} - dn_{\alpha\beta},
dq_\alpha)
\end{equation}
and the assignment of a Deligne class to a $D$-stable isomorphism
class of bundle 2-gerbes is a homomorphism.

To show injectivity of this homomorphism we shall show that the 
Deligne class of a bundle 2-gerbe is trivial only if the bundle 2-gerbe
is $D$-trivial. The corresponding result for $\delta$-trivialisations and
\v{C}ech classes 
has already been given by Stevenson \cite{ste}. This means that if
a bundle 2-gerbe is not $\delta$-trivial then its Deligne class is 
not trivial, so we may assume without loss of generality that our 
bundle 2-gerbe is $\delta$-trivial, but not $D$-trivial. Furthermore by
standard arguments this trivialisation may be given connection and curving 
which are compatible with the trivialisation. In this 
case there exists a $D$-obstruction form $\zeta \in \Omega^3(M)$ 
which is defined by $\pi^*\zeta = \nu - \omega_Q$, where 
$\nu$ is 
the 3-curving of the bundle 2-gerbe and $\omega_Q$ is the curvature of the
trivialisation $Q$. By assumption the $D$-obstruction form is non-trivial,
this implies that $\zeta \notin \Omega^3_0(M)$.
To find the Deligne class of this bundle 2-gerbe the arguments above, 
where the Deligne class of a $D$-trivial bundle 2-gerbe was calculated, 
still apply with the exception of the local 3-curving, which is 
now given by $\nu_\alpha = dq_\alpha + \zeta$. Thus we have
\begin{equation}
\begin{split}
(g_{\alpha\beta\gamma\delta},\kappa_{\alpha\beta\gamma},
\chi_{\alpha\beta},\nu_\alpha) &=
(\delta(\rho)_{\alpha\beta\gamma\delta},\delta(n)_{\alpha\beta\gamma}
+ d\log \rho_{\alpha\beta\gamma}, \delta(q)_{\alpha\beta} - dn_{\alpha\beta},
dq_\alpha + \zeta) \\
&= D(\rho,n,q) + (1,0,0,\zeta) 
\end{split}
\end{equation}
Since $\zeta \notin \Omega^3_0(M)$ the class $(1,0,0,\zeta)$ is not $D$-exact,
so we have shown that a bundle 2-gerbe which is not $D$-trivial cannot 
have a trivial Deligne class.

Finally we need to show that there exists a bundle 2-gerbe with
connection and curvings which 
is represented by any given Deligne class $(g_{\alpha\beta
\gamma\delta},\kappa_{\alpha\beta\gamma},\chi_{\alpha\beta},
\nu_{\alpha})$. Define this bundle 2-gerbe by
\[
\begin{array}{ccccc}
T^* \otimes T & & T & & \\
\downarrow & & \downarrow & & \\
\amalg U_{ij} & \rightrightarrows & \amalg U_{ij} & & \\ 
& & \downarrow & & \\
& & \amalg U_{ij} & \rightrightarrows & \amalg U_i \\
& & & & \downarrow \\
& & & & M
\end{array}
\]
where the projection $\amalg U_{ij} \rightarrow \amalg U_{ij}$ is the 
identity and $T$ is the flat $D$-trivial bundle. This basic structure 
was suggested by Danny Stevenson.
The 3-curving at $m \in U_i$ is  
$\nu_i$. Clearly $s_\alpha^* \nu = \nu_\alpha$. 
The 2-curving on $\amalg U_{ij}$ is $\chi_{ij}$. The local 2-curving is
the $D$-obstruction form for $P_{\alpha \beta}$. In this example $P_{\alpha
\beta}$ is simply the restriction of the trivial bundle gerbe over 
$\amalg U_{ij}$ to $U_{\alpha\beta}$. A trivialisation is given by $T_{\alpha
\beta}$. Since $F_{T_ij} = 0$ the $D$-obstruction form for $P_{\alpha
\beta}$ is $\chi_{\alpha \beta}$. To define the bundle 2-gerbe product
we need to consider the following trivial bundle gerbe 
\[
\begin{array}{ccc}
T^* \otimes T & & T \\
\downarrow & & \downarrow \\
\amalg U_{ijk} & \rightrightarrows & \amalg U_{ijk} \\
& & \downarrow \\
& & \amalg U_{ijk}
\end{array}
\]
The bundle 2-gerbe product is defined by a $D$-trivialisation of this 
bundle gerbe, $J_{ijk}$. We define this to be the trivial bundle on
$\amalg U_{ijk}$ with connection at $m \in U_{ijk}$ given by
$\kappa_{ijk}$. 

Now we would like to find the local connection 1-forms. It might
appear that these would have to be trivial since the bundle 2-gerbe
connection is, however the product also carries information on
the local 1-connections.
We know from 
equation \eqref{defK} that there is an isomorphism of bundles with
connection
\begin{equation}
T_{\alpha\beta} \otimes T_{\beta\gamma} = 
T_{\alpha\gamma} \otimes J_{\alpha\beta\gamma} \otimes K^*_{\alpha\beta\gamma}
\end{equation}
where we have used the fact that the projection is just the identity to
pull back all of these to $U_{\alpha\beta\gamma}$. The $D$-obstruction
form for $K_{\alpha\beta\gamma}$ is now just $\kappa_{\alpha
\beta\gamma}$. This is because the $T$'s have zero connections and flat
trivialisations so their $D$-obstructions are zero. 

Finally we define the associator function on $\amalg U_{ijkl}$ by
$g_{ijkl}$ which satisfies the coherency condition and so we have a bundle 
2-gerbe which, by construction, has Deligne class
$(g_{\alpha\beta
\gamma\delta},\kappa_{\alpha\beta\gamma},\chi_{\alpha\beta},
\nu_{\alpha})$.

\end{proof}

It is sometimes possible, for example when $Y \rightarrow X^{[2]}$ is a 
fibration, to calculate the transition functions of 
a bundle 2-gerbe by an easier method as described by Stevenson \cite{ste3}. 
We shall give a brief outline of this method. It applies when the trivial 
bundle gerbes
$(P_{\alpha\beta},Y_{\alpha\beta},U_{\alpha\beta})$ admit a section 
$\sigma_{\alpha\beta} : U_{\alpha\beta} \rightarrow Y_{\alpha\beta}$.
Recall that unlike bundles, trivial bundle gerbes do not necessarily 
admit a section (see comments after lemma \ref{striv}). If they 
do then we have a map $\sigma_{\alpha\beta\gamma} :
U_{\alpha\beta\gamma} \rightarrow Y^{[2]}_{\alpha\beta\gamma}$ 
given by $(\sigma_{\alpha\gamma},\sigma_{\beta\gamma}\circ
\sigma_{\alpha\beta})$ where $\sigma_{\beta\gamma}\circ
\sigma_{\alpha\beta}$ is the map $Y_{\beta\gamma} \times
Y_{\alpha\beta} \rightarrow Y_{\alpha\gamma}$ which
is implicitly defined by the bundle 2-gerbe product. Use $\sigma_{
\alpha\beta\gamma}$ to pull back the bundle $P \rightarrow Y^{[2]}$
to $P_{\alpha\beta\gamma} \rightarrow U_{\alpha\beta\gamma}$.  
These bundles play the role of $K_{\alpha\beta\gamma}$ in the
general method. We now continue as in the general case by choosing
sections $h_{\alpha\beta\gamma}$ (which are equivalent to trivialisations
for bundles) of $P_{\alpha\beta\gamma}$ which
then satisfy equation \eqref{def4g} (with notation adjusted for sections),
\[
a_{\alpha \beta \gamma \delta}
h_{\beta \gamma \delta}
h_{\alpha \beta \delta} = g_{\alpha \beta \gamma \delta}
h_{\alpha \gamma \delta} h_{\alpha \beta \gamma}
\]
With the presence of sections rather than trivialisations the 
$D$-obstruction forms used to define the full Deligne class may
be replaced by the pullbacks of the relevant connections and curvings
by the sections.

\end{subsubsection}

We conclude our discussion of bundle 2-gerbes with two constructions
involving the cup product in Deligne cohomology which provide 
concrete examples and demonstrate the usefulness of the geometric 
picture of Deligne cohomology which bundle gerbes provide. 

\subsubsection{The Cup Product of Two Bundles}
Let $L$ and $J$ be two bundles over $M$. Then there is 
a bundle 2-gerbe defined by the cup product $L \cup J$. If $L$ and
$J$ both have connections then $L \cup J$ has a 2-curving. In
this case the Deligne class is $(h^{m_{\alpha\beta\gamma}}_{\gamma,\delta},
m_{\alpha\beta\gamma}B_{\gamma},\log (g_{\alpha\beta}) F_B, A_\alpha 
\wedge F_B)$, where the Deligne classes of $L$ and $J$ are $(g_{\alpha
\beta},A_\alpha)$ and $(h_{\alpha\beta},B_\alpha)$ respectively, 
$m_{\alpha\beta\gamma} = \log (g_{\alpha\gamma}) - \log (g_{\beta\gamma})
- \log (g_{\alpha\beta})$ is the Chern class of $L$ and $F_B$ is
the curvature of $J$. We define the bundle 2-gerbe $L \cup J$ by the 
following diagram:
\[
\begin{array}{ccc}
& & g \cup J \\
& & \downarrow \\
& & S^1 \times M \\
& \stackrel{(g,\pi_J)}{\nearrow} & \\
(L \times S^1) \times _\pi (J \times S^1) & \rightrightarrows &
L \times_\pi J \\
& & \downarrow \\
& & M
\end{array}
\]
where the map $(g,\pi_J)$ is defined by $(g,\pi_J)(l,\theta,j,\phi) =
(\theta,\pi_J(j))$. The fibre product bundle $L \times_\pi J$ with
structure group $S^1 \times S^1$ is often written as $L \oplus J$.
We take the bundle 
gerbe product to be the trivial and the associator function to 
be identically 1. The bundle gerbes $P_{\alpha\beta}$ are given by
$g \cup J$ over $(g_{\alpha\beta},m)$. 
We may define sections of these bundle gerbes
and so use the simpler method for calculating the Deligne class. 
Recall that we may write $g \cup J$ as
\[
\begin{array}{ccc}
J^{-n} & &\\
\downarrow & & \\
\R \times \Z \times M & \rightrightarrows & \R \times M \\
& & \downarrow \\
& & S^1 \times M
\end{array}
\]
Define sections $\sigma_{\alpha\beta}$ by
\begin{equation}
\sigma_{\alpha\beta} = (\log (g_{\alpha\beta}),m)
\end{equation}
The section $\sigma_{\alpha\beta\gamma}: U_{\alpha\beta\gamma}
\rightarrow (\R \times \Z \times M)_{\alpha\beta}$ is then given 
by
\begin{equation}
\sigma_{\alpha\beta\gamma} = (\log (g_{\alpha\gamma}), -m_{\alpha\beta\gamma},
m)
\end{equation}
Sections $h_{\alpha\beta\gamma}$ of the bundle $P_{\alpha\beta\gamma}$
are $t_{\gamma}^{m_{\alpha\beta\gamma}}$ where $t_\gamma$
is a local section of $J$. We can now calculate the transition functions
\begin{equation}
\begin{split}
g_{\alpha\beta\gamma\delta} h_{\alpha\gamma\delta}h_{\alpha\beta\gamma} 
&= h_{\beta\gamma\delta}h_{\alpha\beta\delta} \\
g_{\alpha\beta\gamma\delta} t_{\delta}^{m_{\alpha\gamma\delta}}
t_\gamma^{m_{\alpha\beta\gamma}} &= t_{\delta}^{m_{\beta\gamma\delta}}
t_{\delta}^{m_{\alpha\beta\delta}} \\
g_{\alpha\beta\gamma\delta} t_\delta^{m_{\alpha\gamma\delta} +
m_{\alpha\beta\gamma}}h_{\gamma\delta}^{-m_{\alpha\beta\gamma}}
&= t_{\delta}^{m_{\beta\gamma\delta} + m_{\alpha\beta\delta}}\\
g_{\alpha\beta\gamma\delta} &= h_{\gamma\delta}^{m_{\alpha\beta\gamma}}
\end{split}
\end{equation}
If the connection is $-nB$ on $J^{-n}$ then the local connection
forms given by the pullback by $h$ are $m_{\alpha\beta\gamma}B_\gamma$ 
as required. If the 2-curving is $f = rF_B$ on $\R \times M$ then 
$\sigma_{\alpha\beta}^* f = \log (g_{\alpha\beta}) F_B$, and 
the 3-curving is given by $A\wedge B$ on $L \times_\pi J$.

The general structure of this bundle 2-gerbe demonstrates a hierarchy
principle for cup product structures. Recall that the cup product
$f \cup L$ may be constructed in a similar way in terms of the 
cup product bundle, so contained within this cup product bundle
2-gerbe is a cup product bundle gerbe and within that a cup product
bundle. 

A second point to note is that this bundle 2-gerbe to some
extent resembles the associated bundle gerbe for a $G$-bundle.
In this case the $G$-bundle would be the $S^1 \times S^1$ 
bundle $L \oplus J$. The 4-curvature of this bundle 2-gerbe is
given in terms of the curvatures of the two bundles, 
$F_L \wedge F_J$. This is the image in real cohomology of 
the first Pontryagin class of $L\oplus J$, which also be the 
case with a bundle 2-gerbe associated to a $G$-bundle.
Since this structure group is not
simply connected such an associated bundle is not 
actually defined, the obstruction being the fact that the 
tautological bundle gerbe on $S^1 \times S^1$ is not 
well defined. In this particular case we are able to build a similar 
bundle 2-gerbe by replacing the tautological bundle with the
cup product bundle.  

\subsubsection{The Cup Product of a Function and a Bundle Gerbe}
Let $f$ be a $U(1)$-function and let $(P,Y,M)$ be a bundle gerbe
with connection $A$, curving $\eta$ and 3-curvature
$\omega$. The cup product $f\cup P$
has Deligne class 
\[
(g_{\beta\gamma\delta}^{n_{\alpha\beta}},n_{\alpha\beta}A_{\beta\gamma},
n_{\alpha\beta}\eta_{\beta}, \log_\alpha (f)\omega)
\]

The second
realisation of $f\cup L$ in \S \ref{cupsect} suggests that this 
bundle 2-gerbe should be realised geometrically with the following diagram:
\[
\begin{array}{ccc}
P^{-n} & & \\
\downarrow & & \\
\R \times \Z \times M & \rightrightarrows & \R \times M \\
& & \downarrow \\
\mspace{50mu} M & \stackrel{(f,m)}{\longrightarrow} & S^1 \times M
\end{array}
\]
where the fibre over $(r,n,m) \in \R \times \Z \times M$ is the
$n$-fold tensor product of $P$ with the trivial bundle 2-gerbe product
and associator function. The 3-curving is 
$r\omega$, the 2-curving $-n\eta$ and connection $-nA$.
Pulling back by the section $(\log_\alpha (f),n_{\alpha\beta},
m)$ we have a bundle gerbe $P^{-n_{\alpha\beta}}$ over $U_{\alpha\beta}$.
In this case the local constructions may once again be simplified 
however this time the construction is slightly different.

There exist trivialisations $J_{\alpha}$ of $P_{\alpha}$ over
$U_{\alpha}$, so over $U_{ij}$ we have trivialisations
$J_{\beta}^{n_{\alpha\beta}}$ of $P^{-n_{\alpha\beta}}$. 
On double overlaps $J_{\beta} \isom J_{\alpha} \otimes \pi^{-1} 
L_{\alpha\beta}$, for some bundle $L_{\alpha\beta} \rightarrow
U_{\alpha\beta}$, so on triple overlaps the local bundles 
are obtained by comparing $J_{\beta}^{-n_{\alpha\beta}} \otimes
J_{\gamma}^{-n_{\beta\gamma}}$ and $J^{-n_{\alpha\gamma}}_{\gamma}$,
\begin{equation}
\begin{split}
J_{\beta}^{-n_{\alpha\beta}} \otimes
J_{\gamma}^{-n_{\beta\gamma}} \otimes 
J^{n_{\alpha\gamma}}_{\gamma} &=
\pi^{-1} L_{\beta\gamma}^{-n_{\alpha\gamma}} \otimes
J_{\gamma}^{-n_{\alpha\beta}}  \otimes
J_{\gamma}^{-n_{\beta\gamma}} \otimes 
J^{n_{\alpha\gamma}}_{\gamma} \\
&= \pi^{-1} L_{\beta\gamma}^{-n_{\alpha\gamma}}
\end{split}
\end{equation}
Each $L_{\alpha\beta}$ has a section $l_{\alpha\beta}$ over
$U_{\alpha\beta}$, this allows us to find a section $l_{\beta\gamma}^{n_
{\alpha\beta}}$. Moreover the sections $l_{\alpha\beta}$ are the
sections which determine the Dixmier-Douady class of $P$, so they 
satisfy $l_{\alpha\beta} l_{\beta\gamma} = l_{\alpha\gamma} g_{\alpha\beta\gamma}$. Using this it is possible to calculate $\delta(l_{\beta\gamma}^{n_
{\alpha\beta}})_{\alpha\beta\gamma\delta}$ over $U_{\alpha\beta
\gamma\delta}$. This gives the correct transition functions and
the local connections and curvings may also be obtained without
difficulty and agree with what is expected.

\end{section}

\begin{section}{$\Z$-Bundle 0-Gerbes}

Consider once again the bundle gerbe hierarchy as given in table
\ref{t1}. A general method for constructing geometric realisations
may be approached in the following way.
Begin with some geometric representation, $R$, of a Deligne cohomology 
group $H^p$. Build representations of higher dimensional Deligne 
cohomology in the following way 

\[
\begin{array}{ccccccc}
H^p: & & R & & & & \\
\\
& & \vdots & & & & \\
\\
H^{p+1}: & & R & & & & \\
& & \downarrow & & & & \\
& & Y^{[2]} & \rightrightarrows & Y  & & \\
& & & & \downarrow & & \\
& & & & M & & \\
\\
& & \vdots & & & & \\
\\
H^{p+2}: & & R & & & & \\
& & \downarrow & & & & \\
& & Y^{[2]} & \rightrightarrows & Y & & \\
& & & & \downarrow & & \\
& & & & X^{[2]} & \rightrightarrows & X \\
& & & & & &  \downarrow \\
& & & & & &  M 
\end{array}
\]
The examples which we have already dealt with are where $R$ is either 
a function or a bundle. Furthermore we have shown that when we start 
with a function the next object in the hierarchy is a bundle 0-gerbe which
is equivalent to a bundle. Continuing this method gives the basic structure of
bundle gerbes and bundle 2-gerbes. Attempts to generalise to higher degree meet
difficulties due to the increasingly complicated nature of product
structures and related associativity conditions. 
In this section we wish to address the 
question of whether there is a starting point in the hierarchy 
below $U(1)$ functions. This leads to consideration of $\Z$-bundle
0-gerbes.

Since $U(1)$-functions represent $H^1(M,\Z(1)_D)$ then the only 
lower object in the hierarchy would be a representative of 
$H^0(M,\Z(0)_D)$, that is, the set of $\Z$ valued functions on $M$.
Using these as a basis for a hierarchy we get a new family of objects, 
$\Z$-bundle $n$-gerbes. A {\it $\Z$-bundle 0-gerbe} is defined as a 
bundle 0-gerbe $(\lambda,Y,M)$ where the function $\lambda$ takes values
in $\Z$.
\begin{proposition} 
The group of stable isomorphism classes of $\Z$-bundle 0-gerbes
is isomorphic to $H^1(M,\Z(0)_D)$.
\end{proposition}
\begin{proof}
The arguments of Proposition \ref{p1} still apply in this case
giving an isomorphism with $H^1(M,\Z)$. Furthermore
$H^1(M,\Z) \isom H^1(M,\Z(0)_D)$ since
$\Z(0)_D = \Z$.
\end{proof}
We wish to find an analogy with the equivalence of bundle 0-gerbes 
and bundles and of bundle gerbes and gerbes. This leads us to 
consider the Deligne cohomology group $H^1(M,\Z(1)_D)$ and the 
isomorphisms
\[
H^1(M,\Z(1)_D) \isom H^1(M,\Z(1) \rightarrow \ur) \isom H^0(M,\uu).
\]
To relate $\Z$-bundle 0-gerbes to our usual 
geometric representation of degree 1 Deligne cohomology,
$U(1)$-functions, we shall require some extra structure.
\begin{definition}
Let $(\lambda, Y,M)$ be a $\Z$-bundle 0-gerbe. A {\it $\Z$-curving} 
on $(\lambda,Y,M)$ is
a map $f:Y \rightarrow \R$ which satisfies $\delta(f) = \lambda$.
\end{definition}
\begin{proposition}
For each $\Z$-bundle 0-gerbe there exists a $\Z$-curving which
is unique up to the pull back of a globally defined $\R$-valued
function on the base.
\end{proposition}
\begin{proof}
Let $(\lambda,Y,M)$ be a $\Z$-bundle 0-gerbe. Choose an open 
cover $\{U_\alpha \}$ of $M$. Let 
$f_\alpha : \pi^{-1}(U_\alpha) \rightarrow
\Z$ be the family of functions defined by
\[
f_{\alpha}(y) = \lambda(s_\alpha (\pi(y)),y).
\]
Let $\{\phi_\alpha \}$ be a partition of unity on $M$ and let 
$f:Y \rightarrow \R$ be defined by
\[
f(y) = \sum_{\alpha} \phi_\alpha (\pi(y)) f_\alpha (y)
\]
Let $(y_1,y_2) \in Y^{[2]}$ with $\pi(y_1)=\pi(y_2)=m$. Then
\begin{eqnarray*}
\delta(f)(y_1,y_2) &=& f(y_2) - f(y_1) \\
&=& \sum_{\alpha} [ \phi_\alpha(\pi(y_2))f_\alpha(y_2) - 
\phi_\alpha(\pi(y_1))f_\alpha(y_1)] \\
&=& \sum_{\alpha} \phi_\alpha (m) [\lambda(s_\alpha(m),y_2) 
- \lambda(s_\alpha(m),y_1)] \\
&=& \sum_\alpha \phi_\alpha (m) \lambda(y_1,y_2) \\
&=& \lambda(y_1,y_2)
\end{eqnarray*}
Thus we see that $f$ is a $\Z$-curving for $(\lambda,Y,M)$. 

Suppose there exists another $\Z$-curving, $g$. Then 
$\delta(f-g) = 0$ so $f-g$ descends to a function on $M$. 
\end{proof}

The correspondence between $U(1)$-bundle 0-gerbes and $U(1)$-bundles
also applies with $U(1)$ replaced by $\Z$, so we may replace
stable isomorphism classes of bundle 0-gerbes with isomorphism
classes of $\Z$-bundles.

\end{section}

\begin{section}{$B^pS^1$-Bundles}
The correspondence between $B^pS^1$-bundles and 
Deligne cohomology was established by Gajer \cite{Gaj}. This 
was proven abstractly using sheaf theoretical arguments and also 
given in terms of explicit classifying maps using the 
bar resolution to obtain a realisation of the classifying 
spaces.
We shall show that this correspondence is suggested by
consideration of the classifying theory of bundles together
with our discussion of $\Z$-bundle 0-gerbes and
the bundle gerbe hierarchy. 

We recall some well known results on the classification of
bundles (see, for example, \cite{hus} or \cite{dup}). 
Let $P_G \rightarrow M$ be a principal $G$-bundle.
There exists a $G$-bundle $EG \rightarrow BG$ such
that $P_G = \psi^{-1} EG$ where $\psi$, called the classifying map, 
is unique up to
homotopy. The space $BG$ is called the universal classifying space and
$EG \rightarrow BG$ is called the universal $G$-bundle. 
The classifying bundle is a $G$ bundle with
a contractible total space which is unique up to homotopy 
equivalence.

The results of the previous section may be interpreted in terms of 
classifying theory. Since $\R \rightarrow S^1$ defines a $\Z$ bundle and 
$\R$ is contractible then $B\Z = S^1$.
The equivalence between $\Z$ bundles and homotopy classes of
$S^1$ valued functions
corresponds to the classifying theory of $\Z$ bundles. Given
a $\Z$-bundle 0-gerbe with $\Z$-curving $(\lambda,Y,M;f)$ we
define an $S^1$-function by $\hat{f}(m) = \exp f(y)$, where
$y \in \pi^{-1}(m)$. This is independent of the choice of $y$ since
for $y, y' \in \pi^{-1}(m)$ 
\begin{eqnarray*}
(\exp f(y))(\exp f(y'))^{-1} &=& \exp(f(y) - f(y')) \\
&=& \exp(\lambda(y',y)) \\
&=& 1
\end{eqnarray*}
since $\lambda(y',y) \in \Z$.
\end{section}
The existence of a $\Z$-curving up to a global $\R$-function is
equivalent to the existence of a classifying map up to 
homotopy. 

We now have the following equivalences of geometric realisations of 
Deligne cohomology:\\
\begin{center}
\begin{tabular}{ccc}
$\Z$-bundle 0-gerbes & $\longleftrightarrow$ & $B\Z$-functions \\
bundle 0-gerbes & $\longleftrightarrow$ & $BS^1$-functions 
\end{tabular}
\end{center}
where the left hand side consists of stable isomorphism classes and
the right hand side are homotopy classes of maps.  The case of 
higher dimensional objects is dealt with by the following 
\begin{proposition}\cite{Gaj}\label{Bp}
The group $H^p(M,\Z)$ is isomorphic to the group of isomorphism
classes of smooth principal $B^{p-2}S^1$-bundles over $M$.
\end{proposition}
Smoothness of classifying spaces is defined in terms of 
a differentiable space structure. For more details see \cite{Gaj} or
\cite{mos}. 
In general the iterated classifying spaces $B^pG$
are only defined if $G$ is Abelian. When this is the case each of
the spaces $B^pG$ is also an Abelian group.
Consider this result for low values of $p$. When $p = 2$ we have
the usual correspondence between $S^1$-bundles and $H^2(M,\Z)$ given
by the Chern class. When $p=3$ we have $H^3(M,\Z)$ and $BS^1$-bundles.
Our usual geometric realisation of $H^3(M,\Z)$ is stable isomorphism
classes of bundle gerbes so the following result is not surprising,
\begin{proposition}\cite{must}
The set of all stable isomorphism classes of bundle gerbes on $
M$ is in 
bijective correspondence with the set of all isomorphism classes
of $BS^1$ bundles on $M$. 
\end{proposition}

Since $BS^1$ is an Abelian group we could replace $BS^1$-bundles with
$BS^1$-bundle 0-gerbes. The bundle gerbe corresponding to
a $BS^1$ bundle which is referred to in the proposition is the 
lifting bundle gerbe
\[
\begin{array}{ccc}
& & ES^1 \\
& & \downarrow \\
& & BS^1 \\
& \nearrow & \\
P_{BS^1}^{[2]} & \rightrightarrows & P_{BS^1} \\
& & \downarrow \\
& & M 
\end{array}
\]
Given a principal $BS^1$-bundle on $M$ it is possible to construct a 
classifying map $M \rightarrow BBS^1$. This implies
that $BBS^1$ is a classifying space for bundle gerbes 
\cite{must}. This gives a series of equivalent 
realisations
\begin{center}
\begin{tabular}{ccccc}
bundle gerbes & $\longleftrightarrow$ & $BS^1$-bundles 
& $\longleftrightarrow$ & $BBS^1$-functions
\end{tabular}
\end{center} 
Note the similarity with the bundle gerbe hierarchy.

Now consider the case $p=4$. Since there is an isomorphism
between $H^4(M,\Z)$ and the stable isomorphism class of
bundle 2-gerbes on $M$ then proposition \ref{Bp} implies that
there is an isomorphism between stable isomorphism classes of
bundle 2-gerbe and isomorphism classes of $BBS^1$ bundles. 
The cases which we have already considered suggest that 
the bundle 2-gerbe corresponding to a 
$BBS^1$ bundle $(P_{BBS^1},M)$ should be the following bundle 2-gerbe
associated to a $BBS^1$-bundle:
\[
\begin{array}{ccc}
& & ES^1 \\
& & \downarrow \\
& & BS^1 \\
& \nearrow & \\
{EBS^1}^{[2]} & \rightrightarrows & EBS^1 \\
& & \downarrow \\
& & BBS^1 \\
& \nearrow & \\
P_{BBS^1}^{[2]} & \rightrightarrows & P_{BBS^1} \\
& & \downarrow \\
& & M
\end{array}
\] 
Since the relevant local data may not be easily calculated then
to prove this we would need to consider the theory of iterated classifying
spaces and the constructions of Gajer \cite{Gaj} in much more detail. 
This leads us away from our principal concerns here and so we have not done
this.
 
Other realisations are obtained by noting that $BBS^1$ bundles
are classified by $BBBS^1$ functions, and that $BS^1$ bundle gerbes 
are equivalent to $BBS^1$ bundles , so we have the following series of 
realisations:
\begin{center}
\begin{tabular}{ccccccc}
bundle 2-gerbes & $\leftrightarrow$ & $BS^1$-bundle gerbes
& $\leftrightarrow$ & $B^{2}S^1$ bundles & $\leftrightarrow$
& $B^{3}S^1$ functions
\end{tabular}
\end{center}
where it is to be understood that we are dealing the 
appropriate equivalence classes in each case, that is, stable isomorphism
for bundle gerbes, isomorphism for bundles and homotopy equivalence for 
functions.


\begin{section}{Comparing the Various Realisations}
We now present a table comparing the various geometric realisations
of Deligne cohomology which we have discussed. We include for completeness
the differential characters of Cheeger and Simons (\cite{chsi}, \cite{bry}). 
Since the relationship between these and bundle gerbes is closely related
to the theory of holonomy we postpone a definition and further discussion
until Chapter 7.
\begin{landscape}
\begin{table}
\begin{center}
\caption{Geometric Realisations of Deligne Cohomology}
\begin{tabular}{llllll}
$H^0(M,\Z(0)_D)$ & $\Z$-functions & &&& \\
$H^0(M,\Z(1)_D)$ & constant $\Z$-functions & &&& \\
\\
$H^1(M,\Z(1)_D)$ & $\Z$-bundle 0-gerbes &
$S^1$-functions & & & deg 1 differential characters \\
$H^1(M,\Z(2)_D)$ & Flat $\Z$-bundle 0-gerbes & constant $S^1$-functions &
& \\
\\
$H^2(M,\Z(1)_D)$ & $\Z$-bundle gerbes & bundle 0-gerbes & $BS^1$-functions
 & bundles & \\
$H^2(M,\Z(2)_D)$ & $\Z$-bundle gerbes&
bundle 0-gerbes &  & bundles&
deg 2 differential characters \\
& with curving & with connection & & with connection& \\
$H^2(M,\Z(3)_D)$ & flat $\Z$-bundle gerbes &
flat bundle 0-gerbes & & flat bundles & \\
\\
$H^3(M,\Z(1)_D)$ & $\Z$-bundle 2-gerbes & bundle gerbes &
$BS^1$-bundles  & gerbes & \\
$H^3(M,\Z(3)_D)$ & $\Z$-bundle 2-gerbes& 
bundle gerbes & $BS^1$-bundles&
 gerbes & deg 3 differential characters \\
& with 2-curving & with curving & with connection & with curving & \\
$H^3(M,\Z(4)_D)$ & flat $\Z$-bundle 2-gerbes &
flat bundle gerbes & flat $BS^1$-bundles & flat gerbes & \\
\\
$H^4(M,\Z(1)_D)$ & & bundle 2-gerbes & $BBS^1$-bundles & 2-gerbes &
\\
$H^4(M,\Z(4)_D)$ & & bundle 2-gerbes & & 2-gerbes &
deg 4 differential characters \\
& & with 2-curving & & with 2-curving & \\
$H^4(M,\Z(5)_D)$ & & flat bundle 2-gerbes & & flat 2-gerbes &
\end{tabular}
\end{center}
\end{table}
\end{landscape}
\end{section}

%% file: chapter5.tex
\chapter{Holonomy and Transgression}

In this chapter we consider the generalisation of the holonomy of
a bundle around a loop to bundle gerbes. We show how this relates to
local formulae which define a transgression map in Deligne cohomology.
The key to this generalisation is to consider holonomy as a property
of a Deligne class. We have already defined the flat holonomy of a 
flat Deligne class. The holonomy of a general class in $H^p(M,\D^p)$ 
associated with a map $\psi: X \rightarrow M$ for some closed $p$-manifold $X$
is defined to be the flat holonomy of the pullback of the class to $X$.
We shall see how this approach relates to the usual construction of
holonomy for bundles, and then go on to consider holonomy for
bundle gerbes, bundle 2-gerbes and general Deligne classes.

\section{Holonomy of $U(1)$-Bundles}              

We review the holonomy of principal $U(1)$-bundles with an emphasis
on Deligne cohomology which is useful for generalisation to 
bundle gerbes.

Recall that flat bundles and bundle 0-gerbes have a flat holonomy which
is a class in $H^1(M,U(1))$.
It is useful to review the equations which define the Deligne cohomology 
class in general and in the particular cases of flat and trivial bundles.
\\

The Deligne class $(\un{g},\un{A})$ satisfies 
\begin{eqnarray}
d\log g_{\alpha \beta} = A_\beta - A_\alpha
\end{eqnarray}
If it is flat then we can find $U(1)$-valued functions satisfying
$d\log a_\alpha = A_\alpha$ and we have
\begin{equation}
c_{\alpha \beta} = g^{-1}_{\alpha \beta} a^{-1}_\alpha a_\beta
\end{equation}
The functions $c_{\alpha \beta}$ are constant and define 
the flat holonomy class.
If the Deligne class has a trivialisation $\un{h}$ then
\begin{eqnarray}
g_{\alpha \beta} &=& h_{\alpha}^{-1} h_{\beta} \label{defh} \\
d\log h_\alpha - A_\alpha &=& d\log h_\beta - A_\beta \label{Ah}
\end{eqnarray}
\begin{proposition}\cite{bry}\label{cg}
The assignment of a flat holonomy to a flat bundle gives an
isomorphism $H^1(M,\D^p) = H^1(M,U(1))$ for $p > 1$.
\end{proposition}
\begin{proof}
Suppose we have a flat bundle represented by a Deligne 
class $(g_{\alpha \beta},A_{\alpha})$ with flat
holonomy $c_{\alpha \beta}$. We first show that the 
class in $H^1(M,U(1))$ is independent of the choices of
$a_\alpha$. Choose $a'_\alpha$ satisfying $d\log a'_\alpha = A_\alpha$.
Then $a'_\alpha = a_\alpha + K_\alpha$ where $K_\alpha $ are $U(1)$-valued
constants.
Thus $c'_{\alpha \beta} = c_{\alpha \beta} + \delta(K)_{\alpha \beta}$
and so $\underline{c'}$ and $\underline{c}$ define  the same class in
$H^1(M,U(1))$. 

Clearly the map $(\un{g},\un{A}) \mapsto \un{c}$ is a homomorphism.
Let $(\delta(\un{h}),d\log\un{h})$ represent
 a trivial class in $H^1(M,\D^p)$. The corresponding
flat holonomy is given by $\delta(\un{h}) \delta(\un{a})^{-1} =
\delta(\un{h}. \un{a}^{-1})$. Since $d\log\un{h} = d\log\un{a}$ this
represents a trivial class in $H^1(M,U(1))$. 

Given a class $\un{c} \in H^1(M,U(1))$ define a class in
$H^1(M,\D^p)$ by $(-\un{c},0)$. Observe that if $\un{c}$
is the flat holonomy of $(\un{g},\un{A})$ then the two
Deligne classes $(-\un{c},0)$ and $(\un{g},\un{A})$ differ
by a trivial class $(\delta(\un{a}),d\log\un{a})$. Therefore the map
$(\un{g},\un{A}) \mapsto \un{c}$ is onto. Also it is 
clear that a trivial class in $H^1(M,U(1))$ leads to a 
trivial class in $H^1(M,\D^p)$. Therefore we have an 
isomorphism.
\end{proof}

Mostly we shall be interested in the flat holonomy of a bundle 
over $S^1$. All bundles with connection over $S^1$ are flat so they all 
have a flat
holonomy $\un{c} \in H^1(S^1,U(1)) = S^1$. We shall demonstrate how
to calculate this element of $S^1$ for a given bundle with connection.
Let $(g_{\alpha \beta},A_\alpha)$ represent a flat bundle on $S^1$ with
flat holonomy $c_{\alpha \beta}$. Since $H^2(S^1,\Z)=0$ there exists 
a trivialisation $\delta(\un{h}) = \un{g}$. Using this the flat holonomy
becomes $c_{\alpha \beta} = h_\alpha \cdot h_\beta^{-1} \cdot 
a_\alpha a_\beta^{-1}$. By considering $\log (c_{\alpha \beta})$ as 
a representative of a class in the \v{C}ech cohomology of $S^1$ we can use the 
following diagram to calculate the isomorphism with de Rham cohomology:
\[
\begin{CD}
A_\alpha - d\log(h_\alpha)  \\
@AAdA \\
\log(a_\alpha) - \log(h_\alpha) @>\delta>>
\log (c_{\alpha \beta}) 
\end{CD}
\]
Thus the 1-form $A_\alpha - d\log(h_\alpha)$ of equation \eqref{Ah}
is the de Rham representative
of the flat holonomy. It is globally defined on $S^1$ and is well 
defined modulo $2\pi$-integral forms since 
the original \v{C}ech class was defined modulo $\Z(1)$. To evaluate 
it as an element, $H(c_{\alpha \beta})$, of $S^1$ we integrate,
\begin{equation}\label{Hc}
H (c_{\alpha \beta}) = \exp \int_{S^1}  A_\alpha -d\log(h_\alpha) .
\end{equation}
Since the flat holonomy class is isomorphic to the Deligne class 
we should be able to write \eqref{Hc} as a function of 
the Deligne class, $H(g_{\alpha \beta}, A_\alpha)$. To do this
we would like to separate the two terms in the integral into 
separate integrals however they are not independently defined globally so
this is not possible. We shall have to break up the integral into
a sum of integrals on intervals where $d\log(h_\alpha)$ and $A_\alpha$
are defined, to do this we use the method used by Gawedski \cite{gaw}.

At the moment we have a Deligne class in terms of some open cover
of $S^1$, denoted by subscripts $\alpha$ and $\beta$.
Let $t$ be a triangulation of $S^1$ consisting of edges, $e$, and
vertices, $v$ such that each edge is contained wholly within $U_\alpha$ for
at least one $\alpha$. Such a triangulation is said to be subordinate to
the open cover and is guaranteed to exist since compactness implies
the existence of a Lebesgue number \cite[p179]{mun}, we simply
triangulate the circle such that all edges have length which
is less than this number and therefore are contained within a set in
the open cover. 
We can express $S^1$ as a sum $\sum e$ over
all $e \in t$, so this should allow us to break up the integral of 
the global 1-form $A_\alpha -d\log(h_\alpha)$ into a sum over 
these edges, but we need to choose an open set that covers 
each edge first. 
For each $e$ let $U_{\rho(e)}$ be an element of the open cover of $S^1$ such
that $e \subset U_{\rho(e)}$. Here $\rho:t \rightarrow\A$ is an
{\it index map} from the triangulation to the index set for the 
open cover of $M$
\footnote{Gawedski did not use index maps explicitly though they were implicit
in his construction. They were used in this
context by Brylinski \cite{bry} and the terminology appears to be due to 
Gomi and Terashima \cite{gote2}}
. We can now 
split up the terms in the 
integral,
\begin{eqnarray*}
H(c_{\alpha \beta}) &=& \exp \sum_e [\int_e A_{\rho(e)} - \int_e
d\log(h_{\rho(e)})] \\
&=& \exp [ \sum_e \int_e A_{\rho(e)} +
\sum_{v,e} \log(h^{-1}_{\rho(e)})(v)] 
\end{eqnarray*}
where we use the convention that $\underset{v,e}{\sum}$ represents a 
sum over all edges and all vertices bounding each edge such that the 
sign is reversed for vertices which inherit the opposite orientation
to the corresponding edge. This means that for each vertex there are
two terms with opposite sign, one each for each of the edges bounded
by that vertex. Observe that the following equality 
follows from \eqref{defh},
\begin{equation}
\sum_{e,v} \log(h^{-1}_{\rho(e)})(v) = \sum_{e,v} \log(g_{\rho(e)\rho(v)})(v)
- \log(h_{\rho(v)})(v)
\end{equation}
Furthermore the second term on the right hand side is equal to zero since
each vertex bounds exactly two edges which give two equal terms with
opposite signs in the summation.
The flat holonomy is now
\begin{eqnarray}
H(c_{\alpha \beta}) 
&=&   \prod_e \exp \int_e A_{\rho(e)} \cdot
 \prod_{e,v} g_{\rho(v)\rho(e)}(v) \label{HgA} \\
&=& H(g_{\alpha \beta},A_{\alpha}) \nonumber
\end{eqnarray}
This construction is independent of the choice of triangulation. Suppose
we choose another triangulation, $\hat{t}$. Since this 
triangulation must also be subordinate to the open cover we may
assume without loss of generality that $\rho(\hat{e}) = \rho(e)$. Denote the
flat holonomy corresponding to $\hat{t}$ by $\hat{H}$, and denote the 
two components $\exp \hat{H}_g$ and $\exp \hat{H}_A$. For this
calculation is advantageous to expand the sum over the pair
$e,v$ as a sum over $v$ in the following way,
\begin{eqnarray}
\sum_{v,e} \log (g_{\rho(e)\rho(v)})(v) &=&
\sum_v \log (g_{\rho(e^+(v))\rho(v)})(v) - \log (g_{\rho(e^-(v))\rho(v)})
\\
&=& \sum_v \log (g_{\rho(e^+(v))\rho(e^-(v))})(v)
\end{eqnarray}
where $e^+(v)$ (resp. $e^-(v)$) is the edge bounded by $v$ such that it 
inherits a positive (negative) orientation. Now the difference between the 
terms corresponding to the two triangulations is 
\[
H_g - \hat{H}_g = \sum_v \log (g_{\rho(e^+(v))\rho(e^-(v))})(v) -
\sum_{\hat{v}} \log (g_{\rho(e^+(\hat{v})\rho(e^-(\hat{v}))})(\hat{v})
\]
Since both triangulations are subordinate to the open cover we 
may consider pairs $(v,\hat{v})$ which are the unique vertices from
each triangulation that lie within a particular double intersection 
of open sets. We can replace both summations in the expression above
by a summation over such pairings. Furthermore given such a pairing 
we have $\rho(e^+(v)) = \rho(e^+(\hat{v}))$ and $\rho(e^-(v)) = 
\rho(e^-(\hat{v}))$.
The difference now becomes 
\begin{eqnarray*}
H_g - \hat{H}_g &=& \sum_{(v,\hat{v})} \log (g_{\rho(e^+(v))\rho(e^-(v))})(v) -
\log (g_{\rho(e^+(v))\rho(e^-(v))})(\hat{v}) \\
&=& \sum_{(v,\hat{v})} \int_{v - \hat{v}} d\log g_{\rho(e^+(v))\rho(e^-(v))} \\
&=& \sum_{(v,\hat{v})} \int_{v - \hat{v}} A_{\rho(e^-(v))} - A_{\rho(e^+(v))}
\end{eqnarray*}
Now consider the difference 
\[
H_A - \hat{H}_A = \sum_e \int_e A_{\rho(e)} - \sum_{\hat{e}} \int_{\hat{e}}
A_{\rho(\hat{e})}
\]
As with the vertices we can pair the edges $(e,\hat{e})$ such that
$\rho(e) = \rho(\hat{e})$ and replace both sums with a sum over these pairings
to get
\begin{eqnarray*}
H_A - \hat{H}_A &=& \sum_{(e,\hat{e})} \int_e A_{\rho(e)} - \int_{\hat{e}}
A_{\rho(e)} \\
&=& \sum_{(e,\hat{e})} \int_{e - \hat{e}} A_{\rho(e)}
\end{eqnarray*}
Each difference $e - \hat{e}$ consists of two components (in terms of
vertices), $e^+ - \hat{e}^+$ and $e^- - \hat{e}^-$. Using
this to split up the integral into two terms we get
\begin{eqnarray*}
H_A - \hat{H}_A &=& \sum_{(e,\hat{e})} \int_{e^+ - \hat{e}^+}
A_{\rho(e)} + \int_{e^- - \hat{e}^-}
A_{\rho(e)} \\ &=& 
\sum_{(v,\hat{v})} \int_{v - \hat{v}} A_{\rho(e^+(v))} - A_{\rho(e^-(v))}
\end{eqnarray*}
where we have changed to a summation over vertices and used the fact that
$v-\hat{v}$ is equal to one component each from $e^+(v) - \hat{e}^+(v)$
and $\hat{e}^-(v) - e^-(v)$. This term is the opposite of $H_g -
\hat{H}_g$ therefore $H = \hat{H}$. 

Since the 1-forms $A_{\rho(e)} - d\log h_{\rho(e)}$ are global then 
the integral defining the holonomy must be independent of the choice
of index map $\rho$. This implies that \eqref{HgA} should also be independent 
of the choice of $\rho$. 
This may be easily verified. Suppose we have two such choices,
$\rho_0$ and $\rho_1$. Then the difference is given by
\begin{equation}
\begin{split}
\prod_e \exp \int_e A_{\rho_1(e)} - A_{\rho_0(e)} 
\cdot \prod_{v,e} g_{\rho_1(e)\rho_1(v)}g^{-1}_{\rho_0(e)\rho_1(e)} 
&= \prod_e \exp \int_e d\log g_{\rho_0(e)\rho_1(e)} \\
&  \mspace{60mu}  
\cdot \prod_{v,e} g_{\rho_1(e)\rho_1(v)}g^{-1}_{\rho_0(e)\rho_1(e)} \\
&= \prod_{v,e} g_{\rho_0(e)\rho_1(e)}g_{\rho_1(e)\rho_1(v)}
g^{-1}_{\rho_0(e)\rho_0(v)}\\
&= \prod_{v,e} g_{\rho_0(v)\rho_1(v)} \\
&= 1
\end{split}
\end{equation}
since for each $v$ there are two identical terms with opposite 
signs corresponding to the two edges which share $v$ as a bounding
vertex. Note that the global version is not explicitly independent 
of the choice of trivialisation, $h$, however we may deduce this 
from the explicit independence of the local version \eqref{HgA}. 
It may also be calculated directly, this calculation is quite
similar to the one described above.

We have defined an element of $S^1$ associated with every isomorphism
class of flat bundle over $S^1$. It is given by equation
\eqref{HgA} and is 
well defined. We would like now to show how this
relates to the usual concept of the holonomy of a bundle with 
connection around a loop. Let $(L,M;A)$ be a bundle with connection (not
necessarily flat) over $M$. Let $\gamma$ be a loop in $M$, that is,
$\gamma$ is a smooth map $S^1 \rightarrow M$. Use $\gamma$ to pull $L$ back
to $S^1$. Let $H(\gamma^{-1}(L;A))$ be the flat holonomy of the pull 
back bundle. We define this to be the holonomy of $(L;A)$ around
$\gamma$. 

We would like to give an explicit formula for the holonomy. 
To do these we need to examine the Deligne class of a pull
back bundle. Once we have this we can apply equation \eqref{HgA}.

Suppose we have a map $N 
\stackrel{\phi}{\rightarrow} M$ between compact manifolds. Let $\{U_\alpha\}_
{\alpha\in A}$ be a good cover on $M$. The set $A$ is finite since $M$ is 
compact. There is a cover $\{V_{\phi(\alpha)}
=\phi^{-1}(U_\alpha)\}_{\alpha \in A}$ on
$N$ called the induced cover. If we have a bundle with connection 
$(L;A)$ on $M$, then we can calculate the Deligne class of 
$(\phi^{-1}L,N;\phi^*A)$ in terms of the induced cover. 
\begin{lemma}
Let $(g_{\alpha \beta},A_{\alpha})$ be the Deligne class of 
$(L,M;A)$. Then the Deligne class of $(\phi^{-1}L,N;\phi^*A)$ with 
respect to the induced cover is
$(g_{\phi(\alpha)\phi(\beta)},A_{\phi(\alpha)})$
where 
\begin{eqnarray*}
g_{\phi(\alpha)\phi(\beta)}(n) &=& g_{\alpha \beta}(\phi(n)) 
\qquad \mbox{and}\\
A_{\phi(\alpha)} &=& \phi^* A_{\alpha}
\end{eqnarray*}
\end{lemma}
Putting this together with \eqref{HgA} we get
\begin{proposition}\cite{bry}\cite{gaw}\label{holprop}
The holonomy of a bundle with Deligne class
$(\un{g},\un{A})$ around a loop
$\gamma$ is given by 
\begin{equation}\label{hloop}
H((\un{g},\un{A});\gamma) =
 \prod_e \exp \int_e \gamma^*A_{\rho(e)}  \cdot
 \prod_{v, e} g_{\rho(e)\rho(v)}(\gamma(v))
\end{equation}
\end{proposition}
Now recall the usual definition of the holonomy of a bundle
\begin{definition}\label{holdef}
Let $(L,M)$ be a bundle with connection $A$. Any path in $M$ has a 
unique lift through each element of the fibre over the starting point
which is horizontal with respect to $A$. In particular each
loop $\gamma$ has a unique horizontal lift $\tilde{\gamma}$ which
defines an automorphism of the fibre over $\gamma(0)$. The 
{\it holonomy of the connection $A$ around $\gamma$} 
is the element of $S^1$ defined by
$\tilde{\gamma}(1) = \tilde{\gamma}(0)\cdot H(\gamma)$.
\end{definition}

\begin{proposition}\cite{gaw}
The holonomy of Proposition \ref{holprop} is the same as the holonomy of
definition \ref{holdef}.
\end{proposition}

We shall relate these two concepts of holonomy  
by considering {\it parallel transport}. Given a path $\mu \in \Map(I,M)$
the horizontal lift $\tilde{\mu}$ defines a morphism of fibres
$P_{\mu(0)} \rightarrow P_{\mu(1)}$. This is called parallel transport. 
Two paths $\mu$ and $\mu'$ such that $\mu(1) = \mu'(0)$ may be
composed and the horizontal lift of the composition defines
a composition of parallel transports. If we consider a loop as 
a composition of a number of paths then the holonomy is 
defined by the composition of the parallel transports along
each path. By breaking up the loop into components 
$\gamma([t_i,t_{i+1}])$ over which $P$
admits sections $s_i$ then there is an explicit formula for 
parallel transport over each component:
\begin{equation}
s_i (\gamma(t_i))  \mapsto s_i(\gamma(t_{i+1})) 
\exp (\int_{t_i}^{t_{i+1}} s_i^* A )
\end{equation}
Composition then gives a product of terms which combine to give the
local formula \eqref{hloop}. This suggests that we have used a rather long
and complicated method for calculating the holonomy of a bundle, 
however it turns out that our method is useful as it generalises to 
bundle gerbes and to higher degrees. In addition to this it allowed us
to demonstrate certain features of the higher theory in a relatively
simple setting.


\section{Holonomy of Bundle Gerbes}
To define the holonomy of a bundle gerbe we follow the procedure used in the
previous section. The standard technique for 
deriving a formula for holonomy of a bundle 
(as described at the end of the previous section) 
cannot be used here for two main
reasons. One is that it turns out that a bundle gerbe has a holonomy
over a surface rather than a loop, so we cannot just choose a direction 
to integrate around as is the case with a loop. 
Secondly it is not clear what a horizontal lift or parallel transport
map would be
in this situation. This motivates us to define the holonomy of a 
bundle gerbe by first considering the holonomy of a Deligne 
class corresponding to a bundle gerbe.

The Deligne class $(\un{g},\un{A},\un{\eta})$ of a bundle gerbe satisfies
\begin{eqnarray}
d\log g_{\alpha \beta \gamma} &=& -A_{\beta \gamma} + A_{\alpha \gamma}
-A_{\alpha \beta} \\
dA_{\alpha \beta} &=& \eta_\beta - \eta_\alpha
\end{eqnarray}
If the bundle gerbe is flat then 
\begin{eqnarray}
\eta_\alpha &=& dB_\alpha \\
A_{\alpha \beta} &=& B_\beta - B_\alpha + d\log a_{\alpha \beta} \\
c_{\alpha \beta \gamma} &=& g^{-1}_{\alpha \beta \gamma} a^{-1}_{\beta \gamma}
a_{\alpha \gamma} a_{\alpha \beta}^{-1}
\end{eqnarray}
and $c_{\alpha \beta \gamma}$ is the flat holonomy class. If the bundle 
gerbe has trivialisation $\un{h}$ then 
\begin{eqnarray}
g_{\alpha \beta \gamma} &=& h_{\beta \gamma} h_{\alpha \gamma}^{-1}
h_{\alpha \beta} \\
d\log h_{\alpha \beta} &=& -A_{\alpha \beta} + k_\beta - k_\alpha \\
\eta_\alpha -dk_\alpha &=& \eta_\beta - dk_\beta
\end{eqnarray}

Now consider the particular case of a bundle gerbe 
over $\Sigma$, a 2-manifold without boundary. In this case the 
bundle gerbe is not only flat, but also trivial. 
The \v{C}ech-de Rham isomorphism is given by the following diagram:
\[
\begin{CD}
 \eta_\alpha - dk_\alpha\\
@AAdA \\
B_\alpha - k_\alpha @>\delta>> 
0 \\
& & @AA-dA \\
& &- \log a_{\alpha \beta} - \log h_{\alpha \beta}
@>\delta>> \log c_{\alpha \beta \gamma} 
\end{CD}
\]
Thus the globally defined 2-form $\eta-dk$ is the 
de Rham representative of the flat holonomy of the 
bundle gerbe. Since $H^2(\Sigma,U(1)) = U(1)$ we may 
evaluate this class as an element of the circle by 
integrating over the surface $\Sigma$ and taking the exponential.

Thus in terms of bundle gerbes holonomy is defined in the 
following way. 
\begin{definition}\cite{camimu}
Let $(P,Y,M;A,\eta)$ be a bundle gerbe with 
connection and curving and let $\psi:\Sigma \rightarrow M$ be a map
of a surface into $M$. The {\it holonomy of $(P,Y,M;A,\eta)$ over
$\Sigma$} is the flat holonomy of $\psi^*P$.
\end{definition}

To see that this is well defined consider that when we pull back the 
bundle gerbe $P$ to $\Sigma$ using $\psi$ the resulting  bundle gerbe 
has an induced curving which we denote $\psi^* \eta$ and for 
dimensional reasons has a trivialisation $L$. Denote the 
curvature of this trivialisation (given some connection which is 
compatible with the bundle gerbe connection) by $F_L$. The 
2-form $\psi^* \eta - F_L$ descends to $\Sigma$ and its integral
over $\Sigma$ defines the flat holonomy which is an element of 
$H^2(\Sigma,U(1)) = U(1)$. This is independent of the 
choice of trivialisation since a different choice just changes $F_L$ by a 
closed 2-form which descends to $\Sigma$. We shall also
see this when we calculate a formula for the holonomy which is 
explicitly independent of this choice.

\begin{proposition}
The holonomy of a bundle gerbe with Deligne class $(\un{g},\un{A},\un{\eta})$ 
on $M$
over a surface $\psi: \Sigma \rightarrow M$ is given by the following 
formula of Gawedski \cite{gaw}:
\begin{equation} \label{holbg}
H((\un{g},\un{A},\un{\eta});\psi) =  \prod_b \exp \int_b \psi^*\eta_{\rho(b)} \cdot
\prod_{e,b} \exp \int_e \psi^*A_{\rho(b)\rho(e)} \cdot
\prod_{v,e,b}  g_{\rho(b)\rho(e)\rho(v)}(\psi(v))
\end{equation}
\end{proposition}
\begin{proof}
To evaluate the holonomy in terms of the original Deligne 
class we shall need to triangulate $\Sigma$. This 
triangulation, $t$, will consist of vertices, $v$, edges, $e$, and faces $b$
and is required to be subordinate to the open cover $\{ U_\alpha \}_{
\alpha \in \A}$. Thus there exists an index map 
$\rho :t \rightarrow \A$ such that
$b \subset U_{\rho(b)}$, $e \subset U_{\rho(e)}$ and $v\subset U_{\rho(v)}$
for all $b,e,v \in t$.
The integral over $\Sigma$ can be broken up into a sum of integrals over
$b$,
\begin{eqnarray*}
H(c_{\alpha \beta \gamma}) &=& \exp \sum_b \int_b (\eta - dk) \\
&=& \exp \sum_b \left(\int_b \eta_{\rho(b)} + \int_b -dk_{\rho(b)}\right) 
\end{eqnarray*}
Applying Stokes' theorem to the second term gives
\[
H(c_{\alpha \beta \gamma}) = \exp \left( \sum_b \int_b \eta_{\rho(b)} + 
\sum_b \int_{\partial b}
-k_{\rho(b)}\right)
\]
In the second term we have a sum $\sum_b \int_{\partial b}$.
If we break $\partial b$ into a sum of edges we can write this as
$\underset{e,b}{\sum} \int_e$ where the convention is that the sum is over all 
faces and all edges bounding each face, and the integral is given the 
corresponding induced orientation.
\begin{equation*}
\begin{split}
H(c_{\alpha \beta \gamma}) &= \exp \left(\sum_b \int_b \eta_{\rho(b)} +
\sum_{e,b} \int_e -k_{\rho(b)} \right)\\
&= \exp \left(\sum_b \int_b \eta_{\rho(b)} +
\sum_{e,b} \int_e \left( A_{\rho(b)\rho(e)} + d\log(h_
{\rho(b)\rho(e)}) - k_{\rho(e)} \right) \right)\\
&= \exp \left(\sum_b \int_b \eta_{\rho(b)} +
\sum_{e,b} \int_e ( A_{\rho(b)\rho(e)}
+ \int_{\partial e} \log(h_{\rho(b)\rho(e)}) 
\right)\\
&= \exp \left(\sum_b \int_b \eta_{\rho(b)} +
\sum_{e,b} \int_e A_{\rho(b)\rho(e)} +
\sum_{v,e,b} \log(h_{\rho(b)\rho(e)}(v)) \right) \\
&= \exp (\sum_b \int_b \eta_{\rho(b)} +
\sum_{e,b} \int_e A_{\rho(b)\rho(e)} +
\sum_{v,e,b} \log (g_{\rho(b)\rho(e)\rho(v)}(v)) \\
& \mspace{195.0mu}  - \log (h_{\rho(e)\rho(v)}
(v)) + \log (h_{\rho(b)\rho(v)}(v)) ) \\
&= \exp \left(\sum_b \int_b \eta_{\rho(b)} +
\sum_{e,b} \int_e A_{\rho(b)\rho(e)} +
\sum_{v,e,b} \log (g_{\rho(b)\rho(e)\rho(v)}(v)) \right)
\end{split}
\end{equation*}
We have claimed in this calculation that certain terms cancel out. Let
$I(e)$ denote a term depending only on $e$, $I(e,b)$ a term depending only
on $e$ and $b$ and so on. Then we have used the following results:
\begin{eqnarray}
\sum_{e,b} I(e) &=& 0 \\
\sum_{v,e,b} I(v,e) &=& 0 \\
\sum_{v,e,b} I(v,b) &=& 0
\end{eqnarray}
The first two are true because for each edge there are exactly two faces 
with that edge as boundary and they have opposite induced orientations.
The third is true since given a face and a vertex of that face there are
exactly two edges which bound the face and have the vertex as a boundary
component. Furthermore the vertex inherits opposite orientations from 
each of these edges. Note that the first two results would no longer 
hold if we triangulate a surface with boundary. We shall deal with this
situation in the next chapter.

We now have a formula for the holonomy of a flat bundle gerbe in terms
of its Deligne class,
\begin{equation} \label{bghol}
H(\un{g},\un{A},\un{\eta}) =  \prod_b \exp \int_b \eta_{\rho(b)} \cdot
\prod_{e,b} \exp \int_e A_{\rho(b)\rho(e)} \cdot
\prod_{v,e,b}  g_{\rho(b)\rho(e)\rho(v)}(v)
\end{equation}
As in the previous section this formula may be adapted to
define the holonomy of a general bundle gerbe with curving, $(P,Y,M;A,\eta)$, 
associated with a smooth map of surface into $M$, $\psi:\Sigma \rightarrow M$.

This leads us to the required formula
\begin{equation*} 
H((\un{g},\un{A},\un{\eta});\psi) =  \prod_b \exp \int_b \psi^*\eta_{\rho(b)} \cdot
\prod_{e,b} \exp \int_e \psi^*A_{\rho(b)\rho(e)} \cdot
\prod_{v,e,b}  g_{\rho(b)\rho(e)\rho(v)}(\psi(v))
\end{equation*}
\end{proof}

\section{Holonomy of Bundle 2-Gerbes}
For the case of a bundle 2-gerbe we must first establish the notation
associated with the flat holonomy and with trivialisations.

The Deligne class $(\un{g},\un{A},\un{\eta},\un{\nu})$ of a bundle 
2-gerbe satisfies the following equations:
\begin{eqnarray}
d\log g_{\alpha \beta \gamma \delta} &=& A_{\beta \gamma \delta}
- A_{\alpha \gamma \delta} + A_{\alpha \beta \delta} - A_{\alpha \beta
\gamma} \\
dA_{\alpha \beta \gamma} &=& -\eta_{\beta \gamma} + \eta_{\alpha \gamma} -
\eta_{\alpha \beta} \label{b2g2} \\
d\eta_{\alpha \beta} &=& \nu_\beta - \nu_\alpha
\end{eqnarray}
If we assume that the bundle 2-gerbe is flat then we have the following
set of equations
\begin{eqnarray}
\nu_\alpha &=& dq_{\alpha} \\
\eta_{\alpha \beta} &=& q_\beta - q_\alpha + dB_{\alpha \beta} \\
A_{\alpha \beta \gamma} &=& -B_{\beta \gamma} + B_{\alpha \gamma}
- B_{\alpha \beta} + d\log a_{\alpha \beta \gamma} \\
c_{\alpha \beta \gamma \delta} &=& g^{-1}_{\alpha \beta \gamma \delta} 
a_{\beta \gamma \delta} a^{-1}_{\alpha \gamma \delta} a_{\alpha \beta
\gamma} a^{-1}_{\alpha \beta \gamma}
\end{eqnarray}
The constants $c_{\alpha \beta \gamma \delta}$ define the flat holonomy 
class.

If we have a bundle 2-gerbe with trivialisation $\un{h}$ then we have the 
following:
\begin{eqnarray}
g_{\alpha \beta \gamma \delta} &=& h_{\beta \gamma \delta}
h_{\alpha \gamma \delta}^{-1} h_{\alpha \beta \gamma} 
h_{\alpha \beta \gamma}^{-1} \\
d\log h_{\alpha \beta \gamma} &=& A_{\alpha \beta \gamma} - k_{\beta \gamma}
+ k_{\alpha \gamma} - k_{\alpha \beta} \label{tb2g2}\\
\eta_{\alpha \beta} &=& -dk_{\alpha \beta} + j_\beta - j_\alpha \\
\nu_\alpha - dj_\alpha &=& \nu_\beta - dj_\beta
\end{eqnarray}
If we have a bundle 2-gerbe over a 3-manifold without boundary, $X$,  
then it is 
both flat and trivial. In this case we have a \v{C}ech - de Rham isomorphism 
as described by the following diagram:
\[
\begin{CD}
\nu_\alpha -  dj_\alpha\\
@AAdA \\
q_\alpha - j_\alpha @>\delta>> 0 \\
& & @AA-dA \\
& &  - k_{\alpha \beta} -B_{\alpha \beta} @>\delta>> 0 \\
& & & & @AAdA \\
& & & &  \log a_{\alpha \beta \gamma} - \log h_{\alpha \beta \gamma}
@>\delta>> \log c_{\alpha \beta \gamma \delta}
\end{CD}
\]
This tells us that the flat holonomy may be realised as an element 
of $S^1$ by the following formula
\begin{equation}
H(c_{\alpha \beta \gamma \delta}) = \exp \int_{X}\nu_\alpha - dj_\alpha
\end{equation}

This suggests the following
\begin{definition}
Let $(P,M;A,\eta,\nu)$ be a bundle 2-gerbe with connection and curvings. The 
{\it holonomy of $(P,M;A,\eta,\nu)$ over a 
closed 3-manifold $X$} with $\psi:X\rightarrow M$, is the flat holonomy 
of $\psi^* P$.
\end{definition}

Over
$X$ the bundle 2-gerbe $\psi^*P$ is trivial. We choose a trivialisation with
connection and curving. The 3-form defined by the difference between the 
3-curving induced by the pullback and the 3-curvature of the 
trivialisation may be integrated over $X$ to define the holonomy. 

Once again to find a corresponding formula in terms of the Deligne class
we shall need a triangulation, $t$, of $X$ which is subordinate to the 
open cover used to define the Deligne class. This triangulation
consists of tetrahedrons, faces, edges and vertices which are denoted by
$w$, $b$, $e$ and $v$ respectively. As usual we choose an index map $\rho$
with respect to the triangulation $t$ and the open cover of $M$.

Replacing the integral over $X$ with a sum of integrals over $w$,
\begin{equation*}
\begin{split}
H(c_{\alpha \beta \gamma \delta}) &= exp \sum_w \int_w (\nu_{\rho(w)} 
- d j_{\rho(w)}) \\
&= \exp \sum_w \int_w  \nu_{\rho(w)} + \int_{\partial w} -j_{\rho(w)} \quad
, \\
\exp \sum_w \int_{\partial w} -j_{\rho(w)} &= \exp \sum_{b,w} \int_b 
-j_{\rho(w)} \\
&= \exp \sum_{b,w} \int_b \eta_{\rho(w)\rho(b)} + dk_{\rho(w)\rho(b)}
- j_\rho(b) \\
&= \exp \sum_{b,w} \int_b \eta_{\rho(w)\rho(b)} + \int_{\partial b}
k_{\rho(w)\rho(b)}
\end{split}
\end{equation*}
where $\underset{b,w}{\sum} \int_b -j_{\rho(b)} = 0$ since each face bounds 
exactly
two tetrahedrons with opposite orientations. 
\begin{equation*}
\begin{split} 
\exp \sum_{b,w} \int_{\partial b} k_{\rho(w)\rho(b)} &=
\exp \sum_{e,b,w} \int_e k_{\rho(w)\rho(b)} \\
&= \exp \sum_{e,b,w} \int_e A_{\rho(w)\rho(b)\rho(e)} - d\log h_{\rho(w)
\rho(b)\rho(e)} + k_{\rho(w)\rho(e)} - k_{\rho(b)\rho(e)} \\
&= \exp \sum_{e,b,w} \int_e A_{\rho(w)\rho(b)\rho(e)} - \int_{\partial e}
\log h_{\rho(w)\rho(b)\rho(e)}
\end{split}
\end{equation*}
where $\underset{e,b,w}{\sum} \int_e k_{\rho(w)\rho(e)} - k_{\rho(b)\rho(e)}
= 0 $ since each edge of a particular tetrahedron bounds exactly two
faces of that tetrahedron and each edge of a particular face is an
edge of exactly two tetrahedrons and in both cases the corresponding
orientations are opposite. Finally,
\begin{equation*}
\begin{split}
\exp \sum_{e,b,w} \int_{\partial e} -\log h_{\rho(w)\rho(b)\rho(e)} &=
\exp \sum_{v,e,b,w} -\log h_{\rho(w)\rho(b)\rho(e)}(v) \\
&= \exp \sum_{v,e,b,w}  \log g_{\rho(w)\rho(b)\rho(e)\rho(v)}(v)
- \log h_{\rho(w)\rho(b)\rho(v)}(v) \\& 
\mspace{115.0mu}+ \log h_{\rho(w)\rho(e)\rho(v)}(v)
- \log h_{\rho(b)\rho(e)\rho(v)}(v) \\
&= \exp \sum_{v,e,b,w} \log g_{\rho(w)\rho(b)\rho(e)\rho(v)}(v)
\end{split}
\end{equation*}
where once again we get cancellation of terms due to opposite 
contributions as we sum over missing indices.

Collecting these results we have
\begin{equation}
\begin{split}
H(\un{g},\un{A},\un{\eta},\un{\nu}) &= \prod_w \exp \int_w   \nu_{\rho(w)}
\cdot \prod_{b,w} \exp \int_b \eta_{\rho(w)\rho(b)} \cdot \prod_{e,b,w}
\exp \int_e A_{\rho(w)\rho(b)\rho(e)} \\ & \mspace{300.0mu} 
\cdot \prod_{v,e,b,w} 
g_{\rho(w)\rho(b)\rho(e)\rho(v)}(v)
\end{split}
\end{equation}
and the corresponding formula for an embedding of a closed 3-manifold
$\psi: X \rightarrow M$ is 
\begin{equation}
\begin{split}
H((\un{g},\un{A},\un{\eta},\un{\nu});\psi) 
&= \prod_w \exp \int_w   \psi^* \nu_{\rho(w)}
\cdot \prod_{b,w} \exp \int_b \psi^* \eta_{\rho(w)\rho(b)} \cdot \prod_{e,b,w}
\exp \int_e \psi^* A_{\rho(w)\rho(b)\rho(e)} \\ & \mspace{330.0mu} 
\cdot \prod_{v,e,b,w} 
g_{\rho(w)\rho(b)\rho(e)\rho(v)}(\psi(v))
\end{split}
\end{equation} 

\section{A General Holonomy Formula} \label{genhol}
Using the results from the previous sections we can find a formula for
the holonomy of a class in $H^{p}(M,\D^p)$ associated with an 
embedding of a closed $p$-manifold $X$. Since we do not necessarily 
have a geometric realisation of this Deligne class in general, here 
holonomy is not meant in the traditional sense. It is defined purely
in terms of the Deligne class, specifically it is the flat holonomy
class of the pullback Deligne class on $X$, evaluated over $X$ as an element
of $S^1$. This formula gives a particular example of the even more
general transgression formula given by Gomi and Terashima (\cite{gote},
\cite{gote2}). The key feature of our derivation is that it clearly 
generalises the geometric notion of holonomy as we have defined it in the
low degree cases.  

\begin{definition}
Denote a Deligne class on $X$ by $(\un{g},\un{A}^1, \ldots , \un{A}^p)$. 
Think of this class as the pull back of a class on $M$. It 
is flat and trivial so there exists a cochain 
$(\un{h},\un{B}^1, \ldots , \un{B}^{p-1})$ such that 
\begin{equation}\label{flatriv}
\begin{split}
\un{g} &= \delta (\un{h}) \\
\un{A}^q &= \delta(\un{B}^q) + (-1)^{p-q} d\un{B}^{q-1} \\
\delta(\un{A}^p - d\un{B}^{p-1}) &=0 
\end{split}
\end{equation} 
The {\it holonomy of the Deligne class} is defined by
\begin{equation}
\exp \int_X \un{A}^p - d\un{B}^{p-1} 
\end{equation}
\end{definition}
This expression is not satisfactory since it depends explicitly on $\un{B}$.
To deal with this we triangulate $X$ with $t: |K| \rightarrow M$, where
$K$ is a $p$-dimensional simplicial complex, and let $\rho$ be an
index map for this triangulation. 
In terms of the triangulation the holonomy is 
\begin{equation}
\exp \left[ \sum_{\sigma^p} \int_{\sigma^p} A^p_{\rho(\sigma^p)}
+ \sum_{\sigma^p} \int_{\sigma^p} -dB^{p-1}_{\rho(\sigma^p)} \right]
\end{equation}
Consider the second term:
\begin{equation}
\sum_{\sigma^p} \int_{\sigma^p} -dB^{p-1}_{\rho(\sigma^p)}
= \sum_{\sigma^p} \int_{\partial \sigma^p} -B^{p-1}_
{\rho(\sigma^p)}
\end{equation}
In this expression we may express the combination
of the sum and the integral in terms of flags of simplices:
\begin{equation} 
\begin{split}
\sum_{\sigma^p} \int_{\partial \sigma^p} &= \sum_{\sigma^p} 
\sum_{\sigma^{p-1} \subset \sigma^p} \int_{\sigma^{p-1}} \\
&\equiv \sum_{\un{\sigma}^{p-1}} \int_{\sigma^{p-1}} 
\end{split}
\end{equation} 
where we have defined a new notation $\un{\sigma}$. In 
general this denotes a flag of simplices,
\begin{equation}
\un{\sigma}^q = \{ (\sigma^q, \sigma^{q+1}, \ldots , \sigma^p) |
\sigma^q \subset \cdots \subset \sigma^p \}
\end{equation} 
All subsimplices inherit relative orientations.
A similar notation was used in \cite{gote2} to generalise transgression
formulae.

Returning to the holonomy formula, we now have
\begin{equation}
\sum_{\un{\sigma}^{p-1}} \int_{\sigma^{p-1}} -B^{p-1}_
{\rho(\sigma^p)}
\end{equation}
Now use equation \eqref{flatriv},
\begin{equation}
\delta(\un{B}^q)  = \un{A}^q - (-1)^{p-q} d\un{B}^{q-1}
\end{equation}
to get 
\begin{equation}
-B^{p-1}_{\rho(\sigma^p)} = -B^{p-1}_{\rho(\sigma^{p-1})} +
A^{p-1}_{\rho(\sigma^{p})\rho(\sigma^{p-1})} - 
dB^{p-2}_{\rho(\sigma^{p})\rho(\sigma^{p-1})} 
\end{equation}
Using the fact that each $(p-1)$-face in the simplicial 
complex bounds exactly two $p$-faces we have
\begin{equation}
\sum_{\un{\sigma}^{p-1}} \int_{\sigma^{p-1}} 
-B^{p-1}_{\rho(\sigma^{p-1})} = 0
\end{equation}
since the two terms inherit opposite orientations from
$\sigma^p$. Thus 
\begin{equation}
\sum_{\un{\sigma}^{p-1}} \int_{\sigma^{p-1}} -B^{p-1}_{\rho(\sigma^p)}
= \sum_{\un{\sigma}^{p-1}} \int_{\sigma^{p-1}}
A^{p-1}_{\rho(\sigma^{p})\rho(\sigma^{p-1})} + 
dB^{p-2}_{\rho(\sigma^{p})\rho(\sigma^{p-1})} 
\end{equation}
The next step would be to extract the $A^{p-1}$ term for the final 
answer and proceed as above to deal with the $dB^{p-2}$ term.
This suggests an inductive approach with respect to $k = p-q$.
\begin{lemma}\label{dblemma}
For every $q$ such that $1 \leq q \leq p$ 
\begin{equation}
\sum_{\un{\sigma}^q} \int_{\sigma^q} dB^{q-1}_{\rho(\sigma^p)\ldots\rho
(\sigma^q)} = \sum_{\un{\sigma}^{q-1}} \int_{\sigma^{q-1}} (-1)^{p-q+1}
A^{q-1}_{\rho(\sigma^{p})\ldots \rho(\sigma^{q-1})} - 
dB^{q-2}_{\rho(\sigma^{p})\ldots 
\rho(\sigma^{q-1})}
\end{equation}
where we use the conventions $\un{A}^0 = \log \un{g}$, 
$\un{B}^0 = \log\un{h}$ and $\un{B}^{-1} = 0$.
\end{lemma}
\begin{proof}
We have already proved the particular case $p=q$.
More generally 
\begin{equation}
\begin{split}
\sum_{\un{\sigma}^q} \int_{\sigma^q}  dB^{q-1}_{\rho(\sigma^p)\ldots\rho
(\sigma^q)} &= \sum_{\un{\sigma}^q}\int_{\partial\sigma^q}
B^{q-1}_{\rho(\sigma^p)\ldots\rho(\sigma^q)} \\
&= \sum_{\un{\sigma}^q}\sum_{\sigma^{q-1} \subset \sigma^q}\int_
{\sigma^{q-1}} B^{q-1}_{\rho(\sigma^p)\ldots\rho(\sigma^q)} \\
&= \sum_{\un{\sigma}^{q-1}} \int_
{\sigma^{q-1}} B^{q-1}_{\rho(\sigma^p)\ldots\rho(\sigma^q)}
\end{split}
\end{equation}
Next we claim that
\begin{equation} \label{deltaclaim}
\sum_{\un{\sigma}^{q-1}} \int_
{\sigma^{q-1}} B^{q-1}_{\rho(\sigma^p)\ldots\rho(\sigma^q)}
= \sum_{\un{\sigma}^{q-1}} \int_{\sigma^{q-1}} (-1)^{p-q+1}(\delta B^{q-1})
_{\rho(\sigma^{p}) \ldots \rho(\sigma^{q-1})}
\end{equation}
The right hand side consists of all terms of the form
\begin{equation}
\sum_{\un{\sigma}^{q-1}} \int_
{\sigma^{q-1}} (-1)^{p-q+1} 
B^{q-1}_{\rho(\sigma^p)\ldots \widehat{\rho(\sigma^k)}
\ldots \rho(\sigma^q)}
\end{equation}
for all $q-1 \leq k \leq p$ and where the hat symbol denotes that a 
subscript should be omitted. The case $k = q-1$ corresponds to the left hand
side of \eqref{deltaclaim}. \\ Now consider $q-1 < k < p$.
Suppose in the summation we have a flag $(\sigma^{q-1},\ldots,\sigma^p)$,
with the summand depending on all simplices in the flag except for 
$\sigma^k$. This leads to a number of identical terms corresponding to
all flags which agree in all degrees except for $k$. There can only be 
two such flags. This is because such flags must satisfy 
\begin{eqnarray}
\sigma^k &\subset \sigma^{k+1} \\
\sigma^{k-1} &\subset \sigma^k
\end{eqnarray}
This means that $\sigma^k$ is defined by $k+1$ of the $k+2$ vertices of
$\sigma^{k+1}$ and $\sigma^{k-1}$ is defined by $k$ of these. Since
$\sigma^{k+1}$ and $\sigma^{k-1}$ are fixed then there are only two choices
for $\sigma^k$ as there are two vertices in $\sigma^{k+1}$ which are not in
$\sigma^{k-1}$. Furthermore the two possible choice of flags will lead to
opposite induced orientations of $\sigma^{q-1}$. The induced orientations 
are derived from the orientation of $\sigma^p$. The orientations of all
the simplices from $\sigma^p$ to $\sigma^{k+1}$ must be the same since they
are all identical. The two choices for $\sigma^k$ must give opposite 
orientations for $\sigma^{k-1}$. This condition is equivalent to the
basic result $\partial^2 =0$ for the boundary operator 
in the theory of simplicial complexes. 
From $\sigma^{k-1}$ down to $\sigma^{q-1}$ all
of the simplices are equal so there can be no further change in the 
relative orientations of the two choices.

Finally we consider the case $k=p$. In this case we once again have only two
choices of flag corresponding to the two choices of orientation and these
contribute terms of opposite sign. This proves the claim.

The lemma now follows from equation \eqref{flatriv}.

\end{proof}
This lemma leads to the following
\begin{proposition}\label{propgh}
For all $p \geq 1$ the holonomy of the Deligne class 
$(\un{g},\un{A}^1, \ldots , \un{A}^p)$
is given by the following formula:
\begin{equation}\label{ghform}
\exp \int_X \un{A}^p - d\un{B}^{p-1} =
\exp \sum_{n=0}^p \sum_{\un{\sigma}^{p-n}} \int_{\sigma^{p-n}} 
 A^{p-n}_{\rho(\sigma^{p})\ldots
\rho(\sigma^{p-n})}
\end{equation} 
\end{proposition}
As before we let $\un{A}^0 = \log \un{g}$. 
\begin{proof}
It is easily verified that the formulae 
obtained in the previous sections of this chapter prove the result for
$p = 1$, 2 and 3. To prove the more general case we use the following 
intermediate result:
\begin{equation}\begin{split}\label{induction}
\exp \int_X \un{A}^p - d\un{B}^{p-1} &=
\exp (\sum_{n=0}^k \sum_{\un{\sigma}^{p-n}} \int_{\sigma^{p-n}} 
 A^{p-n}_{\rho(\sigma^{p})\ldots \rho(\sigma^{p-n})}
)\\
& \mspace{100.0mu} \cdot \exp
\sum_{\un{\sigma}^{p-k}} \int_{\sigma^{p-k}} (-1)^{k+1}
dB^{p-k-1}_{\rho(\sigma^{p})\ldots \rho(\sigma^{p-k})}
\end{split}
\end{equation}
For $k=0$ this is simply rewriting the integral over $X$ in terms of the
triangulation. We prove the general case, $0< k \leq p$ by induction.
Suppose \eqref{induction} is true for some $k < p$. Applying
Lemma \ref{dblemma} to the $dB$ term gives
\begin{equation}\label{k+1}
\begin{split}
\exp
\sum_{\un{\sigma}^{p-k}} \int_{\sigma^{p-k}} (-1)^{k+1}
dB^{p-k-1}_{\rho(\sigma^{p})\ldots \rho(\sigma^{p-k})} &=
\exp \sum_{\un{\sigma}^{p-k-1}} \int_{\sigma^{p-k-1}}
(-1)^{k+1}(-1)^{k+1} A^{p-k-1}_{\rho(\sigma^{p})\ldots
\rho(\sigma^{p-k-1})} \\
& \mspace{70mu} - (-1)^{k+1} dB^
{p-k-2}_{\rho(\sigma^{p})\ldots \rho(\sigma^{p-k-1})} \\
&\mspace{-60mu} = \exp \sum_{\un{\sigma}^{p-(k+1)}} \int_{\sigma^{p-(k+1)}}
A^{p-(k+1)}_{\rho(\sigma^{p})\ldots
\rho(\sigma^{p-(k+1)})} \\
& \mspace{70mu} + (-1)^{(k+1)+1} dB^
{p-(k+1)-1}_{\rho(\sigma^{p})\ldots \rho(\sigma^{p-(k+1)})} \\
\end{split}
\end{equation}
Substituting \eqref{k+1} back into \eqref{induction} gives
\begin{equation}\begin{split}
\exp \int_X \un{A}^p - d\un{B}^{p-1} &=
\exp (\sum_{n=0}^{k+1} \sum_{\un{\sigma}^{p-n}} \int_{\sigma^{p-n}} 
A^{p-n}_{\rho(\sigma^{p})\ldots \rho(\sigma^{p-n})}
)\\
& \mspace{70.0mu} \cdot \exp
\sum_{\un{\sigma}^{p-(k+1)}} \int_{\sigma^{p-(k+1)}} 
(-1)^{(k+1)+1}
dB^{p-(k+1)-1}_{\rho(\sigma^{p})\ldots \rho(\sigma^{p-(k+1)})}
\end{split}
\end{equation}
thus the statement is true for $k+1$ and therefore by induction is true for 
all $1 \leq k \leq p$. In particular the case $k = p$ is equivalent to
the statement of the proposition since $\un{B}^{-1} = 0$, 
thus this is sufficient to prove the proposition.

\end{proof}

\section{Transgression for Closed Manifolds}\label{trgn}
Consider the constructions of the previous sections of this chapter. In 
each case we start with a bundle $(n-1)$-gerbe with curving ($n$ = 1, 2 or 3).
Then we construct
an element of $S^1$ corresponding to a smooth mapping of a closed manifold
of dimension $n$. Furthermore for $n > 3$ we can carry out this construction
purely in terms of the Deligne class. We would like to consider 
the holonomy as a smooth function
on the infinite dimensional manifold $\Map(X,M)$. We give this mapping
space the compact-open smooth topology \cite[p34]{hir}. 
Since the holonomy is defined in terms of sums, integrals and pull backs 
it will define a smooth, continuous function on $\Map(X,M)$. 
To see that it defines a class in Deligne cohomology consider the 
following open cover of the mapping space:
\begin{definition}
Let $\U \equiv \{ U_\alpha \}_{\alpha \in \A}$ be an open cover of $M$.
Let $t$ be a triangulation of $X$ consisting of simplices $\sigma$ and 
suppose we have an index map 
$\rho: t \rightarrow \A$. Then the set $V_{(t,\rho)}$ is
defined by
\begin{equation}
V_{(t,\rho)} = \{ \phi \in \Map(X,M) | \quad 
\phi(\sigma
) \subset U_{\rho(\sigma)} \}
\end{equation}
Denote open cover defined by these sets by $\V$
\end{definition}
These sets are open in the compact-open smooth topology since they are 
made up of smooth maps of simplices (which
are compact) into open sets in $M$.
Following \cite{gaw} we use $\V$ as our open cover of $\Map(X,M)$. 
We have already used this cover to calculate the holonomy,
so we may think of the holonomy as a collection of $S^1$ functions
defined on open sets in $\V$, that is, a cochain in 
$C^0(\Map(X, M),\un{U(1)})$. The fact that our construction was independent
of the choice of the pair $(t,\rho)$ implies that this cochain is
actually a cocycle in $H^0(\Map(X,M),\un{U(1)})$.
Following \cite{gaw} and \cite{bry} we define the {\it transgression} 
homomorphism 
$\tau_X : H^{n}(M,\D^n) 
\rightarrow H^0(\Map(X, M),\un{U(1)})$.

This homomorphism has been interpreted (\cite{bry},\cite{gote2}) as 
a composition of an evaluation map 
\[
ev^* : H^{n}(M,\D^n) \rightarrow H^{n}(\Map(X,M) \times 
X,\D^n) 
\]
and a fibre integration map
\[
\int_X : H^{n}(\Map(X,M) \times 
X,\D^n) \rightarrow H^0(\Map(X,M),\un{U(1)}).
\]
This homomorphism is compatible with the corresponding map on curvatures,
that is, if the curvature of the Deligne class on $M$ is $\omega$ then
the curvature of the transgressed class on $\Map(X,M)$ is 
$\int_X ev^* \omega$. 

To see that this agrees with our constructions of the preceding
sections suppose that
\mbox{$(\un{g},\un{A^1},\ldots,\un{A^n}) \in
H^{n}(M,\D^n)$}. Pulling back by the evaluation map gives the class
$(ev^*\un{g},ev^*\un{A^1},\ldots,ev^*\un{A^n})$. 
The pull back of the evaluation map gives a 
homomorphism in cohomology.
Restricted to a fixed $\psi \in \Map(X,M)$ this class is equal to
$(\psi^*\un{g},\psi^*\un{A^1},\ldots,\psi^*\un{A^n})$ which represents
a flat bundle $(n-1)$-gerbe on $X$. The fibre integration map
evaluates the flat holonomy for each value of $\psi$. 
It was proven in \cite{gote2} that the fibre integration map
is also a homomorphism. 

In conclusion, we have developed the geometric notion of holonomy from
the familiar case of line bundles to the case of bundle gerbes and 
bundle 2-gerbes. The generalisation was guided by the consideration of  
holonomy 
as a property of the Deligne class, specifically as the evaluation of 
the flat holonomy class of the pullback of the Deligne class to a 
closed manifold of appropriate dimension. The relationship between these
cohomological and geometric concepts was demonstrated. 
As a property of Deligne 
cohomology holonomy could be extended to higher degree classes and also 
considered as an example of the more general notion of a transgression 
homomorphism.

%% file: chapter6.tex
\chapter{Parallel Transport and Transgression with Boundary}

In this chapter we investigate what happens to the constructions of
the previous chapter when we consider manifolds with boundary. This
leads to generalised notions of parallel transport. These results may also
be viewed in terms of an extension of the transgression homomorphism to
manifolds with boundary.

We shall see that parallel transport may be thought of as a section of
a trivial bundle. 
Let $ev_t : \P M \rightarrow M$ be the evaluation
map that takes $\mu \in \P M$ to $\mu(t) \in M$. 
The the parallel
transport map, which is a map between the fibres over 
$\mu(0)$ and $\mu(1)$ of a bundle $L$ may be thought of as an 
element of $L^*_{\mu(0)} \otimes L_{\mu(1)}$. This is 
the same as a section of the bundle $(ev^*_0L)^* \otimes (ev_1^* L)$
on $\P M$. We shall see that such a section arises from the extension
of holonomy from loops to paths.

This approach to parallel transport will lead to a similar interpretation
in the case of bundle gerbes. Here we consider a surface with boundary 
made up of loops. The parallel transport is now defined by pulling 
back a bundle on the loop space to the space $\Map(\Sigma,M)$ of maps
of the surface into $M$. The holonomy of a closed surface generalises to 
give a section of this bundle. For the example of the cylinder this 
construction gives a map between fibres over the two end loops, a
situation similar to parallel transport for bundles, however there
is no problem considering surfaces with different topologies. This 
construction will give a geometric interpretation to Gawedski's 
results on holonomy of classes in $H^2(M,\D^2)$ over surfaces
with boundary \cite{gaw}.

\section{Parallel Transport for Bundles} \label{PTbun}
We now attempt to calculate the holonomy of a bundle over a path in the 
same way in which we calculated the holonomy over a loop. Once again
it may seem that we are using a long and unnecessarily complicated
method, however there are good reasons for this. Firstly it is not
an unreasonable assumption that a method of describing holonomy which
generalises to higher cases should be a good starting point for 
generalising parallel transport. It turns out that there are a number
of different ways of approaching this, and since these appear 
in the literature it is worthwhile seeing how they arise in this context 
and how they relate to each other. We give as much detail as possible at 
the level of bundles since the key features of the theory are 
present, but relatively easy to deal with compared with the bundle gerbe
case.
 
We consider
a path as a smooth map $\mu: I \rightarrow M$ where $I$ is the unit
interval $[0,1] \in \R$. The path space $\P M$ is the space which 
consists of all such maps, we give it the compact-open smooth 
topology (see \S \ref{trgn}).

As in the previous case the pull back of any
bundle to $I$ is flat and trivial. The flat holonomy class is 
an element of $H^1(I,U(1)) = 0$, thus we cannot evaluate a holonomy in the
same sense as the case without boundary. If we cannot define holonomy then
we would like to define something which is as close as possible to holonomy.
This turns out to be parallel transport. For motivation let us consider
the case of a principal bundle. Parallel transport assigns to each
path $\mu$ a $U(1)$-equivariant
map between the fibres over $\mu(0)$ and $\mu(1)$. When we have a loop
this gives an equivariant map from a fibre to itself which is of the
form $p \mapsto pz$ for some $z\in U(1)$ which is the holonomy. In 
general the parallel transport map takes $p$ to $\tilde{\mu}(1)$ where
$\tilde{\mu}$ is a horizontal lift of $\mu$ satisfying $\tilde{\mu}(0)
=p$.
Another view is that parallel
transport satisfies the condition that given any two paths which 
may be joined to form a loop then composing the respective parallel 
transports gives the holonomy. In our case the holonomy is given
by $\exp \int_{\gamma} \chi$ where $\chi$ is a $D$-obstruction form
given by $A_\alpha - d\log h_\alpha$. We have shown that this is equivalent
to a formula $H(\un{g},\un{A})$ in terms of the Deligne class.
These two definitions of the function on the loop space lead to two 
equivalent ways of defining a function on the path space which satisfies the
required criteria.

Given two paths $\mu_1, \mu_2$ with the same endpoints we may define a 
loop, by convention we define this loop associated with a pair $(\mu_1,\mu_2)$
to be the composition $\mu_1 \star \mu_2^{-1}$.
We may define a map on the path space by $H_B(\mu) = \exp \int_\mu \chi$, so
$H_B(\mu_1)H_B^{-1}(\mu_2) = H(\mu_1 \star
\mu_2^{-1})$. It is important to remember that in this case $\chi$ is no longer
a $D$-obstruction form so we cannot be sure that the construction is 
independent of the choice of the trivialisation $h$. 
This is because in the construction
of holonomy the value at a particular loop is given by the flat holonomy
which is a property of flat bundle 0-gerbes. The triviality is only used to 
express it as a differential form, which turns out to be the 
$D$-obstruction, a property
of trivial bundles. In fact we find that the map $H_B$ does depend on the 
choice of trivialisation:
\begin{equation} \label{trivfn}
\begin{split}
H_B(\mu ) &= \exp \sum_e \int_e \mu^*A_{\rho(e)} \cdot \prod_{v,e}
g_{\rho(e) \rho(v)}(\mu(v)) \cdot \prod_{v,e} h^{-1}_{\rho(v)}(\mu(v)) \\
&= \exp \sum_e \int_e \mu^*A_{\rho(e)} \cdot \prod_{v,e}
g_{\rho(e) \rho(v)}(\mu(v)) \cdot h_{\rho(v_0)}(\mu(0))h^{-1}_{\rho(v_1)}
(\mu(1))
\end{split}
\end{equation}
where $v_0$ and $v_1$ are the endpoint vertices of the triangulation of $I$.
The final term fails to cancel this time because of contributions from
$\mu(0)$ and $\mu(1)$. This expression is 
independent of the choice of $\rho$, but the dependence on a choice of
trivialisation causes difficulties. We would prefer to have
a $\rho$ dependence instead, this could be achieved by restricting to a 
particular
$\rho$ and then choosing a trivialisation, however this is rather technical.
It is easy enough to choose a trivialisation over each element of the 
triangulation using the canonical trivialisation over an element of 
the open cover which was described in \S\ref{trivial}. The problem is
that these trivialisations then need to be glued together in some way. 
This is possible but will become even more complicated in higher degrees
(see for example the proof of Proposition 6.5.1 in \cite{bry}), so is
not suitable for our purposes.

Instead, let us turn to another expression for the holonomy, the
transgression formula of proposition \ref{holprop},
\begin{equation}\label{hmu}
H_{(t_0,\rho_0)}(\mu) =
 \prod_e \exp \int_e \mu^*A_{\rho_0(e)}  \cdot
 \prod_{v, e} g_{\rho_0(e)\rho_0(v)}(\mu(v))
\end{equation}
It is easily shown that on paths with boundary this function
is  not independent of the choice of $\rho_0$ so it
is not globally defined on the path space. 
This dependence on the triangulation means that we have
to be careful about describing the open covers
on the loop space and path space. Given a decomposition of a loop into
two paths with the same boundary where the loop is considered as 
an element of some $V_{(t_0,\rho_0)} \subset \L M$ there is an
inherited triangulation and corresponding element of an open cover of 
$\P M$ for each path. In general the two paths 
do not lie in the same open set in terms of this cover, however it turns
out that there is a more appropriate cover on the path space in this 
situation. Suppose that
when $\gamma \in V_{(t_0,\rho_0)}$ is split into paths 
$\mu_1$ and $\mu_2$ the induced open sets on $\P M$ are
$W_{(t_1,\rho_1)}$ and
$W_{(t_2,\rho_2)}$ respectively. Since they are both induced from the
same cover on the loop space they must satisfy the condition 
$\rho_1(v) = \rho_2(v)$ for $v\in \partial \mu_1 (= \partial \mu_2)$.
In this case we have 
\begin{equation}\label{hpath}
H(\mu_1 \star \mu_2^{-1}) = H_{(t_1,\rho_1)}(\mu_1)H_{(t_2,\rho_2)}^{-1}(\mu_2)
\end{equation}
where we have just broken up the expression for the holonomy \eqref{hloop}
into the parts corresponding to each path. 

Since we are using the same formula for both cases we use $H$ for
for both loops and paths. The distinction should be clear in all cases
from the argument and the fact that the function on the path space
has a local dependence indicated by a subscript.


Consider now the case where we are given two paths which share a boundary.
When is equation \eqref{hpath} satisfied? The fact that $H$ is independent of
the open cover suggests that this equation is satisfied whenever the 
two covers on the path space combine to form a cover on the loop space, that
is, precisely whenever $\rho_1(v) = \rho_2(v)$ on the boundary of the paths.
Suppose that $\mu_2$ also lies in the open set $W_{(t_3,\rho_3)}$ and 
$\rho_3(v) = \rho_2(v)$ on the boundary. Then using \eqref{hpath}
\begin{equation}
\begin{split}
H_{(t_2,\rho_2)}(\mu_2) &= H^{-1}(\mu_1 \star \mu_2^{-1}) H_{(t_1,\rho_1)}
(\mu_1) \\
&= H^{-1}_{(t_1,\rho_1)}(\mu_1)H_{(t_3,\rho_3)}(\mu_2)H_{(t_1,\rho_1)}
(\mu_1) \\
&= H_{(t_3,\rho_3)}(\mu_2) 
\end{split}
\end{equation}
This result may also be seen by considering what happens when 
we change $\rho_0(e)$ to $\rho_0'(e)$ for any $e$ in the triangulation for 
$\mu$ in the formula for $H_0(\mu)$. The two expressions for $h_0(\mu)$ agree
except on the terms corresponding to $e$ where the difference is
\begin{equation}
\int_e (A_{\rho(e)} - A_{\rho'(e)}) \cdot \prod_{v \subset \partial e} 
g_{\rho(e)\rho(v)}g^{-1}_{\rho'(e)\rho'(v)} 
= \prod_{v \subset \partial e} g_{\rho'(e)\rho(e)} 
g_{\rho(e)\rho(v)}g^{-1}_{\rho'(e)\rho'(v)} = 1
\end{equation}
Thus we may use a coarser open cover on $\P M$, the cover induced
by the projection to $M \times M$. An open set $\tilde{W}_{\rho_0,t_0}$ in
this cover consists of all paths $\mu$ with triangulation $t_0$ such that
the endpoints of $\mu$ are vertices $v_0$ and $v_1$ satisfying
$\mu(0) \in U_{\rho_0(v_0)}$ and $\mu(1) \in U_{\rho_0(v_1)}$. In terms
of this cover we shall denote the functions by $H_0$, $H_1$ and so on.
We now have a set of locally defined functions on $\P M$ from which
we may recover the holonomy in the required manner. It is natural to
ask what the obstruction is to these defining a global function. Another
way of viewing this is that given a set of local functions, on overlaps we
may define transition functions for a trivial bundle, with trivialisations
(or equivalently sections) defined by the original functions. Using a 
similar calculation to those performed in the previous chapter we find
\begin{equation}\label{tauL}
(H^{-1}_{(t_0,\rho_0)}H_{(t_1,\rho_1)})(\mu)
= \prod_{v,\partial e} g_{\rho_0(v)\rho_1(v)}(\mu(v))
\end{equation}
Since the right hand side depends only on the boundary of $\mu$ it may 
be written as $r^*G_{(t_0,\rho_0)(t_1,\rho_1)}$ where $r$ is the restriction
to the boundary $\partial \mu$ and $G_{(t_0,\rho_0)(t_1,\rho_1)}$ are defined
on $\Map(\partial I,M)$. By applying $d\log$ to $H_{(t_0,\rho)}$ it is possible
to obtain a formula for local connection 1-forms. 
To distinguish the differential on $\P M$ from $d$ on 
$M$ we shall denoted it $\tilde{d}$ (this notation will carry over to
other spaces of maps into $M$ as well). Let $\xi \in T_{\mu}(\P M)$.
\begin{equation}
\begin{split}
(\tilde{d}\log H_0) (\xi) 
&= \sum_e \int_e \mu^* (d\iota_\xi A_{\rho_0(e)} + \iota_\xi d A_{\rho_0(e)})
+ \sum_{v,e} \iota_\xi(v) d\log g_{\rho_0(e)\rho_0(v)}(\mu(v)) \\
&= \sum_{v,\partial e} \mu^* \iota_{\xi(v)} A_{\rho_0(v)}
\end{split}
\end{equation}
Since this depends only on the boundary we may write it as $r^* B_0$
where $B_0$ is a local one form on $\Map(\partial I,M)$.
We denote the trivial bundle with connection on $\P M$ by $D(\tau_I L)$ 
(so $\tau_I L$ is the trivialisation $H_0$) and
the bundle with connection represented by local data $(G_{01},B_0)$ is
denoted by $\tau_{\partial I} L$.

We would like to consider a more global description of the bundles 
$D(\tau_I L)$ and 
$\tau_{\partial I} L$. This is given by the following diagram:
\[
\begin{array}{ccccc}
& &  S^1 \quad & & D(\tau_I L)\\
& \stackrel{H}{\nearrow} & & \swarrow &\\
\L M & \rightrightarrows & \P M & & \\
& & \quad \downarrow \partial & & \\ 
& & M \times M & & 
\end{array}
\]
Here we are claiming that the bundle 0-gerbe described locally by the 
construction of $G_{01}$ is the same
as that represented by the above diagram. We shall consider this
as a particular example of a general result. Let $(\lambda,Y,M)$ be
a bundle 0-gerbe and let $H_\alpha$ be locally defined functions on
$\pi^{-1}(U_\alpha) \subset Y$ such that $H_{\alpha}(y_1)H^{-1}_\alpha(y_2)
= \lambda(y_1,y_2)$. Then the transition functions of the bundle
0-gerbe are given by $g_{\alpha\beta}(m) = H^{-1}_\alpha(y)H_\beta(y)$ for
any $y \in \pi^{-1}(m)$. Observe that 
\begin{equation}
\begin{split}
\lambda(s_\alpha(m),s_\beta(m)) &= H_\alpha(s_\alpha)
H^{-1}_\alpha(s_\beta) \\
&= H_\alpha(s_\alpha) H^{-1}_\beta(s_\beta) H_\beta(s_\beta)
H^{-1}_\alpha(s_\beta) \\
&= H_\alpha(s_\alpha) H^{-1}_\beta(s_\beta) g_{\alpha\beta}(m) \\
&= \delta (s_\alpha^* H_\alpha)(m)g_{\alpha\beta}(m)
\end{split}
\end{equation}
The bundle 0-gerbe defined by this diagram is $\tau_{\partial I} L$ and
as we have seen it is $D$-stably isomorphic to $m_0^{-1}L^*
\otimes m_1^{-1}L$ where $m_0$ and $m_1$ are the projections of the 
two components of $M\times M$ onto $M$. 
This is related to the the holonomy reconstruction 
theorem as described in \cite{mapi}, which says that a bundle with
connection may be reconstructed, up to isomorphism, from the function
on the loop space defined by holonomy. We shall consider this theorem and 
bundle gerbe generalisations of it in the next chapter.

The bundle $\tau_I L$ is canonically trivial, with local trivialisation
functions $H_{\rho}$ such that $H_{\rho}(\mu_1)H_{\rho}^{-1}(\mu_2) 
=  H  (\mu_1 \star \mu^{-1}_2)$. It is also the canonically trivial
bundle obtained by pulling back $\tau_{\partial I} L$ to $\P M$ 
with the map $\partial$, the restriction to the boundary. 
The connection on $\tau_I L$ is given locally on $\P M$ by
$\tilde{d} \int_\mu \chi$.

There is another approach to this problem. 
Following Hitchin \cite{hit} we may consider the
space of trivialisations of $L$ over $\Map(I,M)$, 
which we shall denote $\triv_I$. Each element of $\triv_I$ has a   particular
path associated with it, giving a projection map onto $\Map(I,M)$.
Using this we may calculate the function $\delta(H): {\triv_I}^{[2]}
\rightarrow S^1$ in the following way,
\begin{equation}\label{deltaH}
\delta(H)(H^0_{(t,\rho)},H^1_{(t,\rho)};\mu) =
\sum_{v,e} \log H^1_{\rho(v)}(\mu(v)) - \log H^0_{\rho(v)}(\mu(v))
\end{equation}
Recall that any two trivialisations of a bundle differ by a function.
Equation \eqref{deltaH} gives the `holonomy' of the function
defined by trivialisations $h_\rho^0$ and $h_\rho^1$ over $\partial \mu$. 
Thus once again we see that this bundle 0-gerbe is pulled back from
$\Map(\partial I,M)$. Over $\Map(\partial I,M)$ we have a space
of trivialisations $\triv_{\partial I}$ and \eqref{deltaH} gives
a function on ${\triv_{\partial I}}^{[2]}$ which defines the bundle 
0-gerbe. Note that this function is no longer of the form $\delta(H)$ 
since $H$ is not defined on $\triv_{\partial I}$. By forming 
the bundle corresponding to this bundle gerbe we obtain the lower
dimensional version of the moduli space of flat trivialisations (\cite{hit})
as the total space.
The difficulty with this approach is having to deal with the space of
trivialisations. The space of trivialisations of a bundle is the 
infinite dimensional set $\Map(M,U(1))$, however the space 
of trivialisations of a bundle gerbe is a collection of 
all line bundles on $M$ which is not a set and would have
to be considered in terms of category theory.   

In \cite{gote} it was proved that the local transgression formulae that
we have described above correspond to the usual
notion of parallel transport for bundles. We would like to describe
how the local functions $H_{(t_0,\rho_0)}$ on the
loop space lead to the parallel transport map.

Recall that given a bundle with connection $L \rightarrow M$ and a path 
$\mu$ in $M$,
parallel transport is a $U(1)$-equivariant map of fibres $PT : L_{\mu(0)} 
\rightarrow L_{\mu(1)}$ which
is defined by the unique lifting of $\mu$ to $L$ which is horizontal
with respect to the connection. Transgression defines a trivialisation of
$r^{-1} L = L_{\mu(0)}^* \otimes L_{\mu(1)}$. This means we have a global
section on $\Map(I,M)$ which assigns an element of $L_{\mu(0)}^* \otimes 
L_{\mu(1)}$ to each path, but this may also be interpreted as a $U(1)$-
equivariant map 
$L_{\mu(0)} \rightarrow L_{\mu(1)}$ which defines the parallel transport.


\section{Loop Transgression of Bundle Gerbes}\label{loopbg}

We consider what happens to the holonomy of a bundle gerbe over $M$ 
when we allow surfaces with boundary.
We obtain a section of a trivial line bundle over
$\Map(\Sigma^\partial, M)$ which gives a generalisation of parallel transport.
Furthermore this is the pullback via restriction to the 
boundary of a possibly non-trivial
line bundle over $\Map(\partial \Sigma ,M)$. This may be related to a 
line bundle over the loop space. We derive a local
formula which gives a transgression homomorphism 
$\tau_{S^1} : H^2(M,\D^2) \rightarrow H^1(\L M,\D^1)$.
Throughout we abbreviate $\Map(X,M)$ by $XM$ for 
various manifolds $X$.

We already have a formula \eqref{holbg}
which gives an $S^1$-function
over the space of smooth mappings of closed surfaces into $M$. This corresponds
to pulling back the bundle gerbe and evaluating the flat holonomy as 
an element of $H^2(\Sigma,U(1)) = S^1$. In the case with boundary
we proceed as in the case of parallel transport. 

Starting with $\exp \int_\Sigma \eta - dk$ and following the 
same procedure as for closed surfaces we find that the terms
involving the trivialisation $(\un{h},\un{k})$ do not 
cancel out on the boundary components. Given a choice $\rho_0$ 
we could use the canonical trivialisation over each $U_{\rho_0(v)}$ and
$U_{\rho_0(e)}$. The problem here, unlike in the previous case where 
boundary components were just points, is that each boundary component is 
a loop and will
consist of a number of edges and vertices, so to define a trivialisation
over the whole loop it would be necessary to glue together each of these
in some way. This is possible (see \cite{bry} for a description of this 
in the gerbe case), however it is not a method that will be suitable
for generalisation to higher degree as it becomes very complicated. Instead 
we turn to the second method that was developed in the previous 
section. 

We define a function
on $\Sigma^\partial M$ by $H_{(t_0,\rho_0)}(\psi)$,
\begin{equation} 
H_{(t_0,\rho_0)}((\un{g},\un{A},\un{\eta});\psi) =  
\prod_b \exp \int_b \psi^*\eta_{\rho_0(b)} \cdot
\prod_{e,b} \exp \int_e \psi^*A_{\rho_0(b)\rho_0(e)} \cdot
\prod_{v,e,b}  g_{\rho_0(b)\rho_0(e)\rho_0(v)}(\psi(v))
\end{equation} 
the usual holonomy formula which is well-defined 
but not globally defined on $\Sigma^\partial M$. This is the 
approach taken by Gawedski \cite{gaw}.
Clearly for two surfaces with the same boundary this `trivialises' the 
holonomy function, but it is not a proper trivialisation since it
is not globally defined. These local functions define a $D$-trivial bundle
$D(\tau_{\Sigma^\partial} P) \rightarrow \Sigma^\partial M$.
As in the previous case we can define a bundle 
0-gerbe $\tau_{\partial \Sigma} P$ on the space of mappings of the boundary,
\begin{equation}\label{trdiag}
\begin{array}{ccccc}
& & S^1 & & D(\tau_{\Sigma^\partial} P) \\
& \stackrel{H}{\nearrow} & & \swarrow & \\
\Sigma M & \rightrightarrows & \Sigma^\partial M & & \\
& & \downarrow & & \\
& & \partial \Sigma M & & 
\end{array}
\end{equation}
The transition functions of $D(\tau_{\Sigma^\partial})$ descend to
$G_{(t_0,\rho_0)(t_1,\rho_1)}$ on $\partial \Sigma M$, the transition 
functions of $\tau_{\partial \Sigma} P$. We may calculate these explicitly,
\begin{equation}
\begin{split}
G_{(t_0,\rho_0)(t_1,\rho_1)} &= H^{-1}_{(t_0,\rho_0)} H_{(t_1,\rho_1)} \\
&= \exp (\sum_b \int_b (\eta_{\rho_1(b)} - \eta_{\rho_0(b)})+ \sum_{e,b} 
(A_{\rho_1(b)\rho_1(e)} - A_{\rho_0(b)\rho_0(e)})) \\ & \mspace{250mu}
\prod_{v,e,b} g_{\rho_1(b)\rho_1(e)\rho_1(v)}g^{-1}
_{\rho_0(b)\rho_0(e)\rho_0(v)} \\
&= \exp (\sum_b \int_b dA_{\rho_0(b)\rho_1(b)}+ \sum_{e,b} 
(A_{\rho_1(b)\rho_1(e)} - A_{\rho_0(b)\rho_0(e)})) \\ & \mspace{250mu}
\prod_{v,e,b} g_{\rho_1(b)\rho_1(e)\rho_1(v)}g^{-1}
_{\rho_0(b)\rho_0(e)\rho_0(v)} \\
&= \exp \sum_{e,b} (A_{\rho_0(b)\rho_1(b)} + A_{\rho_1(b)\rho_1(e)} - 
A_{\rho_0(b)\rho_0(e)}) \prod_{v,e,b} g_{\rho_1(b)\rho_1(e)\rho_1(v)}g^{-1}
_{\rho_0(b)\rho_0(e)\rho_0(v)} \\
&= \exp \sum_{e,b} (A_{\rho_0(e)\rho_1(e)} - d\log g_{\rho_0(b)\rho_1(b)
\rho_1(e)} + d\log g_{\rho_0(b)\rho_0(e)\rho_1(e)}) \\ & \mspace{250mu}
\prod_{v,e,b} g_{\rho_1(b)\rho_1(e)\rho_1(v)}g^{-1}
_{\rho_0(b)\rho_0(e)\rho_0(v)} \\
&= \exp (\sum_{e,b} A_{\rho_0(e)\rho_1(e)}) \prod_{v,e,b} 
g^{-1}_{\rho_0(b)\rho_1(b)\rho_1(e)} g_{\rho_0(b)\rho_0(e)\rho_1(e)} 
g_{\rho_1(b)\rho_1(e)\rho_1(v)}g^{-1}
_{\rho_0(b)\rho_0(e)\rho_0(v)} \\
&= \exp (\sum_{e,b} A_{\rho_0(e)\rho_1(e)}) \prod_{v,e,b}
g^{-1}_{\rho_0(e)\rho_0(v)\rho_1(v)}g_{\rho_0(e)\rho_1(e)\rho_1(v)}
\end{split}
\end{equation}
where the last step involves repeated applications of the cocycle condition
on $g$ and the elimination of terms depending only on $b$ and $v$. All 
interior terms will cancel due to the summation over $b$, leaving only 
edges in the boundary which we denote $\partial e$. In terms
of $\partial \psi$, the restriction of the map $\psi$ to the boundary, we have
\begin{equation}\label{bgtls}
\begin{split}
G_{(t_0,\rho_0)(t_1,\rho_1)}(\partial \psi) &= 
\exp \sum_{\partial e} \int_e \pa \psi^* 
A_{\rho_0(e)\rho_1(e)} 
\cdot
\prod_{v,\partial e} g^{-1}_{\rho_0(e)\rho_0(v)\rho_1(v)}
g_{\rho_0(e)\rho_1(e)\rho_1(v)}(\pa\psi(v)) \\
&= \exp \sum_{\partial e} \int_e \pa\psi^* 
A_{\rho_0(e)\rho_1(e)} 
\cdot \prod_{v,\partial e} g^{-1}_{\rho_1(e)\rho_0(v)\rho_1(v)}
g_{\rho_0(e)\rho_1(e)\rho_0(v)}(\pa\psi(v))
\end{split}
\end{equation}
where the second line is an alternate formula obtained using the cocycle
condition which is equivalent to that given by Brylinski \cite{bry}.

It is possible to relate this directly with the approach of 
Hitchin \cite{hit} where the holonomy of a gerbe over a loop is 
defined in terms of the holonomy of a bundle defined by the
difference between two trivialisations. In this case we have 
trivialisations $(\un{k}^0,\un{h}^0)$ and $(\un{k}^1,\un{h}^1)$
defined over loops in $V_{(t_0,\rho_0)}$ and $V_{(t_1,\rho_1)}$ 
respectively. In this case transition functions would be
\begin{equation}
\hat{G}_{(t_0,\rho_0)(t_1,\rho_1)} = \exp \sum_e \int_e
(k^0_{\rho(e)} - k^1_{\rho(e)}) \cdot \prod_{v,e} 
h^0_{\rho(e)\rho(v)} {h^1}^{-1}_{\rho(e)\rho(v)}
\end{equation} 
where $(t,\rho)$ is any index of the open cover of the loop
space which is defined on $V_{(t_0,\rho_0)(t_1,\rho_1)}$. The
transition functions are independent of this choice since they
are defined as a holonomy. In the following calculation we shall choose this 
to be $(t_0,\rho_0)$. We now directly compare $G_{(t_0,\rho_0)
(t_1,\rho_1)}$ and $\hat{G}_{(t_0,\rho_0)(t_1,\rho_1)}$:
\begin{equation}
\begin{split}
G_{(t_0,\rho_0)(t_1,\rho_1)}\hat{G}^{-1}_{(t_0,\rho_0)(t_1,\rho_1)}
&= \exp \sum_e \int_e (A_{\rho_0(e)\rho_1(e)} - k^0_{\rho_0(e)} 
+ k^1_{\rho_0(e)} ) \\
& \mspace{20mu} \cdot \prod_{v,e} g^{-1}_{\rho_0(e)\rho_0(v)\rho_1(v)}
g_{\rho_0(e)\rho_1(e)\rho_1(v)} {h^0}^{-1}_{\rho_0(e)\rho_0(v)}
h^1_{\rho_0(e)\rho_0(v)} \\
&= \exp \sum_e \int_e (A_{\rho_0(e)\rho_1(e)} - k^0_{\rho_0(e)} 
+ k^1_{\rho_0(e)} ) \\
& \mspace{80mu} \cdot \prod_{v,e} h^1_{\rho_0(e)\rho_1(e)} 
h^1_{\rho_1(e)\rho_1(v)} 
{h^0}^{-1}_{\rho_0(e)\rho_0(v)} \\
&= \exp \sum_e \int_e (k^1_{\rho_1(e)} - k^0_{\rho_0(e)}) \cdot 
\prod_{v,e} h^1_{\rho_1(e)\rho_1(v)} {h^0}^{-1}_{\rho_0(e)\rho_0(v)}
\end{split}
\end{equation}
This is a trivial Deligne class so $\un{G}$ and $\hat{\un{G}}$
define isomorphic bundles.\\

To get the local connection 1-form on $V_{(t_0,\rho_0)}$ we calculate
$\tilde{d}\log H_{(t_0,\rho_0)}$. 
\begin{equation}
\begin{split}
\tilde{d}\log H_{(t_0,\rho_0)}(\xi)
&= \sum_b \int_b (\i_\xi d\eta_{\rho_0(b)} + d \iota_\xi \eta_{\rho_0(b)}) 
+ \sum_{e,b}
\int_e (\iota_\xi dA_{\rho_0(b)\rho_0(e)} + d\iota_\xi
A_{\rho_0(b)\rho_0(e)}) \\
&\mspace{340mu} + \sum_{v,e,b} (\iota_\xi d\log g_{\rho_0(b)
\rho_0(e)\rho_0(v)})\\
&= \sum_b \int_b \i_\xi \omega + 
\sum_{e,b} \int_e \iota_\xi (\eta_{\rho_0(b)} - \eta_{\rho_0(b)} +
\eta_{\rho_0(e)}) \\
& \mspace{120mu} + \sum_{v,e,b} \iota_\xi (A_{\rho_0(b)\rho_0(e)} - 
A_{\rho_0(b)\rho_0(e)} + A_{\rho_0(b)\rho_0(v)} - A_{\rho_0(e)\rho_0(v)}) \\
&=  (\int_\Sigma ev^*\omega)(\xi) + 
\sum_{\partial e} \int_e \iota_\xi \eta_{\rho_0(e)} +  
\sum_{v,\partial e} -\iota_\xi A_{\rho_0(e)\rho_0(v)}
\end{split}
\end{equation}
The 1-form $\tilde{d}\log H_{(t_0,\rho_0)}$ does not descend to the boundary
however $\tilde{d}\log H_{(t_0,\rho_0)} - \int_\Sigma ev^*\omega$ is
an equivalent connection on the trivial bundle $D(\tau_{\Sigma^\partial})$
which does descend, so the local connection forms are given by
\begin{equation}
B_{(t_0,\rho_0)} = \sum_{\partial e} \int_e \iota_\xi \eta_{\rho_0(e)} +  
\sum_{v,\partial e} -\iota_\xi A_{\rho_0(e)\rho_0(v)}
\end{equation}
If the map $\pa\psi$ consists of a number of components $\pa\psi_i$ then 
\begin{equation}
G_{(t_0,\rho_0)(t_1,\rho_1)}(\partial \psi) =
\prod_i G_{(t_0,\rho)(t_1,\rho_1)}(\pa\psi_i)
\end{equation}
and
\begin{equation}
B_{(t_0,\rho_0)}(\xi) = \sum_i B_{(t_0,\rho_0)}(\xi_i)
\end{equation}
thus if we consider maps, $\partial_i$, of the boundary components into the 
free loop
space $LM$ then the bundle described by $\un{G}$ is a product of bundles
over each component loop, $\bigotimes_i \partial_i^{-1} L$,
where $L$ is the line bundle over the loop space defined by the
transgression formula \eqref{bgtls}. The bundle was described in this way 
in \cite{bry}. 

Let us summarise what we know about transgression of bundle gerbes. 
Given a bundle gerbe with connection and curving, $P$, over $M$ we obtain
a trivial bundle $D(\tau_{\Sigma^\partial}P)$ over $\Map(\Sigma^\partial,M)$.
A section of this trivial bundle defines a generalised notion of 
parallel transport.
Also we have seen that $D(\tau_{\Sigma^\partial}P)$ is the pullback of a 
bundle on the boundary, $\tau_{\partial \Sigma}P$. This bundle may be defined
in terms of the transgression to the loop space which is a bundle
$\tau_{S^1} P$.

It follows that whenever $M$ is 1-connected, we have a simple geometric picture
of the bundle over the loop space (as a bundle 0-gerbe),
\[
\begin{array}{ccc}
& & S^1 \\
& \stackrel{H}{\nearrow} & \\
S^2(M) & \rightrightarrows & D^2(M) \\
& & \downarrow \\
& & \L M 
\end{array}
\]
where $S^2(M)$ is the space of smooth maps of the 2-sphere into $M$
and $D^2(M)$ is the space of smooth maps of the 2-disc into $M$.
We know from our previous discussion that the transition functions
of this bundle gerbe are $G_{(t_0,\rho_0)(t_1,\rho_1)}$. 
Previously we considered the bundle 
0-gerbe defined by the holonomy function on $\Sigma M$ where $\Sigma$ is 
any surface on $M$ however in this case spheres are sufficient to obtain all
possible loops in $M$ since every loop bounds a disc. 
For more general $M$ there exist loops which do not bound a disc, 
in which case there is no simple geometric picture of the bundle on
the loop space, however, the local transition functions are still
well defined. There is an analogy between this
situation and the holonomy of bundles over $M$. When $M$ is 1-connected
the holonomy may be defined in simple terms as an integral of the 
curvature over some surface, however more generally if we want to
define the holonomy as a function on $M$ we use a local formula.

\section{A General Formula for Parallel Transport}
In this section we take a similar approach to that taken in \S \ref{genhol}.
The formulae obtained are examples of the more general fibre
integration formulae of Gomi and Terashima (\cite{gote2}\cite{gote}).
The difference with our approach is that rather than starting with a 
formula and proving that it satisfies the 
requirements for parallel transport, we are deriving formulae by first
generalising the notion of parallel transport.
Suppose we have a class in $H^{p}(M,\D^p)$ and a $p$-manifold, 
$X^\partial$, with boundary $\partial X$. Then we may define a set 
of local functions $H_{(t_0,\rho_0)}$ on $X^\partial M$ such that given 
any two mappings $\psi_1, \psi_2$ which have the same boundary then 
$H^{-1}_{(t_0,\rho_0)}(\psi_1)H_{(t_0,\rho_0)}(\psi_2) =
H(\psi_1 \# \psi_2)$ where $H$ is the holonomy function on $XM$, 
$\#$ is the connected sum which replaces the composition for paths and 
loops and 
the open cover is induced by the restriction to $\partial X$.  
We may define transition functions $G_{(t_0,\rho_0)(t_1,\rho_1)} =
H^{-1}_{(t_0,\rho_0)}H_{(t_1,\rho_1)}$ which define a trivial bundle on
$X^\partial M$ and descend to define a possibly non-trivial bundle on 
$\partial X M$. We calculate the formula for these transition functions 
using our general formula for holonomy \eqref{ghform}.
\begin{equation}\label{gentranf}
G_{(t_0,\rho_0)(t_1,\rho_1)} = \exp \sum_{n=0}^p \sum_{\un{\sigma}^{p-n}}
\int_{\sigma^{p-n}}  [
A^{p-n}_{\rho_1(\sigma^{p}) \ldots \rho_1(\sigma^{p-n})} -
A^{p-n}_{\rho_0(\sigma^{p}) \ldots \rho_0(\sigma^{p-n})}]
\end{equation}
We would like to express this in a form which explicitly relies only
on the restriction to the boundary. Consider the first term in the summation
over $n$,
\begin{equation}\label{k=1}
\begin{split}
\sum_{\un{\sigma}^p} \int_{\sigma^p} [A^p_{\rho_1(\sigma^p)} -
A^p_{\rho_0(\sigma^p)}] &=
\sum_{\un{\sigma}^p} \int_{\sigma^p}  \delta(A^p)_{\rho_0(\sigma^p)
\rho_1(\sigma^p)} \\
&= \sum_{\un{\sigma}^p} \int_{\sigma^p}  dA^{p-1}_{\rho_0(\sigma^p)
\rho_1(\sigma^p)} \\
&= \sum_{\un{\sigma}^{p-1}} \int_{\sigma^{p-1}} A^{p-1}_{\rho_0(\sigma^p)
\rho_1(\sigma^p)} 
\end{split}
\end{equation}
Combining this with the term corresponding to $n=1$ gives
\begin{equation}
\begin{split}
\sum_{\un{\sigma}^{p-1}} \int_{\sigma^{p-1}} A^{p-1}_{\rho_0(\sigma^p)
\rho_1(\sigma^p)} + A^{p-1}_{\rho_1(\sigma^{p})\rho_1(\sigma^{p-1})} 
- A^{p-1}_{\rho_0(\sigma^{p})\rho_0(\sigma^{p-1})}&=\\
& \mspace{-400mu}
\sum_{\un{\sigma}^{p-1}} \int_{\sigma^{p-1}} A^{p-1}_{\rho_0(\sigma
^{p-1})\rho_1(\sigma^{p-1})}
 + \delta(A^{p-1})_{
\rho_0(\sigma^p)\rho_1(\sigma^p)\rho_1(\sigma^{p-1})} -
\delta(A^{p-1})_{\rho_0(\sigma^p)\rho_0(\sigma^{p-1})
\rho_1(\sigma^{p-1})}
\end{split}
\end{equation}
We would now like to iterate this process. First we shall prove a formula
which will simplify some of the terms,
\begin{lemma}\label{indlem}
For $1 \leq n \leq p$
\begin{equation}
\begin{split}
\sum_{r=0}^{n} (-1)^{r} \delta(A^{p-n})_{\rho_0(\sigma^{p})\ldots
\rho_0(\sigma^{p-r})\rho_1(\sigma^{p-r})\ldots \rho_1(\sigma^{p-n})}
&= \\& \mspace{-330mu}
(-1)^{n+1}[\sum_{r=0}^{n-1} (-1)^{r} A^{p-n}_{\rho_0(\sigma^{p})\ldots
\rho_0(\sigma^{p-r})\rho_1(\sigma^{p-r})\ldots\rho_1(\sigma^{p-n+1})}] \\
& \mspace{-290mu}
+ \sum_{r=1}^{n} (-1)^{r} A^{p-n}_{\rho_0(\sigma^{p-1})\dots
\rho_0(\sigma^{p-r})\rho_1(\sigma^{p-r})\ldots\rho_1(\sigma^{p-n})}
+ A^{p-n}_{\rho_1(\sigma^{p-n})\ldots\rho_1(\sigma^p)}
\\& \mspace{60mu}
 - A^{p-n}_{\rho_0(\sigma^{p-n})\ldots\rho_0(\sigma^p)} + I
\end{split}
\end{equation}
where $I$ consists of a sum of terms in which one subscript
$\rho_i(\sigma^m)$ is omitted, for $p-n < m < p$ and $i \in
\{0,1\}$.
These are exactly the terms which cancel
out under the sum $\underset{\un{\sigma}^{p-n}}{\sum} \int_{\sigma^{p-n}} I$.

\end{lemma}

\begin{proof}
Most of the terms on the left hand side are absorbed into the 
term $I$ on the right hand side. Consider the remaining terms. There are three
distinct types:
\begin{enumerate}
\item Those for which the subscript $\rho_1(\sigma^{p-n})$ is omitted. There is
one term of this type for each $r$ in the summation. In each case a factor
of $(-1)^{n+1}$ is introduced by the definition of $\delta$.\\
\item Those for which the subscript $\rho_0(\sigma^{p})$ is omitted. There is
one term of this type for each $r$ in the summation. \\
\item Terms for which the subscripts include $\sigma^m$ for all 
$p-n \leq m \leq p$. This is only possible if the subscript omitted
by $\delta$ is either $\rho_0(\sigma^{p-r})$ or $\rho_1(\sigma^{p-r})$.
Note that the term obtained by omitting $\rho_0(\sigma^{p-r})$ is the
same as that obtained by omitting the $\rho_1(\sigma^{p-r-1})$ subscript
in the $(r-1)$-st term
of the summation over $r$. In both of these cases we are eliminating the
subscript immediately after $\rho_0(\sigma^{p-r-1})$ so the signs arising
from the $\delta$ map are equal, however the two terms get opposite signs 
from the factor $(-1)^{r}$, and hence they cancel out. \\
There are two terms 
which do not cancel out in this way. One is obtained by omitting the 
subscript $\rho_0(\sigma^{p})$ from the $r=0$ term in the summation.
The factor $(-1)^{r}$ is equal to 1 and the coefficient from $\delta$
is also 1 since we are omitting the first term.\\
The other term which does not cancel out is the one obtained by omitting the 
subscript $\rho_1(\sigma^{p-n})$ from the $r=n$ term.  
\end{enumerate}
\end{proof}

\begin{proposition}
The formula for the transition functions given in equation \eqref{gentranf}
is equivalent to
\begin{equation}\label{propind}
\begin{split}
\exp \sum_{n=1}^{k-1} \sum_{\un{\sigma}^{p-n}} \int_{\sigma^{p-n}}
-
\sum_{r=1}^n (-1)^{r}A^{p-n}_{\rho_0(\sigma^{p-1})
\ldots\rho_0(\sigma^{p-r})\rho_1(\sigma^{p-r})\ldots\rho_1(\sigma^{p-n})
} &\\& \mspace{-505mu} \cdot \exp \sum_{\un{\sigma}^{p-k}} 
\int_{\sigma^{p-k}}  
(-1)^{k+1}[\sum_{r=0}^{k-1}(-1)^{r}
A^{p-k}_{\rho_0(\sigma^{p})\ldots\rho_0(\sigma^{p-r})\rho_1(\sigma^{p-r})
\ldots\rho_1(\sigma^{p-k+1})}] \\
&\mspace{-425mu}\cdot 
\exp \sum_{n=k}^{p} \sum_{\un{\sigma}^{p-n}} \int_{\sigma^{p-n}}
[A^{p-n}_{\rho_1(\sigma^{p})\ldots \rho_1(\sigma
^{p-n})} - A^{p-n}_{\rho_0(\sigma^{p})\ldots\rho_0(\sigma^{p-n})}]
\end{split}
\end{equation}
for each $k$ such that $1 \leq k \leq p$.
\end{proposition}
\begin{proof}
Equation \eqref{k=1} shows that the result is true for $k=1$. We now proceed 
by induction on $k$. Assume the result is true for $k=l$. Define terms
$A(n)$, $B(n)$ and $C(n)$ such that equation \eqref{propind} becomes
\begin{equation}
\exp [ \sum_{n=1}^{k-1}A(n) + B(k) + \sum_{n=k}^p C(n)]
\end{equation}
To prove that the equation holds for $k=l+1$ we need to show
\begin{equation}\label{ll+1}
\exp [ \sum_{n=1}^{l-1}A(n) + B(l) + \sum_{n=l}^p C(n)] 
= \exp [ \sum_{n=1}^{l}A(n) + B(l+1) + \sum_{n=l+1}^p C(n)]
\end{equation}
Clearly if $B(l) + C(l) = A(l) + B(l+1)$ then
the result holds since all other terms in \eqref{ll+1} are identical.
Consider
\begin{equation}
\begin{split}
B(l) + C(l) &= \sum_{\un{\sigma}^{p-l}} \int_{\sigma^{p-l}}
(-1)^{l+1}\sum_{r=0}^{l-1}(-1)^{r}A^{p-l}_{
\rho_0(\sigma^{p})\ldots\rho_0(\sigma^{p-r})\rho_1(\sigma^{p-r})
\ldots\rho_1(\sigma^{p-l+1})}\\& \mspace{70mu}
+ \sum_{\un{\sigma}^{p-l}}\int_{\sigma^{p-l}}
[A^{p-l}_{\rho_1(\sigma^{p})\ldots\rho_1(\sigma^{p-l})} - A^{p-l}_
{\rho_0(\sigma^{p})\ldots\rho_0(\sigma^{p-l})}]\\
&= \sum_{\un{\sigma}^{p-l}} \int_{\sigma^{p-l}}
\sum_{r=0}^l (-1)^{r}\delta(A^{p-l})_
{\rho_0(\sigma^{p})\ldots\rho_0(\sigma^{p-r})\rho_1(\sigma^{p-r})
\ldots\rho_1(\sigma^{p-l})} \\& \mspace{10mu}
- \sum_{\un{\sigma}^{p-l}}\int_{\sigma^{p-l}}
\sum_{r=1}^l (-1)^{r}A^{p-l}_{\rho_0(\sigma^{p-1})\ldots
\rho_0(\sigma^{p-r})\rho_1(\sigma^{p-r})\ldots\rho_1(\sigma^{p-l})} \\
\intertext{by Lemma \ref{indlem},}
&= \sum_{\un{\sigma}^{p-l}} \int_{\sigma^{p-l}}
\sum_{r=0}^l (-1)^{r}(-1)^l dA^{p-l-1}_
{\rho_0(\sigma^{p})\ldots\rho_0(\sigma^{p-r})\rho_1(\sigma^{p-r})
\ldots\rho_1(\sigma^{p-l})} \\& \mspace{10mu}
- \sum_{\un{\sigma}^{p-l}}\int_{\sigma^{p-l}}
\sum_{r=1}^l (-1)^{r}A^{p-l}_{\rho_0(\sigma^{p-1})\ldots
\rho_0(\sigma^{p-r})\rho_1(\sigma^{p-r})\ldots\rho_1(\sigma^{p-l})} \\
\intertext{since $\delta(\un{A}^{p-n}) = (-1)^nd\un{A}^{p-n-1}$ for a 
Deligne class,}
&= \sum_{\un{\sigma}^{p-l-1}} 
\int_{\sigma^{p-l-1}} (-1)^l
\sum_{r=0}^l (-1)^{r}A^{p-l-1}_
{\rho_0(\sigma^{p})\ldots\rho_0(\sigma^{p-r})\rho_1(\sigma^{p-r})
\ldots\rho_1(\sigma^{p-l})} \\&
\mspace{10mu}
- \sum_{\un{\sigma}^{p-l}}\int_{\sigma^{p-l}}
\sum_{r=1}^l (-1)^{r}A^{p-l}_{\rho_0(\sigma^{p-l})\ldots
\rho_0(\sigma^{p-r})\rho_1(\sigma^{p-r})\ldots\rho_1(\sigma^{p-1})} \\
&= B(l+1) + A(l)
\end{split}
\end{equation}
\end{proof}
This proposition gives a general formula for parallel transport
by considering the case $k=p$. In this case the 
transition functions are given by
\begin{equation}
\begin{split}
G_{(t_0,\rho_0)(t_1,\rho_1)} =& 
\exp \sum_{n=1}^{p-1} \sum_{\un{\sigma}^{p-n}} \int_{\sigma^{p-n}}
-
\sum_{r=1}^n (-1)^{r}A^{p-n}_{\rho_0(\sigma^{p-1})
\ldots\rho_0(\sigma^{p-r})\rho_1(\sigma^{p-r})\ldots\rho_1(\sigma^{p-n})
}\\
& \cdot \exp (B(p) + C(p))
\end{split}
\end{equation}
It is not difficult to show that $B(p) + C(p) = A(p)$. This follows from
the proof that $B(l) + C(l) = A(l) + B(l+1)$. Recall that the $B(l+1)$ term
appears after applying Stokes' theorem to a $dA^{p-l-1}$ term which in
turn comes from applying the Deligne cocycle condition to a $\delta(A^{p-l})$
term. In this case we have $\delta(A^0) = 0$ so the $B(l+1)$ term is absent,
leaving the required result. 

Thus we have
\begin{corollary}\cite{gote2,gote}
The transition functions for the parallel transport bundle associated 
with a Deligne $(p+1)$-class over $M$ and a smooth map of a $p$-manifold
with boundary into $M$ are given by
\begin{equation}\label{genG}
G_{(t_0,\rho_0)(t_1,\rho_1)} = 
\exp \sum_{n=1}^{p} \sum_{\un{\sigma}^{p-n}} \int_{\sigma^{p-n}}
\sum_{r=1}^n (-1)^{r+1}A^{p-n}_{\rho_0(\sigma^{p-1})
\ldots\rho_0(\sigma^{p-r})\rho_1(\sigma^{p-r})\ldots\rho_1(\sigma^{p-n})
}
\end{equation}
\end{corollary}
We now derive a formula for the connection. Once again we start with
the local functions
\begin{equation}
H_{(t_0,\rho_0)} = \exp \sum_{n=0}^p \sum_{\un{\sigma}^{p-n}}
\int_{\sigma^{p-n}} A^{p-n}_{
\rho_0(\sigma^{p}) \ldots \rho_0(\sigma^{p-n})}
\end{equation}
We  apply $\tilde{d}\log$ and evaluate at $\xi \in T(\Map(\partial\Sigma,M))$
to get
\begin{equation}\label{genB}
\begin{split}
\tilde{d}\log(H_{(t_0,\rho_0)})(\xi) - (\int_\Sigma ev^* \omega)(\xi)
=& \sum_{n=0}^p \sum_{\un{\sigma}^{p-n}}
\int_{\sigma^{p-n}} [d\iota_\xi A^{p-n}_{
\rho_0(\sigma^{p}) \ldots \rho_0(\sigma^{p-n})} + \iota_\xi d A^{p-n}_{
\rho_0(\sigma^{p}) \ldots \rho_0(\sigma^{p-n})}] \\
=& \sum_{n=0}^p [\sum_{\un{\sigma}^{p-n-1}}
\int_{\sigma^{p-n-1}} \iota_\xi A^{p-n}_{
\rho_0(\sigma^{p}) \ldots \rho_0(\sigma^{p-n})} \\& \mspace{60mu}
+ \sum_{\un{\sigma}^{p-n}}
\int_{\sigma^{p-n}} (-1)^{n+1} \iota_\xi
\delta(A^{p-n+1})_{\rho_0(\sigma^{p}) 
\ldots \rho_0(\sigma^{p-n})}]\\
=& \sum_{n=0}^{p-1}  \sum_{\un{\sigma}^{p-n-1}}
\int_{\sigma^{p-n-1}} 
\iota_\xi A^{p-n}_{
\rho_0(\sigma^{p}) \ldots \rho_0(\sigma^{p-n})} \\ & \mspace{35mu}
+ \sum_{n=1}^p  \sum_{\un{\sigma}^{p-n}}
\int_{\sigma^{p-n}} 
(-1)^{n+1} \iota_\xi \delta(A^{p-n+1})_{\rho_0(\sigma^{p}) 
\ldots \rho_0(\sigma^{p-n})}\\
\intertext{where we have used $\iota_\xi \un{A}^0 = 0$ 
and $\un{A}^{p+1}=0$,}
=& \sum_{n=0}^{p-1}  \sum_{\un{\sigma}^{p-n-1}}
\int_{\sigma^{p-n-1}} \iota_\xi A^{p-n}_{
\rho_0(\sigma^{p}) \ldots \rho_0(\sigma^{p-n})} \\& \mspace{140mu}
+ (-1)^{n} \iota_\xi 
\delta(A^{p-n})_{\rho_0(\sigma^p)
\ldots \rho_0(\sigma^{p-n-1})}
\end{split}
\end{equation}
Only two terms from the $\delta$ part survive under the sum over
$\un{\sigma}^{p-n-1}$. These are the ones which omit the subscripts 
$\rho_0(\sigma^{p})$  and $\rho_0(\sigma^{p-n-1})$ respectively. Consider
the first of these. Since it is the first term in the $\delta$ expansion
it is positive, thus we have $(-1)^{n}
\iota_\xi A^{p-n}_{\rho_0(\sigma^{p-1})\ldots\rho_0(\sigma^{p-n-1})}$.
The term omitting the subscript $\rho_0(\sigma^p)$ in
the $\delta$ expansion will have the additional coefficient 
$(-1)^{n+1}$ which makes equal to the negative of the first term in
\eqref{genB}.
Thus we have shown that the connection on the bundle over 
$\Map(\partial \Sigma,M)$ representing parallel transport
is given by
\begin{equation}
B_{(t_0,\rho_0)}(\xi) = \sum_{n=0}^{p-1}  \sum_{\un{\sigma}^{p-n-1}}
\int_{\sigma^{p-n-1}} (-1)^n
\iota_\xi A^{p-n}_{\rho_0(\sigma^{p-1})\ldots\rho_0(\sigma^{p-n-1})}
\end{equation}
\begin{example}
Using these general formulae we can calculate the transition functions
and local connections corresponding to the parallel transport of a bundle
2-gerbe on $M$ over a smooth map of $3$-manifold with boundary, $X$ into $M$. 
Let
the bundle 2-gerbe be represented by the Deligne class $(\un{g},\un{A},
\un{\eta},\un{\nu})$, let $X$ be triangulated by 3-faces $w$, 2-faces $b$, 
edges $e$ and vertices $v$. Denote the map $X \rightarrow M$ by $\psi$.
Then the transition functions are given by
\begin{equation}
\begin{split}
G_{(t_0,\rho_0)(t_1,\rho_1)}(\psi) &=
\exp \left( 
\sum_{b,w} \int_b \psi^*\eta_{\rho_0(b)\rho_1(b)} + \sum_{e,b,w} \int_e
\psi^*(A_{\rho_0(b)\rho_1(b)\rho_1(e)} - A_{\rho_0(b)\rho_0(e)\rho_1(e)} )
\right) \\
& \mspace{65mu}
\cdot \prod_{v,e,b,w} g_{\rho_0(b)\rho_1(b)\rho_1(e)\rho_1(v)}
g^{-1}_{\rho_0(b)\rho_0(e)\rho_1(e)\rho_1(v)}
g_{\rho_0(b)\rho_0(e)\rho_0(v)\rho_1(v)}(\psi(v))
\end{split}
\end{equation} 
and the local connections are
\begin{equation}
B_{(t_0,\rho_0)}(\xi) = \sum_{b,w}\int_b \iota_\xi \nu_{\rho_0(b)}
- \sum_{e,b,w} \int_e \iota_\xi \eta_{\rho_0(b)\rho_0(e)}
+ \sum_{v,e,b,w} \iota_\xi A_{\rho_0(b)\rho_0(e)\rho_0(v)}
\end{equation}

In the case where $M$ is 2-connected then we may, by analogy with 
$\tau_{S^1}P$ in the bundle gerbe case, define a bundle 0-gerbe 
$\tau_{S^2} P$ which may be represented in the 
following way:
\[
\begin{array}{ccc}
&& S^1 \\
& \stackrel{H}{\nearrow} & \\
S^3(M) & \rightrightarrows & D^3(M) \\
& & \downarrow \\
& & S^2(M)
\end{array}
\]
Gomi and Terashima \cite{gote2,gote} suggest that it would be of interest to
find geometric realisations of transgression in higher degrees. This 
construction gives such a realisation in terms of bundle gerbes.
\end{example}

\section{Loop Transgression of Bundle 2-Gerbes}\label{ltb2g}
Recall that so far we have considered a generalised form of parallel 
transport. This has involved transgression of a bundle gerbe to a 
bundle 0-gerbe (or, equivalently, a bundle) over the
loop space and the transgression of a bundle 2-gerbe to a bundle 
0-gerbe on the space $\partial X M$ of smooth maps of boundaries of
3-manifolds in $M$. 
As in the case of loop transgression of a bundle gerbe the fact that the
transgression formula may be broken up into a product of factors
over each boundary component implies that this bundle 0-gerbe 
may be realised as a product $\bigotimes_i \partial_i^{-1} L$ where 
$L$ is a bundle 0-gerbe on $\Sigma M$ which is defined locally 
by the transgression formula.


We may now proceed as in the case of loop transgression of bundle gerbes.
To do this we apply the hierarchy principle, replacing functions with 
bundles. This means that we want to find a bundle on $\Sigma^{\partial}
M$ that locally trivialises $\tau_{\partial X} P$. To do this simply
apply the formula for $G_{(t_0,\rho_0)(t_1,\rho_1)}$. If this fails to
be consistent on triple intersections then it will define a 
trivial bundle gerbe. This is quite a long calculation involving the 
repeated application of the various cocycle conditions which
define the Deligne class of the bundle 2-gerbe. We start with the 
highest term, $\eta$ and apply the appropriate cocycle condition and then
use Stokes' Theorem to move down to the next level. This process is repeated 
at each level with terms such as $\sum A$ and $\prod g$ used as an abbreviation
of all of the terms at the other levels. We also write $G_{01}$ for
$G_{(t_0,\rho_0)(t_1,\rho_1)}$ and so on.
\begin{equation}
\begin{split}
G_{01}G_{12}
G^{-1}_{02} &= \exp (
\sum_b \int_b \eta_{\rho_0(b)\rho_1(b)} + \eta_{\rho_1(b)\rho_2(b)}
-\eta_{\rho_0(b)\rho_2(b)}) \cdot \exp \sum A \cdot \prod g \\
& = \exp (\sum_{e,b} \int_e -A_{\rho_0(b)\rho_1(b)\rho_2(b)} +
A_{\rho_0(b)\rho_1(b)\rho_1(e)}
-A_{\rho_0(b)\rho_0(e)\rho_1(e)} + A_{\rho_1(b) \rho_2(b)\rho_2(e)}
\\& \mspace{150mu} - A_{\rho_1(b)\rho_1(e)\rho_2(e)} 
- A_{\rho_0(b) \rho_2(b)\rho_2(e)}
+ A_{\rho_0(b)\rho_0(e)\rho_2(e)}) \cdot \prod g \\
& = \exp (\sum_{e,b} \int_e -A_{\rho_0(e)\rho_1(e)\rho_2(e)} + 
d\log g_{\rho_0(b)\rho_1(b)\rho_2(b)\rho_1(e)}
+ d\log g_{\rho_1(b)\rho_2(b)\rho_1(e)\rho_2(e)} \\& \mspace{150mu}
- d\log g_{\rho_0(b)\rho_2(b)\rho_1(e)\rho_2(e)}
+ d\log g_{\rho_0(b)\rho_0(e)\rho_1(e)\rho_2(e)})
\cdot \prod g \\
&=  \exp (\sum_{e,b} \int_e -A_{\rho_0(e)\rho_1(e)\rho_2(e)})
\cdot \prod_{v,e,b} g^{-1}_{\rho_0(e)\rho_0(v)\rho_1(v)\rho_2(v)}
g_{\rho_0(e)\rho_1(e)\rho_1(v)\rho_2(v)} \\
& \mspace{420mu} g^{-1}_{\rho_0(e)\rho_1(e)
\rho_2(e)\rho_2(v)}(v)
\end{split}
\end{equation}
All of the interior terms will cancel in the sum over $b$ so these
transition functions descend to $\partial \Sigma M$. We also
have a canonical choice of connection on this bundle gerbe which 
is given by the $D$-trivial local connection forms on $\Sigma^\partial M$,
\begin{equation}
(B_1 - B_0 - \tilde{d}\log G_{01})(\xi) =
\sum_{e,b} \int_e -\iota_\xi \eta_{\rho_0(e)\rho_1(e)} + \sum_{v,e,b}
\iota_\xi(A_{\rho_0(e)\rho_1(e)\rho_1(v)} - A_{\rho_0(e)\rho_0(v)\rho_1(v)})
\end{equation}
The details here are similar to those of previous calculations and have
been omitted.

The canonical choice of curving is defined by
\begin{equation}
\begin{split}
\tilde{d}B_0(\xi,\mu) &= \sum_b \int_b \i_\mu\i_\xi d\nu_{\rho_0(b)}
-d\i_\mu\i_\xi \nu_{\rho_0(b)} + \sum_{e,b} \int_e -\i_\mu\i_\xi
d\eta_{\rho_0(b)\rho_0(e)} + d\i_\mu \i_\xi \eta_{\rho_0(b)\rho_0(e)} \\
& \mspace{365mu} +
\sum_{v,e,b} \i_\mu\i_\xi dA_{\rho_0(b)\rho_0(e)\rho_0(v)} \\
&= \sum_b \int_b \i_\mu\i_\xi \omega + \sum_{e,b} \int_e -\i_\mu\i_\xi 
\nu_{\rho_0(e)} + \sum_{v,e,b} -\i_\mu\i_\xi \eta_{\rho_0(e)\rho_0(v)}
\end{split}
\end{equation}
The 2-form $\tilde{d}B_0 - \int_\Sigma ev^* \omega$ descends to give
local curving 2-forms 
\begin{equation}
\zeta_0(\xi,\mu) = 
\sum_{e,b} \int_e
- \i_\mu\i_\xi \nu_{\rho_0(e)} + 
\sum_{v,e,b} -\i_\mu \i_\xi \eta_{\rho_0(e)\rho_0(v)}
\end{equation}

This local bundle gerbe data splits into a product of terms over each
component of $\partial \Sigma M$ and defines a bundle gerbe over the loop
space. This situation is more complicated than the usual concept of 
parallel transport. We shall give a direct comparison for the example
of a cylinder. For a bundle gerbe there is a line bundle $L$ on the loop
space, so associated with each boundary loop of the cylinder is 
a fibre of this bundle. Over the cylinder there is a trivial
bundle with fibres given by the product of the fibres of $L$ at the 
end loops. A trivial bundle has a section, which in this case defines 
a $U(1)$-equivariant map between the fibres at the end loops. If, on the
other hand, we start with a bundle 2-gerbe then we have seen that a 
bundle gerbe on the loop space is obtained. Thus associated with each boundary
loop is the fibre of a bundle gerbe, which is a $U(1)$-groupoid. 
Associated with a cylinder is a trivial bundle gerbe however this does not
necessarily define a section. 

This is similar to what happens when we consider a bundle gerbe over a 
path. Starting with the line bundle on the loop space obtained by 
transgression we obtain a trivialisation of a bundle gerbe on $\P M$.
The trivial bundle gerbe is the pull back of a bundle gerbe on 
$\Map(\partial I ,M)$. The fibre over each component of 
$\Map(\partial I,M)$ is simply the fibre of the original bundle
gerbe over the relevant point in $M$. This example leads to holonomy
reconstruction which shall be discussed in the next chapter, detailed
calculations shall be given there. 


%

%% file: chapter7.tex
\chapter{Further Results on Holonomy and Transgression}

We briefly discuss the terminology used in this chapter.
Strictly speaking holonomy is the $U(1)$-valued map corresponding to
a bundle $n$-gerbe over a closed $(n+1)$ manifold. Corresponding
to $(n+1)$-manifolds there are sections of trivial bundles which
generalise parallel transport. Corresponding to closed $n$-manifolds
there is the transgression bundle, named after a homomorphism in
Deligne cohomology. Holonomy and parallel transport may also both be
viewed as arising from such transgression maps so we sometimes use
the term transgression to apply to all of these cases. It is also 
possible to view all cases as a generalisation of the notion of 
holonomy, sometimes the term holonomy is used in this context.

We present some basic properties of holonomy and transgression. We 
consider the holonomy of some of the examples we have encountered and 
discuss consequences for the theory of holonomy reconstruction. We
also consider gauge invariance properties which are relevant to 
applications.

\section{Some General Results} \label{genres}
In Chapter 5 we considered generalisations of holonomy to bundle gerbes and
higher objects. 
Here we present some basic properties of bundle holonomy and show that they
also apply to bundle gerbe holonomy. We then go on to see how these results
apply to the more general notion of transgression which was discussed in
Chapter 6.  
Throughout the final two chapters
$\hol$ will denote the holonomy map as defined in Chapter 5, with
the bundle $n$-gerbe ($n$ = 0,1,2), connections and curvings and
the closed $n+1$ manifold indicated where necessary.

\subsubsection{Functoriality}
Let $P \rightarrow M$ be a bundle with connection $A$ and suppose we have
a map
$\phi : N \rightarrow M$. Then 
\begin{equation}
\hol(\phi^{-1} P; \phi^* A) = \phi^* \hol(P;A)
\end{equation}
as functions on $LN$.

Now suppose we have a bundle gerbe $(P,Y,M;A,f)$ and map
$\phi : N \rightarrow M$. Then we can form the pullback bundle gerbe 
over $N$ and use this to define a holonomy function on $\Map(\Sigma,N)$. 
For an element $\psi$ of $\Map(\Sigma,N)$ the holonomy function is defined by
the evaluation of the flat holonomy class of the bundle gerbe $\psi^{-1}(
\phi^{-1} P)$. This is the same as the bundle gerbe over $\Sigma$ obtained
by the pullback of $P$ by $\phi \circ \psi \in \Map(\Sigma,M)$. The flat 
holonomy class of this bundle gerbe is the holonomy function on 
$\Map(\Sigma,M)$ evaluated at $\phi \circ \psi$, so it follows that
\begin{equation}
\hol(\phi^{-1}P,\phi^{-1}Y,N;\phi^*A,\phi^*f) = \phi^* \hol(P,Y,M)
\end{equation}
as functions on $\Map(\Sigma,N)$. Furthermore a similar result clearly 
holds for bundle 2-gerbes and maps of 3-manifolds.

Now consider the transgression of the bundle gerbe $\phi^{-1} P$ to the loop
space $\L N$. The map $\phi$ induces a map $\L_\phi : \L N \rightarrow \L M$ 
such that given any $\mu \in \L N$ the map $L_\phi \mu \in \L M$ is defined by
\begin{equation}
\L_\phi \mu (\theta) = \phi (\mu (\theta))
\end{equation}

Now consider the transgression of $P$ to the loop space $\L M$. This gives
a bundle 0-gerbe which may be pulled back to $\L N$ by the map $\L_\phi$.
We claim that this is stably isomorphic to the bundle 0-gerbe over $\L N$ 
mentioned above. First consider the case where there are no non-trivial 
loops and there is a geometric picture of this bundle 0-gerbe.
At each level we have a map between surfaces induced
by $\phi$. In particular given $\psi \in \Sigma N$ then we can map
to the surface given pointwise by $\phi \circ \psi$ in $M$. The holonomy
map on $\Sigma N$ gives the flat holonomy class of $\psi^{-1} \phi^{-1} P$
which is clearly equal to the holonomy map on $\Sigma M$ corresponding
to the surface $\phi \circ \psi$. Furthermore since the connection 
is derived from the holonomy function then we have a $D$-stable isomorphism.
This result also extends to the bundle gerbe over the loop space associated
with a bundle 2-gerbe since once again we have an induced map between
smooth maps of manifolds at each level which is invariant under the
holonomy map. Since the transgression formulae are obtained by this 
construction we claim that the functoriality applies to transgression
in the general case, for example, functoriality of the holonomy
of a bundle gerbe extends also to the local functions on surfaces
with boundary and the transition functions are then defined in terms
of these.

\subsubsection{Orientation}
We consider the effect of a change of orientation of the manifold over 
which the holonomy is evaluated.
In the case without boundary then it is 
clear that the holonomy function changes sign under a change of orientation
for bundles, bundle gerbes and bundle 2-gerbes. Recall that the function 
at a particular point of the relevant mapping space is the evaluation of
the flat holonomy class of the corresponding pullback. This evaluation 
involves integrating the $D$-obstruction form,
so a change of orientation will reverse the sign of this integral and 
lead to an overall reversal of the sign of the holonomy function.

This argument easily extends to the case with boundary. Our usual 
approach here is to proceed as in the case without boundary to obtain a 
function which locally trivialises the holonomy function on the fibre product.
As above, a reverse in orientation leads to a reverse in sign of this function.
The transition functions and 
local connections of the transgression bundle are derived from the local
expressions for this function, so they too have signs reversed. Thus the
bundle 0-gerbes with connection  obtained by transgression 
change to their duals under a reverse of orientation of the embedded manifold.
In the case of a bundle 2-gerbe over a surface with boundary then we have
a bundle 0-gerbe which locally trivialises the bundle 0-gerbe on $\Sigma M$.
Since this has a change of orientation then so do the trivialisations and
thus the bundle gerbe on the loop space has a change of orientation.

\subsubsection{Multiplicativity}
Suppose $\Sigma_1 \sqcup \Sigma_2$ is a disjoint union of closed 
$n$-manifolds, for $n = 1,2,3$. If P is a bundle $(n-1)$-gerbe then 
\begin{equation}
\hol_{\Sigma_1 \sqcup \Sigma_2}(P) = \hol_{\Sigma_1}(P)\cdot\hol_{\Sigma_2}(P)
\end{equation}
The left hand side is equal to $\exp \int_{\Sigma_1 \sqcup \Sigma_2} \chi$ 
where $\chi$ is the $D$-obstruction form for the flat, trivial bundle 
$(n-1)$-gerbe $(\Sigma_1 \sqcup \Sigma_2)^{-1}P$. This integral separates into 
the two holonomies on the right hand side.\\
In the case with boundary then we find a similar relation however the 
functions are only defined locally.
Taking $D$ of this to get the bundle on the boundary, $L$ we have
\[
L_{\Sigma_1 \sqcup \Sigma_2} = L_{\Sigma_1} \otimes L_{\Sigma_2}
\]
an isomorphism of bundle with connection. Furthermore since the sections
of the pull backs to the mapping space of the manifold with boundary
are given by the local function these also satisfy an additivity property.

Suppose we have a disjoint union of 
a closed manifold, $\Sigma_1$, and one with boundary, $\Sigma_2$. In this 
case the union is a manifold with boundary so there is a function $H_\rho$ on 
$\Sigma^\partial M$ leading to a bundle on the restriction to the 
boundary, however since the $H_\rho$ term corresponding to $\Sigma_1$ has no
boundary then it is globally defined and so does not contribute to the 
transition function and thus the 
bundle does not depend on $\Sigma_1$ and is simply
$L_{\Sigma_2}$. The difference from the case where we just have $\Sigma_2$
is the section of the trivial bundle defined by $H_\rho$. The term corresponding
to $\Sigma_1$ gives a different choice of trivialisation of the bundle $D(H_\rho)$
than if we just have the $\Sigma_2$ term.

These results are easily extended to bundle 2-gerbes. The functions, bundles
and bundle gerbes obtained by the various transgressions satisfy the
obvious additivity conditions. In the situation where we have one
component with boundary and one closed then the result of transgression
is the transgression corresponding to the component with boundary, with 
the closed component contributing to the canonical section/trivialisation. 

\subsubsection{Gluing}
Let $\Sigma_1$ and $\Sigma_2$ be two $n$ manifolds where $n = 1,2$ or 3.
Suppose that $\partial \Sigma_1 = X_1 \cup X$ and $\partial \Sigma_2 =
X_2 \cup X$. In this case we may glue $\Sigma_1$ and $\Sigma_2$ along 
$X$ to get a new manifold $\Sigma$ with $\partial \Sigma = X_1 \cup -X_2$.
When evaluating the holonomy of $\Sigma$ with respect to 
a bundle $(n-1)$-gerbe we get functions $H_{(t_0,\rho_0)}$ which 
separate into a boundary component and an
interior component which descends to the mapping space of the boundary. 
We claim
that this will be equal to the product of $H_{\Sigma_1(t_0,\rho_0)}$ and 
$-H_{\Sigma_2(t_0,\rho_0)}$. 
The change in sign is due to the need to reverse orientation
for gluing. Clearly the claim holds for the interior components since these
are not affected by gluing. Since the boundary components are summations over
a triangulation then we may break these down to get $H^{\partial}_{X_1(t_0,
\rho_0)} +
H^{\partial}_{X(t_0,\rho_0)} - H^{\partial}_{X(t_0,\rho_0)} - 
H^{\partial}_{X_2(t_0,\rho_0)}$ which is equal to the
boundary term for $\Sigma$. Thus we have
\begin{equation}
\tau_{\partial \Sigma} P = \tau_{\partial X_1} P  \otimes (\tau_{\partial 
X_2}P)^*
\end{equation}

Now consider the special case where 
we have two $n$-manifolds $\Sigma_1$ and $\Sigma_2$, such that
$\partial \Sigma_1 = - \partial \Sigma_2$. Then it is possible to glue 
them together to obtain a closed manifold $\Sigma$. This time all of the 
boundary terms cancel leaving only the interior terms. Putting these together
we get the holonomy function on $\Map(\Sigma,M)$. This function corresponds
to the product of the two sections over $\Map(\Sigma,M)$. 

In a similar way, given a bundle 2-gerbe $P$ and surfaces $\Sigma_1$ and 
$\Sigma_2$ then the same results apply in terms of bundle gerbes.

\section{The Holonomy of the Tautological Bundle Gerbe}
We begin by considering the holonomy of the tautological bundle 0-gerbe.
Recall that given a closed, $2\pi$-integral 2-form, $F$, on a 1-connected
manifold, $M$, this is defined by the 
diagram
\[
\begin{array}{ ccc}
&& U(1)     \\
&\stackrel{\rho}{\nearrow} &     \\
   \L_0M & \rightrightarrows & \P_0M \\
   & & \downarrow \\
   & & M \\
\end{array}
\]
where the map $\rho$ is defined by
\begin{equation}\label{tbhol}
\rho(\gamma) = \exp \int_\Sigma F
\end{equation}
and $\Sigma$ is any surface which is bounded by $\gamma$. 
The connectedness requirement ensures
that we can choose such a $\Sigma$ and
we have already established that $\rho$ is independent of the choice of
$\Sigma$.

Now consider the holonomy of this bundle 0-gerbe. Since the curvature is
$F$ then we know that the holonomy must also satisfy the condition on 
$\rho$ in \eqref{tbhol}. This condition completely characterises the function
and thus the holonomy of the tautological bundle 0-gerbe is the function
$\rho: \L_0 (M) \rightarrow U(1)$. This is a rather trivial fact however
it will serve as an indication of what to expect as we move up the bundle
gerbe hierarchy.

Now consider the tautological bundle gerbe, with the tautological bundle over
$\L_0(M)$ considered as a bundle 0-gerbe,
\begin{equation}\label{tautbg}
\begin{array}{ccccc}
& & U(1) & & \\
& \stackrel{\rho}{\nearrow} & & & \\
S^2( M)  & \rightrightarrows & D^2( M) & & \\
& & \downarrow & & \\
& & \L_0 M  & \rightrightarrows & \P_0 M  \\
& & & & \downarrow \\
& & & & M
\end{array}
\end{equation}
By our connectedness assumption on $M$ all closed 2-surfaces in $M$ are of the
form $\partial X$ for some 3-manifold $X$.
The holonomy satisfies
$H(\partial X) = \exp \int_X \omega = \rho(\partial X)$. Since by the usual 
arguments this is independent of the choice of $X$ we can conclude that the
holonomy function on $S^2(M)$ is equal to $\rho$.

Recall that we can also transgress a bundle gerbe 
to get a bundle over $\L M$. We wish to compare this with the 
bundle over $\L_0 M$ in the tautological bundle gerbe described above.
Recall Lemma \eqref{indbp} tells us that the tautological construction is
independent of the choice of base point, so without loss of generality
we may replace $\L M$ with $\L_0 M$. Now the fact that $\rho$ is 
equal to the holonomy map shows that the bundle 0-gerbe obtained by
transgression of the tautological bundle gerbe is the same as the
bundle 0-gerbe
on the loop space which is used to define the tautological bundle gerbe in
\eqref{tautbg}.
This implies the following
\begin{proposition} \label{propLG}
The transgression to the loop space $\L G$ of the tautological 
bundle gerbe over a compact, simply connected, semi-simple Lie group
$G$ with curvature 3-form $\Tr <g^{-1}dg \wedge [ g^{-1}dg \wedge 
g^{-1}dg]>$ is the bundle associated with the central extension
$\widetilde{\L G} \rightarrow \L G$ (see example \ref{egLG}).
\end{proposition}

Next consider the tautological bundle 2-gerbe,
\begin{equation}\label{taut2}
\begin{array}{ccccccc}
& & S^1 & & & & \\
& \stackrel{\rho}{\nearrow} & & & & & \\
S^3(M) & \rightrightarrows & D^3( M) & & & & \\
& & \downarrow & & & & \\
& & S^2(M) & \rightrightarrows & D^2(M) & & \\
& & & & \downarrow & & \\
& & & & \L_0 M & \rightrightarrows & \P_0 M \\
& & & & & & \downarrow \\
& & & & & & M 
\end{array}
\end{equation}
where $M$ is 3-connected. By the assumption on $M$, every element of $D^3(M)$ 
is
the boundary of a 4-manifold, $W$ and thus the holonomy function is uniquely 
determined by the property $H(\partial W) = \exp \int_W \Theta$ where 
$\Theta$ is the 4-curvature. This is exactly the definition of the function 
$\rho$, so it is clear that the transgression bundle $\tau_{S^2}$ is the
bundle over $S^2( M)$ in \eqref{taut2} and the transgression
bundle gerbe $\tau_{\L M}$ is the bundle gerbe over 
$\L M$ in \eqref{taut2}.

\section{Holonomy Reconstruction}\label{DCHR}


In this section we consider the theory of holonomy reconstruction 
in the bundle gerbe context.

Let us consider the implications of our connectedness assumptions in the
tautological case. This will lead to an understanding of the more general
case. Let $M$ be 1-connected so that the tautological
bundle 0-gerbe is defined. Then we also have $H_1(M) = 0$. Recall that 
there is an exact sequence defined by the map of a Deligne class
to its curvature,

\[
0 \rightarrow H^1(M,U(1)) \rightarrow H^1(M,\D^1) \rightarrow
A^2_0(M) \rightarrow 0 
\]
First we use the Universal Coefficient Theorem for cohomology (\cite{botu});
\begin{equation}
H^q(X,G) \isom \Hom (H_q(X),G) \oplus \Ext(H_{q-1},G)
\end{equation}
Since we have $H_1(M) = 0$ and $\Ext(H_0(M),U(1)) = 0$ then
\begin{equation}
H^1(M,U(1)) \isom \Hom (0,U(1)) \isom 0
\end{equation}
Together with $\eqref{gaj}$ this gives the result
\begin{equation}
H^1(M,\D^1) \isom A^2_0(M)
\end{equation}
If we think of this Deligne class as a bundle then this tells
us that when $\pi_1(M) = 0$ all bundle gerbes with connection on $M$ are
completely determined by their curvature. The tautological construction
gives us a bundle with connection over a 1-connected base $M$ 
which has a particular curvature, so this tells us that the tautological
bundle is in fact the unique (up to isomorphism) bundle with connection
satisfying these requirements. Rather than constructing a bundle from
its curvature we shall construct it from its holonomy function. In
the tautological case the holonomy function is completely determined by
the curvature, so this approach does make sense as a generalisation of the
tautological case. 

We now relax the
requirement that $M$ be 1-connected. Let us start with a bundle 0-gerbe with
connection defined by a class in $H^1(M,\D^1)$. There is now no map
$\rho$, however in the previous section we found that for the tautological
bundle $\rho$ is equal to the holonomy, which does exist in the general
case. An explicit construction of the Deligne class from the holonomy is 
described in \cite{mapi}. We shall construct a bundle 0-gerbe using a similar
approach. 

We use the tautological construction replacing $\rho$ with the 
holonomy, giving the following bundle 0-gerbe:
\begin{equation}\label{recon0}
\begin{array}{ccc}
& & U(1) \\
& \stackrel{H}{\nearrow} &  \\
\L_0M & \rightrightarrows & \P_0M \\
& & \downarrow \\
& & M
\end{array}
\end{equation}
We define the connection to be $ \int_I ev^* F$ as in the tautological case.
The inverse of the holonomy map follows from the relationship
between transition functions and the flat holonomy class, see the proof
of proposition \ref{cg} for an example of this.

Now we outline the method for finding local expressions of a bundle with
connection from \cite{mapi}. The first step is to define a section
over
each $U_i \subset M$ which gives a path in $U_i$ for each which ends at
$m_i \in U_i$. 
These are composed with paths $p_i$ which connect them
to the base point and so define local sections of the path fibration. 
The sections over $U_i$ and $U_j$ are then composed and the holonomy is 
evaluated over the resulting loop. Observe that this gives precisely the
transition functions of the bundle 0-gerbe \eqref{recon0}. In particular
note that the inverse is present since the direction of
the path is positive on the lift over $U_i$ and negative over
$U_j$. Furthermore it
is shown that the resulting bundle is independent of the choice of base point.

Next we show that our connection is the same as that obtained in \cite{mapi}.
Since our connection is defined on the path space we need to pull it back 
using a local connection for comparison. Let $s_i$ be the local section over
$U_i$ and let $v$ be a vector at $x \in U_i$ which may be represented by 
$q'(0)$ for some curve $q \subset U_i$. The vector $s_{i*} X$ is a vector 
in the path space which is based at the path $s_i(x)$. It is given by
$\frac{d}{dk}s_i(q(k))|_{k=0}$, where for each $k$ the path $s_i(q_k)$
is from $m_i$, the centre of $U_i$, to $q(k)$. Strictly speaking the paths
actually begin at $m_0$, but since $s_i(q(k))$ is constant between $m_0$ and
$m_i$ we are only interested in what happens between $m_i$ and $q(k)$.
The local form of our connection evaluated at $v$ is $\int_I 
\iota_{s_{i*v}} ev^* F =  \int_I F(s_{i*v}(t),dt)$ where 
$s_{i*v}(t)$ is the vector on $U_i$ given by the element of the
vector field on the path $s_i(x)$ at $s_i(x)(t)$.

In \cite{mapi} the connection is defined by first taking the holonomy
of a particular loop associated with a vector $v$ as above. Omitting the 
trivial part again, this loop consists of three components, first the 
lift $s_i(x)$, then the curve $q$ and finally the lift $s_i(q(k))$ in the 
reverse direction. The local one form evaluated at $v$ is then defined
by taking minus the $k$ derivative of log of this holonomy at 
$k=0$. Note that
the loop around which the holonomy is taken corresponds to the boundary
of a surface, $\Sigma_k$, defined by the collection of all of the loops
$s_i(q(s))$ from 0 to $k$. This means that we can express this holonomy
as an integral of the curvature, $\exp \int_{\Sigma_k} F$. The integral
over $\Sigma_k$ may be parametrised by $I \times [0,k]$ where $I$ parametrises
the curves from $m_i$ to $q(s)$ and $[0,k]$ parametrises the curves. 
We can now
calculate the connection:
\begin{equation}
\frac{d}{dk} \int_I \int_0^k F(s_i(q'(s))(t),dt) |_{k=0}
= \int_I ev^* F(s_{i*v}(t),dt)
\end{equation} 
as required.

Now consider a general bundle 0-gerbe with connection $(g,Y,M;A)$. 
If we calculate the holonomy and then use it to reconstruct the bundle 0-gerbe
then we have a new bundle 0-gerbe $(H,\P_0 M,M;B)$ which is 
$D$-stably isomorphic
to $(g,Y,M;A)$. This may be shown explicitly by considering the 
product bundle 0-gerbe 
\begin{equation}
\begin{array}{ccc}
& & U(1) \\
& \stackrel{g^{-1}H}{\nearrow} & \\
Y^{[2]} \times_\pi \L_0 M & \rightrightarrows & Y \times_\pi \P_0 M \\
& & \downarrow \\
& & M
\end{array}
\end{equation}
The transition functions for the bundle 0-gerbe are given by
\begin{equation}
g^{-1}_{\alpha\beta} \cdot \exp \int_{\mu_\alpha \star \mu_\beta^{-1}} 
(A_i - d\log h_i)
\end{equation}
where $\mu_\alpha$ and $\mu_\beta$ are paths given by the local sections
of the path fibration at $m \in U_{\alpha \beta}$. Normally we deal with 
the term in the integral by breaking it down into a sum of edges in a 
triangulation of the loop. Using the sections of the path space described
in \cite{mapi} we can break down $\mu_\alpha \star \mu^{-1}_\beta$ into
a sum of four components. First there is a path from $m_0$ to $m_\alpha$, the 
``centre'' of $U_\alpha$. This is independent of $m$ and so contributes
a constant factor $K_\alpha$ to the transition function. Next there is a path, 
$\tilde{\mu}_\alpha$, 
from $m_\alpha$ to $m$. This is contained within $U_\alpha$ so the integral
corresponding to this component is $\int_{\tilde{\mu}_\alpha} A_\alpha - 
d\log h_\alpha$. The remaining components are from $m$ to $m_\beta$ and
from $m_\beta$ to $m_0$ and contribute similar terms to give
\begin{equation}
g^{-1}_{\alpha \beta}(m)\cdot K_\alpha K^{-1}_\beta \cdot \exp 
\int_{\tilde{\mu}_\beta} (-A_\beta + d\log h_\beta) \cdot 
\exp \int_{\tilde{\mu}_\alpha} (A_\alpha - d\log h_\alpha)
\end{equation}
Applying Stokes' theorem to the $h$ terms gives
\begin{equation}
g^{-1}_{\alpha \beta}(m)\cdot K_\alpha K^{-1}_\beta \cdot
\exp (\int_{\tilde{\mu}_\beta} - A_\beta + \int_{\tilde{\mu}_\alpha} A_\alpha)
\cdot h^{-1}_\beta(m_\beta)h_\beta(m)h^{-1}_\alpha(m)h_\alpha(m_\alpha)
\end{equation}
We may now cancel out $g^{-1}_{\alpha \beta}(m)$ with
$h_\beta(m)h^{-1}_\alpha(m)$ and incorporate the other $h$ terms into the
constants since they don't depend on $m$. This leaves
\begin{equation}
\delta ( K_\alpha^{-1} \cdot \exp (\int _{\tilde{\mu}_\alpha} -A_\alpha))
\end{equation}
proving that the two bundle 0-gerbes are stably isomorphic. To see that this
extends to a $D$-trivialisation it remains to observe that applying
$d\log$ gives the pull back of the connection form, $-A + \int_I ev^* F$ by 
a local section.

Using these results we now have a canonical representative of the 
$D$-stable isomorphism class of any bundle 0-gerbe with connection which
we shall call the {\it holonomy representative}. Given any bundle 0-gerbe
with connection this is obtained by taking the holonomy and then reconstructing
a bundle 0-gerbe.

We now consider the case of bundle gerbes. We assume for now that $M$ is
1-connected. We shall postpone discussion of this requirement until the
next section. 
Given a function representing 
the holonomy of a bundle gerbe 
we reconstruct the bundle gerbe in the following way
\begin{equation}
\begin{array}{ccccc}
& & U(1) & & \\
& \stackrel{H}{\nearrow} & & & \\
S^2(M) & \rightrightarrows & D^2(M) & & \\
& & \downarrow & & \\
& & \L_0 M & \rightrightarrows & \P_0 M \\
& & & & \downarrow \\
& & & & M 
\end{array}
\end{equation} 
Over $U_{ij}$ we have the sections $s_i$ and
$s_j$. Use these to define the pullback bundle 0-gerbe
\begin{equation}
\begin{array}{ccc}
& & U(1) \\
& \stackrel{H_{ij}}{\nearrow} & \\
S^2(M)_{ij} & \rightrightarrows & D^2(M)_{ij} \\
& & \downarrow \\
& & U_{ij}
\end{array}
\end{equation}
where elements of $D^2(M)_{ij}$ which lie in the fibre over 
$m$ are surfaces bounded
by the 1-cycle $s_i(m) \star s_j(m)^{-1}$. Since $U_{ij}$ is contractible
this bundle 0-gerbe is trivial. Let the trivialisation be defined by
a function $e_{ij} : D^2(M)_{ij} \rightarrow U(1)$ which is a homomorphism
with respect to the gluing of surfaces and is equal to $H_{ij}$ for surfaces
with no boundary. The function $e_{ij}$ plays the same role as the section 
$\sigma_{ij}$ in the usual construction of the transition functions of a 
bundle gerbe. In a similar way we may define $e_{jk}$ and $e_{ik}$ and
form the products of these over $U_{ijk}$.
Consider the product $e_{ij} e_{jk} e_{ik}$. This is defined on the fibre 
product $D^2(M)_{ij} \times_\pi D^2(M)_{ij} \times_\pi D^2(M)_{ij}$. 
Consider an a general element of this space. First you glue together 
a surface with boundary $s_i(m) \star s_j(m)^{-1}$ and one with
boundary $s_j(m) \star s_k(m)^{-1}$ This forms a new surface with boundary 
$s_i(m) \star s_k(m)^{-1}$. Finally we glue this to another surface with
the same boundary, thus giving a surface with no boundary. Using the
gluing properties this means that the result of the product
is the holonomy around this surface. We claim that this is equivalent 
to the construction in \cite{mapi}.

The bundle gerbe connection on $D^2(M)$ is $A = \int_\Sigma ev^* \omega$ where
$\Sigma \in D^2(M)$.
To get a local expression we consider again the bundle 0-gerbe on $U_{ij}$
obtained by pulling back with $(s_i,s_j)$. To get a local formula we 
pull back using a section $s_{ij}: U_{ij} \rightarrow D^2(M)_{ij}$. This 
section takes $x \in U_{ij}$ to a surface $\Sigma_{ij}$ with boundary 
$s_i \star s_j^{-1}$. Given a vector $v \in T_x(M)$, the vector $s_{ij*}v$ is
a vector field on $s_{ij}(x)$. Suppose we parameterise $\Sigma_{ij}$ with
$s$ and $t$. 
Then the pull back is $\int_{\Sigma_{ij}} \omega(s_{ij*}v(s,t),ds,dt)$. The 
construction in \cite{mapi} involves taking the holonomy around the
surface $\Sigma_{ij}$. We write this as $\int_{X_{ij}} \omega$ where $X_{ij}$
is the 3-manifold bounded by $\Sigma_{ij}$ which is defined by the family
of surfaces given by lifting the curve, $q$, defining $v$ with $s_{ij}$. Taking
the derivative with respect to $k$, the parameter giving the endpoint of $q$, 
we get
\begin{equation}
{\frac{d}{dk}}\int_{\Sigma_{ij}} \int_{u=0}^k \omega(s_{ij}(q'(u))(s,t),ds,dt)
= \int_{\Sigma_{ij}} \omega(s_{ij*}v(s,t),ds,dt)
\end{equation}
and hence our connection is the same as that in \cite{mapi}.

The curving may be dealt with in a similar way. We start with the 
2-form on $\P_0 M$, $\int_I ev^* \omega$ and pull this back to 
$U_i$ with $s_i$, and then evaluate at a pair of vectors $v=q'(0)$ and 
$w=r'(0)$.
This gives $\int_I \omega(s_{i*v}(t),s_{i*w}(t),dt)$. In $\cite{mapi}$
the approach is to take the holonomy over a surface $\Sigma_i$ which 
is defined in the following way. Locally the vectors $v$ and $w$ are
extended to vector fields defined by commuting flows. These flows 
consist of families of curves, $q_m$ and $r_m$ where $q_m'(0)$ and
$r_m'(0)$ give the elements in the respective vector fields at the 
point $m$. Thus we have $q_x'(0) = v$ and $r_x'(0)=w$.
Furthermore these flows commute, that is, $q_{r(l)}(k) = r_{q(k)}(l)$.
The vector $s_{i*}v$ is defined by a path of paths which forms a surface
bounded by $s_i(x)$, $q_x$ and $s_i(q_x(k))$. Similarly associated with 
$s_{i*}w$ is a surface bounded by $s_i(x)$, $r_x$ and $s_i(r_x(l))$. To
get a closed surface we make similar constructions at the point
$q_{r(l)}(k) = r_{q(k)}(l)$. This gives a surface which is a cone
from $m_i$ to the surface $l_{q_k,r_l}$ (to use the notation of \cite{mapi})
bounded by $q_x, r_x, q_{r(l)}$ and $r_{q(k)}$ with appropriate orientations.
This cone defines the surface $\Sigma_i$ and we denote the enclosed volume
$X_i$. The holonomy over $\Sigma_i$ may be expressed as $\int_{X_i} \omega$.
We may parametrise $X_i$ by $I \times [0,k] \times [0,l]$ where the last 
two give a parametrisation of the discs in the cross section of the cone
and $I$ parametrises the length. Taking the partial derivatives in $k$ and
$l$ then gives
\begin{equation}
\frac{\partial}{\partial k \partial l} \int_I \int_0^l \int_0^k
\omega(s_i(q_x'(s))(t),s_i(r_x'(u))(t),dt) =
\int_I \omega(s_{i*v}(t),s_{i*w}(t),dt)
\end{equation}
Thus our definition of a bundle gerbe reconstructed from holonomy
agrees with the definition of a gerbe given in \cite{mapi}.

\section{Reconstruction via Transgression}
The techniques of the previous section do not easily extend to the case of 
bundle 2-gerbes since the construction of the local data from sections 
becomes quite complicated. Another problem is that we have only been 
able to deal with base manifolds $M$ which are 1-connected. Instead we use
transgression formulae to approach holonomy reconstruction.

Recall that in \S \ref{PTbun} we noted that the parallel transport of a 
bundle 0-gerbe gives a bundle over $M$ which
is $D$-stably isomorphic to the original bundle. This gives us an alternative
way of calculating local data using transgression formulae. Let us consider
the bundle gerbe case first, assuming for now that $M$ is 1-connected. 
Given a bundle gerbe on $M$ with corresponding
holonomy map $H: S^2(M) \rightarrow U(1)$ we wish to show that the 
following bundle gerbe, with connection
$\int_{D^2} ev^* \omega$ and curving
$\int_I ev^* \omega$ is $D$-stably isomorphic to the original one with 
curvature $\omega$:
\[
\begin{array}{ccccc}
& & U(1) & & \\
& \stackrel{H}{\nearrow} & & & \\
S^2(M) & \rightrightarrows & D^2 (M) & & \\
& & \downarrow & & \\
& & \L_0 M & \rightrightarrows & \P_0 M \\
& & & & \downarrow \\
& & & & M 
\end{array}
\]
We do this by using transgression to find the local data for this bundle 
gerbe. First we consider only the transition functions. The transition 
functions for the bundle 0-gerbe $(H,D^2(M),\L_0 M)$ are those obtained by 
transgression to the loop space, $G_{(t_0,\rho_0)(t_1,\rho_1)}$. We
define functions on $\P_0 M$ by the same formula, and denote these by 
$h_{(t_i,\rho_i)(t_j,\rho_j)}$ or just $h_{ij}$. These functions satisfy 
\begin{equation}
h_{ij}(\mu_1)h^{-1}_{ij}(\mu_2) = G_{01}(\mu_1 \star \mu_2^{-1})
\end{equation}
whenever $\rho_i$ and $\rho_j$ agree with $\rho_0$ and $\rho_1$ on the 
boundary. Thus by similar arguments to those used in the bundle 0-gerbe 
case we may use the open cover on $\P_0 M$ which is induced by the projection
to $M$. Next consider what happens on $\pi^{-1}(U_{\alpha \beta \gamma})$ 
by calculating $h_{\alpha \beta}h_{\beta \gamma}h_{\alpha \gamma}^{-1}$.
Let $(\rho_\alpha,t_\alpha)$ denote any choice of $(\rho,t)$ such that 
$\rho_\alpha(v) = \alpha$ where $v$ is the endpoint of the path. Then 
we have 
\begin{equation}
\begin{split}
h_{\alpha\beta}h_{\beta\gamma}h^{-1}_{\alpha\gamma} &=
\exp \sum_{e} \int_e A_{\rho_\alpha(e)\rho_\beta(e)} + A_{\rho_\beta(e)
\rho_\gamma(e)} - A_{\rho_\alpha(e)\rho_\gamma(e)} \cdot
\prod_{v,e} g^{-1}_{\rho_\alpha(e)\rho_\alpha(v)\rho_\beta(v)}\\ & 
\mspace{27mu}
g_{\rho_\alpha(e)\rho_\beta(e)\rho_\beta(v)}  g^{-1}_{\rho_\beta(e)
\rho_\beta(v)\rho_\gamma(v)} g_{\rho_\beta(e)\rho_\gamma(e)\rho_\gamma(v)}
g_{\rho_\alpha(e)\rho_\alpha(v)\rho_\gamma(v)}g^{-1}_{\rho_\alpha(e)
\rho_\gamma(e)\rho_\gamma(v)}(v) \\
&= \prod_{v,e} g^{-1}_{\rho_\alpha(e)\rho_\beta(e)\rho_\gamma(e)}
g^{-1}_{\rho_\alpha(e)\rho_\alpha(v)\rho_\beta(v)}
g_{\rho_\alpha(e)\rho_\beta(e)\rho_\beta(v)}  g^{-1}_{\rho_\beta(e)
\rho_\beta(v)\rho_\gamma(v)} g_{\rho_\beta(e)\rho_\gamma(e)\rho_\gamma(v)}
\\ & \mspace{330mu}
g_{\rho_\alpha(e)\rho_\alpha(v)\rho_\gamma(v)}g^{-1}_{\rho_\alpha(e)
\rho_\gamma(e)\rho_\gamma(v)}(v) \\
&= \prod_{v,e} g^{-1}_{\rho_\alpha(v)\rho_\beta(v)\rho_\gamma(v)}(v)
\end{split}
\end{equation}
where the last line is obtained by repeated application of the 
cocycle identity on $\un{g}$. It is not difficult to see that this 
descends to $M$, so we have $h^{-1}_{\alpha\beta}h^{-1}_{\beta\gamma}
h_{\alpha\gamma}(\mu) = g_{\alpha\beta\gamma}(\pi(\mu))$. We need to show
that these are transition functions for the bundle gerbe described above.
So far we have
\begin{eqnarray}
G_{01}(\mu_1^{-1} \star \mu_2) = h^{-1}_{\alpha \beta}(\mu_1)h_{\alpha\beta}
(\mu_2) \\
h^{-1}_{\alpha\beta}h^{-1}_{\beta\gamma}
h_{\alpha\gamma}(\mu_1) = g_{\alpha\beta\gamma}(\pi(\mu_2))
\end{eqnarray}
Since we wish to avoid using local sections we shall prove directly that 
this gives the obstruction to this bundle gerbe being trivial. First 
suppose that there exists a bundle gerbe trivialisation. This means 
that there exist functions $q_{\alpha \beta}$ such that
\begin{eqnarray}
G_{01}(\mu_1^{-1} \star \mu_2) = q_{\alpha \beta}^{-1}(\mu_1)q_{\alpha\beta}
(\mu_2) \\
q_{\alpha\beta}q_{\beta\gamma}
q^{-1}_{\alpha\gamma}(\mu_1) = 1
\end{eqnarray}
Consider the functions $h^{-1}_{\alpha\beta}q_{\alpha\beta}$ on 
$\P_0 M$. Since $h^{-1}_{\alpha\beta}q_{\alpha\beta}(\mu_1)
h_{\alpha\beta}q^{-1}_{\alpha\beta}(\mu_2) = G_{01}G^{-1}_{01} = 1$ these
functions descend to $M$. On $M$ we have 
$\delta(h^{-1}q^{-1})_{\alpha\beta\gamma} = 
\delta(h^{-1})_{\alpha\beta\gamma} = g_{\alpha\beta\gamma}$, thus
$\un{g}$ is a trivial cocycle. Conversely if $\un{g}$ is trivial 
then let $g_{\alpha\beta\gamma} = g_{\alpha\beta}g_{\beta\gamma}g^{-1}_{
\alpha\gamma}$. On $\P_0 M$ the functions $h^{-1}_{\alpha \beta}(\mu)
g_{\alpha\beta}(\pi(\mu))$ are globally defined and are transition 
functions for a bundle which trivialises the bundle defined by $\un{G}$.

Next we show that the connection may also be reconstructed by this 
method. We need to show that the connection $\int_{D^2} ev^* \omega$ 
corresponds to the original connection $A$. 
Transgression gives a formula for $B_{t_0,\rho_0}$, which are
local connection 1-forms on $\Omega^1(M)$. These are defined by
\begin{equation}
\pi^* B_0 = \tilde{d}\log h_0 - \int_{D^2} ev^* \omega
\end{equation}
Note that the term $\tilde{d}\log h_0$ is trivial when considering this as a 
connection on the bundle over $\L_0 M$, so we see that the one
forms $-B_0$ are local representatives for the connection.
Over $\P_0 M$ we have 1-forms $-k_\alpha$ induced by the extension of $B_0$ from
loops to paths. These satisfy $\delta(k_\alpha) = B_0$ whenever $\rho_0(v) =
\alpha$ where $v$ is the endpoint of the path. The bundle gerbe connection
is trivial if these form a connection on $\P_0 M$, that is, if 
$k_\alpha - k_\beta = \tilde{d}\log h_{\alpha \beta}$. 
If these are not equal then
they differ by a 1-form which descends to $M$ which is the local
representative $A_{\alpha \beta}$ of the bundle gerbe connection.

In terms of transgression formulae we have
\begin{eqnarray}
k_\alpha(\xi) &=& \sum_e \int_e \iota_\xi f_{\rho_\alpha(e)} + \sum_{v,e}
-\iota_\xi A_{\rho_\alpha(e)\rho_\alpha(v)}\\
\tilde{d}\log h_{\alpha\beta} (\xi) &=& \sum_e \int_e 
\iota_\xi(f_{\rho_\beta(e)}
- f_{\rho_\alpha(e)}) \\
&&  + \sum_{v,e} \iota_\xi(A_{\rho_\alpha(e)\rho_\alpha(v)}
+ A_{\rho_\alpha(v)\rho_\beta(v)} - A_{\rho_\beta(e)\rho_\beta(v)}) \\
(\tilde{d}\log h_{\alpha \beta} - k_\beta + k_\alpha)(\xi) &=&
\sum_{v,e} \iota_\xi A_{\rho_\alpha(v)\rho_\beta(v)} \\
&=&  A_{\alpha \beta}(\xi)
\end{eqnarray}
Thus we have reconstructed the local representative of the original 
connection. 

The local curving is given by $f - \tilde{d}k_\alpha$, the extent to which
the curving fails to be a curvature for the connection $k_\alpha$. 
\begin{equation}
\begin{split}
\tilde{d}k_\alpha (\xi,\nu) &= \tilde{d}(k_\alpha(\nu))(\xi) - 
\tilde{d}(k_\alpha(\xi))(\nu)
- \iota_{[\xi,\nu]}k_\alpha \\
&= \sum_e \int_e d\iota_\xi \iota_\nu f_{\rho(e)} 
+ \iota_\xi d \i_\nu f_{\rho(e)}
- d\i_\nu \i_\xi f_{\rho(e)} - \i_\nu d\i_\xi f_{\rho(e)} - \i_{[\xi,\nu]}
f_{\rho(e)} \\
& \mspace{130mu} + \sum_{v,e} -\iota_\xi d \i_\nu A_{\rho(e)\rho(v)} + 
\iota_\nu d \i_\xi A_{\rho(e)\rho(v)} + \iota_{[\xi,\nu]} A_{\rho(e)\rho(v)}
\\
&= \sum_e \int_e d\iota_\xi \iota_\nu f_{\rho(e)} + \L_\xi \i_\nu f_{\rho(e)} -
d \iota_\xi \i_\nu f_{\rho(e)} - \i_\nu \L_\xi f_{\rho(e)} +
\i_\nu \i_\eta \omega - \i_{[\xi,\nu]}f_{\rho(e)} \\& \mspace{115mu} +
\sum_{v,e} -\i_\nu\i_\eta f_{\rho(e)} - \L_\xi \i_\nu A_{\rho(e)\rho(v)}
+ \i_\nu \L_\xi A_{\rho(e)\rho(v)} - \i_\nu \i_\xi f_{\rho(v)} 
\\ & \mspace{365mu} + \i_\nu \i_\xi f_{\rho(e)} 
+ \iota_{[\xi,\nu]} A_{\rho(e)\rho(v)} \\
&= \sum_e \int_e \i_\nu \i_\eta \omega - \sum_{v,e} \i_\nu \i_\xi f_{\rho(v)}
\\
&= (\int_I ev^* \omega)(\xi,\nu) - \pi^* f_{\alpha}(\xi,\nu) 
\end{split}
\end{equation}
therefore the local curving is given by $f_\alpha$ as required.

In trying to deal with the case where $M$ is not 1-connected we still have the
problem that $D^2(M) \rightarrow \L_0 M$ is not well defined. 
It has been noted \cite{mapi} that in 
this case the holonomy map may not be used to reconstruct the bundle gerbe,
instead reconstruction is given in terms of a parallel transport structure.
We may think of this as equivalent to a transgression
line bundle on the loop space. If we start with a bundle gerbe $P$ then
there is a transgression bundle $L$ on the loop space regardless of whether 
$M$ is 1-connected. Using this we define a bundle gerbe 
\[
\begin{array}{ccc}
L & & \\
\downarrow & & \\
\L_0 M & \rightrightarrows & \P_0 M \\
& & \downarrow \\
& & M
\end{array}
\]
with curving defined once again by $\int_I ev^* \omega$.
The same arguments used in the previous case apply to prove that this is 
equivalent to the original bundle gerbe.

\section{Reconstruction of Bundle 2-Gerbes}
Let $P$ be a bundle 2-gerbe on $M$, a 2-connected manifold, with
holonomy function $H: S^3(M) \rightarrow U(1)$ and curvature $\Theta$.
We may apply the techniques of the previous section to prove that
the original bundle 2-gerbe may be reconstructed from the holonomy 
using the following diagram:
\[
\begin{array}{ccccccc}
& & S^1 & & & & \\
& \stackrel{H}{\nearrow} & & & & & \\
S^3(M) & \rightrightarrows & D^3(M) & & & & \\
& & \downarrow & & & & \\
& & S^2(M) & \rightrightarrows & D^2(M) & & \\
& & & & \downarrow & & \\
& & & & \L_0 M & \rightrightarrows & \P_0 M \\
& & & & & & \downarrow \\
& & & & & & M 
\end{array}
\]
with connection, 2-curving and 3-curving given by
\begin{eqnarray}
A &=& \int_{D^3} ev^* \Theta \\
\eta &=& \int_{D^2} ev^* \Theta \\
\nu &=& \int_I ev^* \Theta
\end{eqnarray}
We know from the previous results on reconstruction that the bundle gerbe
over $\L_0 M$ is that obtained by transgression in \S \ref{ltb2g}. 
We recall the local data on $\L_0 M$:
\begin{equation*}
\begin{split}
G_{012} &=
 \exp (\sum_{e,b} \int_e -A_{\rho_0(e)\rho_1(e)\rho_2(e)})
\cdot \prod_{v,e,b} g^{-1}_{\rho_0(e)\rho_0(v)\rho_1(v)\rho_2(v)}
g_{\rho_0(e)\rho_1(e)\rho_1(v)\rho_2(v)} \\
& \mspace{420mu} g^{-1}_{\rho_0(e)\rho_1(e)
\rho_2(e)\rho_2(v)}(v) \\
B_{01} &= \sum_{e,b} \int_e -\iota_\xi \eta_{\rho_0(e)\rho_1(e)} + \sum_{v,e,b}
\iota_\xi(A_{\rho_0(e)\rho_1(e)\rho_1(v)} - A_{\rho_0(e)\rho_0(v)\rho_1(v)})
\\
\zeta_0 &= \sum_{e,b} \int_e
- \i_\nu \i_\xi \nu_{\rho_0(e)} + 
\sum_{v,e,b} -\i_\nu \i_\xi \eta_{\rho_0(e)\rho_0(v)}
\end{split}
\end{equation*}
We extend these to $\P_0 M$ where they are locally defined. As usual they
are independent of the choice of $\rho$ up to the choice on the boundary
so we may express them in terms of the open cover induced from the 
base, $(G_{\alpha\beta\gamma}, B_{\alpha\gamma},\zeta_{\alpha})$. 
The $D$-trivial bundle 2-gerbe obtained from applying $-D$ to 
this local data descends to
$M$ to give the local data for the bundle 2-gerbe described in the
diagram above. We now calculate this data to show that it is the same as 
that for the original bundle 2-gerbe with holonomy $H$.
\begin{equation}
\begin{split}
G^{-1}_{\beta\gamma\delta}G_{\alpha\gamma\delta}G^{-1}_{\alpha\beta\delta}
G_{\alpha\beta\gamma} &= \exp (\sum_e \int_e A_{\rho_{\beta}(e)\rho_{
\gamma}(e)\rho_\delta(e)} - A_{\rho_\alpha(e)\rho_\gamma(e)\rho_\delta(e)}
+ A_{\rho_\alpha(e)\rho_\beta(e)\rho_\delta(e)} \\
& \mspace{300mu} - A_{\rho_{\alpha}(e)\rho_{\beta}(e)\rho_{\gamma}(e)}) 
\cdot \prod g \\
&= \prod_{v,e} g_{\rho_\alpha(e)\rho_\beta(e)\rho_\gamma(e)\rho_\delta(e)}
   \delta(g_{\rho_\alpha(e)\rho_\alpha(v)\rho_\beta(v)\rho_\gamma(v)}
g^{-1}_{\rho_\alpha(e)\rho_\beta(e)\rho_\beta(v)\rho_\gamma(v)}\\ 
& \mspace{320mu}
g_{\rho_\alpha(e)\rho_\beta(e)\rho_\gamma(e)\rho_\gamma(v)})_{
\alpha\beta\gamma\delta}\\
&= \prod_{v,e} g_{\rho_\alpha(v)\rho_\beta(v)\rho_\gamma(v)\rho_
\delta(v)}\\
&= g_{\alpha\beta\gamma\delta}
\end{split}
\end{equation}  

\begin{equation}
\begin{split}
(-\tilde{d}\log G_{\alpha\beta\gamma} - \delta(B)_{\alpha\beta\gamma})(\xi) 
&= \sum_e \i_\xi dA_{\rho_\alpha
(e)\rho_\beta(e)\rho_\gamma(e)} + d\i_\xi A_{\rho_\alpha(e)\rho_\beta(e)
\rho_\gamma(e)} + \sum_{v,e} g  \\
& \mspace{35mu} +\sum_e \int_e \i_\xi \eta_{\rho_{\beta}(e)\rho_
\gamma(e)}- \i_\xi \eta_{\rho_{\alpha}(e)\rho_{\gamma}(e)} + 
\i_\xi \eta_{\rho_\alpha(e)\rho_\beta(e)} + \sum_{v,e} A \\
&= \sum_{v,e} \i_\xi (A_{\rho_\alpha(e)\rho_\beta(e)\rho_\gamma(e)}
+ d\log g_{\rho_\alpha(e)\rho_\alpha(v)\rho_\beta(v)\rho_\gamma(v)}\\
& \mspace{60mu} 
- d\log g_{\rho_\alpha(e)\rho_\beta(e)\rho_\beta(v)\rho_\gamma(v)}
+ d\log g_{\rho_\alpha(e)\rho_\beta(e)\rho_\gamma(e)\rho_\gamma(v)}\\
& \mspace{160mu}
- \delta(A_{\rho_\alpha(e)\rho_\beta(e)\rho_\beta(v)} -
A_{\rho_\alpha(e)\rho_\alpha(v)\rho_\beta(v)})_{\alpha\beta\gamma}) \\
&= \sum_{v,e} \i_\xi A_{\rho_\alpha(v)\rho_\beta(v)\rho_\gamma(v)}
\end{split}
\end{equation}

\begin{equation}
\begin{split}
(\tilde{d}B_{\alpha \beta} - \delta(\zeta)_{\alpha\beta})(\xi,\nu)
&= \sum_e \int_e -\i_\xi d\i_\nu \eta_{\rho_\alpha(e)\rho_\beta(e)}
- d\i_\xi \i_\nu \eta_{\rho_\alpha(e)\rho_\beta(e)}
+\i_\nu d\i_\xi \eta_{\rho_\alpha(e)\rho_\beta(e)}\\ & \mspace{135mu}
+ d\i_\nu \i_\xi \eta_{\rho_\alpha(e)\rho_\beta(e)} 
+ \i_{[\xi,\nu]}\eta_{\rho_\alpha(e)\rho_\beta(e)} + \sum_{v,e} A \\
& \mspace{175mu} 
- \sum_e \int_e -\i_\nu \i_\xi \nu_{\rho_\beta(e)} 
- \i_\nu \i_\xi 
\nu_{\rho_\alpha(e)} + \sum_{v,e} \eta \\
&= \sum_{v,e} \i_\nu\i_\xi \eta_{\rho_\alpha(e)\rho_\beta(e)} + \i_\nu \i_\xi
dA_{\rho_\alpha(e)\rho_\beta(e)\rho_\beta(v)}  
- \i_\nu \i_\xi
dA_{\rho_\alpha(e)\rho_\alpha(v)\rho_\beta(v)} \\
&\mspace{256mu} + \i_\nu\i_\xi 
\eta_{\rho_\beta(e)\rho_\beta(v)} - \i_\nu\i_\xi \eta_{\rho_\alpha(e)\rho_\alpha(v)}
\\
&= \sum_{v,e} \i_\nu \i_\xi \eta_{\rho_\alpha(v)\rho_\beta(v)}
\end{split}
\end{equation}

\begin{equation}
\begin{split}
-\tilde{d} \zeta_\alpha (\xi,\nu,\mu) &= \sum_e \int_e \i_\mu\i_\nu\i_\xi
d\nu_{\rho_\alpha(e)} + d \i_\mu\i_\nu\i_\xi \nu_{\rho_\alpha(e)}
+ \sum_{v,e} +\i_\nu\i_\mu\i_\xi d\eta_{\rho_\alpha(e)\rho_\alpha(v)}\\
&= \sum_e \int_e \i_\mu\i_\nu\i_\xi \Theta + \sum_{v,e}
\i_\mu\i_\nu\i_\xi \nu_{\rho_\alpha(v)}
\end{split}
\end{equation}

If $M$ is only 1-connected then we may reconstruct the bundle 2-gerbe 
from the transgression bundle on $S^2(M)$. If $M$ is not 
1-connected then the bundle 2-gerbe may only be reconstructed from
the transgression bundle gerbe on the loop space. These are both 
proven using the same calculations as above.

We would like to briefly comment on {\it Cheeger-Simons Differential
Characters} and their relation to holonomy reconstruction. Let $Z_p(M)$ 
denote the group of smooth singular $p$-cycles on $M$.
A degree $p$ differential character (\cite{chsi}, \cite{bry}) 
is a $U(1)$-valued homomorphism
on $Z_{p-1}(M)$, $c$,  together with a $p$-form, $\alpha$, on $M$ which 
satisfy the condition 
\begin{equation}\label{defdc}
c(\partial \gamma) = \exp \int_\gamma \alpha
\end{equation}
for $\gamma \in Z_p(M)$. These are classified by Deligne cohomology so
there must be a one to one correspondence between degree 3 differential 
characters
and $D$-stable isomorphism classes of bundle gerbes with connection
and curving. The holonomy map satisfies equation \eqref{defdc} on
$p$-manifolds which are the boundary of a $(p+1)$-manifold, and the
additivity property is similar to the homomorphism property of
differential characters. We have seen that unlike differential characters
a holonomy map is
only sufficient to reconstruct the Deligne class under certain
assumptions regarding the topology of the base. The difference appears to 
correspond to the distinction between smooth mapping spaces
and simplicial complexes. Differential characters may be 
useful in further investigation of holonomy and Deligne cohomology 
however given that they are not needed for our applications
we have not studied this in any further detail.


\section{Geometric Transgression}
In this section we consider some examples where we are able to give a 
geometric interpretation of transgression. These include 
lifting bundle gerbes and the bundle 2-gerbe associated
to a principal bundle.

We begin with the more general case where rather than a fibre bundle over
$M$ we only have a fibration, $Y \rightarrow M$. 
\begin{proposition}\label{transfib}
Let $Y \rightarrow M$ be a fibration with $M$ 1-connected. 
Then the transgression of a bundle gerbe
$(P,Y,M)$ to $LM$ is the bundle 0-gerbe described by the following diagram:  
\[
\begin{array}{ccc}
& & S^1 \\
& \stackrel{\hol(P)}{\nearrow} &  \\
LY^{[2]} & \rightrightarrows & LY \\
 & & \downarrow \\
& & LM
\end{array}
\]
The function $\hol(P)$ is evaluated with respect to the bundle gerbe 
connection $A$ which is also a bundle connection on $P$.
The connection 1-form on $LY$ is given by $\int_{S^1} ev^* f$ where 
$f$ is the curving 2-form on $Y$.
\end{proposition}
This is a well-defined bundle 0-gerbe since the cocycle condition
is satisfied due to the gluing property of holonomy.
\begin{proof}
We need to show that the bundle 0-gerbe described above is equivalent to 
\[
\begin{array}{ccc}
& & S^1 \\
& \nearrow & \\
S^2(M) & \rightrightarrows & D^2(M) \\
& & \downarrow \\
& & \L M 
\end{array}
\] 
We do this by showing that the following bundle 0-gerbe is trivial:
\[
\begin{array}{ccc}
& & S^1 \\
& \stackrel{\Lambda}{\nearrow} & \\
S^2(M) \times_\pi LY^{[2]} & \rightrightarrows & D^2(M) \times LY \\
& & \downarrow \\
& & \L M
\end{array}
\]
where the map $\Lambda: S^2(M) \times_\pi LY^{[2]} \rightarrow S^1$
is defined by 
\begin{equation}
\Lambda(\Sigma_1,\Sigma_2,\mu_1,\mu_2) = 
\hol^{-1}_{\Sigma_1 * \Sigma_2}(P,Y,M) \cdot 
\hol_{(\mu_1 , \mu_2)}(P,Y^{[2]})
\end{equation}
where $\Sigma_1,\Sigma_2 \in S^2(M)$ and $\mu_1,\mu_2 \in \L Y$. 
Note that the first factor is 
a bundle gerbe holonomy and the second is a bundle 0-gerbe holonomy.

Define a trivialisation of $\Lambda$ by
\begin{equation}
l(\Sigma,\mu) = \exp \int_\Sigma (-f + dk_J) \cdot \hol_\mu(J)
\end{equation}
where $J \rightarrow Y_\Sigma$ is a trivialisation of the bundle 
gerbe $P$ over $\Sigma$. The second factor is well defined
since $\mu$ is a lift of $\gamma = \partial \Sigma$ to $Y$ (this is where the
assumption that $Y\rightarrow M$ is a fibration is required).
This is independent of the choice of
trivialisation since the difference is 
\begin{equation}
\begin{split}
\exp \int_\Sigma (dk_J - dk_{J'}) \cdot \hol_\mu (J^* \otimes J')
&= \hol^{-1}_{\partial_\Sigma} (L) \cdot \hol_{\mu} (\pi^{-1} L) \\
&= \hol^{-1}_\gamma (L) \cdot \hol_\gamma (L)
\end{split}
\end{equation}
where $\gamma = \partial \Sigma = \pi(\mu)$ and $\pi^{-1} L = J^* \otimes
J'$. 

If we let $J$ be a trivialisation over $\Sigma_1 \# \Sigma_2$ with
restrictions $J_1$ and $J_2$ to $\Sigma_1$ and $\Sigma_2$ respectively
then we have
\begin{equation}
\begin{split}
l^{-1}(\Sigma_1,\mu_1)l(\Sigma_2,\mu_2) &= 
\exp \int_{\Sigma_1 \# \Sigma_2} (-f + dk_J) \cdot \hol_{\mu_1}(J_1)
\cdot \hol_{\mu_2}(J_2) \\
&= \hol^{-1}_{\Sigma_1 \# \Sigma_2} (P,Y,M) \cdot \hol_{(\mu_1,\mu_2)}
(P,Y^{[2]})
\end{split}
\end{equation}

We now show that the connection of the theorem is a bundle 0-gerbe 
connection.  
\begin{eqnarray*}
\delta( \int_{S^1} ev^* f ) &=& \int_{S^1} ev^* \delta(f) \\
&=& \int_{S^1} ev^* F
\end{eqnarray*}
\begin{eqnarray*}
d\log (\hol(L)) &=& d (\int_{S^1} ev^* A) \\
&=& \int_{S^1} ev^* F.
\end{eqnarray*}

\end{proof}
We would like to use this result to examine the holonomy of a 
lifting bundle gerbe (see \S \ref{lifting}) when $M$ is 1-connected,
\[
\begin{array}{ccc}
& & \hat{G} \\
& & \downarrow \\
& & G \\
& \stackrel{\rho}{\nearrow} & \\
P_G^{[2]} & \rightrightarrows & P_G \\
& & \downarrow \\
& & M 
\end{array}
\]
If we assume that this bundle gerbe has a connection and curving then 
immediately we see that the transgression to the loop space is given
by the following diagram:
\[
\begin{array}{ccc}
& & S^1 \\
& \stackrel{\hol(\rho^* \hat{G})}{\nearrow} & \\
\L P_G^{[2]} & \rightrightarrows & \L P_G \\
& & \downarrow \\
& & \L M
\end{array}
\]
It is tempting to apply the functorial property of holonomy here 
however the bundle gerbe connection may not be equal to
the pullback of the connection on $\hat{G}$ since, in general, this
does not give a bundle gerbe connection however there is a 1-form
$\epsilon$ such that $\phi^* A - \epsilon$ is a connection. Some
explicit calculations of such 1-forms are given in \cite{must2}. In terms of
holonomy we have
\begin{equation}
\begin{split}
\hol(\phi^{-1} \hat{G};\phi^*A - \epsilon) &= 
\hol(\phi^{-1} \hat{G};\phi^*A)\cdot I(\epsilon)\\
&= \phi^* \hol(\hat{G};A) \cdot I(\epsilon)
\end{split}
\end{equation}
where $I: \L P^{[2]} \rightarrow S^1$ is defined by 
\begin{equation}
I_\gamma(\epsilon) = \exp \int_{S^1} \gamma^* \epsilon
\end{equation}
This is well defined since $\epsilon$ is the difference between two choices
of connection and thus descends to $P^{[2]}$.
The connection on this bundle 0-gerbe depends on the curving chosen. Unlike
the bundle gerbe connection there is no canonical choice.

Now we turn to the case of bundle 2-gerbes. 
\begin{proposition}
The transgression to the mapping space $S^2(M)$ of a bundle 2-gerbe
$(P,Y,X,M;A,\eta,\nu)$ such that $Y \rightarrow M$ is a fibration with
simply connected fibres and $M$ is 2-connected is given by the
following diagram:
\[
\begin{array}{ccc}
& & S^1 \\
& \stackrel{\hol(P;\eta,\nu)}{\nearrow} & \\
S^2(X^{[2]}) & \rightrightarrows & S^2(X) \\
& & \downarrow \\
& & S^2(M)
\end{array}
\]
The connection is given by $\int_\Sigma ev^* \nu$.
\end{proposition}
\begin{proof}
This is proven using an argument that is similar to that
used for proposition \ref{transfib}. The stable 
isomorphism is given by a function $l: D^3(M) \times_\pi S^2( X)
\rightarrow U(1)$ which is defined by
\begin{equation}
l(\Delta,\Sigma) = \exp \int_\Delta (-\nu + dj_R) \cdot \hol_\Sigma(R)
\end{equation}
where $R$ is a trivialisation of the bundle 2-gerbe over
$\Delta \in  D^3(M)$, which defines a bundle gerbe over $X_\Delta$ and
the second term is the holonomy of this bundle gerbe over
a closed surface. We require the connected condition on the fibres
so that lifts of surfaces to $X$ are well defined. The proof now follows 
that of Proposition \ref{transfib}.
\end{proof}
Now we would like to apply this to the bundle
2-gerbe of a principal $G$-bundle \cite{ste} where 
$G$ is simply connected. Recall that this bundle 
2-gerbe is defined in the following way:
\[
\begin{array}{ccc}
& & (R,Y) \\
& & \Downarrow \\
& & G \\
& \stackrel{\rho}{\nearrow} & \\
P_G^{[2]} & \rightrightarrows & P_G \\
& & \downarrow \\
& & M
\end{array}
\]
where $P_G \rightarrow M$ is a principal $G$-bundle, 
$\rho: P_G^{[2]} \rightarrow 
G$ is the usual map to the group element which acts on $p_2$ to give $p_1$ and
$(R,Y)$ is a bundle gerbe over $G$.
On the pullback bundle gerbe $\rho^{-1}R$ the pull back connection may be used 
as the bundle 2-gerbe connection, however the curving may not be
a bundle 2-gerbe 2-curving. It is a result of Stevenson \cite{ste} that
given a curving on a bundle gerbe over $Y^{[2]}$ it is possible
to make it into a bundle 2-gerbe 2-curving. This involves subtraction of
$\pi^* \epsilon$ where $\epsilon$ is some 2-form on $Y^{[2]}$. We will provide
a detailed calculation of such an $\epsilon$ for a specific example in 
\S \ref{CSa}.
\begin{proposition}\label{tranapr}
Let $(\rho^{-1}R,\rho^{-1}Y,P_G,M)$ be a bundle 2-gerbe associated with a 
principal $G$-bundle 
$P_G \rightarrow M$ over a 2-connected base $M$  
and a bundle gerbe with connection and 
curving $(R,Y,G;\eta,A)$
Then the transgression to $S^2(M)$ is the bundle
0-gerbe given  by the following diagram:
\begin{equation}
\begin{array}{ccc}
& & S^1 \\
& \stackrel{\rho^*\hol(R;A,\eta)\cdot I(\epsilon)}{\nearrow} & \\
S^2(P_G^{[2]}) & \rightrightarrows & S^2(P_G) \\
& & \downarrow \\
& & S^2(M)
\end{array}
\end{equation}
where $I(\epsilon)$ is the $S^1$-valued function on 
$S^2(P_G^{[2]})$ defined by
$I(\epsilon)(\psi) = \exp \int \psi^* \epsilon$.
\end{proposition}
\begin{proof}
This follows easily from the previous discussion and the bundle gerbe case.
\end{proof}

\section{Gauge Transformations}
We would like to define gauge transformations for bundle gerbes and consider
their effect on holonomy. Since a gauge transformation of a bundle gerbe is
basically a stable isomorphism of a bundle gerbe with itself it is no 
surprise that it turns out to be invariant under holonomy however it 
is of interest to see how this invariance arises in the bundle gerbe 
context. This will be of interest in subsequent applications of bundle gerbe
theory.

\subsubsection{Bundles and Bundle 0-Gerbes}
We recall some basic facts about gauge transformations of 
$U(1)$-bundles.
\begin{definition}
A {\it gauge transformation} of a principal $G$-Bundle is an
automorphism of the total space which covers the identity
on the base space.
\end{definition}
The automorphism property guarantees that a gauge transformation
preserves fibres, hence for
any gauge transformation $\phi : P \rightarrow P$
we have a map $g_{\phi} : P \rightarrow G$ defined by
\[
\phi(p) = pg_\phi (p).
\]
Since $\phi(pg) = \phi(p)g$ then we have $g_\phi (pg) = g^{-1}g_\phi(p)g$.
We are interested in the case where $G = U(1)$, so this becomes
$g_\phi (pg) = g_\phi(p)$, that is, $g_\phi$ is constant on 
fibres, so it induces a map $\hat{g}_\phi : M \rightarrow G$. 

Suppose we have a connection $A$ on the $U(1)$-bundle $P
\rightarrow M$. Pulling back by the gauge transformation $\phi$
gives
\begin{equation}\label{phi^*A}
\phi^* A = A + g^{-1}_\phi dg_\phi.
\end{equation}
To generalise to the bundle gerbe case we first consider bundle 0-gerbes.
Recall that to each bundle $P \rightarrow M$ we can associate the 
following bundle 0-gerbe
which has the same Deligne class:
\begin{equation}
\begin{array}{ccc}
& & U(1) \\
& \stackrel{\rho}{\nearrow} & \\
P^{[2]} & \rightrightarrows & P \\
& & \downarrow \\
& & M 
\end{array}
\end{equation}
where the map $\rho : P^{[2]} \rightarrow S^1$ is defined by
$\rho(p_1,p_2) = g$ where $p_2 = p_1 g$. Equivalently we can identify
$P^{[2]}$ with $P  \times S^1$ in which case $\rho$ is simply the 
identity map on $S^1$. 

Given a general bundle 0-gerbe $(\lambda,Y,M)$ a gauge transformation should 
obviously be a map $\phi:Y \rightarrow Y$
such that $\pi \circ \phi = \pi$ however the condition $\phi(pg) = 
\phi(p)g$ cannot be used here as in general we don't have a group 
action on Y. Instead we need a condition on the map $\phi^{[2]}: Y^{[2]}
\rightarrow Y^{[2]}$ and $\lambda : Y^{[2]} \rightarrow S^1$. 
Consider again the bundle 0-gerbe $(\rho,P,M)$. Applying $\phi^{[2]}$ to
$P^{[2]}$ we have
\begin{eqnarray*}
\rho(\phi^{[2]}(p_1,p_2)) &=& \rho(\phi(p_1),\phi(p_2)) \\
&=& \rho(\phi(p_1),\phi(p_1 \rho(p_1,p_2))) \\
&=& \rho(\phi(p_1),\phi(p_1) \rho(p_1,p_2)) \\
&=& \rho(p_1,p_2)
\end{eqnarray*}
This suggests the following 
\begin{definition}
Let $(\lambda,Y,M)$ be a bundle 0-gerbe. A {\it gauge transformation}
of $(\lambda,Y,M)$ is a smooth map $\phi: Y \rightarrow Y$ which satisfies
the following conditions:
\begin{eqnarray}
\pi \circ \phi &= \pi \\
\lambda \circ \phi^{[2]} &= \lambda
\end{eqnarray}
\end{definition}
Let $A$ be a bundle 0-gerbe connection 1-form on $Y$. 
This means that $A$ satisfies
the equation $\delta (A) = d\log \lambda$. The gauge transformation $\phi$ 
may be used to construct a map $(1,\phi): Y \rightarrow Y^{[2]}$. Use
this map to pull back $\delta (A)$,
\begin{equation}\label{gaugeA}
\begin{split}
(1,\phi)^* \delta (A) &= (1,\phi)^* (\pi_2^* - \pi_1^*) A \\
&= (\pi_2 \circ (1,\phi))^* A - (\pi_1 \circ (1,\phi))^* A \\
&= \phi^* A - A
\end{split}
\end{equation}
thus we have
\begin{equation}\label{phi^*A2}
\phi^* A = A + (1,\phi)^* \delta(A) = A + (1,\phi)^* d\log \lambda
\end{equation}
We would like to compare this result with equation \eqref{phi^*A}.
By definition $\rho (p,\phi(p)) = g_\phi (p)$, so this immediately
shows the equivalence of the two equations.\\
The function $(1,\phi)$ may also be used to define an $S^1$-function 
$(1,\phi)^*\lambda$ on $Y$.
\begin{equation}
\begin{split}
\delta ((1,\phi)^*\lambda)(y_1,y_2) &= \lambda^{-1}(y_1,\phi(y_1))
\lambda(y_2,\phi(y_2)) \\
&= \lambda(\phi(y_1),y_1)\lambda(y_2,\phi(y_2)) \\
&= \lambda(\phi(y_1),\phi(y_2))\lambda(\phi(y_2),y_1)\lambda(y_2,\phi(y_2)) \\
&= \lambda(\phi(y_1),\phi(y_2))\lambda(y_2,y_1) \\
&= \lambda(y_1,y_2)\lambda^{-1}(y_1,y_2) \\
&= 1
\end{split}
\end{equation}
therefore the function $(1,\phi)^* \lambda$ descends to $M$. Note that 
once again this function has similar properties to $g_\phi$. We shall
denote the function on $M$ by $\lambda_\phi$. 

Finally we calculate the holonomy of the bundle 0-gerbe $(\lambda,Y,M)$
with respect to the transformed connection $\phi^*A$. To do this we
need the $D$-obstruction form. In general for bundle 0-gerbes this
is given by $A_\alpha - d\log h_\alpha$ where $A_\alpha$ are the 
local connection forms and $h_\alpha$ is a trivialisation. In this
case to get the connection we use a local section of $Y \rightarrow M$ 
to pull back the right hand side of \eqref{phi^*A2} giving 
\begin{equation}
s_\alpha^*(A + \pi^* d\log \lambda_\phi) = A_\alpha + d\log \lambda_\phi
\end{equation} 
The holonomy is then given by
\begin{equation}
\begin{split}
H((\un{g},\un{\phi^* A});\gamma) &= 
\exp \int_\gamma A_\alpha - d\log h_\alpha + d\log \lambda_\phi \\
&= H((\un{g},\un{A});\gamma) \cdot \exp \int_\gamma d\log \lambda_\phi \\
&= H((\un{g},\un{A});\gamma)
\end{split}
\end{equation}
So the gauge transformation leaves the holonomy unchanged. This result is
not so surprising if we consider that the difference of the 
Deligne classes $(\un{g},\un{A})$ and $(\un{g},\un{\phi^*A})$ is
$(1,\un{d\log \lambda_\phi}) = D(\un{\lambda_\phi})$. 

Now we consider parallel transport. If $\mu$ is an open path in $M$ then
recall that the holonomy function now depends on the choice of trivialisation,
\begin{equation}
\begin{split}
H((\un{h},\un{\phi^*A});\mu) &=
\exp \int_\mu A_\alpha - d\log h_\alpha + d\log \lambda_\phi \\
&= H((\un{h},\un{A});\mu) \cdot \lambda^{-1}_\phi(\mu(0)) \lambda_\phi
(\mu(1))
\end{split}
\end{equation}
Thus the gauge transformation contributes an extra term to the Deligne 
cochain on $\Map(I;M)$ obtained by transgression. If we apply $D$ to this
cochain to get a bundle then the extra term will cancel out as it doesn't
depends on the choice of trivialisation. The local connections on this 
bundle will pick up an extra term $\iota_{\xi(\mu(1))}d\log \lambda_\phi
- \iota_{\xi(\mu(0))}d\log \lambda_\phi$, however, just as on the original
bundle this difference is $D$-trivial.

\subsubsection{Bundle Gerbes and Bundle 2-gerbes}\label{gaugebb2}
We shall extend the concept of gauge transformation to bundle gerbes
and bundle 2-gerbes. We shall refer to these collectively as bundle
$n$-gerbes where it is to be understood that $n = 0, 1$ or 2.
\begin{definition}
Let $(P,Y,M)$ be a bundle gerbe. A {\it gauge transformation} of 
$(P,Y,M)$ is a smooth map $\phi: Y \rightarrow Y$ which
satisfies the following conditions:
\begin{eqnarray}
\pi \circ \phi &=& \pi \\
\phi^{[2]*}(P,Y^{[2]}) &=& (P,Y^{[2]})
\end{eqnarray}
with the second condition involving a choice of isomorphism of bundles, 
$\tilde{\phi}$, over $Y^{[2]}$ which preserves the bundle gerbe product.
\end{definition}
Similarly,
\begin{definition}
Let $(P,Y,X,M)$ be a bundle 2-gerbe. A {\it gauge transformation} of 
$(P,Y,X,M)$ is a smooth map $\phi: X \rightarrow X$ which
satisfies the following conditions:
\begin{eqnarray}
\pi \circ \phi &=& \pi \\
\phi^{[2]*}(P,Y,X^{[2]}) &=& (P,Y,X^{[2]})
\end{eqnarray}
with the second condition being a choice of stable isomorphism of bundle 
gerbes,
$\tilde{\phi}$  
over $X^{[2]}$, which preserves the bundle 2-gerbe product.
\end{definition}
The next step is too see how the various connections and curvings transform.
We begin with the 2-curving for a bundle gerbe. This situation is 
similar to that of the bundle 0-gerbe connection. 
If we denote the curving by $\eta$ then following \eqref{gaugeA} we have
\begin{equation}
(1,\phi)^* \delta(\eta) = \phi^*\eta - \eta
\end{equation}
so
\begin{equation}\label{gaugef}
\phi^* \eta = \eta + (1,\phi)^* F
\end{equation}
where $F$ is the curvature of the bundle $P \rightarrow Y^{[2]}$. 
Next we consider the bundle $(1,\phi)^{-1}P$ over $Y$. 
\begin{equation}
\begin{split}
\delta((1,\phi)^{-1}P)_{(y_1,y_2)} &= P^*_{(y_1,\phi(y_1))} \otimes 
P_{(y_2,\phi(y_2))} \\
&= P_{(\phi(y_1),\phi(y_2))} \otimes P_{(\phi(y_2),y_1)} \otimes
P_{(y_2,\phi(y_2))} \\
&= P_{(\phi(y_1),\phi(y_2))} \otimes P^*_{(y_1,y_2)}
\end{split}
\end{equation}
Thus we see that $(1,\phi)^{-1}P$ is a trivialisation of the bundle gerbe 
$P^* \otimes \phi^{[2]*}P$. More importantly observe that we have an isomorphism
of bundles,
\begin{equation}
\begin{split}
\delta((1,\phi)^{-1}P) &= \phi^{[2]*}P \otimes P^* \\
&= P \otimes P^* 
\end{split}
\end{equation}
therefore $(1,\phi)^{-1}P$ descends to a bundle on $M$ which we shall call 
$P_\phi$.\\
There is a connection $\nabla_{(1,\phi)^{-1}P}$ on $(1,\phi)^{-1}P$ which satisfies
\begin{equation}
\delta (\nabla_{(1,\phi)^{-1}P}) = \nabla_{\phi^{[2]*}P} \otimes \nabla^*_{P}
\end{equation}
On the right hand sides we have two choices of connection on isomorphic
bundles, so they differ by a 1-form $\alpha$ on $Y^{[2]}$ such that
$d\alpha = \phi^{[2]*}F - F$. Over $Y$ we may compare the connections
on $\pi^{-1} P_\phi$ and $(1,\phi)^{-1}P$. These differ by a 1-form $\beta$ on $Y$.
Furthermore since $\delta(\nabla_{\pi^{-1} P_\phi}) = 0$ we have $\delta(\beta) =
\alpha$. The curvatures of these two bundles are related by 
$\pi^*F_\phi + d\beta = (1,\phi)^*F$. We may now express \eqref{gaugef} in
terms of $F_\phi$,
\begin{equation}
\phi^* \eta = \eta + \pi^* F_\phi + d\beta
\end{equation}
We also have the relationship between the connections, where $\phi^* A$ 
refers
to the induced connection on $\phi^{[2]*} P$,
\begin{equation}
\phi^*A = A + \pi_P^* \delta (\beta)
\end{equation}
We may interpret these results in terms of $D$-obstructions. The bundle
gerbe $\phi^{[2]*} P \otimes P$ has a trivialisation $(1,\phi)^{-1} P$ which
is not necessarily a $D$-trivialisation. The obstruction is given by
a 2-form $\chi$ on $M$ such that locally $\pi^*\chi = \eta - F_L$
where $F_L$ is the curvature of the connection on the trivialisation.
In this case $F_L = d\beta$, so the $D$-obstruction form is $F_\phi$.
This means that the $D$-obstruction form is trivial and so the two
bundle gerbes $\phi^{[2]*}P$ and $P$ are $D$-stably isomorphic. Thus we 
expect the holonomy to be invariant, as we shall see from explicit
calculations.

We can now work out the holonomy corresponding to the new connection
and curving by considering the local formula. We may assume without 
loss of generality that the transition functions are equal, since introducing
an extra factor from the stable isomorphism will be cancelled by an 
additional term in the connection. 

Substituting $\phi^* \eta$ and $\phi^* A$ into the local formula for the
holonomy of a bundle gerbe over a closed surface $\Sigma$ gives
\begin{equation} 
\begin{split}
H(\phi^{-1}P;\phi^*A,\phi^*\eta) &= H(p;A,\eta) \cdot \exp \left( \sum_b \int_b 
F_\phi + d\beta_{\rho(b)} \sum_{e,b} \int_e \beta_{\rho(e)} - \beta_{\rho(b)}
\right) \\
&=H(p;A,\eta) \cdot \exp  \sum_b \int_b  F_\phi \\
&= H(p;A,\eta) \cdot \exp \int_{\Sigma} F_\phi
\end{split}
\end{equation} 
where the $\beta$ terms cancel due to Stokes' theorem and the
usual combinatorial arguments. Since $F_\phi$ is a curvature we have 
\begin{equation}
H(\phi^{-1}P;\phi^*A,\phi^*\eta) = H(p;A,\eta)
\end{equation}
so the holonomy is an invariant of the gauge transformation of a 
bundle gerbe. 

In the case where $\Sigma$ has boundary, there are extra terms in the
function $H$ on the space of trivialisations,
\begin{equation}
\exp \int_{\Sigma} F_\phi + \sum_{e,b} \int_e \beta_{\rho(e)}
\end{equation}
These are independent of the
choice of trivialisation so they do not affect the transition functions
on the bundle over $\partial \Sigma$. The local connections of this 
bundle gain an extra term $\int_{\partial \Sigma} \iota_\xi F_\phi$ which
is $D$-trivial. 

For bundle 2-gerbes the situation is very similar. Let $(P,Y,X,M;A,\eta,\nu)$
be a bundle 2-gerbe and let $\phi$ be a gauge transformation. There is 
an isomorphism of bundle 2-gerbes
\begin{equation}
\phi^{[2]*} P \otimes P^* = \delta((1,\phi)^{-1} P)
\end{equation}
by the same arguments as in the lower cases. We may give the trivialisation
a connection and curving which are compatible with $\delta$. The 3-curvings
satisfy 
\begin{equation}
\phi^* \nu = \nu + (1,\phi)^*\omega
\end{equation} 
where $\omega$ is the three curvature of the bundle gerbe $(P,Y,X^{[2]})$
which satisfies $\omega = \delta(\nu)$. By similar arguments to the 
bundle gerbe case above, the bundle gerbe $(1,\phi)^{-1}P$ 
descends to a bundle gerbe $P_\phi$ on $M$ with 3-curvature $\omega_\phi$
which satisfies 
\begin{equation}
 (1,\phi)^*\omega = \pi^* \omega_\phi + d\beta
\end{equation}
where $\beta$ is the 2-curving of the trivialisation.
Thus we see that $\omega_\phi$ is the $D$-obstruction. Since it is a curvature
then the $D$-obstruction is trivial and the holonomy of bundle 2-gerbes is
invariant under transgression, though as in the previous terms there will
be a different choice of section of the bundle on the mapping space in the
case with boundary.

\subsubsection{$G$-Gauge Transformations} 
We have seen examples of bundle gerbes and bundle 2-gerbes, namely
lifting bundle gerbes and the bundle 2-gerbe associated with a principal
$G$-bundle, which have the following general form:
\begin{equation}
\begin{array}{ccc}
& & Q \\
& & \downarrow / \Downarrow\\
& & G \\
& \stackrel{\rho}{\nearrow} & \\
P_G \times G & \rightrightarrows & P_G \\
& & \downarrow  \\
& & M
\end{array}
\end{equation}
where $\downarrow / \Downarrow$ indicates that $Q$ may be a bundle or 
bundle gerbe and $\rho$ is projection 
onto the second factor of $P \times G$ which may also be thought 
of as the map $P^{[2]} \rightarrow G$ defined by $p_2 = p_1 \rho(p_1,p_2)$.
We shall refer to a gauge transformation of the $G$-bundle 
$P \rightarrow M$ as a {\it $G$-gauge transformation}.
\begin{proposition}\label{Ggauge}
Let $Q$ be a bundle (2-)gerbe as described above. A $G$-gauge transformation 
of $P$ defines a gauge transformation of $Q$.
\end{proposition}
\begin{proof}
Let $\psi: P \rightarrow P$ be a $G$-gauge transformation. By definition 
we have $\pi \circ \psi = \pi$ so we need only verify that there is 
an isomorphism of line bundles (or stable isomorphism of bundle gerbes)
$\psi^{[2]*} (\rho^{-1} Q) = \rho^{-1}Q$. On $P^{[2]}$ we have
$\psi^{[2]}(p_1,p_2) = \psi^{[2]}(p_1,p_1g_{12})
= (\psi(p_1),\psi(p_1g_{12})) = (\psi(p_1),\psi(p_1)g_{12})$. So on
$P \times G$ we have $\psi^{[2]}(p,g) = (\psi(p),g)$ and 
$\rho ( \psi^{[2]}(p,g)) = g = \rho(p,g)$, therefore 
$\psi^{[2]*} (\rho^{-1} Q) = \rho^{-1}Q$ as required.
\end{proof}

%% file: chapter8.tex
\chapter{Applications}

We consider some applications of the various constructions in bundle gerbe 
theory which we have discussed to physics.

\section{The Wess-Zumino-Witten Action}\label{WZWact}
We shall review the bundle gerbe model of the Wess-Zumino-Witten (WZW) 
action as described in \cite{camimu}. This example
serves as motivation for the use of bundle gerbe holonomy
to study topological actions. Furthermore the WZW theory plays a role
in the discussion of Chern-Simons theory which follows. 

Following \cite{camimu} the WZW action is defined as a function on the
space of maps from a Riemann surface $\Sigma$ to a compact Lie group $G$,
which we denote by $\Sigma G$. This function is defined by the equation
\begin{equation}\label{defwzw}
WZW(\psi) = \exp \int_X \hat{\psi}^* \omega
\end{equation}
where $X$ is a 3-manifold with boundary $\Sigma$, $\hat{\psi}$ is an
extension of $\psi \in \Sigma M$ to $XM$ and $\omega$ is a closed 3-form
which generates the integral cohomology of $G$. This is well defined 
as long as such a $\hat{\psi}$ exists, for example if $G$ is simply connected. 
In this case \eqref{defwzw} is the holonomy of the 
tautological bundle gerbe with curvature $\omega$. 
When such a $\hat{\psi}$ does not exist we may replace \eqref{defwzw} with the 
holonomy of 
any bundle gerbe with curvature $\omega$, though, as is observed
in \cite{diwi} where similar constructions are made using differential
characters, this bundle gerbe is not uniquely determined by $\omega$. Any two
choices differ by a flat bundle gerbe which is classified by $H^2(G,U(1))$.
To eliminate this ambiguity the action must be defined in terms of
a full Deligne class rather than just the Dixmier-Douady class. This
leads to 
\begin{definition}
Let $\alpha \in H^2(G,\D^2)$ be 
a Deligne class. The WZW action evaluated on a map $\psi: \Sigma \rightarrow 
G$ is the holonomy of the class $\alpha$, that is, the flat holonomy of 
$\psi^* \alpha$.
\end{definition} 
If we represent $\alpha$ by a bundle gerbe with connection
and curving $(P,Y,G;A,\eta)$ then the action may be written as \cite{camimu}
\[
\exp \int_\Sigma \psi^* \eta - F_L
\]
where $F_L$ is the curvature of a trivialisation of $\psi^*P$. In this way 
each bundle gerbe with connection and curving over $G$ defines a WZW action.
This general form of the WZW action was given in terms of a transgression
formula by Gawedski \cite{gaw}.

From \cite{fre} and \cite{gaw} we see that when we attempt to define this 
action for surfaces with boundary 
we need to consider line bundles over the boundary maps, just as with
holonomy. 
Recall
that we start with the holonomy function on $\Sigma G$, which for each 
$\phi \in \Sigma G$ may be thought of as the evaluation of the
corresponding flat holonomy class $\chi$. This function is 
extended to surfaces with boundary in such a way that it is 
multiplicative with respect to unions so that given 
two surfaces with the same boundary the product is equal to the
holonomy of the combined surface. It turns out that such a function can
in general only be defined locally on $\Sigma^\partial G$. These local
functions are used to define a trivial bundle with connection
which pulls back to $\partial \Sigma G$ to give a possibly non-trivial
bundle with connection. The local data corresponding to this 
line bundle as derived in Chapter 6 agrees with the formulae given by 
Gawedski \cite{gaw}.


In the case where $G$ is simply connected then 
the theory of transgression of tautological bundle gerbes tells us that the 
bundle on the loop space is just the bundle over $\L_0 G$ in the
definition of the bundle gerbe. This is of the form
\[
\begin{array}{ccc}
& & S^1 \\
& \stackrel{\rho}{\nearrow} & \\
S^2(G)  & \rightrightarrows & D^2(G) \\
& & \downarrow \\
& & \L_0 G 
\end{array}
\]

Furthermore if $G$ has an integral bilinear form $<.,.>$ then we may
write $\omega = -\frac{1}{6}<\theta\wedge[\theta\wedge\theta]>$ where 
$\theta$ is the Maurer-Cartan form and where we have used the same 
normalisation as \cite{fre}. When $G$ is simple we have $H^3(G,\Z) = \Z$
and it is generated by this $\omega$.
We know the connection and curvature of $L$ since these come from the
corresponding objects on the tautological bundle gerbe. The connection
on $D^2( G)$ is $\int_{D} ev^* \omega$ and the
curvature is $\int_{S^1} ev^* \omega$. We would like to 
get more concrete expressions of a similar nature to those in \cite{fre}.
To do this recall that for the Maurer-Cartan form, $\theta$, we
have $\iota_\xi \theta = \xi$.
Consider the connection evaluated at $\Xi \in T_\phi (D^2( G))$,
\begin{equation}
\begin{split}
\int_{D} \iota_\Xi \omega &= -\frac{1}{6} \int_{D}
\iota_\Xi <\phi^{-1}d\phi \wedge [\phi^{-1}d\phi \wedge \phi^{-1}d\phi]> \\
&= -\frac{1}{2} \int_{D} <\Xi \wedge [\phi^{-1}d\phi\wedge \phi^{-1}d\phi]>
\end{split}
\end{equation}
For the curvature evaluated at vectors $\xi_1,\xi_2 \in T_\gamma (\L_0 G)$ 
we have
\begin{equation}
\begin{split}
\int_{S^1} \iota_{\xi_1} \iota_{\xi_2} \omega &= -\frac{1}{6}
\int_{S^1} \iota_{\xi_1} \iota_{\xi_2} <\gamma^{-1}d\gamma
\wedge [\gamma^{-1} d\gamma \wedge
\gamma^{-1} d\gamma]> \\
&= -\frac{1}{2}
\int_{S^1} \iota_{\xi_1} <\xi_2 \wedge [\gamma^{-1} d\gamma \wedge
\gamma^{-1} d\gamma]>\\ 
&= - \int_{S^1} <[\xi_1,\xi_2]\wedge \gamma^{-1}d\gamma>
\end{split}
\end{equation}
Proposition \ref{propLG} implies that the hermitian
lines over the loop space defined by the transgression of the
WZW bundle gerbe give the standard central extension of the loop
group \cite{fre}.


\section{The Chern-Simons Action}\label{CSa}

Our description of basic Chern-Simons (CS) theory follows Freed \cite{fre} and
Dijkgraaf-Witten \cite{diwi}. We show that there is a bundle gerbe 
interpretation of the cases they deal with and see that it is 
useful for generalisation to more general theories. We show that the 
bundle gerbe CS theory reproduces the expected results when 
restricted to specific cases (usually relying on restriction of the 
group $G$ such that it satisfies particular properties).

Let $G$ be a compact Lie group and $X$ an
oriented 3-manifold. Let $P_G \rightarrow X$ be a principal
$G$-bundle with connection 1-form $A$. Define
a 3-form on $P_G$, called the {\it Chern-Simons form} by
\begin{equation}
CS(A) = \Tr(A \wedge dA) + \frac{2}{3}\Tr(A \wedge  A \wedge A )
\end{equation}
If the bundle $P_G \rightarrow X$ is trivial, with section $s$, then the 
{\it Chern-Simons action} associated a 3-manifold $X$, is defined by 
\begin{equation}
\exp \int_X s^* CS(A)
\end{equation}
Ideally this should be independent of the choice of section. A change of
section is given by a gauge transformation $\phi:P_G\rightarrow P_G$, or 
alternatively $g_\phi :P_G \rightarrow G$. Under such a gauge 
transformation the CS form transforms as follows,
\begin{equation}
\phi^* CS(A)= CS(A) + d\Tr(g_\phi^{-1}Ag_\phi\wedge g_\phi^{-1}dg_\phi)
- \frac{1}{3} \Tr(g_\phi^{-1}dg_\phi)^3
\end{equation}
This suggests that for the action to be independent of the choice 
of section we should require that the trace be normalised to
make $\frac{1}{3} \Tr(g_\phi^{-1}dg_\phi)^3$ a $2\pi$-integral form.

If $A$ extends over a 4-manifold $W$ such that $\partial W = X$ then
the action is 
\begin{equation}
\exp \int_W \Tr(F \wedge F)
\end{equation}
Note that if the bundle is non-trivial then this definition of the 
action still makes sense as long as the bundle and connection extend over
$W$.

This situation very closely parallels the problem of defining the WZW action 
for general $G$ (\S \ref{WZWact})
and the problem of generalising the tautological bundle gerbe 
to get holonomy reconstruction (\S \ref{DCHR}). 
This suggests that a general definition of the CS action may be obtained 
by considering it as the holonomy of a bundle 2-gerbe. Furthermore the 
dependence on a principal $G$-bundle with connection suggests that we are
interested in particular in a bundle 2-gerbe associated with a $G$-bundle.

\subsubsection{The Chern-Simons Bundle 2-Gerbe}
This construction is based on 
the bundle 2-gerbe associated with the principal bundle 
$P_G$. Thus we require that $G$ is connected, simply connected and
simple. A 2-gerbe of a similar nature was described in \cite{brmc1}.
Our basic geometric structure is given in the following diagram
\[
\begin{array}{ccc}
\rho^{-1} Q[\omega] & &  \\
\Downarrow & &  \\
P_G^{[2]} & \rightrightarrows & P_G  \\
& & \downarrow \\
& & X  
\end{array}
\]
where $\rho : P_G^{[2]} \rightarrow G$ satisfies $p\rho(p,q) = q$ and
$Q[\omega]$ is a tautological bundle gerbe associated with a 
curvature 3-form $\omega$ on $G$. 
Following \cite{brmc1} we let
$\omega = k \Tr (g^{-1}dg \wedge g^{-1}dg \wedge g^{-1}dg)$. Here the 
trace replaces the more general bilinear form we considered in the case
of WZW theory.
In this case we can set 
$\beta = 3k \Tr(g_1^{-1}dg_1 \wedge dg_2g_2^{-1})$ and we have 
$\delta(\omega) = d\beta$ and
$\delta(\beta) = 0$. 
Consider the following diagram:
\[
\begin{array}{ccc}
P_G^{[3]} & \rightarrow & G \times G \\
\downarrow & & \downarrow \\
P_G^{[2]} & \rightarrow & G
\end{array}
\]
We have $\omega \in \Omega^3(G)$ and $\beta \in \Omega^2(G \times G)$. We can pull these
back to $\rho^* \omega \in \Omega^3(P_G^{[2]})$ and $\rho^* \beta \in 
\Omega^2(P_G^{[3]})$.
Now we have $ \delta \rho^* \beta = \rho^* \delta  \beta = 0$ therefore there exists 
$\epsilon \in \Omega^2(P_G{[2]})$ such that $\delta \epsilon = \rho^* \beta$. On 
$\Omega^3(P_G^{[2]})$ we have the equation
\[
\delta (\rho^* \omega  - d\epsilon) = \rho^* (\delta \omega - d \beta) = 0
\]
so we may define $\alpha \in \Omega^3(P_G)$ such that
\[
\delta \alpha = \rho^* \omega - d\epsilon.
\]

When we pull back the tautological bundle gerbe by $\rho$ to $P_G^{[2]}$, 
the curvature pulls back to $\rho^*\omega$, however this is
not adequate as the 3-curvature on our bundle 2-gerbe
(meaning the 3-curvature of the bundle gerbe
$\rho^{-1}Q[\omega]$) since this is
not $\delta$-exact for $\delta : P_G \rightarrow P_G^{[2]}$. As we have shown above,
subtraction of $d \epsilon$ will result in $\delta$-exactness. 
This is a specific example of the general $\epsilon$ referred to in
proposition \ref{tranapr}.
The 2-curving of the 
bundle gerbe is given by 
\begin{eqnarray*}
d\tilde{\eta} &=& \pi^* \rho^* \omega - \pi^* d\epsilon \\
&=& \rho^* \pi^* \omega - d \pi^* \epsilon \\
&=& \rho^* d\eta - d \pi^* \epsilon 
\end{eqnarray*}
where $\eta$ is the curving of the tautological bundle gerbe. So let
$\tilde{\eta} = \eta - \pi^* \epsilon$. 

Now we will find solutions for $\epsilon$ and $\alpha$. To do this
we will identify $P_G \times G$ and $P_G^{[2]}$ via the map 
$(p,g) \mapsto (p,pg)$. Similarly we have a map from $P_G \times G \times G$
to $P_G^{[3]}$ given by $(p,g_1,g_2) \mapsto (p, pg_1,pg_1g_2)$. This will 
change the $\delta$ maps. We want the following diagram to commute
\[
\begin{array}{ccc}
P_G^{[3]} & \rightarrow & P_G \times G \times G \\
\downarrow & & \downarrow \\
P_G^{[2]} & \rightarrow & P_G \times G
\end{array}
\]
For each map $\pi_i$ of $\delta$ we will have a diagram which
shows what the induced map from $P_G \times G \times G$ to $P_G \times G$
should be. For $\pi_0$ we have
\[
\begin{array}{ccc}
(p,pg_1,pg_1g_2) & \rightarrow & (p,g_1,g_2) \\
\downarrow & & \downarrow \\
(pg_1,pg_1g_2) & \rightarrow & (pg_1, g_2)
\end{array}
\]
For $\pi_1$,
\[
\begin{array}{ccc}
(p,pg_1,pg_1g_2) & \rightarrow & (p,g_1,g_2) \\
\downarrow & & \downarrow \\
(p,pg_1g_2) & \rightarrow & (p, g_1g_2)
\end{array}
\]
For $\pi_2$,
\[
\begin{array}{ccc}
(p,pg_1,pg_1g_2) & \rightarrow & (p,g_1,g_2) \\
\downarrow & & \downarrow \\
(p,pg_1) & \rightarrow & (p, g_1)
\end{array}
\]
These diagrams give us the following equations:
\begin{eqnarray*}
\pi_0(p,g_1,g_2) &=& (pg_1,g_2) \\
\pi_1(p,g_1,g_2) &=& (p,g_1g_2) \\
\pi_2(p,g_1,g_2) &=& (p,g_1) 
\end{eqnarray*}
\begin{lemma}
Let $\epsilon \in \Omega^2(P_G \times G)$ be defined by $\epsilon = 3k\Tr
(A\wedge \dg)$
where $A$ is a connection for $P_G \rightarrow X$.
Then $\delta \epsilon = \beta$ and thus $\delta(\rho^* \omega - d\epsilon)
=0$.
\end{lemma}
\begin{proof}
We will omit $3k$ since it appears in all expressions.
We need to evaluate $\delta \epsilon = \pi_0^* \epsilon - \pi_1^* \epsilon
+ \pi_2^* \epsilon$. There are three types of pullback map that we will need.
Let $(Z,\xi) \in TP_G \times TG$. Then $p_* Z = Z$ and $g_* \xi = \xi$. This leaves
$pg_* Z$. Applying the chain rule gives the result $pg_* Z = R_{g*} Z + g^{\#}$ where
$R_g$ is the right action of $g$ and $g^{\#}$ is the fundamental field of $g$. We can
now use the two defining properties of $A$, which are $R_{g*} A = g^{-1} A g$ and
$A(g^{\#}) = g$. Now we can write
\begin{eqnarray*}
\pi_0^* \epsilon &=& \Tr
((g_1^{-1} A g_1 - g_1^{-1}dg_1 )\wedge dg_2 g_2^{-1}) \\
\pi_1^* \epsilon &=& \Tr
(A \wedge d(g_1g_2) (g_1g_2)^{-1}) \\
&=& \Tr (A \wedge (dg_1 g_1^{-1} + g_1 dg_2 g_2^{-1} g_1^{-1})) \\
\pi_2^* \epsilon &=& \Tr (A \wedge dg_1 g_1^{-1}) 
\end{eqnarray*}
Putting these together we have 
\[
\delta \epsilon = \Tr (g_1^{-1}dg_1 \wedge dg_2 g_2^{-1}) 
\]
where Ad invariance of the trace has been used to eliminate the other terms.
\end{proof}

Now we can write down an expression for $\rho ^* \omega - d\epsilon$.
The map $\rho: P_G \times G \rightarrow G$ is defined by $(p,g) \mapsto
g$, so $\rho ^* \omega = \omega$. Applying $d$ to $\epsilon$ yields
\begin{eqnarray*} 
d \Tr (A \wedge \dg) &=& \Tr (dA \wedge \dg) - \Tr (A \wedge d (\dg )) \\
&=& \Tr (dA \wedge \dg) - \
Tr(A \wedge \dg \wedge \dg)
\end{eqnarray*}
Thus we have 
\begin{equation*}
\begin{split}
\rho ^* \omega - d\epsilon &= k \Tr (g^{-1}dg \wedge g^{-1}dg \wedge g^{-1}dg)
- 3k \Tr (dA \wedge \dg) \\
&  \quad + 3k \Tr (A \wedge \dg \wedge \dg).
\end{split}
\end{equation*}

\begin{proposition}
Let $\alpha \in \Omega^3(P_G)$ be defined by the Chern-Simons form
\[
\Tr (A \wedge dA) + \frac{2}{3}\Tr (A \wedge  A \wedge A ).
\]
Then 
\[
\delta( -3k \alpha) = \rho ^* \omega - d\epsilon.
\]
and $\alpha$ is a 3-curving for the bundle 2-gerbe described above.
\end{proposition}
\begin{proof}

The map $\delta$ is given by $\pi_2^* - \pi_1^*$ where
$\pi_1 : (p,g) \mapsto p$ and $\pi_2 : (p,g) \mapsto pg$.
First we will calculate $\pi_2^* Tr(A \wedge dA)$. Recall
that $pg_* Z = R_{g*} Z + g^{\#}$. Throughout the following calculations 
we will make use of the Ad-invariance and the cyclic symmetry of
the trace. We also omit the symbol $\wedge$.
\begin{eqnarray*}
pg^* \Tr (A \wedge dA) &=& \Tr ((g^{-1} A g + g^{-1} dg) 
\wedge d(g^{-1} A g + g^{-1} dg)) \\
&=& \Tr ((g^{-1} A g + g^{-1} dg) \\
&& \quad \wedge
(-g^{-1} \dg  A g + g^{-1} dA g - g^{-1} A dg - g^{-1} \dg dg)) \\
&=& - \Tr (A A\dg) + \Tr (AdA) - \Tr (A A\dg) \\
&& \quad - \Tr (A\dg\dg)
-\Tr (A \dg \dg) + \Tr (dA \dg)  \\
&& \qquad - \Tr (A \dg \dg) - \Tr (\dg\dg\dg) \\
&=& \Tr (AdA) - 2\Tr (AA\dg) + \Tr (dA\dg) \\
&& \quad - 3 \Tr(A\dg\dg) - \Tr (\dg\dg\dg)
\end{eqnarray*}
\begin{eqnarray*}
pg^* \Tr (A \wedge A \wedge A) &=& \Tr ((g^{-1}Ag + g^{-1}dg)\wedge
(g^{-1}Ag + g^{-1}dg) \\
&& \quad \wedge(g^{-1}Ag + g^{-1}dg)) \\
&=& \Tr (AAA) + 3\Tr (AA\dg) + 3\Tr (A\dg\dg) \\
&& \quad + \Tr (\dg\dg\dg)
\end{eqnarray*}
\begin{equation*}
p^* \Tr (A\wedge dA) = \Tr (AdA)
\end{equation*}
\begin{equation*}
p^* \Tr (A\wedge A\wedge A) = \Tr (AAA)
\end{equation*}
Putting all of this together we get
\begin{eqnarray*} 
&&(pg^* - p^*)\Tr (A\wedge dA + \frac{2}{3}A\wedge A\wedge A) = \\
&& \qquad \Tr (AdA)- 2\Tr (AA\dg) + \Tr (dA\dg) \\
&& \qquad- 3\Tr (A\dg\dg) - \Tr (\dg\dg\dg) + \frac{2}{3}Tr(AAA) \\
&&\qquad + 2\Tr (AA\dg) + 2\Tr (A\dg\dg) + \frac{2}{3}\Tr (\dg\dg\dg) \\
&& \qquad - \Tr (AdA) - \frac{2}{3}\Tr (AAA)
\end{eqnarray*}
Collecting terms gives
\[
\Tr (dA\dg) - \Tr (A\dg\dg) - \frac{1}{3}\Tr (\dg\dg\dg)
\]
which is the desired result.
\end{proof}

We call this bundle 2-gerbe with connection and curvings the 
{\it Chern-Simons bundle 2-gerbe}. Its holonomy 
satisfies the properties of the 
Chern-Simons action, this shows
\begin{proposition}
The Chern-Simons action associated with a principal $G$-bundle, where
$G$ is connected, simply connected and simple, 
may be realised as the holonomy of
the Chern-Simons bundle 2-gerbe over a closed 3-manifold.
\end{proposition}
The real purpose of using bundle gerbe theory to approach this problem is
to understand how it generalises when we relax the requirements on 
$G$ in the proposition above. So far we have removed only the requirement 
that the bundle $P$ be trivial in the original definition of the action. 
This is possible because defining the holonomy only requires that the bundle
2-gerbe be trivial, which is always true over a 3-manifold. The existence
of a section of $P$ implies this triviality however there exist trivial
bundle 2-gerbes for which such a section does not exist. 

Suppose we wish to allow $G$ to be only semi-simple instead of simple. In this
case we can still define the CS form using the Killing form on the 
Lie algebra. The difference with the simple case is that we may have
$H^3(G,\Z) \neq \Z$. We can still define the CS-bundle 2-gerbe as above,
the only difference is that the possible bundle gerbes over $G$ in the 
construction do not necessarily account for all bundle gerbes over $G$.
To allow non semi-simple groups we may replace the trace with an 
invariant quadratic polynomial on the Lie algebra, as in \cite{diwi}.
If $G$ is not simply connected then the tautological bundle gerbe on $G$
must be replaced with a general bundle gerbe with curvature $\Tr (g^{-1}
dg)^3$. This may require additional data (such as a full Deligne class) since
the bundle gerbe is no longer determined completely by its curvature.

Using this interpretation we may consider further 
aspects of CS theory in terms of the theory of holonomy of bundle
2-gerbes.

\subsubsection{Chern-Simons Lines and Gauge Invariance}
For the purposes of this section we shall follow \cite{fre} and
set $\omega = -\frac{1}{3}\Tr(g^{-1}dg)^3$ and $\epsilon
= -Tr(Adgg^{-1})$ so that the curving is precisely $CS(A)$.

It is a standard fact (\cite{fre},\cite{diwi}) that given a 
3-manifold, $X$, with non-empty boundary $\partial X$, 
the CS action cannot be defined as 
a function, rather it is a section of a line bundle called a
{\it Chern-Simons line}. This is, of course, exactly what we 
would expect since we have interpreted the CS action as the
holonomy of a bundle 2-gerbe. We shall give arguments as to 
why the line bundle corresponding to the transgression of 
a bundle 2-gerbe as described in Chapter 6 above is the 
same as the CS lines described in \cite{fre} and \cite{diwi}.

Recall that Proposition \ref{tranapr} tells us that 
when $M$ is 2-connected and $G$ is simply connected
the transgression of 
the bundle 2-gerbe associated with a principal
bundle is described by the following diagram:
\[
\begin{array}{ccc}
& & S^1 \\
& \stackrel{\rho^*\hol( Q[\omega])\cdot I(\epsilon)}{\nearrow} & \\
S^2 (P_G) \times S^2( G) & \rightrightarrows & S^2 (P_G) \\
& & \downarrow \\
& & S^2( M)
\end{array}
\]
For purposes of comparison it is useful to describe the fibres of the 
line bundle which may be obtained from this
bundle 0-gerbe. Over $\Sigma \in S^2(M)$ 
the fibre consists of elements of an equivalence class 
$[\tilde{\Sigma},\theta]$, where $\tilde{\Sigma}$ is a lift of
$\Sigma$ to $P_G$ and $\theta \in S^1$, and the equivalence is given
by $[\tilde{\Sigma}_1,\hol_{(\tilde{\Sigma}_1 , \tilde{\Sigma}_2)}
(\rho^{-1} Q[\omega])] \sim [\tilde{\Sigma}_1 , 1]$. 
The trivial bundle over $X$, where $\partial X = \Sigma$
is obtained by pulling back this bundle using the restriction to the 
boundary. A trivialisation is given by the extension of
the holonomy function on $S^3(M)$.

We now recall some earlier results on gauge transformations which are
relevant here. By Proposition \ref{Ggauge} any $G$-gauge transformation
on $P_G \rightarrow X$ is also a bundle 2-gerbe gauge transformation, so we
may apply the results of \S \ref{gaugebb2} to examine gauge invariance of the
CS action. In the case of a closed 3-manifold the gauge invariance of
the holonomy implies the same for the CS action. In the case with boundary
there are two additional terms in the section of the trivial
line bundle,
\begin{equation}
\int_X \omega_\phi - \sum_{b,w} \int_b \beta_{\rho(b)}
\end{equation}
Recall that $\omega_\phi$ is the 3-curvature of the bundle gerbe 
$(1,\phi)^{-1} \rho^{-1} Q[\omega]$. Since $\rho \circ (1,\phi) = g_\phi$ this is 
just $\Tr(g^{-1}_\phi dg_\phi)^3$. Recall also that the 2-form $\beta$ 
arises from the failure of the curvature of $(1,\phi)^{-1}P_G$ to descend
to a curvature on M. In this case the former is 
$(1,\phi)^* (\rho^* \omega - d\epsilon) = g_\phi^* \omega -
d(1,\phi)^* \epsilon$. Observe that $\delta(g_\phi^* \omega) 
= \Ad_g \omega - \omega = 0$ by the invariance of the trace, so this
part descends and $\beta = (1,\phi)^* \epsilon$. For $\epsilon =
\Tr (A dg g^{-1})$ we have $(1,\phi)^* \epsilon =
\Tr(A  dg_\phi g_\phi^{-1})$, so the section changes by 
\begin{equation}
\exp (\int_X -\frac{1}{3}\Tr(g^{-1}_\phi dg_\phi)^3 + \sum_{b,w} \int_b 
\Tr(A_{\rho(b)}  dg_\phi g_\phi^{-1}))
\end{equation}
where we have used the fact that $\delta(dg_\phi g_\phi^{-1})=0$, so
only $A$ need be expressed in local form. If the $G$-bundle $P_G \rightarrow 
M$ is trivial with section $s$ then we recover proposition 2.10 of 
\cite{fre}, where the section changes by
\begin{equation} \label{freedline}
\exp(\int_X -\frac{1}{3}\Tr(g^{-1}_\phi dg_\phi)^3
+ \int_{\partial X} \Tr(s^*A dg_\phi g_\phi^{-1}))
\end{equation}
We may now compare our construction of the CS lines with that of Freed 
\cite{fre}. Suppose the bundle $P_G \rightarrow Y$ is trivial where 
$Y$ is a closed 2-manifold. We think of $Y$ as the image in $M$ of 
an element of 
$\Map(\Sigma,M)$. Each choice of a section $Y \rightarrow P_G$ gives
a lift $\tilde{Y}$. Let $s_1$ and $s_2$ be two such choices with 
corresponding lifts $\tilde{Y}_1$ and $\tilde{Y}_2$. There is a 
$G$-gauge transformation, $\phi$, which gives the difference between
these two sections. The pair $(\tilde{Y}_1,\tilde{Y}_2) \in (\Sigma P_G)^{[2]}$
is then equivalent to $(\tilde{Y}_1,g_\phi(\tilde{Y}_1)) \in \Sigma P_G
\times \Sigma G$. Thus we have $\rho(\tilde{Y}_1,\tilde{Y}_2) = 
g_\phi(\tilde{Y}_1)$, though it should be kept in mind that the $\tilde{Y}_2$ 
dependence is contained in the definition of $\phi$. If we let $Y = \partial X$
then the equivalence relation in the definition of the line bundle 
is given by the function \eqref{freedline}, which is used in an analogous 
way in the construction of the line in \cite{fre}.

Viewing CS theory from a bundle gerbe point of view it is no surprise that the 
WZW action arises when we apply gauge transformations. The CS bundle 2-gerbe 
includes a bundle gerbe over $G$ with curvature $Tr(g^{-1}dg)^3$ in its 
definition, this is the bundle gerbe which produces the most common form
of the WZW action (that is, the one obtained when $G$ is simple) via 
its holonomy. That the holonomy of this bundle gerbe should be relevant
here follows from the results on the effects of gauge transformations
on holonomy. 

\subsubsection{Relationship with the Central Extension of the Loop Group}
In the previous section we considered the transgression of the 
Chern-Simons bundle 2-gerbe to a line bundle on $S^2( M)$. 
We have also seen (\S \ref{ltb2g}) that it is possible to 
transgress a bundle 2-gerbe to a bundle gerbe on the loop space.
We are interested here in the case where $G$ is simply connected, so the
only bundle gerbe (up to D-stable isomorphism) is the tautological one.
In the case of the CS bundle 2-gerbe we get the following bundle gerbe:
\[
\begin{array}{ccc}
& & \tau_{S^1} (Q[\omega]) \\
& & \downarrow \\
& & \L G \\
& \stackrel{\L\rho}{\nearrow} & \\
\L P \times \L  G & \rightrightarrows & \L  P \\
& & \downarrow \\
& & \L M 
\end{array}
\]
where $\tau_{S^1} (Q[\omega])$ is the loop space transgression of the tautological
bundle gerbe on $G$. Recall (see proposition \ref{propLG}) that this 
transgression is isomorphic to the bundle corresponding to the 
central extension $\widetilde{\L G} \rightarrow \L G$, so we have
\[
\begin{array}{ccc}
& & \widetilde{\L G} \\
& & \downarrow \\
& & \L G \\
& \stackrel{\L\rho}{\nearrow} & \\
\L P_G \times \L G & \rightrightarrows & \L P_G \\
& & \downarrow \\
& & \L M 
\end{array}
\]
This is the lifting bundle gerbe which describes the obstruction to 
lifting the structure group of the bundle $\L P_G \rightarrow \L M$ from
$\L G$ to $\widetilde{\L G}$. This result is given in terms of gerbes in
\cite{brmc1} and \cite{gom}.\\

In conclusion, we have seen that the standard Chern-Simons action may be
interpreted as the holonomy of a bundle 2-gerbe. Just as with
WZW theory, this allows us to understand the failure of
the action to be well defined in the general case, that is, when there
is no section of the $G$-bundle or it cannot be defined as 
an integral of a 4-curvature. These features are key characteristics
of holonomy. A number of other features of the theory have been explained 
in bundle gerbe terms. In \cite{diwi} more general theories are discussed
in terms of general WZW theories. In terms
of bundle gerbes these could be interpreted as a generalisation of the
associated bundle gerbe to a case where the bundle gerbe on $G$ is 
not tautological (an example that was similar to this, $L\cup J$ was 
described in section \ref{bun2ger}). Even more generally 
differential characters are used, since these correspond to 
classes in Deligne cohomology this suggests that bundle 2-gerbes can play the
same role.

\section{D-Branes and Anomaly Cancellation}

In \cite{cajomu} it is shown how to use bundle gerbes to cancel anomalies
in $D$-brane theory. Here we concentrate on the local aspects of this
approach as an application of the holonomy of bundle gerbes.

The basic situation described in \cite{cajomu} is that we have 
actions which are functions associated with maps of a surface with
boundary, $\Sigma$, into a manifold $M$ with submanifold $Q$ such
that $\partial \Sigma \subset Q$. The submanifold $Q$ is referred to as the 
brane. The action turns out to be a 
section of a trivial line bundle on $\Sigma M$. This should come as 
no surprise by now since we have seen examples of actions which behave
as holonomies and this is precisely the behaviour we would expect from
the holonomy of a bundle gerbe. This failure of the action to be
a well defined function is called the anomaly. Anomaly cancellation 
involves the introduction of an extra term (or terms) such that 
together they define a function. Our approach is guided by the knowledge
that two trivialisations differ by a global function, so to 
cancel the anomaly we need to find another trivialisation of the 
bundle on $\Sigma M$. The general technique for doing this is as 
follows. Recall that if we transgress a bundle gerbe to the loop space
then we can obtain the trivial line bundle over $\Sigma M$ by 
pulling back the line bundle on the loop space with the map
$\partial: \Sigma M \rightarrow \L M$, which is induced from the restriction
to the boundary. Suppose we have a term in the action which corresponds
to a section of the trivial bundle over $\Sigma$ corresponding
to the transgression, $L \rightarrow \L M$,  of a bundle gerbe, $P$ on $M$. 
If we can find another bundle $L'$ on the loop space which is isomorphic to 
$L$ then the product $L \otimes L^{'*}$ will be trivial, and this 
trivialisation will induce a trivialisation of $\partial^{-1}L 
\otimes \partial^{-1} L^{'*}$. Thus the combination of the 
usual trivialisation of the pull back of $L'$ to $\Sigma M$ and 
and the trivialisation of the product bundle on $\L M$ will cancel the 
anomaly. Furthermore the functoriality of transgression tells us that
a suitable bundle $L'$ may be found via the transgression to $\L M$ of
a bundle gerbe $P'$ with the same Dixmier-Douady class as $P$. This
bundle gerbe $P'$ is known as a $B$-field in the physics literature
and the requirement that $\dd (P) = \dd (P')$ leads to a natural 
division of the anomaly cancellation problem into three distinct 
cases.

First we consider the situation described by Freed and Witten \cite{frwi}
as interpreted in \cite{cajomu}. In this case the first term in the
action is derived from the transgression of a torsion bundle gerbe,
that is, a bundle gerbe with a torsion Dixmier-Douady class.
The Deligne class of this bundle gerbe is determined by the 
second Stiefel-Whitney class, $w_2 \in H^2(M,\Z_2)$, 
of the normal bundle to $Q$. This class determines a Deligne class 
$(w_{\alpha\beta\gamma}, 0, 0)$ 
where $w_{\alpha\beta\gamma} \in H^2(M,U(1))$ is induced by the inclusion
$\Z_2 \subset U(1)$. Let $P_{w_2}$ be a bundle gerbe which is classified
by this Deligne class.

The $B$-field is defined as a triple $(g_{\alpha\beta\gamma}, 
k_{\alpha \beta}, B_\alpha)$ which defines a Deligne cohomology class and 
hence a $D$-stable isomorphism class of bundle gerbes (note that 
$B$ is a 2-form field). Let $P_B$ be 
a representative of this class. In this case 
we specify that the Dixmier-Douady class of the $B$-field is equal to 
that of the torsion bundle gerbe described above.
Thus the two
transgression bundles on the loop space are isomorphic and there 
exists a section which may be pulled back to $\Sigma M$ to cancel 
the anomaly. We wish to get a local expression for this term.

The product $P^*_w \otimes P_B$ is represented locally by the 
Deligne class $(g_{\alpha\beta\gamma}w^{-1}_{\alpha\beta\gamma},
k_{\alpha\beta},B_{\alpha})$. The local formula for the transition functions
of the 
transgression to the loop space is obtained by 
applying equation \eqref{bgtls},
\begin{equation}\label{dbr1}
G_{01} = \exp \sum_{e} \int_e k_{\rho_0(e)\rho_1(e)} \cdot
\prod_{v,e} g^{-1}_{\rho_0(e)\rho_0(v)\rho_1(v)}
g_{\rho_0(e)\rho_1(e)\rho_1(v)}w_{\rho_0(e)\rho_0(v)\rho_1(v)}
w^{-1}_{\rho_0(e)\rho_1(e)\rho_1(v)}(v)
\end{equation}
Since $P_w$ and $P_B$ have the same Dixmier-Douady classes then
by functoriality of the transgression (see \S \ref{genres}) the line 
bundles $\tau(P_w)$ and $\tau(P_B)$ have the same Chern class, 
so there exists a trivialisation $(h_{\alpha \beta},A_\alpha)$ satisfying
\begin{eqnarray}
g_{\alpha \beta \gamma}w^{-1}_{\alpha\beta\gamma} = 
h_{\beta\gamma}h_{\alpha\gamma}^{-1}h_{\alpha\beta} \\
k_{\alpha\beta} = -d\log h_{\alpha\beta} + A_\alpha - A_\beta
\end{eqnarray}
The pair $(\un{h},\un{A})$ defines an $A$-field \cite{frwi}.
Substituting into \eqref{dbr1} and using the usual combinatorial arguments
we get
\begin{equation}
G_{01} = \exp \sum_e \int_e (A_{\rho_1(e)} - A_{\rho_0(e)}) \cdot
\prod_{v,e} h_{\rho_1(e)\rho_1(v)}h^{-1}_{\rho_0(e)\rho_0(v)}(v) 
\end{equation}
and thus we have local functions 
\begin{equation}
\Gamma_0 = \exp \sum_e \int_e A_{\rho_0(e)} \cdot
\prod_{v,e} h_{\rho_0(e)\rho_0(v)}(v) 
\end{equation}
satisfying
$\Gamma^{-1}_0\Gamma_1 = G_{01}$. These local functions may be pulled back
to give local functions (or equivalently sections of a trivial bundle) on
$\Sigma M$ and cancel the anomaly.

We would like to indicate how this local picture relates to the global
version given in \cite{cajomu}. Denote the transgressions of
$P_w$ and $P_B$ by $L_w$ and $L_B$ respectively. The original 
term in the action from which the anomaly arises is the Pfaffian of the
Dirac operator on the world sheet, denoted Pfaff, which is a section of
$L_w$. We refer to \cite{frwi} for further details since this 
term does not arise from bundle gerbe considerations. 

The bundle $L_B$ corresponds to the following bundle 0-gerbe,
\[
\begin{array}{ccc}
& & S^1 \\
& \stackrel{\hol(P_B)}{\nearrow} & \\
S^2(M) & \rightrightarrows & D^2(M) \\
& & \downarrow \\
& & \L M 
\end{array}
\]
so the bundle $\partial^{-1} L_B$ is given by
\[
\begin{array}{ccc}
& & S^1 \\
& \stackrel{\hol(P_B)}{\nearrow} & \\
\Sigma M \times_\pi S^2(M) & \rightrightarrows & \Sigma M \times_\pi D^2(M) \\
& & \downarrow \\
& & \Sigma M 
\end{array}
\]  
The section of
$\partial^{-1} L_B$ may be defined by $\phi_B(\Sigma,\sigma) = \hol(P_B;\Sigma
\# \sigma)$ where $\sigma \in D(M)$ satisfies $\partial \sigma = \partial
\Sigma$. The gluing property of holonomy ensures that this is a bundle
0-gerbe trivialisation on $\Sigma M \times_\pi D(M)$, so it defines a 
section of the line bundle $\partial^{-1} L_B$. 

To get a local expression let $\chi_B$ be a $D$-obstruction form
for $P_B$ over $\Sigma$ corresponding to a trivialisation $T_B$, then
\begin{equation}
\begin{split}
\hol(P_B;\Sigma \# \sigma) &= \exp \int_{\Sigma \# \sigma}
\chi_B \\
&= H_{int} (B;\Sigma) H^{-1}_{int}(B;\sigma) H_\partial(T_B;\partial \Sigma)
H_\partial^{-1} (T_B ; \Sigma) \\ 
&= H_{int} (B;\Sigma) H^{-1}_{int}(B;\sigma)
\end{split}
\end{equation}
where we have used the fact that the local expression for 
$\exp \int_{\Sigma \# \sigma}
\chi_B$ splits into terms on the interior of $\Sigma \# \sigma$, $H_{\int}$,
which are the extension of the holonomy formula on closed manifolds,
and terms on the boundary, $H_\partial$, which depend on a choice of
trivialisations. See Chapter 5 for further details. 

The section corresponding to the trivialisation of
$L_w \otimes L_B$ may be defined as a function $\lambda$ on
$D(Q)$ such that $\lambda (\sigma_2) = \lambda (\sigma_1)
\hol(\sigma_1 \# \sigma_2)$,
so it is a trivialisation of the transgression bundle 0-gerbe
\[
\begin{array}{ccc}
& & S^1 \\
& \stackrel{\hol(P_B\otimes P^*_w)}{\nearrow} & \\
S^2(Q) & \rightrightarrows & D^2(Q) \\
& & \downarrow \\
& & \L (Q)
\end{array}
\]
Note that this bundle is only defined on $Q$ since it is only on 
the brane that the Dixmier-Douady classes of $P_w$ and $P_B$ agree.
The standard section of this bundle is obtained by extending the
holonomy function to discs. Let $A$ be a trivialisation of 
$P_w^* \otimes P_B$ (the trivialisation defined by
the $A$-field). Then in terms of the corresponding
$D$-obstruction form $\chi_A$ we have
\begin{equation}
\begin{split}
\lambda_A (\sigma) &= \exp \int_\sigma \chi_A \\
&= H_{int}(B-w;\sigma) H_\partial(A;\Sigma) \\
&= H_{int}(B;\sigma) H_{int}^{-1}(w;\sigma) H_\partial(A;\Sigma)
\end{split}
\end{equation}
Now we combine $\phi_B$ and $\lambda_A$ to get
\begin{equation}\label{notcanc}
H_{int} (B;\Sigma) 
H_{int}^{-1}(w;\sigma) H_\partial(A;\Sigma)
\end{equation}
In this context the anomaly corresponds to $\sigma$ dependence, so
while we have cancelled some $\sigma$ terms there is still one left.
This is because we have not yet incorporated the section of $\partial^{-1}
L_w$. 
Consider this as a bundle 0-gerbe,
\begin{equation}
\begin{array}{ccc}
& & S^1 \\
& \stackrel{\hol(P_w)}{\nearrow} & \\
\Sigma_Q (M) \times S^2(Q) & \rightrightarrows & \Sigma_Q(M) \times_\pi D^2(Q) \\
& & \downarrow \\
& & \Sigma_Q(M) 
\end{array}
\end{equation}
Observe that in this case we cannot use the same approach that we used to
define the section $\phi_B$ since $P_w$ is only defined on $Q$ and in general
elements of $\Sigma_Q$ may not lie entirely in $Q$. Due to the 
definition of $P_w$ is turns out that there is a section of this bundle
called Pfaff \cite{frwi}. Given any such section we may find a 
$\C$-valued function (since the section may vanish) on 
$\Sigma_Q \times_\pi D(Q)$ via the corresponding local functions, $p_0$.
This function is defined by $\pi^* p_0 (\Sigma) H_{int}(w;\sigma)$.
It is easily verified that this is a section and is a globally defined 
function since the local dependence of the two terms cancels. Thus 
when we incorporate this term into \eqref{notcanc} the anomaly is cancelled
and we are left with terms derived from the Pfaffian, the $B$-field and
the $A$-field.

The second case, which appears in \cite{kap}, involves a $B$-field which
has a different Dixmier-Douady class to $P_w$, but it is still
required to be torsion. In this case the line bundle obtained by transgressing
$P^*_w \otimes P_B$ is no longer trivial so we need some further 
structure to cancel the anomaly. Before we introduce this we would like
to see the nature of the obstruction from a local point of view. The 
following arguments follow \cite{kap} closely. 
We have a torsion bundle gerbe $P_{(w,B)} = P^*_w \otimes P_B$
with Deligne class $(\un{g}\un{w}^{-1},\un{k},\un{B})$. Since the
image of $H^3(M,\Z)$ in $H^3(M,\R)$ is zero then the curvature is
exact, so denote it by $d\tilde{B}$. We now have a series of equations
\begin{eqnarray}
dB_\alpha &=& d\tilde{B} \\
B_\alpha - dm_\alpha &=& B_\beta - dm_\beta \\
k_{\alpha\beta} &=& m_\beta - m_\alpha + d\log q_{\alpha\beta} \\ 
g_{\alpha\beta\gamma}w^{-1}_{\alpha\beta\gamma} &=& 
q_{\beta\gamma}q_{\alpha\gamma}^{-1}q_{\alpha\beta}\zeta_{\alpha\beta\gamma}
\end{eqnarray}
where $m_\alpha$ are 1-forms, $q_{\alpha\beta}$ are $U(1)$-valued 
functions and $\zeta_{\alpha\beta\gamma}$ are $U(1)$-valued constants.
These constants correspond to the torsion class which measures the
obstruction to the equality of $\dd(P_w)$ and $\dd(P_B)$. 
Since $ g_{\alpha\beta\gamma}w^{-1}_{\alpha\beta\gamma}$ 
represents a torsion bundle gerbe it admits a bundle gerbe 
module, so from a local point of view there exist matrix
valued functions $\lambda_{\alpha\beta}$ satisfying
\begin{equation} 
g_{\alpha\beta\gamma}w^{-1}_{\alpha\beta\gamma}
= \lambda_{\beta\gamma}\lambda^{-1}_{\alpha\gamma}\lambda_{\beta\gamma}
\end{equation}
and so we have a sort of `non-Abelian trivialisation' of $\un{\zeta}$,
\begin{equation}
\zeta_{\alpha\beta\gamma} = \lambda_{\beta\gamma}q^{-1}_{\beta\gamma}
\lambda^{-1}_{\alpha\gamma}q_{\alpha\gamma}
\lambda_{\beta\gamma}q^{-1}_{\alpha\gamma}
\end{equation}
where it is assumed that all scalar functions are multiplied by the 
unit matrix of the appropriate dimension so that this expression 
makes sense.
We may view this in terms of a more general problem: if we are given a 
bundle gerbe with a trivialisation then we may find 
a trivialisation of the transgression bundle on the 
loop space, so if we have a bundle gerbe module represented 
locally by $(\lambda_{\alpha\beta}, \theta_\alpha)$ then we want to 
know to what extent we can use this 
to trivialise the bundle on the loop space. The answer is that in
general we cannot trivialise the bundle, this would violate
functoriality since the original bundle gerbe is non trivial, however
we can find a $\C$-valued function which `trivialises' it. The
distinction is analogous to that between a non-vanishing 
section of a line bundle and a section in general. It is a result of
Kapustin \cite{kap} that this is section is given by the trace of the
holonomy of the bundle gerbe module \footnote{Kapustin dealt with
Azumaya modules which have the same local data as bundle gerbe modules.}. 
We may realise this locally in terms of the holonomy of a non-Abelian 
bundle \cite{cajomu}. 
Over a disc the bundle gerbe which acts on the module is
trivial, choose a trivialisation $J$. Let $J$ be represented locally
by the pair $(K_\alpha,j_{\alpha\beta})$. The bundle $E \otimes J^*$ then 
descends to the disc. The trace of the holonomy of this bundle can be 
calculated over the boundary of the disc. To eliminate the $J$ dependence
we must introduce another term, $\exp \int_\sigma \chi_J$, where $\chi_J$ is
a $D$-obstruction form for the bundle gerbe $\zeta$ and trivialisation $J$.
It is easily shown \cite{cajomu} that this defines a section, as a $\C$ valued
function on $D(Q)$. 

To examine the anomaly cancellation from a local point of view we must
be careful as we cannot use the usual holonomy formula in the non-Abelian
case. When the boundary loop $\partial \sigma$ is triangulated the
holonomy breaks down in to an ordered product of parallel transport
terms along edges and jumps at vertices. Following Kapustin \cite{kap}
we denote the parallel transport for the bundle with connection
$\theta - K$ along the edge $e$ by $\hol_e(\theta_{\rho(e)} - K_{\rho_e})$.
The jumps are given by terms of the form $\lambda_{\rho(e)\rho(v)}j^{-1}_
{\rho(e)\rho(v)}$. The trace of the holonomy is then given by
\begin{equation}
\Tr [ \hol_{e_0}(\theta - K)_{\rho(e_0)} \cdot (\lambda j^{-1})_{\rho(
e_0)\rho(v_{01})}\cdot (j \lambda^{-1})_{\rho(v_{01})\rho(e_{1})}
\hol_{e_1}(\theta - K)_{\rho(e_1)} \ldots ]
\end{equation}
The Abelian parts may be pulled out of the trace leaving the trace of
holonomy term of Kapustin, which we denote simply by $\Tr\hol(\theta;\partial
\sigma)$.
The Abelian terms maybe be then dealt with by the usual combinatorial
methods to give the term $H^{-1}_\partial(J;\sigma)$. Combining all terms
corresponding to the bundle gerbe $\zeta$ and module $E$ we now have
\begin{equation}
\Tr\hol(\theta)  H^{-1}_\partial(J;\partial \sigma) H_{int}(\zeta;\sigma)
H_\partial(J;\partial \sigma) = \Tr\hol(\theta;\partial \sigma) 
H_{int}(\zeta;\sigma)
\end{equation}
The contributions from the $B$-field and the torsion class $w$ are
as in the previous case,
\begin{equation}
H_{int}(B;\Sigma)H^{-1}_{int}(B;\sigma)\pfa H^{-1}_{int}(w;\sigma)
\end{equation}
The $A$-field now trivialises $P^*_w \otimes P_B \otimes P^*_\zeta$, so
the corresponding terms are the opposites of all of the 
$H_{int}$ terms in the previous expressions as well as $H_\partial(A;
\partial{\Sigma})$. Thus combining all terms and using $\partial \sigma
= \partial \Sigma$ we gain a combination of terms,
\begin{equation}
H_{int}(B;\Sigma)\cdot\pfa\cdot\Tr\hol(\theta;\partial \Sigma)\cdot 
H_\partial(A;
\partial{\Sigma})
\end{equation}
which is independent of
$\sigma$ and so the anomaly is cancelled.

The third case is where the $B$-field is non-torsion, so the 
class $\zeta$ is non-torsion and so does not represent a bundle
gerbe which admits a bundle gerbe module. To get around this 
it is possible to define bundle gerbe modules with infinite
dimensional fibres \cite{bcmms} which are acted on by
non-torsion bundle gerbes. These may be used to define a 
$\Tr\hol$ term which cancels the anomaly \cite{cajomu}.
The details of this approach are not particularly relevant here, however
we make note of it since it shows that the bundle gerbe theoretical
approach of the simpler cases described above leads to a way of
dealing with the general case.\\

{\bf $C$-Fields.} It seems likely that bundle gerbes could be useful in other
string theory applications. In particular it has been noted (\cite{sha},
\cite{sek}) that $C$-fields in five-brane theory may be represented 
locally by the following data:
\begin{eqnarray}
C_{\alpha} - C_\beta &=& dB_{\alpha\beta} \\
B_{\alpha\beta} + B_{\beta\gamma} + B_{\gamma \alpha} &=& d A_{\alpha 
\beta \gamma} \\
A_{\beta \gamma \delta} - A_{\alpha\gamma\delta} + A_{\alpha \beta \gamma}
- A_{\alpha \beta\gamma} &=& d\log h_{\alpha\beta\gamma\delta} \\
\delta h_{\alpha\beta\gamma\delta} &=& 1
\end{eqnarray}
This data defines a class in $H^3(M,\D^3)$ or an equivalence
class of bundle 2-gerbes. The actions which are defined using 
$C$-fields are not the holonomy of this bundle 2-gerbe since
they are usually defined in seven or eleven dimensions rather than 
three (\cite{wit},\cite{fhmm}). These actions are higher dimensional
generalisations of Chern-Simons theory, and while we do not have 
a theory of higher bundle gerbes that applies in such dimensions 
the actions may still be interpreted in terms of Deligne cohomology.
If the curvature of the $C$-field is $G$ then the seven dimensional 
Chern-Simons term is defined on a 7-manifold $M$ by its extension to
an 8-dimensional manifold $X$ as
\begin{equation}
CS_7(C) = \exp \int_X G \wedge G
\end{equation}
We may think of this as the holonomy of a Deligne class in $H^7(M,\D^7)$
with curvature $G \wedge G$. Such a Deligne class may be constructed
via a cup product. Let $[C] \in H^3(M,\D^3)$ be the Deligne class 
of the $C$-field. Then $[C] \cup [C]$ is a class in $H^7(M,\D^7)$
\footnote{Recall that the cup product of two classes in $H^p(M,\D^p)$ 
gives a class in $H^{2p+1}(M,\D^{2p+1})$ since the cup product is 
actually defined on $H^{p+1}(M,\Z(p+1)_D) \isom H^p(M,\D^p)$.}
with
curvature $G \wedge G$. The action may then be defined without the 
extension $X$ as the holonomy of this class over $M$. If the local
3-curving forms $C_\alpha$ are actually globally defined (corresponding
to $G$ being de-Rham trivial and the $C$-field representing a torsion
bundle 2-gerbe) then this may be expressed as
\begin{equation}
CS_7(C) = \exp \int_X C \wedge G
\end{equation}  
In the general case it would be necessary to use the formula for the 
cup product (definition \ref{defcp}) and to substitute the resulting
Deligne class into the general formula for holonomy given by
proposition \ref{propgh}. Given a 6-manifold $W$ then it
is possible to construct a line bundle by transgression, a 
local formula for the transition functions would be
given by equation \eqref{genG}.

An 11-dimensional Chern-Simons theory may be defined in a similar way.
This time the Deligne class is given by the triple cup product 
$[C]\cup[C]\cup[C]$ so that the curvature is $G\wedge G \wedge G$. The
holonomy is defined as an integral of this curvature over a 12-manifold,
an integral of $C\wedge G \wedge G$ over an 11-manifold or
more generally by a transgression formula. There is a transgression
line bundle obtained by considering the holonomy over 10-manifolds.

These observations give only a starting point for a bundle gerbe analysis
of $C$-fields and 5-brane theories. We have not analysed anomaly
cancellation for this theory however it is possible that our
approach to anomaly cancellation in the $D$-brane case could
also apply here to some extent.

\section{Axiomatic Topological Quantum Field Theory}

We would like to relate the properties of bundle gerbe holonomy and
transgression to the axioms of topological quantum field theory (TQFT) 
(\cite{ati},\cite{qui}). This arises from the relationship between holonomy
and topological actions that has been demonstrated in the previous sections,
however it should be noted that we have considered only classical actions.  
It is possible to proceed to topological quantum field 
theories using the technique of path integration (see \cite{wit2} 
for a discussion of the case of Chern-Simons theory),
however this is not generally well defined and will not be discussed here.
The axiomatic definition of TQFT is of interest however since TQFTs may be
derived from classical
theories satisfying similar axioms. 

As additional motivation we cite some relevant literature.
It has been noted that the line bundle obtained via transgression
of a gerbe \cite{brmc1} and the Chern-Simons lines \cite{fre}, 
which we have derived via transgression satisfy certain axioms which are
closely related to those of TQFT. 
This approach to
quantum Chern-Simons theory has been taken by Freed \cite{fre2} in the
case where the gauge is group is finite since in this case the path
integral reduces to a finite sum which is well defined. In this instance the
properties of the classical theory carry over to give the axiomatic 
properties of the quantum theory. Another approach is to consider 
{\it homotopy quantum field theories} \cite{butuwi} which in certain cases
are closely related to gerbes. 
Generalisations relating to higher categories have 
also 
been considered \cite{brtu}. Also Segal \cite{seg} considers an axiomatic 
approach to $B$-fields in string theory, following the axiomatic definition 
of conformal field theory (CFT). Thus the link between gerbes, topological 
field theories and axiomatic definition of such theories has arisen in a 
number of different ways.

We shall consider first the axioms given by Atiyah and then
examine the extent to which they relate to bundle gerbe theory.
\begin{definition}\cite{ati}\label{ati}
A {\it topological quantum field theory} (TQFT) in dimension $d$ defined over
a ground ring $\Lambda$, consists of a finitely generated 
$\Lambda$-module $Z(\Sigma)$ associated to each oriented closed
smooth $d$-manifold $\Sigma$, and an element $Z(X) \in Z(\partial X)$ 
associated to each oriented smooth $(d+1)$-manifold $X$. These are
required to satisfy the following axioms:
\begin{enumerate}
\item $Z$ is functorial with respect to orientation preserving 
diffeomorphism of $\Sigma$ and $X$, 
\item $Z$ is involutory, that is, reversing orientation of the manifold
gives the dual module, 
\item $Z$ is multiplicative under disjoint unions and gluing of 
manifolds.
\end{enumerate}
\end{definition}
We also note some further explanation from \cite{ati} about each of
these axioms.
\begin{itemize}
\item {\bf Functoriality.} Let $\phi: \Sigma \rightarrow \Sigma'$ be 
an orientation preserving diffeomorphism. Then there is an isomorphism
of modules $Z(\phi): Z(\Sigma) \rightarrow Z(\Sigma')$ such that
\mbox{$Z(\psi \circ \phi) = Z(\psi) \circ Z(\phi)$} where $\psi: 
\Sigma' \rightarrow \Sigma''$. When $\phi$ extends to an
orientation preserving diffeomorphism $X \rightarrow X'$, with
$\partial X = \Sigma$, $\partial X' = \Sigma'$ then the
isomorphism $Z(\phi)$ maps $Z(X)$ to $Z(X')$.
\item {\bf Involution.} In general a reverse in orientation gives a 
`dual' module. When $\Lambda$ is a field then a reverse of orientation
gives dual vector spaces. We need not be concerned here with details 
of the general case.
\item {\bf Multiplication.} In the case where $\Sigma$ and $\Sigma'$ are
disjoint we have $Z(\Sigma \cup \Sigma') = Z(\Sigma) \otimes Z(\Sigma')$.
If $X$ has boundary $\Sigma_1 \cup \Sigma_2$ and we cut $X$ along
$\Sigma_3$ to get two components such that $\partial X_1 =
\Sigma_1 \cup \Sigma_3$ and $\partial X_2 = \Sigma_2 \cup \Sigma_3$ then
\begin{equation}\label{pair}
Z(X) = <Z(X_1),Z(X_2)>
\end{equation}
which is defined to be
the natural pairing 
\begin{equation} 
Z(\Sigma_1)\otimes Z(\Sigma_3) \otimes 
Z(\Sigma_3)^* \otimes Z(\Sigma_2) \rightarrow Z(\Sigma_1) \otimes
Z(\Sigma_2)
\end{equation}

If we set $Z(\emptyset_d)= \Lambda$, where $\emptyset_d$ denotes the empty
$d$-manifold then we may extend this to the case where $X$ is closed and 
may be cut along $\Sigma$ to make $X = X_1 \cup_\Sigma X_2$. In this
case we also get \eqref{pair} however this time it represents a pairing
\begin{equation}
Z(\Sigma) \otimes Z(\Sigma)^* \rightarrow \Lambda
\end{equation}
If $\emptyset_{d+1}$ is considered as the empty $(d+1)$-manifold then we let 
$Z(\emptyset_{d+1}) = 1$.
\end{itemize}

Our model for relating these axioms to bundle gerbe theory will be the 
case where the module is a vector space defined as the fibre 
of a vector bundle. To stay consistent with the preceding work we 
shall allow $\Lambda = U(1)$, so instead of a vector space we have a
principal $U(1)$-space. This should be thought of in the same terms as the 
equivalence between principal bundles and associated vector bundles.
Holonomy and transgression will not define TQFTs in this sense, however we
have given the axioms in this form since they are well known in the literature.
Instead we consider the ``classical'' TQFTs which Quinn \cite{qui} uses 
in a study of Chern-Simons theory in terms of axiomatic TQFT. These theories
differ from the TQFTs defined above in that there is extra data associated
with the manifolds on which the theory is defined. It is required that 
there is a topology on this extra data, for example it may consist of 
a space of mappings. To be specific, if the holonomy of a $U(1)$ 
bundle $L \rightarrow M$ with connection $A$
was considered to be a theory associated with a closed 1-manifold $\Gamma$ 
(a disjoint union of loops)
then
the problem is that the theory does not just depend on $\Gamma$
itself, it also depends on the map of the $\Gamma$ into $M$. Recall
that to define the holonomy we pull $L$ back using this map, so the
extra data could be considered either as an appropriate equivalence 
class of maps of
$\gamma$ into $M$ or alternatively as a space of isomorphism classes of
bundles with connection over $\Gamma$. Since all such bundles are trivial then
this is actually a space of gauge equivalent connections which arises in the
path integral. Note that these theories are not to be confused with
classical theories which take values in a field (for example $\R$) 
and which differ significantly from the quantum theory in that the
multiplicative property involves a scalar product rather than 
a tensor product. See \cite{qui} for a discussion of the importance of this
difference.

Allowing for extra data as discussed above, the following examples 
all satisfy the axioms by the
results discussed in \S \ref{genres}. \\

{\bf \noindent Holonomy and Parallel Transport of Line Bundles.} 
In terms of definition \ref{ati} we are dealing with a 0-dimensional
theory where we think of a 0-dimensional manifold as a point. Let $(L,M)$ be
a line bundle with connection, then associated with any point $m \in M$ 
we have a group defined by $L_m$, the fibre of $L$ at $m$. Given
a path $\mu$ in $M$ with $\partial \mu = \{m_0,m_1 \}$ then 
there is an element $Z(\mu) \in Z(\partial \mu) = L_{m_0}^* \otimes
L_{m_1}$ which is defined by parallel transport. Given a closed loop
$\gamma$ then there is an element $Z(\gamma) \in Z(\emptyset) = U(1)$
which is defined by the holonomy around $\gamma$.\\

{\bf \noindent Holonomy and Parallel Transport of Bundle Gerbes.} These define
1-dimensional theories. Let $(P,Y,M)$ be a bundle gerbe with connection 
and curving. For any $\gamma \in \L M$ let $Z(\gamma)$ be the fibre
of the transgression bundle $L = \tau_{S^1} P$ at $\gamma$. Given
a 2-manifold $\Sigma$ with boundary $\partial \Sigma = 
\bigcup_i \gamma_i$ then there is an element $Z(\Sigma) \in
Z(\partial \Sigma) = \bigotimes_i L^{\sigma(i)}_{\gamma_i}$, where
$\sigma(i)$ is the orientation of the boundary component, defined by
the section of $\partial^{-1}L \rightarrow \Sigma M$ which is derived
from the holonomy. For a closed 2-manifold the element of $Z(\emptyset)
= U(1)$ is defined by the holonomy.\\

{\bf \noindent Holonomy and Parallel Transport of Bundle 2-Gerbes.} These
define 2-dimensional theories in precisely the same way as the previous
two examples so we omit details.\\

The bundle gerbe hierarchy and the properties of holonomy and 
transgression imply the existence of more general theories where 
$Z(\Sigma)$ is no longer a vector space but which
essentially satisfy the same axioms. We have used fibres of bundles 
in the place of the modules $Z(\Sigma)$, so the next step in the hierarchy
is to use the fibre of a bundle gerbe. Such a fibre is a $U(1)$-groupoid
in the sense of \cite{mur}. 
We shall review this construction here, the 
important point being that all operations on modules which
are required in definition \ref{ati} have analogous constructions in the 
groupoid setting. 
\begin{definition}
A {\it $U(1)$-groupoid} with base $X$ 
is a principal $U(1)$-bundle $P \rightarrow X^2$
which has a product which is a bundle morphism covering the map
$((x,y),(y,z)) \rightarrow (x,z)$. This product is required to be
associative.
\end{definition}
We may denote this groupoid by the pair $(P,X)$.
A $U(1)$-groupoid has an identity, given by a section of $P$ over the 
diagonal $(x,x) \in X^2$, and an inverse which is given by taking the 
dual bundle $P^* \rightarrow X^2$. The existence of the identity and 
inverse is implied by the definition (see \cite{mur} for details).
A morphism of $U(1)$-groupoids is
a morphism of $U(1)$-bundles which respects the product structure.
It is clear from the definition of a bundle gerbe that the fibre over
a point in the base has the structure of a $U(1)$-groupoid. 
Given a point $m$ in the base then the objects are all $y\in Y$ such
that $\pi(y) = m$ and the morphism between two objects $y_0$ and
$y_1$ is given by $P_{(y_0,y_1)}$. Composition of morphisms is given by
the groupoid product. The tensor product of two groupoids is defined by 
analogy with the tensor product of bundle gerbes. Given groupoids
defined by $(P,X)$ and $(Q,Y)$ the product is the 
groupoid defined by the tensor product bundle $P\otimes Q \rightarrow
X^{[2]} \times Y^{[2]}$. A trivialisation of a $U(1)$-groupoid $(P,X)$ is a 
$U(1)$-bundle $L \rightarrow X$ such that there is a bundle isomorphism 
$P_{(x_1,x_2)} \isom L_{x_1}^* \otimes L_{x_2}$.

We now consider theories for which
the modules $Z(\Sigma)$
are
replaced by groupoids and elements of modules $Z(M) \in Z(\partial M)$ 
are replaced by trivialisations. The reason for this is that when $Z(\partial
M)$ is the fibre of a bundle then an element of the fibre is determined by
a section which is equivalent to a trivialisation. All operations involved in
the axioms are replaced by those described above. If we are dealing with
a $d$-dimensional theory 
then we define $Z(\emptyset_d)$ to be the groupoid with one object, that is
$(P,x)$ where $x$ is a single point, which may also be viewed as a fibre of 
a bundle given by the trivial morphism $P_{(x,x)}$. 
This ensures consistency of the multiplicative property. Similar ideas have 
been explored by Freed \cite{fre2} using actions which take values in 
torsors.\\


{\bf \noindent Fibres of a bundle gerbe}. 
Let $(P,Y,M)$ be a bundle gerbe with connection
and curving. We use this to define a 0-dimensional groupoid theory. 
For any point $m \in M$ let $Z(m)$ be the fibre of the bundle gerbe
over $m$. Given a path $\mu$ from $m_0$ to $m_1$ there is a trivialisation
of $P^*_{m_0} \otimes P_{m_1}$ defined by the extension of the loop
space transgression to paths. Recall that when using the transgression 
approach to holonomy reconstruction the transition functions of the
bundle over the loop space extend to a trivialisation of a bundle gerbe
on $\P M$ which is isomorphic to the pull back of the 
original bundle gerbe by the boundary restriction map. Given a closed
loop $\gamma$ we have a fibre of the transgression bundle $L_\gamma$ 
which we consider as an element of the trivial groupoid.\\

{\bf \noindent The loop space transgression of a bundle 2-gerbe} If we
have a bundle 2-gerbe $(P,Y,X,M)$ then we may define a 1-dimensional
theory by applying the previous example to the bundle gerbe on $\L M$ 
which is obtained by transgression (see \S \ref{ltb2g}).

In theory this approach could be extended to an even more abstract setting 
by moving further up the bundle gerbe hierarchy.
By considering fibres of a bundle 2-gerbe over a point one would obtain
a theory where the modules are replaced by 2-groupoids.\\

Finally we comment on the fact the theories
which we have described here all correspond
to theories involving modules (or groupoids) 
which are one dimensional vector spaces. This
is because there is not currently a satisfactory theory of non-Abelian 
bundle gerbes, so we only have generalisations of line bundles and
not vector bundles of higher rank.